\documentclass{book}
\usepackage{latexsym}
\usepackage{amsmath,amssymb}
\usepackage{txfonts}
\usepackage[dvips]{graphicx, color}
\begin{document}
\title{\bf {\LARGE Arthur's Periods, 
Regularized Integrals and 
Refined Structures of Non-Abelian~$L$-Functions}}
\author{\bf \large{WENG Lin}}
\date{May 22, 2005}
\maketitle
In 1970's, in order to establish a trace formula for reductive groups, 
Arthur introduced an analytic truncation. It plays
a key role in the study of theories of automorphic forms, $L$-functions 
and representations.

Historically, the simplest form of such a truncation was first used by 
Rankin and Selberg respectively along with the line of
 what we now call the Rankin-Selberg method. The essential point is that 
classical Eisensetin series for the principal 
 congruence group $SL_2(\mathbb Z)$ is not integrable over the 
corresponding fundamental domain, due to the existence of 
the constant term. It is to remove such a constant term, which is only of 
slow growth, in Rankin-Selberg method, we need to use 
 an analytic truncation so that after truncating the constant term, 
the (remaining) truncated Eisenstein series becomes rapidly 
 decreasing, and hence integrable over the original fundamental domain.

Such a method was generalized by Langlands in his \lq extraordinary' 
study of Eisenstein series in general in 1960's. There 
with more general type of reductive groups, various parabolic subgroups 
play a central role via the so-called parabolic 
reduction. In terms of truncations, what is best known now is  the 
so-called Langlands' combinatorial
lemma. (However, in Langlands' version of reduction theory,
somehow only maximal parabolic subgroups are in the 
central stage.) 

All this may be viewed as  predecessors of Arthur's analytic truncation. 
To put such truncations in a simplest form, we may 
say that Arthur's truncation is the most powerful and fruitful method 
in getting the integrability of Eisenstein series: various 
level of parabolic subgroups are very-well oriented in a beautiful 
way as they should be so as to yield a neat integrability 
statement (for truncated Eisenstein series).

Besides making truncated automorphic forms become rapidly decreasing, 
key properties of Arthur's analytic truncation
are that (1) it is self-adjoint;  and (2) it is idempotent. (Hence
it is an orthogonal projection). 
Consequently, integrations over fundamental domains of 
truncated automotphic forms are equal to integrations of the original 
automorphic forms over certain compact subsets of  the 
original fundamental domain obtained by applying Arthur's analytic 
turncation to the constant function one. (See the main 
text for a more precise statement.)

On the other hand, there is also a very natural way to intrinsically 
obtain  arithmetically meaningful compact subsets of  
fundamental domains. It has its roots in Mumford's intersection stability 
for vector bundles over algebraic curves: 
The compact subsets are obtained from the so-called semi-stable lattices, 
a concept  we independently introduce when defining  
and studying non-abelian zeta functions for number fields. (We note here 
also that it was Sthuler who first introduces the 
stability concept for lattices in literature.) Later on, motivated by our 
understandings of these non-abelian zeta functions, in 
particular, after revealing the intrinsic relation of these zeta functions 
with Epstein type zeta functions, we further introduce 
more general non-abelian $L$-functions for global fields: 
Simply put, non-abelian 
$L$-functions are integrations of Esiensetin series associated 
with $L^2$-automorphic forms over  certain compact subsets obtained  
geometrically from truncations of the corresponding 
fundamental domains via staibility of lattices. 

Being naturally defined, these new $L$-functions admit basic properties 
such as that (1) they are well-defined and admit
unique meromorphic continuation; (2) they satisfy functional equations; 
and (3) their singularities may be described in a 
resonable way, based on Langlands' theory of Eisesnetin series. More 
surprisinly, albert special, we show that (4$)_2$ the 
analogue of the Riemann Hypothesis for rank two non-abelian zeta 
functions of all number fields hold, i.e., all zeros of
such functions lie on the critical line $\Re(s)=1/2$, based on a 
clever discovery of Lagarias and Suzuki, who first show that
the rank two non-abelian zeta for the rational field $\mathbb Q$ 
satisfies the analogue of the RH. (We should mention 
here that
in the case of function fields, examples show that modified RH are 
satisfied. Say in the case of rank two non-abelian zeta 
functions of elliptic curves, Siegel zeros do exist. However it is 
expected that even there the modified RH holds for rank two 
zeta functions.)

Parallelly, as we find later, in a study of relative trace formula and 
automorphic forms, Jacquet-Lapid-Rogawski introduce the 
concept of period of automorphic forms. Simply put, these periods are 
integrations of truncated automorphic forms over the 
corresponding fundamental domains.

As said above, due to the properties of analytic truncations being 
self-adjoint and being orthogonal projections, these general 
periods may be understood as integrations of automorphic forms over 
compact subsets obtained from Arthur's analytic 
truncation for the constant function one applied to the fundamental 
domains (i.e., by removing suitable neighborhood around 
cusps in a well-organized way).
In particular, by considering Eisenstein series, a special type of  
automotphic forms (special enough so that nice properties 
such as
meromorphic continuations and functional equations hold while general 
enough so that discrete spectrum may be obtained as 
\lq residues' of Eisenstein series while continuous spectrum may be 
obtained as \lq integrations' of them,) we are led to a 
special type of periods. For our own convenience, we will call these 
special type of period as Eisenstein periods.

One of the central themes of this paper is to relate these Esenstein 
priods of Jacquet-Lai-Lapid-Rogawski with our non-abelian 
$L$-functions. Simply put, we are going to expose the fundamental fact 
that Eisenstein periods correspond to what we call the 
abelian parts of our non-abelian $L$-functions.

As such, the first thing we should expose is an intrinsic relation between 
Arthur's analytic truncation and our geometric 
truncation. The first crue of such a connection is hidden in the work of 
Lafforgue on Langlands' conjectures for function fields.
In his clever work, Lafforgue first exposes a basic relation between 
Harder-Narasimhan's filtration and modified analytic 
truncation and vaguely says that somehow such a modified analytic 
truncation is closely related with Arthur's original 
truncation.

In this paper, we start with a discussion along this line of Lafforgue. 
Roughly speaking, instead of the so-called positive cones
in Arthur's analytic truncation defined by certain positivities 
with respect to the so-called fundamental weights, we use 
different type of cones, defined by positive
roots in a twisted way. (We point out here that it is not proper to
introduce analytic truncation associated with positive chambers, which are
cones defined using positive roots directly.) Such a change proves to be 
essential to our geometrically oriented truncation: it 
beautifully makes extra rooms for creating essential non-abelian 
contributations, which, we believe, should play a key role in studying 
non-abelain arithmetic aspects of global fields.

In this sense, our present study of  
refined structures of non-abelian $L$-functions is 

\noindent
(1) in exposing the fact that non-abelain
$L$-functions consist of an abelian part and an essential non-abelian part;

\noindent
(2)  in establishing a general formula for the abelian part.

However, all these two are proved to be formidable: We only have a partial 
success in studying abelian parts at this stage. 
More precisely, when the original Eisenstein series are defined by cusp 
forms of lower level parabolic subgroups, a fine study 
can be carried out in a very satisfactory way, via the so-called 
Rankin-Selberg methods, or better in its
 mordern version via the discussion of Zagier for moderate growth functions 
and/or regularized intergations over cones of 
 Jacquet-Lapid-Rogawski along the line of Langlands-Arthur's formula on inner 
product of Eisenstein series. 
 In general, we may pave the way of Arthur in his study of inner products of 
Eisenstein series. Even so,
 due to the complication of singularities of Eisensetin series, we expect
 that the situation becomes very very complicated, if it is not out of hand. 
(However, as one may expect, in the so-called 
 relative rank one case, a precise formula can also be established due to 
the fact that then the poles of the associated 
 Eisenstein series are simple poles, so via Langlands' theory of Eisenstein
 system, it is managable along with the line of Osborne 
 and Warner.)
 
 As for the so-called essential non-abelian part, we undetsand that this 
is a totally new field  which should be taken very 
 seriously. Such parts are governed by yet to be discoveried non-abelian 
reciprocity law. (See our Program for Geometric 
 Arithmetic for details.)  On the other hand,  a good example of such a 
non-abelian study is the so-called Kronecker
 limit formula, in which the non-abelian part, i.e., the non-constant part 
in the Fourier expansion conststing of Bessel $K$-
 functions plays a key role. 
 
Now it becomes quite clear that a detailed quantitative
study of Langlands' theory of 
Eisenstein systems for the groups at hands 
are urgently needed. Pioneer works  are Langlands' fundamental work and 
Arthur's work on the inner product of truncated 
Eisenstein series. (See also the beautiful treatment presented by 
Osborne-Warner.) We are going to study this aspect of the theory elsewhere, 
with the works of Jacquet and Moeglin-Waldspurger on $GL_r$ 
and Henry Kim on $G_2$ in mind. 
Instead here, we would like to draw reader's attention to a 
fundamental vanishing result of abelian parts of our 
non-abelain $L$-functions. This vanishing statement is an analogue of a
remarkable result of Jacquet-Lapid-Rogawski on 
Eisenstein periods resulting from certain non-existence of representation 
of groups over finite adeles due to Bernstein: While 
these abelian parts of non-abelian $L$-functions vanish identilcally, they 
do have an explicit expression in terms of certain 
combinations of integrations of partial constant terms. In this way, the 
total vanishing result leads to various non-trivial 
relations of abelian $L$-functions (appeared as coefficients of constant 
terms resulting from various ways of 
 taking parabolic reductions). It is our understanding that such a type 
of natural relations are the principal reasons behind 
 special yet very important relations among multiple  and ordinary zeta 
functions in literature.

The sole aim of this paper is to give an easy access
to our paper on \lq Non-Abelian $L$-Functions for Number Fields'. 
For this purpose, we add two chapters on Arthur's truncation and periods, 
which are simply taken from several related papers of Arthur and 
Jacquet-Lapid-Rogawski respectively. I sincerely hope that, with such a 
detailed treatment, the reader from now on will have no further problem 
to understand the definition and basic properties 
of non-abelian $L$-functions together with their abelian parts.
\vfill
\eject
\vskip 1.0cm
\centerline{\Large Contents}
\vskip 0.45cm
\noindent
Abstract
\vskip 0.5cm
\noindent
{\bf Ch.1. Space of $\mathcal O_K$-Lattices and Geometric Truncation}

\noindent
1.1. Space of $\mathcal O_K$-Lattices

1.1.1. Projective $\mathcal O_K$-modules

1.1.2. $\mathcal O_K$-Lattices

1.1.3. Space of $\mathcal O_K$-Lattices

\noindent
1.2. Geometric Truncation

1.2.1. Semi-Stable Lattices

1.2.2. Canonical Filtration

1.2.3. Compactness

\noindent
1.3. Space of $\mathcal O_K$-Lattices via Special Linear Groups
\vskip 0.5cm
\noindent
{\bf Ch. 2. Geometric Truncation and Analytic Truncation}

\noindent
2.1. Geometric Truncation: Revised

2.1.1. Slopes, Canonical Filtrations
 
2.1.2. Canonical Polygons and Geometric Truncation

\noindent
2.2. Partial Algebraic Truncation

2.2.1. Parabolic Subgroups

2.2.2. Logarithmic Map
 
2.2.3. Roots, Coroots, Weights and Coweights
 
2.2.4. Positive Cone and Positive Chamber
 
2.2.5. Example with $SL_r$
 
2.2.6. An Algebraic Truncations

\noindent
2.3. A Bridge between Algebraic and Analytic Truncations
 
\noindent
2.4. Global Relation between Analytic Truncation and Geometric Truncation
\vskip 0.5cm 
\noindent
{\bf Ch.3. Non-Abelian L-Functions}

\noindent
3.1. Moduli Spaces as Integration Domains

\noindent
3.2. Choice of Eisenstein Series: First Approach to Non-Abelian 
$L$-Function

\noindent
3.3. New Non-Abelian $L$-Functions

\noindent
3.4. Meromorphic Extension and Functional Equations
 
\noindent
3.5. Holomorphicity and Singularities
\vskip 0.5cm
\noindent
{\bf Ch. Arthur's Analytic Truncation}

\noindent
4.1. Langlands' Combinatorial Lemma and Arthur's Partition

4.1.1. Height Functions

4.1.2. Partial Truncation and First Estimation

4.1.3. Langlands' Combinatorial Lemma

4.1.4. Langlands-Arthur's Partition: Reduction Theory

\noindent
4.2. Arthur's Analytic Truncation

4.2.1. Definition

4.2.2. Idempotence

4.2.3. Self-Adjointness

4.2.4. $\Lambda^T\phi$ is Rapidly Decreasing
\vskip 0.5cm
\noindent
{\bf Ch.5. Abelian Part of Non-Abelian L-Functions}

\noindent
5.1. Arthur's  Periods

\noindent
5.2. Periods as Integrals over Truncated Fundamental Domains

\noindent
5.3. Geometrically Oriented Analytic Truncation

\noindent
5.4. Abelian and Essential Non-Abelain Parts
\vskip 0.5cm
\noindent
{\bf Ch. 6. Rankin-Selberg Method: Evaluation of Abelian Part}

\noindent
6.1. Abelian Part of Non-Abelian $L$-Functions: Fomal Calculation
 
6.1.1. From Geometric Truncation to Periods
 
6.1.2. Regularization
 
6.1.3. Berstein's Principle
 
6.1.4. Constant Terms of Eisenstein Series
 
6.1.5. Abelian Part of Non-Abelian $L$-Functions

\noindent
6.2. Regularized Integration

6.2.1. Integration over Cones
 
\noindent
6.3 Positive Cones and Positive Chambers

6.3.1. Modified Truncations
 
6.3.2. JLR's Period of an Automorphic Form
 
6.3.3. Arthur's Period
 
6.3.4. Bernstein Principle

6.3.5. Eisenstein Period
 
\noindent
6.4 Non-Abelian Part of Abelain L-Functions: Open Problems

6.4.1. Cone Corresponding to ${\bold 1}(p_P^g>_Pp)$

6.4.2. Open Problems
\vskip 0.30cm
\noindent
Acknowledgements
\vskip 0.5cm
\noindent
{\bf REFERENCES}
\chapter{Space of $\mathcal O_K$-Lattices and Geometric Truncation}
\section{Space of $\mathcal O_K$-Lattices}
\subsection{Projective $\mathcal O_K$-modules}

Let $K$ be an algebraic number field.
 Denote by $\mathcal O_K$ the ring of integers of $K$. Then
an $\mathcal O_K$-module $M$ is called {\it projective} 
 if there exists an $\mathcal O_K$-module $N$ such that 
$M\oplus N$ is a free $\mathcal O_K$-module. 
Easily,  for a fractional ideal $\frak a$,
 $$P_{\frak a}:=P_{r;\frak a}:=\mathcal O_K^{r-1}\oplus 
\frak a$$ is a rank $r$ projective 
 $\mathcal O_K$-module. The nice thing is that such 
types of projective $\mathcal O_K$-modules, up to
  isomorphism, give all rank $r$ projective 
$\mathcal O_K$-modules. Indeed, we have the following well-known
 \vskip 0.30cm
\noindent
{\bf{\Large Proposition}.} 
 (1)  {\it For fractional ideals $\frak a$ and $\frak b$, 
$P_{r;\frak a}\simeq P_{r;\frak b}$ if and only if 
 $\frak a\simeq \frak b$;}
 
 \noindent
 (2) {\it For a rank $r$  projective $\mathcal O_K$-module $P$, 
there exists a  fractional ideal $\frak a$ such that
 $P\simeq P_{\frak a}.$}
 
The reader can find a complete proof in [FT]. 
\vskip 0.30cm  
Now use the natural inclusion of fractional ideals in $K$ to embed
$P_{r;\frak a}$ into $K^r$. View an element in $K^r$ as a column
vector. In such a way, any $\mathcal O_K$-module
 morphism  $A:P_{r;\frak a}\to P_{r;\frak b}$ may 
be written down as an element in $A\in M_{r\times r}(K)$ so that
the image $A(x)$ of $x$ under $A$ becomes simply the matrix multiplication
 $A\cdot x$. 

\subsection{$\mathcal O_K$-Lattices}

Let $\sigma$ be an Archimedean place of $K$, and  
$K_\sigma$ be the $\sigma$-completion of $K$. It is 
well-known that  $K_\sigma$ is 
either equal to $\mathbb R$, or equal to $\mathbb C$.
 Accordingly, we call $\sigma$ (to be) real or complex, 
write sometimes in terms of $\sigma:\mathbb R$ or 
 $\sigma:\mathbb C$ accordingly.
 
Recall also that a finite dimensional $K_\sigma$-vector 
space $V_\sigma$ is called a {\it metrized 
space} if it is equipped with an inner product.
\vskip 0.20cm
By definition, an $\mathcal O_K$-{\it lattice} $\Lambda$ 
consists of 

\noindent
(1) a projective $\mathcal O_K$-module 
$P=P(\Lambda)$ of finite rank; and 

\noindent
(2) an inner product on the vector space 
$V_\sigma:=P\otimes_{\mathcal O_K}K_\sigma$ for each of the
 Archmidean place $\sigma$ of $K$.

\noindent
Set $V=P\otimes_{\mathbb Z}\mathbb R$, 
then $V=\prod_{\sigma\in S_\infty}V_\sigma$,
  where $S_\infty$ denotes the collection of all (inequivalent) 
Archimedean places of $K$. Indeed, this is a
  direct consequence of the fact that as a $\mathbb Z$-module,
 an $\mathcal O_K$-ideal is of rank $n=r_1+2r_2$ 
  where $n=[K:\mathbb Q]$, $r_1$ denotes the number of real places 
and $r_2$ denotes the number of complex places
   (in $S_\infty$).

\subsection{Space of $\mathcal O_K$-Lattices}

Let $P$ be a rank $r$ projective $\mathcal O_K$-module. Denote by 
$GL(P):=\mathrm{Aut}_{\mathcal O_K}(P)$. Let
 $\widetilde{\mathbf{\Lambda}}:=\widetilde{\mathbf{\Lambda}}(P)$ be 
the space of ($\mathcal O_K$-)lattices $\Lambda$ whose 
 underlying $\mathcal O_K$-module is $P$. For $\sigma\in S_\infty$, 
let $\widetilde{\mathbf{\Lambda}}_\sigma$ be 
 the space of inner products on $V_\sigma$; if a basis is chosen for
 $V_\sigma$ as a real or a complex vector space according to whether
$\sigma$ is real or complex,   
 $\widetilde{\mathbf{\Lambda}}_\sigma$ may be realized as
 an open set of a real or complex vector space. (See the discussion in \S 1.3 
below for details.) We have
  $\widetilde{\mathbf{\Lambda}}=\prod_{\sigma\in S_\infty}
\widetilde{\mathbf{\Lambda}}_\sigma$ and this provides 
  us with a natural topology on $\widetilde{\mathbf{\Lambda}}.$

Consider $GL(P)$ to act on $P$ from the left. 
Given $\Lambda\in \widetilde{\mathbf{\Lambda}}$ and
$u,w\in V_\sigma$, let $\langle u,w\rangle_{\Lambda,\sigma}$ or
$\langle u,w\rangle_{\rho_\Lambda(\sigma)}$
denote the value of the inner product on the 
vectors $u$ and $w$ associated to the lattice $\Lambda$.
Then, if $A\in GL(P)$, we may define a new 
lattice $A\cdot\Lambda $ in 
$\widetilde{\mathbf{\Lambda}}$ by the
following formula $$\langle u,w
\rangle_{A\cdot\Lambda ,\sigma}:=
\langle A^{-1}\cdot u,A^{-1}\cdot w\rangle_{\Lambda,\sigma}.$$ 
This defines an action of $GL(P)$ on 
$\widetilde{\mathbf{\Lambda}}$ from the left. Clearly, 
then the map $v\mapsto Av$ gives an isometry
 $\Lambda\cong A\cdot\Lambda $ of the lattices. 
(By an isometry here, we mean an isomorphism of 
 $\mathcal O_K$-modules for the underlying 
$\mathcal O_K$-modules subjecting the condition that the
 isomorphism also keeps the inner product unchanged.) 
Conversely, suppose that $A:\Lambda_1\cong\Lambda_2$
  is an isometry of $\mathcal O_K$-lattices, each of 
which is in $\widetilde{\mathbf{\Lambda}}$. Then, $A$
   defines an element, also denoted by $A$, of $GL(P)$. 
Clearly $\Lambda_2\cong A\cdot \Lambda_1.$
\vskip 0.20cm 
Therefore, the orbit set $GL(P)\backslash 
\widetilde{\mathbf{\Lambda}}(P)$ can be regarded as the set of isometry
 classes of $\mathcal O_K$-lattices whose underlying 
$\mathcal O_K$-modules are isomorphic to $P$.

\section{Geometric Truncation}
\subsection{Semi-Stable Lattices}

Let $\Lambda$ be an $\mathcal O_K$-lattice with 
underlying $\mathcal O_K$-module $P$. Then any submodule
$P_1\subset P$ can be made into an 
$\mathcal O_K$-lattice by restricting the inner 
product on each $V_\sigma$ 
to the subspace $V_{1,\sigma}:=P_1\otimes_KK_\sigma$. 
Call the resulting $\mathcal O_K$-lattice $\Lambda_1:=\Lambda\cap P_1$
 and write $\Lambda_1\subset \Lambda$. If moreover,
 $P/P_1$ is  projective, we say that $\Lambda_1$
is a {\it sublattice} of $\Lambda$.

The orthogonal projections $\pi_\sigma:V_\sigma\to 
V_{1,\sigma}^\perp$ to the orthogonal complement $V_{1,\sigma}^\perp$
of $V_{1,\sigma}$ in $V_{\sigma}$
provide isomorphisms 
$(P/P_1)\otimes_{\mathcal O_K}K_\sigma\simeq V_{1,\sigma}^\perp$, 
which can be used to make $P/P_1$ into an
 $\mathcal O_K$-lattice. We call this resulting lattice the 
 {\it quotient lattice} of $\Lambda$ by $\Lambda_1$, and denote it by 
$\Lambda/\Lambda_1$. 
\vskip 0.20cm
There is a procedure called {\it restriction of scalars} 
which makes an $\mathcal O_K$-lattice into a standard 
$\mathbb Z$-lattice. Recall that $V=\Lambda
\otimes_{\mathbb Z}\mathbb R=\prod_{\sigma\in S_\infty}V_\sigma$. 
Define an inner product on the real vector space $V$ by
 $$\langle u,w\rangle_\infty:=
\sum_{\sigma:\mathbb R}\langle u_\sigma,w_\sigma\rangle_\sigma
+\sum_{\sigma:\mathbb C}\mathrm{Re}\,\langle u_\sigma,
w_\sigma\rangle_\sigma.$$ Let $\mathrm{Res}_{K/\mathbb Q}\Lambda$ 
denote the $\mathbb Z$-lattice obtained by 
equipped $P$, regarding as a $\mathbb Z$-module, with this 
inner product (at the unique infinite place 
$\infty$ of $\mathbb Q$).

We let $\mathrm{rk}(\Lambda)$ denote the $\mathcal O_K$-module
 rank of $P$ (or of $\Lambda$), and define the 
{\it Lebesgue volume} of $\Lambda$, 
denoted by $\mathrm{Vol}_{\mathrm{Leb}}(\Lambda)$, to be
the (co)volume of the lattice $\mathrm{Res}_{K/\mathbb Q}
\Lambda$ inside its inner product space $V$. 

Clearly, if $P'$ is a submodule of 
finite index in $P$, then 
$$\mathrm{Vol}_{\mathrm{Leb}}(\Lambda')=
[P:P']\mathrm{Vol}_{\mathrm{Leb}}(\Lambda),$$ 
where $\Lambda'=\Lambda\cap P$ is the lattice induced from $P'$.
\vskip 0.20cm
\noindent 
{\bf{\large Examples}}. Take $P=\mathcal O_K$ and for each place 
$\sigma$, let $\{1\}$ be an orthonormal basis of 
$V_\sigma=K_\sigma$, i.e., equipped 
$V_\sigma=\mathbb R$ or $\mathbb C$ with the standard Lebesgue
 measure. This makes $\mathcal O_K$ into an 
$\mathcal O_K$-lattice $\overline{\mathcal O_K}=(\mathcal O_K,\mathbf 1)$ 
in a natural  way. It is a well-known fact, see e.g.,
[L1], that $$\mathrm{Vol}_{\mathrm{Leb}}
\Big(\overline{\mathcal O_K}\Big)=2^{-r_2}\cdot \sqrt{\Delta_F},$$ where
  $\Delta_F$ denotes the absolute value of the discriminant of $K$.

 More generally, take $P=\frak a$ an fractional 
  idea of $K$ and equip the same inner product as above 
on $V_\sigma$. Then $\frak a$ becomes an
   $\mathcal O_K$-lattice $\overline {\frak a}=(\frak a,\mathbf 1)$ 
in a natural 
way with $\mathrm{rk}(\frak a)=1$.  It is a well-known fact, see e.g.,
[Neu], that
$$\mathrm{Vol}_{\mathrm{Leb}}\Big(\overline {\frak a}\Big)=
2^{-r_2}\cdot \Big(N(\frak a)\cdot\sqrt{\Delta_K}\Big),$$ 
where $N(\frak a)$ denote the norm of $\frak a$.
\vskip 0.30cm
Due to the appearence of the factor $2^{-r_2}$, we  also define 
the {\it canonical volume} of $\Lambda$, denoted
 by $\mathrm{Vol}_{\mathrm{can}}(\Lambda)$ or simply by
$\mathrm{Vol}(\Lambda)$, to be 
$2^{r_2\mathrm{rk}(\Lambda)}\mathrm{Vol}_{\mathrm{Leb}}(\Lambda).$ 
 So in 
particular, $$\mathrm{Vol}\Big(\overline 
{\frak a}\Big)=
N(\frak a)\cdot\sqrt{\Delta_K},$$ with $$\mathrm{Vol}
\Big(\overline {\mathcal O_K}\Big)
=\sqrt{\Delta_K}$$ as its special case.

Now we are ready to introduce our first key definition.

\vskip 0.20cm
\noindent
{\bf \Large Definition.} An $\mathcal O_K$ lattice $\Lambda$ is called 
{\it semi-stable} (resp. {\it stable}) if for any proper 
sublattice $\Lambda_1$ of $\Lambda$,
 $$\mathrm{Vol}(\Lambda_1)^{\mathrm{rk}(\Lambda)}\geq\,(\mathrm{resp.}>)
\,\mathrm{Vol}(\Lambda)^{\mathrm{rk}(\Lambda_1)}.$$ 

Clearly the last inequality is equivalent to 
$$\mathrm{Vol}_{\mathrm{Leb}}(\Lambda_1)^{\mathrm{rk}(\Lambda)}\geq
\mathrm{Vol}_{\mathrm{Leb}}
(\Lambda)^{\mathrm{rk}(\Lambda_1)}.$$ So it does not matter which volume, 
the canonical one or the Lebesgue one, we use.
\vskip 0.20cm
\noindent
{\bf Remark.}  The canonical measures has an 
advantage theoretically. For example, we have
the following

\noindent
{\bf {\large Arakelov-Riemann-Roch Formula}:} 
For an  $\mathcal O_K$-lattice $\Lambda$ of rank $r$, 
$$-\log \Big(\mathrm{Vol}(\Lambda)\Big)=
\mathrm{deg}(\Lambda)-
\frac {r}{2}\log\Delta_K.$$ 
(For the reader who does not know the definition of the Arakelov 
degree, she or he may simply take this 
relation as the definition.)

\subsection{Canonical Filtration}

Based on stability, we may introduce a more general 
geometric truncation for the space of lattices. For this, we start with
the following well-known
\vskip 0.30cm
\noindent
{\bf {\large Lemma.}} {\it For a fixed ${\mathcal O}$-lattice 
$\Lambda$,}
$\Big\{\mathrm{Vol}(\Lambda_1):\Lambda_1\subset\Lambda\Big\}\subset 
{\mathbb R}_{\geq 0}$ {\it is  discrete and bounded from
below.}
\vskip 0.20cm
\noindent
{\it Sketch of proof.} It is based on the follows:

\noindent 
1) If a lattice $\Lambda_1$ induced from a submodule of
 the lattice $\Lambda$ has the minimal volume among all
lattices induced from submodules of $\Lambda$ of the same
rank, then $\Lambda_1$ is a sublattice;

\noindent
2) Taking wedge product, we may further assume that the rank is 1. 
It is then clear.
\vskip 0.30cm
As a direct consequence, we have the following

\noindent
{\bf {\large Proposition.}} {\it Let $\Lambda$ be an ${\mathcal O}$-lattice. 
Then}

\noindent
(1) ({\bf Canonical Filtration}) {\it There exists a unique 
filtration of proper sublattices
$$0=\Lambda_0\subset\Lambda_1\subset\cdots \subset\Lambda_s=\Lambda$$ 
such that}
\begin{description}
\item[(i)] {\it for all $i=1,\cdots, s$, $\Lambda_i/\Lambda_{i-1}$ 
is semi-stable}; and

\item[(ii)] {\it for all} $j=1,\cdots s-1$,
$$\Big(\mathrm{Vol}(\Lambda_{j+1}/\Lambda_j)\Big)^{
\mathrm{rk}(\Lambda_j/\Lambda_{j-1})}>
\Big(\mathrm{Vol}(\Lambda_{j}/\Lambda_{j-1})\Big)^{
\mathrm{rk}(\Lambda_{j+1}/\Lambda_{j})};$$
\end{description}
\noindent
(2) ({\bf Jordan-H\"older Filtration}) {\it If moreover, $\Lambda$ 
is semi-stable, then there exists a filtration of proper
sublattices,
$$0=\Lambda^{t+1}\subset\Lambda^t\subset\cdots \subset\Lambda^0=\Lambda$$ 
such that}
\begin{description}
\item [(i)] {\it for all $k=0,\cdots, t$, $\Lambda^k/\Lambda^{k+1}$ 
is stable; and}

\item[(ii)] {\it for all} $l=1,\cdots, t$,
$$\Big(\mathrm{Vol}(\Lambda^{l}/\Lambda^{l+1})\Big)^{
\mathrm{rk}(\Lambda^{l-1}/\Lambda^l)}
=\Big(\mathrm{Vol}(\Lambda^{l-1}/\Lambda^l)\Big)^{
\mathrm{rk}(\Lambda^l/\Lambda^{l+1})}.$$
\end{description}
\noindent
{\it Furthermore, up to isometry, the graded lattice 
${\rm Gr}(\Lambda):=\oplus_{k=0}^t\Lambda^k/\Lambda^{k+1}$ 
is uniquely determined by $\Lambda.$}
\vskip 0.20cm
\noindent
{\it Sketch of Proof.} Existence is clear and the uniqueness of 
Jordan-H\"older graded lattice is fairly standard. Let us look at 
the uniqueness of canonical filtration, an analogue of the well-known
Harder-Narasimhan filtration for vector bundles, a bit carefully: Existence of 
two such filtrations will lead to a contradiction 
by applying  fact that for any two sublattices $\Lambda_1,\, \Lambda_2$  
of $\Lambda$,
 $$\mathrm{Vol}\Big(\Lambda_1/(\Lambda_1\cap\Lambda_2)\Big)\geq
\mathrm{Vol}\Big((\Lambda_1+\Lambda_2)/\Lambda_2\Big).$$

\subsection{Compactness}

For a lattice $\Lambda\in \widetilde{\mathbf{\Lambda}}$ 
with the associated canonical filtration
$$0=\overline\Lambda_0\subset\overline\Lambda_1\subset\cdots 
\subset\overline\Lambda_s=\Lambda$$
define the associated {\it canonical polygone} 
$\overline p_\Lambda:[0,r]\to {\mathbb R}$ by the following conditions:
\begin{description}
\item[(1)] $\overline p_\Lambda(0)=\overline p_\Lambda(r)=0$;

\item[(2)] $\overline p_\Lambda$ is affine over the closed interval
$[\mathrm{rk}\overline\Lambda_i,\mathrm{rk}\overline\Lambda_{i+1}]$; and

\item[(3)] $\overline p_\Lambda(\mathrm{rk}\overline\Lambda_i)
=\mathrm{deg}(\overline\Lambda_i)-\mathrm{rk}(\overline\Lambda_i)
\cdot\frac{\mathrm{deg}(\overline\Lambda)}{r}.$
\end{description}
\vskip 0.30cm
Clearly, the canonical polygone is well-defined on the space 
$\widetilde{\mathbf{\Lambda}}$.

To go further, let us introduce an operation among
 the lattices in $\widetilde{\mathbf{\Lambda}}$. If $T$ 
is a positive real number, then from $\Lambda$, 
we can produce a new $\mathcal O_K$-lattice called
$\Lambda[T]$ by multiplying each of the inner products 
on $\Lambda$, or better, on $\Lambda_\sigma$ for 
$\sigma\in S_\infty$, by $T^2$. 
Obviously, the $[T]$-construction changes volumes of lattices 
in the following way
 $$\mathrm{Vol}(\Lambda[T])
=T^{\mathrm{rk}(\Lambda)\cdot[K:\mathbb Q]}\cdot \mathrm{Vol}(\Lambda).$$
That is to say,  the $[T]$-construction naturally fixes 
a specific volume for a certain family of lattices,
while  does not really change the \lq essential' structures of 
lattices involved.

Accordingly, let  
${\mathbf{\Lambda}}={\mathbf{\Lambda}}(P)$ be the quotient topological
space of $\widetilde{\mathbf{\Lambda}}$ modulo the equivalence relation 
$\Lambda\sim\Lambda[T]$. 
Since canonical polygone is  invariant under the scaling
operation and hence decends to this quotient space ${\mathbf{\Lambda}}$. 

With all this, we are ready to state the following fundamental 
\vskip 0.30cm
\noindent
{\bf{\Large Theorem.}} 
{\it For any fixed  convex  polygone $p:[0,r]\to {\mathbb R}$, 
the subset
$\Big\{[\Lambda]\in {\mathbf{\Lambda}}:\overline p_{\Lambda}\leq p\Big\}$ is 
compact in ${\mathbf{\Lambda}}$.}
\vskip 0.20cm
\noindent
{\it Sketch of Proof.} This is based on the follows:

\noindent
(1) Dirichlet's Unit Theorem (see the section below for details); and 

\noindent
(2) Minkowski's Reduction Theory: Due to the fact that volumes of 
lattices involved are fixed, semi-stability condition implies that 
the first Minkowski successive minimums of these lattices 
admit a natural lower bound away 
from 0 (depending only on $r$). Hence by the standard 
reduction theory, see e.g., Borel [Bo1,2], the subset
$\Big\{[\Lambda]\in {\mathbf{\Lambda}}:\overline p_{\Lambda}\leq p\Big\}$ 
is compact.

\section{Space of $\mathcal O_K$-Lattices via Special Linear Groups}

Recall that an $\mathcal O_K$-lattice 
$\Lambda$ consists of  a underlying projective 
$\mathcal O_K$-module $P$ and a metric structure on the space 
$V=\Lambda\otimes_{\mathbb Z}\mathbb R=
\prod_{\sigma\in S_\infty}V_\sigma.$ Moreover, for 
the projective $\mathcal O_K$-module $P$, in assuming 
that the $\mathcal O_K$-rank of $P$ is $r$, we can identify 
$P$ with one of the $P_i:=P_{r;\frak a_i}:=
\mathcal O_K^{(r-1)}\oplus\frak a_i$, where
 $\frak a_i, i=1,\cdots,h,$ are chosen integral $\mathcal O_K$-ideals 
so that $\Big\{[\frak a_1],\,[\frak a_2],\,\ldots,\,[\frak a_h]\Big\}=CL(K)$
 the class group of $K$. 
In the sequel, we often use $P$ as a running symbol for the $P_i$'s.

With this said,  via the Minkowski embedding 
$K\hookrightarrow \mathbb R^{r_1}\times\mathbb C^{r_2}$,  
we obtain a natural embedding for $P$:
$$P:=\mathcal O_K^{(r-1)}\oplus\frak a\hookrightarrow 
K^{(r)}\hookrightarrow \Big(\mathbb R^{r_1}\times
\mathbb C^{r_2}\Big)^{r}\cong \Big(\mathbb R^r\Big)^{r_1}\times
\Big(\mathbb C^r\Big)^{r_2},$$ which is simply the space 
$V=\Lambda\otimes_{\mathbb Z}\mathbb R$ above. 
As a direct consequence,  our lattice $\Lambda$ then is determined
by a metric structure on $V=\prod_{\sigma\in S_\infty}V_\sigma$, or the same, 
on $\big(\mathbb R^r\big)^{r_1}\times\big(\mathbb C^r\big)^{r_2}$.
Note that all metrized structures on  $\mathbb R^{r}$ 
(resp. on $\mathbb C^{r}$) are parametrized by the quotient space
$GL(r,\mathbb R)/O(r)$ (resp. $GL(r,\mathbb C)/U(r)$).
Consequently, metrized structures on 
$\big(\mathbb R^r\big)^{r_1}\times\big(\mathbb C^r\big)^{r_2}$ 
are parametrized by the space 
$$\Big(GL(r,\mathbb R)\Big/O(r)\Big)^{r_1}\times 
\Big(GL(r,\mathbb C)\Big/U(r)\Big)^{r_2}.$$
 
Motivated by the $[T]$-construction for lattices, naturally
at the group level, we need to shift our discussion from 
general linear group $GL$ to  special linear group $SL$.
For this, let us start with a local discussion on $\mathcal O_K$-lattice 
structures. 
 
First, look at complex places $\tau$. Clearly, by fixing a branch of the 
$n$-th root, we get  natural identifications
  $$\begin{matrix}GL(r,\mathbb C)&\to& SL(r,\mathbb C)\times \mathbb C^*&\to& 
 SL(r,\mathbb C)\times S^1 \times \mathbb R_+^*\\
 g&\mapsto&\Big(\frac{1}{\root{r}\of{\det g}}g,\det g\Big)&\mapsto& 
\Big(\frac{1}{\root{r}\of{\det g}}g,\frac{\det g}{|\det g|},|\det g|\Big)
\end{matrix}$$ and $$U(r)\to SU(r)\times S^1,\qquad  
U\mapsto\Big(\frac{1}{\root{r}\of{\det U}}U,\det U\Big),$$ where 
$SL$ (resp. $SU$) denotes the special linear group
 (resp. the special unitary group) and $S^1$ denotes the 
unit circle $\{z\in\mathbb C: |z|=1\}$ in 
 $\mathbb C^*$. Consequently,  
we have the following natural identification $$GL(r,\mathbb C)\Big/U(r)\cong 
\Big(SL(r,\mathbb C)\Big/SU(r)\Big)\times \mathbb R_+^*.$$

Then, let us trun to real places $\sigma$. This is slightly complicated. 
As an intermediate step, we use  subgroups 
$$GL^+(r,\mathbb R):=
  \Big\{g\in GL(r,\mathbb R):\det g>0\Big\}\quad 
\mathrm{and}\quad O^+(r):=\Big\{A\in O(r,\mathbb R):\det g>0\Big\}.$$ 
Clearly, we have the following relations;

\noindent
(i) $O^+(r)=SO(r)$, the special orthogonal group consisting of 
these $A$'s in $O(r)$ whose determinants are exactly 1; and 

\noindent
(ii) $GL(r,\mathbb R)\Big/O(r)\cong GL^+(r,\mathbb R)\Big/SO(r);$ moreover

\noindent
(iii) There is an identification  
$$GL^+(r,\mathbb R)\to  SL(r,\mathbb R)\times \mathbb R_+^*, \qquad
  g\mapsto  \Big(\frac{1}{\root{r}\of{\det g}}g,\det g\Big).$$ 

\noindent
Consequently,  we have a natural identification
 $$GL(r,\mathbb R)\Big/O(r)\cong \Big(SL(r,\mathbb R)\Big/SO(r)\Big)
\times \mathbb R_+^*.$$ 

\vskip 0.30cm 
Now we are ready to resume our global 
discussion on $\mathcal O_K$-lattices of rank $r$.
From  above, the metrized structures on 
$V=\prod_{\sigma\in S_\infty} V_\sigma\simeq(\mathbb R^r)^{r_1}\times
(\mathbb C^r)^{r_2}$ are  
parametrized by the space $$\bigg(\Big(SL(r,\mathbb R)\Big/SO(r)\Big)^{r_1}
\times \Big(SL(r,\mathbb C)\Big/SU(r)\Big)^{r_2}\bigg)\times 
(\mathbb R_+^*)^{r_1+r_2}.$$ 
Furthermore, when we really work with $\mathcal O_K$-lattice strucures on $P$, 
i.e., with the space $\bold\Lambda=\bold\Lambda(P)$, from 
the above parametrized space of metric 
structures on $V=\prod_{\sigma\in S_\infty} V_\sigma$, we need 
to further factor out $GL(P)$, i.e., the automorphism  group 
$\mathrm{Aut}_{\mathcal O_K}(\mathcal O_K^{(r-1)}\oplus\frak a)$ of 
$\mathcal O_K^{(r-1)}\oplus\frak a$ as $\mathcal O_K$-modules. So our next aim
is to use $SL$ to understand the quotient space
$$GL(P)\backslash\Bigg(\bigg(\Big(SL(r,\mathbb R)\Big/SO(r)\Big)^{r_1}
\times \Big(SL(r,\mathbb C)\Big/SU(r)\Big)^{r_2}\bigg)\times 
(\mathbb R_+^*)^{r_1+r_2}\Bigg).$$
  
As such, naturally, now we want  

\noindent
(a) To study the structure of the group 
$\mathrm{Aut}_{\mathcal O_K}(\mathcal O_K^{(r-1)}\oplus\frak a)$ in terms
of $SL$ and units; 
and 

\noindent
(b) To see how this group acts on the space of metrized structures 
$$\bigg(\Big(SL(r,\mathbb R)/SO(r)\Big)^{r_1}
\times \Big(SL(r,\mathbb C)/SU(r)\Big)^{r_2}\bigg)\times 
(\mathbb R_+^*)^{r_1+r_2}.$$
\vskip 0.20cm
View $\mathrm{Aut}_{\mathcal O_K}(\mathcal O_K^{(r-1)}\oplus\frak a)$ 
as a subgroup of $GL(r,K)$. Then, 
   without too much difficulty we see that 
$$\begin{aligned}~&\mathrm{Aut}_{\mathcal O_K}
(\mathcal O_K^{(r-1)}\oplus\frak a)
=GL(r,\mathcal O_K^{(r-1)}\oplus\frak a)\\
:=&\Bigg\{(a_{ij})\in GL(r,K):\begin{aligned} &a_{rr} \&
a_{ij}\in \mathcal O_K,\\
& a_{ir}\in\frak a,\ a_{rj}\in\frak a^{-1},\end{aligned}  i,j=1,\cdots,r-1;
\qquad \det(a_{ij})\in U_K\Bigg\}.\end{aligned}$$
In other words, $$\mathrm{Aut}_{\mathcal O_K}
(\mathcal O_K^{(r-1)}\oplus\frak a)=\Bigg\{A\in GL(r,K)\cap\begin{pmatrix}
&&&\frak a\\
&\mathcal O_K&&\vdots\\
&&&\frak a\\
\frak a^{-1}&\ldots&\frak a^{-1}&\mathcal O_K\end{pmatrix}:\det A
\in U_K\Bigg\}.$$

To go further,  introduce the subgroup
$\mathrm{Aut}_{\mathcal O_K}^+(\mathcal O_K^{(r-1)}\oplus\frak a)$ 
of $\mathrm{Aut}_{\mathcal O_K}
(\mathcal O_K^{(r-1)}\oplus\frak a)$ consisting of these elements 
whose local determinants at real places are all 
positive. 
Then, by checking directly,  we obtain a natural 
identification of quotient spaces between 
$$\mathrm{Aut}_{\mathcal O_K}(\mathcal O_K^{(r-1)}\oplus\frak a)
\backslash\Big(\big(GL(r,\mathbb R)/O(r)\big)^{r_1}\times 
\big(GL(r,\mathbb C)/U(r)\big)^{r_2}\Big)$$ 
and  
$$\mathrm{Aut}_{\mathcal O_K}^+(\mathcal O_K^{(r-1)}\oplus\frak a)
\backslash 
\Big(\big(GL^+(r,\mathbb R)/O^+(r)\big)^{r_1}\times
\big(GL(r,\mathbb C)/U(r)\big)^{r_2}\Big).$$ 

With all this, we are ready  to shift further to the special 
linear group $SL$. 
It is here that Dirichlet's Unit Theorem plays a key role.
As to be expected, the discussion here is a bit involved, for the 
reason that when dealing with metric 
 structures, locally, the genuine realizations are precisely given by 
the following identifications:
$$\begin{aligned}GL(r,\mathbb R)\Big/O(r)\to& GL^+(r,\mathbb R)\Big/SO(r)\\
\to& \Big(SL(r,\mathbb R)\Big/SO(r)\Big)\times \Big(\mathbb R_+^*
\cdot\mathrm{diag}(1,\cdots,1)\Big)\simeq \Big(SL(r,\mathbb R)/SO(r)\Big)
\times \mathbb R_+^*\end{aligned}$$ via
$$\begin{aligned} [A]&\mapsto[A^+]\\
&\mapsto \Big(\frac{1}{\root{r}\of {\det A^+}}A^+,
\mathrm{diag}\big(\root{r}\of{\det A^+},
 \cdots,\root{r}\of{\det A^+}\big)\Big)\mapsto 
\Big(\frac{1}{\root{r}\of {\det A^+}}A^+,\root{r}\of{\det A^+}\Big),
\end{aligned}$$ 
for real places, and 
$$\begin{aligned}GL(r,\mathbb C)/U(r)\to &\Big(SL(r,\mathbb C)
\times \mathbb C\Big)/\Big(SU(r)\times S^1\Big)\\
\to& 
\Big(SL(r,\mathbb C)/SU(r)\Big)\times \Big(\mathbb R_+^*\cdot
\mathrm{diag}(1,\cdots,1)\Big)\simeq \Big(SL(r,\mathbb C)/SU(r)\Big)
\times \mathbb R_+^*\end{aligned}$$ via
$$\begin{aligned} [A]&\mapsto[A]\\
&\mapsto \Big(\frac{1}{\root{r}\of {\det A}}A,\mathrm{diag}\big(
\root{r}\of{\det A}, \cdots,\root{r}\of{\det A}\big)\Big)\to 
\Big(\frac{1}{\root{r}\of {\det A}}A,\root{r}\of{\det A}\Big),\end{aligned}$$ 
for complex places.
Ideally, we want to have corresponding identifications 
for elements in $\mathrm{Aut}_{\mathcal O_K} 
(\mathcal O_K^{(r-1)}\oplus\frak a)$.   However, this cannot be done
in general, due to the fact that,  the $r$-th roots of a unit in $K$ 
lie only in a finite extension of $K$. 
So suitable modifications have to be made. More precisely, we go as follows:
\vskip 0.30cm
Recall that for a unit $\varepsilon \in U_K$,  

\noindent
(a) $\mathrm{diag}(\varepsilon,\cdots,\varepsilon)\in 
\mathrm{Aut}_{\mathcal O_K}(\mathcal O_K^{(r-1)}\oplus\frak a);$ and

\noindent
(b) $\det\mathrm{diag}(\varepsilon,\cdots,\varepsilon)
=\varepsilon^r\in U_F^r:=\{\varepsilon^r:\varepsilon\in U_K\}.$

\noindent
So to begin with, note that to pass from $GL$ to $SL$ over $K$, 
we need to use the intermediate subgroup $GL^+$. Consequently, 
we introduce a subgroup $U_K^+$ of $U_K$ by setting
$$U_K^+:=\{\varepsilon\in U_K:\varepsilon_\sigma>0,
\forall\sigma\ \mathrm{real}\}$$ so as to get
a well-controlled subgroup $U_K^{r,+}:=U_K^+\cap U_K^r$. Indeed, by 
Dirichlet's Unit Theorem, the quotient group 
$U_K^+/(U_K^+\cap U_K^r)$ is  finite. 
 
With this said, next we use $U_K^+\cap U_K^r$ to decomposite the automorphism
group  $\mathrm{Aut}_{\mathcal O_K}(\mathcal O_K^{(r-1)}\oplus\frak a).$ 
Thus, choose elements $$u_1,\cdots,
 u_{\mu(r,F)}\in U_K^+$$ such that 
$\Big\{[u_1],\cdots, [u_{\mu(r,F)}]\Big\}$ gives a complete 
representatives of the finite quotient group
$U_K^+\big/\Big(U_K^+\cap U_K^r\Big)$, where $\mu(r,K)$ denotes
 the cardinality of the group $U_K^+/(U_K^+\cap 
U_K^r)$. Set also  $$SL(\mathcal O_K^{(r-1)}\oplus\frak a)
:=SL(r,K)\cap GL(\mathcal O_K^{(r-1)}\oplus\frak a).$$  
\vskip 0.20cm
\noindent
{\bf{\large Lemma}:} {\it There exist elements $A_1, \ldots, A_{\mu(r,K)}$ 
in $GL^+(\mathcal O_K^{(r-1)}\oplus\frak a)$ such that}

\noindent 
(i) $\det A_i=u_i$ for $i=1,\ldots,\mu(r,K)$; and

\noindent
(ii) $\Big\{A_1,\cdots,A_{\mu(r,K)}\Big\}$ {\it is a completed representatives 
of the  quotient of
$\mathrm{Aut}_{\mathcal O_K}^+\Big(\mathcal O_K^{(r-1)}\oplus\frak a\Big)$ 
modulo} $SL(\mathcal O_K^{(r-1)}
\oplus\frak a)\times \Big(U_K^{r,+}\cdot\mathrm{diag}(1,\cdots,1)\Big)$.
\vskip 0.30cm
\noindent
That is to say, for  automorphism groups, 

\noindent
(a) $\mathrm{Aut}_{\mathcal O_K}^+(\mathcal O_K^{(r-1)}\oplus\frak a)$
is naturally identified with the disjoint union $$\cup_{i=1}^{\mu(r,K)}
A_i\cdot\Bigg(SL(\mathcal O_K^{(r-1)}\oplus\frak a)\times \Big(U_K^{r,+}
\cdot\mathrm{diag}(1,\cdots,1)\Big)\Bigg);$$ and, consequently,

\noindent 
(b) The  $\mathcal O_K$-lattice 
structures $\bold\Lambda(P)$ on the projective $\mathcal O_K$-module 
$P=\mathcal O_K^{(r-1)}\oplus \frak a$ are parametrized by the disjoint union 
 $$\begin{aligned}\cup_{i=1}^{\mu(r,K)}&A_i\backslash 
\Bigg(
\bigg(SL(\mathcal O_K^{(r-1)}\oplus\frak a)\backslash 
 \Big(\Big(SL(r,\mathbb R)/SO(r)\Big)^{r_1}\times 
\Big(SL(r,\mathbb C)/SU(r)\Big)^{r_2}\Big)
\bigg)\\
&\qquad\times
 \Big(\Big|U_{K}^r\cap U_K^+\Big|\backslash 
\Big(\mathbb R_+^*\Big)^{r_1+r_2}\Big)\Bigg).\end{aligned}$$

\noindent
{\it Proof}. It is a direct consequence of the follows:

\noindent
(1) For all $\varepsilon\in U_{K}^+$,
$\mathrm{diag}(\varepsilon,\cdots, \varepsilon)\in 
\mathrm{Aut}_{\mathcal O_K}^+
(\mathcal O_K^{(r-1)}\oplus\frak a)$ and its determinant 
belongs to $U_K^+\cap U_K^r$;
 
\noindent
(2) For $A\in \mathrm{Aut}_{\mathcal O_K}^+
(\mathcal O_K^{(r-1)}\oplus\frak a)$, by definition, 
$\det A\in U_{K}^+$. completed.
\vskip 0.30cm
Therefore, to understand the space of 
$\mathcal O_K$-lattice structures, beyond the spaces
$SL(r,\mathbb R)/SO(r)$ and $SL(r,\mathbb C)/SU(r)$, 
we further need to study 
 
\noindent
(i) the quotient space 
$\big|U_{K}^r\cap U_K^+\big|\backslash 
 \Big(\mathbb R_+^*\Big)^{r_1+r_2}$, which is more or less standard 
(see e.g. [Neu] or [We4]); and more importantly,

\noindent
(ii) the (modular) space $$SL(\mathcal O_K^{(r-1)}
\oplus\frak a)\backslash \Bigg(\Big(SL(r,\mathbb R)/SO(r)\Big)^{r_1}
\times \Big(SL(r,\mathbb C)/SU(r)\Big)^{r_2}\Bigg).$$ 
\vskip 0.30cm
Now denote by $\widetilde{\mathcal M}_{K,r}(\frak a)$ the moduli space of rank 
$r$ semi-stable $\mathcal O_K$-lattices
with underlying projective module $\mathcal O_K^{(r-1)}\oplus\frak a$. 
For our own convenience, for a set $X$ of (isometry classes of)
lattices, we use the notation $X_{\mathrm{ss}}$ 
to denote the subset of $X$ consisting of lattices
which are semi-stable. 
As such, then we have the following variation of the previous lemma.
\vskip 0.30cm
\noindent
{\bf{\Large Proposition}}. {\it There is a natural identification
between the moduli space $\widetilde{\mathcal M}_{K,r}(\frak a)$ of rank 
$r$ semi-stable $\mathcal O_K$-lattices
on the projective module $\mathcal O_K^{(r-1)}\oplus\frak a$
and the disjoint union of (the ss part of) the quotient spaces}
$$\begin{aligned}\cup_{i=1}^{\mu(r,K)}&A_i\backslash 
\Bigg(
\bigg(SL(\mathcal O_K^{(r-1)}\oplus\frak a)\backslash 
 \Big(\Big(SL(r,\mathbb R)/SO(r)\Big)^{r_1}\times 
\Big(SL(r,\mathbb C)/SU(r)\Big)^{r_2}\Big)
\bigg)_{\mathrm{ss}}\\
&\qquad\times
 \Big(|U_{K}^r\cap U_K^+|\backslash 
(\mathbb R_+^*)^{r_1+r_2}\Big)\Bigg).\end{aligned}$$

\chapter{Geometric Truncation and Analytic Truncation}
\section{Geometric Truncation: Revised}
\subsection{Slopes, Canonical Filtrations}

Following Lafforgue [Laf], we call an abelian category ${\mathcal A}$ 
together with two additive
morphisms $$\mathrm {rk}:{\mathcal A}\to {\mathbb N},
\qquad \mathrm{deg}:{\mathcal A}\to {\mathbb R}$$ a {\it category with slope 
structure}. In particular, for non-zero $A\in {\mathcal A}$, 
\vskip 0.30cm
\noindent
({\bf 1}) define the {\it slope} of $A$ by 
$$\mu(A):={{\mathrm{deg}(A)}\over {\mathrm{rk}A}};$$
\vskip 0.30cm 
\noindent
({\bf 2}) If $0=A_0\subset A_1\subset\cdots\subset A_l=A$ is a filtration of $A$ in 
${\mathcal A}$ with $\mathrm {rk}(A_0)<\mathrm {rk}(A_1)<\cdots<
\mathrm {rk}(A_l)$, then define the {\it associated polygon} 
to be the continuous 
function $[0,\mathrm {rk}A]\to {\mathbb R}$ such that

\noindent
(i) its values at 0 and $\mathrm {rk}(A)$ are 0;

\noindent
(ii)  it is affine on the intervals 
$[\mathrm {rk}(A_{i-1}),\mathrm {rk}(A_i)]$ 
with slope $\mu(A_i/A_{i-1})-\mu(A)$  for all $1\leq i\leq l$;
\vskip 0.30cm
\noindent
({\bf 3}) If $\frak a$ is a collection of subobjects of $A$ in ${\mathcal A}$, then 
$\frak a$ is said to be {\it nice} if

\noindent
(i) $\frak a$ is stable under intersection and finite summation;

\noindent
(ii) $\frak a$ is Noetherian in the sense that every increasing chain 
of elements
in $\frak a$ has a maximal element in $\frak a$;

\noindent
(iii) if $A_1\in\frak a$ then $A_1\not=0$ if and only if 
$\mathrm {rk}(A_1)\not=0$; and

\noindent
(iv) for $A_1,A_2\in\frak a$ with $\mathrm {rk}(A_1)=\mathrm {rk}(A_2)$. 
Then
$A_1\subset A_2$ is proper implies that $\mathrm {deg}(A_1)<
\mathrm {deg}(A_2);$
\vskip 0.30cm
\noindent
({\bf 4}) For any nice $\frak a$, set 
$$\mu^+(A):=\mathrm {sup}\,\Big\{\mu(A_1):A_1\in\frak a,
\mathrm {rk}(A_1)\geq 1\Big\},\ 
\mu^-(A):=\mathrm {inf}\,\Big\{\mu(A/A_1):A_1\in\frak a,
\mathrm {rk}(A_1)< \mathrm {rk}(A) \Big\}.$$ 
Then we say $(A,\frak a)$ is {\it semi-stable} if $\mu^+(A)=\mu(A)=\mu^-(A)$. 
Moreover if $\mathrm {rk}(A)=0$, 
set also $\mu^+(A)=-\infty$ and $\mu^-(A)=+\infty$.
\vskip 0.30cm
\noindent
{\bf{\large Proposition 1}.} ([Laf]) {\it Let ${\mathcal A}$ be a 
category with slope, $A$ an object in ${\mathcal A}$ and $\frak a$ a nice 
family of subobjects of $A$ in ${\mathcal A}$. Then}

\noindent
(1) ({\bf Canonical Filtration}) {\it $A$ admits a unique filtration $0=
\overline A_0\subset \overline A_1\subset\cdots\subset\overline A_l=A$
with elements in $\frak a$  such that}

\noindent
(i) {\it $\overline A_i,0\leq i\leq k$ are maximal in $\frak a$;}

\noindent 
(ii) {\it $\overline A_i/\overline A_{i-1}$ are semi-stable; and}

\noindent 
(iii) $\mu(\overline A_1/\overline A_{0})>
\mu(\overline A_2/\overline A_{1}>\cdots>
\mu(\overline A_k/\overline A_{k-1})$;

\noindent
(2) ({\bf Boundness}) {\it All polygons of filtrations of $A$ with elements 
in $\frak a$ are bounded from above by $\overline p$, where 
 $\overline p:=\overline p^A$ is the associated polygon for the canonical 
filtration in (1);}

\noindent
(3) {\it For any $A_1\in \frak a, \mathrm {rk}(A_1)\geq 1$ implies
$\mu(A_1)\leq \mu(A)+\frac{\overline p(\mathrm {rk}(A_1))}
{\mathrm {rk}(A_1)};$}

\noindent
(4) {\it The polygon $\overline p$ is convex with maximal slope 
$\mu^+(A)-\mu(A)$ 
and minimal slope $\mu^-(A)-\mu(A)$;}

\noindent
(5) {\it If $(A',\frak a')$ is another pair, and $u:A\to A'$ is a homomorphism 
such that $\mathrm {Ker}(u)\in\frak a$ and $\mathrm{Im}(u)\in\frak a'$. Then
$\mu^-(A)\geq\mu^+(A')$ implies that $u=0$.}
\vskip 0.30cm
\noindent
{\it Proof.} This results from a Harder-Narasimhan type filtration 
consideration. A detailed proof may be found at pp. 87-88 in [Laf].
\vskip 0.30cm
As an example, we have the following  

\noindent
{\bf {\large Proposition 2}.} {\it  Let $F$ be a number field. Then}

\noindent
(1) {\it the abelian category of hermitian vector sheaves on 
$\mathrm {Spec}\,{\mathcal O}_F$ together with the natural rank and the 
Arakelov degree is a category with slopes;}

\noindent
(2) {\it For any hermitian vector sheaf $(E,\rho)$, $\frak a$ consisting 
of pairs $(E_1,\rho_1)$ with $E_1$ sub vector sheaves of $E$ and 
$\rho_1$ the restrictions of $\rho$, forms a nice family.}
\vskip 0.30cm
\noindent
{\it Proof.} (1) is obvious, while (2) is a direct consequence of 
the following standard facts: 

\noindent
(i) For a fixed $(E,\rho)$,
$\Big\{\mathrm {deg}(E_1,\rho_1):(E_1,\rho_1)\in \frak a\Big\}$ is discrete 
subset of ${\mathbb R}$; and 

\noindent
(ii) for any two sublattices $\Lambda_1,\, \Lambda_2$  
of $\Lambda$,
 $$\mathrm{Vol}\Big(\Lambda_1/(\Lambda_1\cap\Lambda_2)\Big)\geq
\mathrm{Vol}\Big((\Lambda_1+\Lambda_2)/\Lambda_2\Big).$$
\vskip 0.30cm
Thus in particular we  get the canonical filtration  of Harder-Narasimhan
 type for 
hermitian vector sheaves over $\mathrm {Spec}\,{\mathcal O}_F$.
Recall that hermitian vector sheaves over $\mathrm {Spec}\,{\mathcal O}_F$
are ${\mathcal O}_F$-lattices in $({\mathbb R}^{r_1}\times 
{\mathbb C}^{r_2})^{r=\mathrm{rk}(E)}$ in the language of Arakelov theory:
Say, corresponding $\mathcal O_F$-lattices are 
induced from their $H^0$ via the natural embedding
$F^r\hookrightarrow ({\mathbb R}^{r_1}\times {\mathbb C}^{r_2})^r$ where $r_1$ 
(resp. $r_2$) denotes the real (resp. complex) embeddings of $F$.

\subsection{Canonical Polygons and Geometric Truncation}

Let $F$ be a number field with $\mathcal O_F$ the ring of integers.
Denote by $\mathbb A$ the ring of adeles and $\mathbb A_f$ the ring of 
finite adeles. If $E$ is a vector sheaf of rank $r$ 
over $X:=\mathrm {Spec}({\mathcal O}_F)$,
 i.e, a locally free ${\mathcal O}_F$-sheaf of rank $r$, denote by $E_F$ the 
fiber 
of $E$ at the generic point $\mathrm {Spec}(F)\hookrightarrow 
\mathrm {Spec}({\mathcal O}_F)$ of $X$ 
and for each finite place $v\in S_{\mathrm{fin}}$, set $E_{{\mathcal O}_v}:=
H^0(\mathrm {Spec}{\mathcal O}_{F_v},E)$. Then $E_F$ is an $F$-vector space of 
dimension $r$, and $E_{{\mathcal O}_v}$ is a free ${\mathcal O}_v$-module of 
rank $r$. In particular, we have a canonical isomorphism:
$$\mathrm{can}_v:F_v\otimes_{{\mathcal O}_v}E_{{\mathcal O}_v}\simeq F_v
\otimes_F E_F.$$
 
Now let $E$ be equipped with a basis 
$\alpha_F:F^r\simeq E_F$ of its generic fiber and a basis 
$\alpha_{{\mathcal O}_v}:{\mathcal O}_v^r\simeq E_{{\mathcal O}_v}$ for any 
$v\in S_{\mathrm{fin}}$. Then, elements 
$g_v:=(F_v\otimes_F\alpha_F)^{-1}\circ\mathrm{can}_v
\circ (F_v\otimes_{{\mathcal O}_v}\alpha_{{\mathcal O}_v})\in GL_r(F_v)$ for 
all $v\in S_{\mathrm{fin}}$ define an element 
$g_{\mathbb A}:=(g_v)_{v\in S_{\mathrm{fin}}}$ of 
$GL_r({\mathbb A}_f)$, since for almost every $v$ we have $g_v\in 
GL_r({\mathcal O}_v)$. By this construction, we obtain a bijection from 
the set of isomorphism classes of triples $(E;\alpha_F;
(\alpha_{{\mathcal O}_v})_{v\in S_{\mathrm{fin}}})$ as above onto 
$GL_r({\mathbb A}_f)$. 
Moreover, if $r\in GL_r(F), k\in GL_r({\mathcal O}_F)$ and if this bijection
 maps the triple 
$(E;\alpha_F;(\alpha_{{\mathcal O}_v})_{v\in S_{\mathrm{fin}}})$ onto 
$g_{\mathbb A}$, the same map maps the triple $(E;\alpha_F\circ r^{-1};
(\alpha_{{\mathcal O}_v}\circ k_v)_{v\in S_{\mathrm{fin}}})$ 
onto $rg_{\mathbb A}k$. 
Therefore the above bijection induces a bijection between the set of 
isomorphism classes of vector sheaves of rank $r$ on 
$\mathrm {Spec}{\mathcal O}_F$ and the double coset space
$GL_r(F)\backslash GL_r({\mathbb A}_f)/GL_r({\mathcal O}_F)$.

More generally, let $r=r_1+\cdots+r_s$ be a partition $I=(r_1,\cdots,r_s)$ of 
$r$ and let $P_I$ be the corresponding standard parabolic subgroup of $GL_r$.
(See next section for details.) 
Then we have a natural bijection from the set of isomorphism classes of 
triple $\Big(E_*;\alpha_{*,F}:
(\alpha_{*,{\mathcal O_v}})_{v\in S_{\mathrm{fin}}}\Big)$ onto 
$P_I({\mathbb A}_f)$, where $E_*:=\Big((0)=E_0\subset E_1\subset\cdots\subset
 E_s\Big)$ is a filtration of vector sheaves of rank $\Big(r_1,r_1+r_2,\cdots, 
r_1+r_2+\cdots+r_s=r\Big)$ over $X$, (i.e, each $E_j$ is a 
vector sheaf of rank 
$r_1+r_2+\cdots+r_j$ over $X$ and each quotient $E_j/E_{j-1}$ is torsion 
free,) which is equipped with an isomorphism of filtrations of $F$-vector 
spaces
$$\alpha_{*,F}:\Big((0)=F_0\subset F^{r_1}\subset\cdots\subset F^{r_1+r_2+
\cdots+r_s=r}\Big)\simeq (E_*)_F,$$ and with an isomorphism of filtrations 
of free ${\mathcal O}_v$-modules
$$\alpha_{*,{\mathcal O}_v}:\Big((0)\subset {\mathcal O}_v^{r_1}
\subset\cdots\subset 
{\mathcal O}_v^{r_1+r_2+\cdots+r_s=r}\Big)\simeq (E_*)_{{\mathcal O}_v},$$ 
for every 
$v\in S_{\mathrm{fin}}$. Moreover this bijection induces a bijection between the set of 
isomorphism classes of the filtrations of vector sheaves of rank 
$(r_1,r_1+r_2,\cdots, r_1+r_2+\cdots+r_s=r)$ over $X$ and the double coset 
space $P_I(F)\backslash P_I({\mathbb A}_f)/P_I({\mathcal O}_F)$. The natural 
embedding $P_I({\mathbb A}_f)\hookrightarrow P_I({\mathbb A})$ (resp. the 
canonical projection $P_I({\mathbb A}_f)\to M_I({\mathbb A}_f)\to GL_{r_j}
({\mathbb A}_f)$ for $j=1,\cdots,s$, where $M_I$ denotes the standard Levi 
of $P_I$) admits the modular interpretation
$$\Big(E_*;\alpha_{*,F}:(\alpha_{*,{\mathcal O_v}})_{v\in S_{\mathrm{fin}}}
\Big)\mapsto 
\Big(E_s;\alpha_{s,F}:(\alpha_{s,{\mathcal O_v}})_{v\in S_{\mathrm{fin}}}
\Big)$$ (resp. 
$$\Big(E_*;\alpha_{*,F}:(\alpha_{*,{\mathcal O_v}})_{v\in S_{\mathrm{fin}}}
\Big)\mapsto 
\Big(\mathrm{gr}_j(E_*);\mathrm{gr}_j(\alpha_{*,F}), \mathrm{gr}_j
(\alpha_{*,{\mathcal O_v}})_{v\in S_{\mathrm{fin}}}\Big),$$ where 
$\mathrm{gr}_j(E_*):=E_j/E_{j-1}$, and 
$$\mathrm{gr}_j(\alpha_{*,F}):F^{r_j}\simeq \mathrm{gr}_j(E_*)_F,\quad 
\mathrm{gr}_j(\alpha_{*,{\mathcal O_v}}):{\mathcal O_v}^{r_j}\simeq 
\mathrm{gr}_j(E_*)_{\mathcal O_v},\qquad\forall v\in S_{\mathrm{fin}}$$ 
are induced by $\alpha_{*,F}$ and $\alpha_{*,{\mathcal O_v}}$.)
\vskip 0.30cm
Furtherover, any $g=\Big(g_f;g_\infty\Big)\in GL_r({\mathbb A}_f)\times 
GL_r({\mathbb A}_\infty)=GL_r({\mathbb A})$ gives first a rank $r$ 
vector sheaf $E_g:=E_{g_f}$ on $\mathrm{Spec}({\mathcal O}_F)$, which via the 
embedding $E_F\hookrightarrow \Big({\mathbb R}^{r_1}\times 
{\mathbb C}^{r_2}\Big)^r$ 
gives a discrete subgroup, a  free  ${\mathcal O}_F$-module of rank $r$. In 
such a way, 
$g_\infty=(g_\sigma)$ then induces a natural metric $\rho_g:=\rho_{g_\infty}$  
on $E_g$ by twisting 
the standard one on $({\mathbb R}^{r_1}\times {\mathbb C}^{r_2})^r$ via the 
linear transformation induced from $g_\infty$. (That is to say, for $\sigma$ 
real, the metric is given by $g\cdot g^t$, while for $\sigma$ 
complex, the metric is given by $g\cdot\overline g^t$.) 
 Consequently, 
see e.g, [L1,3], $$\mathrm{deg}\Big(E_g,\rho_g\Big)
=-\log\Big(N(\mathrm{det}g)\Big)$$ 
with $N:GL_1({\mathbb A}_F)={\mathbb I}_F\to {\mathbb R}_{>0}$ the standard 
norm map of the idelic group of $F$.
\vskip 0.30cm
With this, for $g=(g_f;g_\infty)\in GL_r({\mathbb A})$ and a parabolic 
subgroup 
$Q$ of $GL_r$, denote by $E_*^{g;Q}$ the filtration of the vector sheaf 
$E_{g_f}$ induced by the parabolic subgroup $Q$. Then we have a filtration 
of hermitian vector sheaves $(E_*^{g;Q},\rho_*^{g;Q})$ with the hermitian 
metrics $\rho_j^{g;Q}$
on $E_j^{g;Q}$ obtained via the restrictions of $\rho_{g_\infty}$.
\vskip 0.30cm
Now introduce the associated polygon $p_Q^g:[0,r]\to {\mathbb R}$ by the 
following 3 conditions:

\noindent
(i) $p_Q^g(0)=p_Q^g(r)=0$;

\noindent
(ii) $p_Q^g$ is affine on the interval 
$[\mathrm{rk}E_{i-1}^{g;Q},\mathrm{rk}E_{i}^{g;Q}]$; and 

\noindent
(iii) for all indices $i$, 
$$p_Q^g(\mathrm{rk}E_{i}^{g;Q})=\mathrm{deg}(E_{i}^{g;Q},\rho_{i}^{g;Q})-
\mathrm{rk}E_{i}^{g;Q}\cdot \frac{\mathrm{deg}(E_{g},\rho_{g})}{r}.$$
Then by Prop. 2 of 2.1.1, there is a unique convex polygon $\overline p^g$ 
which bounds all $p_Q^g$ from above for all parabolic subgroups $Q$ for 
$GL_r$. 
Moreover there exists a parabolic subgroup $\overline Q^g$ such that 
$p_{{\overline Q}^g}^g=\overline p^g$.

With this, we are now ready to state the following foundamental
\vskip 0.30cm
\noindent
{\bf {\Large Theorem}$'$.} {\it For any fixed polygon 
$p:[0,r]\to {\mathbb R}$ and any $d\in {\mathbb R}$, the subset
$$\Big\{g\in GL_r(F)\backslash GL_r({\mathbb A}):
\mathrm{deg}\,g=d,\overline p^g\leq p\Big\}$$ 
is compact.}
\vskip 0.30cm
\noindent
{\it Proof.} This is a restatement of the classical reduction theory, 
i.e., Thm. 1.2.3. 
Indeed, it consists of two parts: In terms of ${\mathcal O}_F$-lattices, 
by fixing the degree, we get a fixed volume for the free 
${\mathcal O}_F$-lattice of rank $r$ corresponding to $(E_g,\rho_g)$, 
by the Arakelov Riemann-Roch formula. Thus the condition $\overline p^g\leq p$ gives an upper bound for
the volumes of all the sublattices of $(E_g,\rho_g)$ and  a lower bound
of the Minkowski successive minimums.

On the other hand, the reduction theory in terms of adelic language also 
tells us that the fiber of the natural map from $GL_r(F)\backslash 
GL_r({\mathbb A})$ to isomorphism classes of ${\mathcal O}_F$-lattices are all 
compact. (This is in fact generally true for all reductive groups, a 
result due to Borel [Bo1,2]. However, our case
where only $GL_r$ is involved is rather obvious: As seen in Chapter 1,
essentially, the fibers are the compact subgroup 
$GL_r({\mathcal O}_F)\times SO_r({\mathbb R})^{r_1}\times 
SU_r({\mathbb C})^{r_2}$.)
This completes the proof.

In particular, as also being seen in Ch. 1,
after fixing the volume of lattices, we may 
use special linear groups instead of general linear groups.
\vskip 0.30cm 
Similarly yet more generally, for a fixed parabolic subgroup $P$ of $GL_r$ 
and $g\in GL_r({\mathbb A})$, there is a unique maximal element 
$\overline p_P^g$ 
among all $p_Q^g$, where $Q$ runs over all parabolic subgroups of $GL_r$ 
which are contaiend in $P$.  And we have
\vskip 0.30cm
\noindent
{\bf {\Large Theorem}}$'$. 
{\it For any fixed polygon $p:[0,r]\to {\mathbb R}$, $d\in {\mathbb R}$ 
and  any standard parabolic subgroup $P$ of $GL_r$, the subset
$$\Big\{g\in GL_r(F)\backslash GL_r({\mathbb A}):\mathrm{deg}\,g=d,
\overline p^g_P\leq p, p^g_P\geq -p\Big\}$$ is compact.}

\section{Partial Algebraic Truncation}

For our limited purpose of exposing a relation between
geometric truncation and analytic truncation, it is enough to deal with 
special linear groups. However, there is not that much difference in 
working with general reductive groups for our later purposes. 
Hence we start with a general framework.

\subsection{Parabolic Subgroups}

Let $F$ be a number field with $\mathbb A=\mathbb A_F$ the ring of adeles. 
Let $G$ be a connected reductive group defined over $F$. Recall that a 
subgroup $P$ of $G$ is called {\it parabolic} if $G/P$ is a complete algebraic 
variety. Fix a minimal $F$-parabolic subgroup $P_0$ of $G$ with its unipotent 
radical $N_0=N_{P_0}$ and fix a $F$-Levi subgroup $M_0=M_{P_0}$ of $P_0$ so as 
to have a Levi decomposition $P_0=M_0N_0$. An $F$-parabolic subgroup $P$ is 
called {\it standard} if it cotains $P_0$.  For such parabolic subgroups $P$, 
there exists a unique  Levi subgroup $M=M_P$ containing $M_0$ which we 
call the 
{\it standard Levi subgroup} of $P$. Let $N=N_P$ be the unipotent radical.
 Let us agree to use the term parabolic subgroups and Levi subgroups to 
denote standard $F$-parabolic subgroups and standard Levi subgroups 
repectively, unless otherwise is stated. 

Let $P$ be a parabolic subgroup of $G$. Write $T_P$ for the maximal split 
torus in the center of $M_P$ and $T_P'$ for the maximal quotient split 
torus of $M_P$. Set $\tilde{\frak a}_P:=X_*(T_P)\otimes \mathbb R$ and 
denote its real dimension
by $d(P)$, where $X_*(T)$ is the lattice of 1-parameter subgroups in the torus 
$T$. Then it is known that $\tilde{\frak a}_P=X_*(T_P')\otimes \mathbb R$ as 
well. The two descriptions of $\tilde {\frak a}_P$ 
show that if $Q\subset P$ is a 
parabolic subgroup, then there is a canonical injection 
$\tilde {\frak a}_{P}\hookrightarrow \tilde{\frak a}_{Q}$ and a 
natural surjection  
$\tilde {\frak a}_{Q}\twoheadrightarrow \tilde{\frak a}_{P}$.  
We thus obtain a canonical decomposition, for a certain subspace 
$\tilde{\frak a}_Q^{P}$ of $\tilde {\frak a}_{Q}$, 
$$\tilde {\frak a}_{Q}=\tilde{\frak a}_Q^{P}\oplus
\tilde {\frak a}_{P}.$$ In particular, $\tilde{\frak a}_G$ is a summand of 
$\tilde {\frak a}=\tilde{\frak a}_P$ for all $P$. Set 
$\frak a_{P}:= \tilde{\frak a}_{P}/\tilde {\frak a}_{G}$ and 
$\frak a_{Q}^P:= \tilde{\frak a}_Q^{P}/ \tilde{\frak a}_{G}.$ Then we have 
$$\frak a_{Q}=\frak a_Q^{P}\oplus\frak a_{P}$$
and $\frak a_P$ is canonically identified as a subspace of 
$\frak a_Q$. Set $\frak a_0:=\frak a_{P_0}$ and 
$\frak a_0^P=\frak a_{P_0}^P$ then we also have 
$\frak a_0=\frak a_0^P\oplus \frak a_P$ for all $P$. 

(For the inexperience reader, please refer the example on SL below 
to have a good understanding of the general discussion in this and the 
following a few subsections.)

\subsection{Logarithmic Map}

For a real vector space $V$, write $V^*$ its dual space over $\mathbb R$. Then
dually we have the spaces such as $\frak a_0^*, \frak a_P^*, 
\Big(\frak a_0^P\Big)^*$ and hence the decompositions 
$$\frak a_0^*=\Big(\frak a_0^Q\Big)^*\oplus
\Big(\frak a_Q^P\Big)^*\oplus\frak a_P^*.$$ In particular, 
$\frak a_P^*=X(M_P)\otimes \mathbb R$ with $X(M_P)$ the group 
$\mathrm {Hom}_F\Big(M_P,GL(1)\Big)$ i.e., collection of characters on 
$M_P$. It is known that
$\frak a_P^*=X(A_P)\otimes \mathbb R$ where $A_P$ denotes the split 
component of the center of $M_P$. Clearly, if $Q\subset P$, then 
$M_Q\subset M_P$ while $A_P\subset A_Q$. Thus via restriction,
the above two expressions of $\frak a_P^*$ also naturally induce
an injection $\frak a_P^*\hookrightarrow \frak a_Q^*$ and a sujection 
$\frak a_Q^*\twoheadrightarrow
\frak a_P^*$, compactible with the decomposition $\frak a_Q^*=
\Big(\frak a_Q^P\Big)^*\oplus\frak a_P^*.$

Every $\chi=\sum s_i\chi_i$ in $\frak a_{P,\mathbb C}^*:=
\frak a_P^*\otimes\mathbb C$ determines a morphism 
$P(\mathbb A)\to \mathbb C^*$ by $p\mapsto p^\chi:=\prod|\chi_i(p)|^{s_i}.$ 
Consequently, we have a natural logarithmic map $H_P:P(\mathbb A)\to\frak a_P$ 
defined by $$\langle H_P(p),\chi\rangle=p^\chi,\qquad 
\forall \chi\in\frak a_P^*.$$ The 
kernel of $H_P$ is denoted by $P(\mathbb A)^1$ and we set 
$M_P(\mathbb A)^1:=P(\mathbb A)^1\cap M_P(\mathbb A)$.
 
Let also $A_+$ be the set of $a\in A_P(\mathbb A)$ such that 

\noindent
(1) $a_v=1$ for all finite places $v$ of $F$; and 

\noindent
(2) $\chi(a_\sigma)$ 
is a positive number independent of infinite places $\sigma$ of $F$ 
for all $\chi\in X(M_P)$. 

\noindent
Then $M(\mathbb A)=A_+\cdot M(\mathbb A)^1$. 

\subsection{Roots, Coroots, Weights and Coweights}

We now define the standard bases of the above spaces and their duals. Let 
$\Delta_0$ and $\widehat\Delta_0$ be the subsets of simple roots and simple 
weights in $\frak a_0^*$ respectively. (Recall that elements of 
$\widehat\Delta_0$ are non-negative linear combinations of elements in 
$\Delta_0$.) Write $\Delta_0^\vee$ (resp. $\widehat\Delta^\vee_0$) for 
the basis of $\frak a_0$ dual to $\widehat\Delta_0$ (resp. $\Delta_0$). 
Thus being the dual of the collection of simple weights (resp. the 
collective of simple roots), $\Delta_0^\vee$ (resp. $\widehat\Delta_0^\vee$) 
is the set of coroots (resp. coweights).

For every $P$, let $\Delta_P\subset\frak a_0^*$ be the set of non-trivial 
{\it restrictions} of elements of $\Delta_0$ to $\frak a_P$. Denote the dual 
basis of $\Delta_P$ by $\widehat\Delta_P^\vee$.
For each $\alpha\in\Delta_P$, let $\alpha^\vee$ be the projection of 
$\beta^\vee$ to $\frak a_P$, where $\beta$ is the root in $\Delta_0$ 
whose restriction to $\frak a_P$ is $\alpha$. Set 
$\Delta_P^\vee:=\Big\{\alpha^\vee:\alpha\in\Delta_P\Big\}$, and define the 
dual basis of $\Delta_P^\vee$ by $\widehat\Delta_P$. 

More generally, if $Q\subset P$, write $\Delta_Q^P$ to denote the {\it subset} 
$\alpha\in\Delta_Q$ appearing in the action of $T_Q$ in the unipotent radical 
of $Q\cap M_P$. (Indeed, $M_P\cap Q$ is a parabolic subgroup of $M_P$ with 
nilpotent radical $N_Q^P:=N_Q\cap M_P$. Thus $\Delta_Q^P$ is simply the set 
of roots of the parabolic subgroup $(M_P\cap Q,A_Q)$. And one checks that the 
map $P\mapsto \Delta_Q^P$ gives a natural bijection between parabolic 
subgroup $P$ containin $Q$ and subsets of $\Delta_Q$.)
Then $\frak a_P$ is the subspace of $\frak a_Q$ annihilated by 
$\Delta_Q^P$. Denote by $(\widehat\Delta^\vee)_Q^P$ the dual of $\Delta_Q^P$. 
Let $(\Delta_Q^P)^\vee:=\Big\{\alpha^\vee:\alpha\in \Delta_Q^P\Big\}$ 
and denote by $\widehat\Delta_Q^P$ the dual of $(\Delta_Q^P)^\vee$.
\subsection{Positive Cone and Positive Chamber}

Let $Q\subset P$ be two parabolic subgroups of $G$.
We extend the linear functionals in $\Delta_Q^P$ and $\widehat\Delta_Q^P$ to
 elements of the dual space $\frak a_0^*$ by means of the canonical 
projection from $\frak a_0$ to $\frak a_Q^P$ given by the decomposition 
$\frak a_0=\frak a_0^Q\oplus\frak a_Q^P\oplus\frak a_P$. Let $\tau_Q^P$ 
be the characteristic function of the {\it positive chamber} 
$$\Big\{H\in\frak a_0:\langle\alpha,H\rangle>0\ \mathrm{for\ all}\ 
\alpha\in\Delta_Q^P\Big\}=\frak a_0^Q\oplus
\Big\{H\in\frak a_Q^P:\langle\alpha,H\rangle>0\ \mathrm{for\ all}\ 
\alpha\in\Delta_Q^P\Big\}
\oplus\frak a_P$$ and let $\widehat\tau_Q^P$ be the characteristic 
function of the positive cone $$\Big\{H\in\frak a_0:\langle\varpi,H\rangle>0\ 
\mathrm{for\ all}\ \varpi\in\widehat\Delta_Q^P\Big\}=\frak a_0^Q\oplus
\Big\{H\in\frak a_Q^P:\langle\varpi,H\rangle>0\ 
\mathrm{for\ all}\ \varpi\in\widehat\Delta_Q^P\Big\}
\oplus\frak a_P.$$ Note that elements in 
$\widehat\Delta_Q^P$ are non-negative linear combinations of elements in 
$\Delta_Q^P$, we have $$
\widehat\tau_Q^P\geq \tau_Q^P.$$ In particular, we have the following
\vskip 0.30cm
\noindent
{\bf {\large Langlands' Combinatorial Lemma}}. 
(Arthur) {\it If $Q\subset P$ 
are parabolic subgroups, then for all $H\in \frak a_0$, we have
$$\sum_{R:\,Q\subset R\subset P}(-1)^{d(R)-d(P)}\tau_Q^R(H)\cdot 
\widehat\tau_R^P(H)=\delta_{QP}$$ and
$$\sum_{R:\,Q\subset R\subset P}(-1)^{d(Q)-d(R)}\widehat \tau_Q^R(H)\cdot 
\tau_R^P(H)=\delta_{QP}.$$}

\subsection{Example with $SL_r$}

Set $G=SL_r(\mathbb R)$. Let $P_0=B$ be the Borel subgroup of $G$ 
consisting of upper triangle matrices, and $M_0$ be the Levi component 
consisting of diagonal matrices. One checks easily that the unipotent radical 
$N_0$ consists of upper triangle matrices whose diagonal entries are all 
equal to 1 and that (the positive part $A_{0+}$ of) the split component $A_0$ 
consists of diagonal matrices $\mathrm{diag}\Big(a_1,a_2,\cdots,a_r\Big)$ with 
$a_1,a_2,\cdots,a_r\in\mathbb R_{>0}$ and $\prod_{i=1}^r a_i=1$.
Hence $\frak a_0=\Big\{H:=(H_1,H_2,\cdots, H_r)\in\mathbb R^r:
\sum_{i=1}^rH_i=0\Big\}$ and the natural logarithmic map $H_{P_0}:P_0\to 
\frak a_0$ is given by
$g\mapsto H_{P_0}\Big(a(g)\Big)=\log a_1+\log a_2+\cdots+\log a_r$,
where we have used the Iwasawa decomposition $g=nma(g)k$ with $n\in N_0$, 
$m\in M_0^1$, $a(g)=\mathrm{diag}\Big(a_1,
a_2,\cdots,a_r\Big)\in A_{0+}$ and $k\in K=SO(r)$.

Set $e_i\in\frak a_0^*$ such that $e_i(H)=H_i$, then 
$$\Delta_0=\Big\{e_1-e_2,e_2-e_3,\cdots, e_{r-1}-e_r\Big\}.$$ We identify 
$\Delta_0$ with the set $\Delta:=\{1,2,\cdots, r-1\}$ by 
the map $\alpha_i:=e_i-e_{i+1}\mapsto i, i=1,2,\cdots r-1.$

It is well-known that there is a natural bijection between 
 partitions of $r$ 
and  parabolic subgroups of $G$. More precisely, 
the partition $I:=\Big(d_1,d_2,\cdots, d_n\Big)$ of $r$ (so that 
$r=d_1+d_2+\cdots+d_n$ and $d_i\in\mathbb Z_{>0}$) corresponds to the 
parabolic subgroup $P_I$ consisting of upper triangle blocked matrices 
whose diagonal blocked entries are matrices of sizes $d_1,d_2,\cdots,d_n$. 
Then clearly, the unipotent radical $N_I$   
consists of upper triangle blocked matrices whose diagonal blocked entries 
are identity matrices of sizes $d_1,d_2,\cdots,d_n$ , the Levi component $M_I$ 
consists of diagonal blocked matrices whose blocked diagonal entries are
 of sizes $d_1,d_2,\cdots d_r$, and the 
corresponding $A_{I+}$ consits of 
matrices of the form
$\mathrm{diag}\Big(a_1I_{d_1},a_2I_{d_2},\cdots,a_nI_{d_n}\Big)$ 
where $I_{d_i}$ 
denotes the identity matrix of size $d_i$, $a_i>0$ for $i=1,2,\cdots,n$, 
and $\prod_{i=1}^na_i^{d_i}=1$. 
As such, $\frak a_I$ is a subspace of $\frak a_0$ defined by
$$\begin{aligned}\frak a_I:=\Big\{(H_1,H_2,\cdots,H_r)\in\frak a_0:
&H_1=H_2=\cdots=H_{r_1},\\
&H_{r_1+1}=H_{r_1+2}=\cdots=H_{r_2},\\
&\cdots, H_{r_{n-1}+1}=
H_{r_{n-1}+2}=\cdots=H_{r_n}\Big\}\end{aligned}$$ where we set 
$r_i=d_1+d_2+\cdots+d_i$ for $i=1,2,\cdots,n$.
Moreover, with the help of Iwasawa decomposition,
the corresponding logarithmic map may be simply seen to be the map 
$H_I:A_{I+}\to\frak a_I$ defined by $$\mathrm{diag}\Big(a_1I_{d_1},
a_2I_{d_2},\cdots,a_nI_{d_n}\Big)
\mapsto \Big((\log a_1)^{(d_1)},(\log a_2)^{(d_2)},\cdots,
(\log a_n)^{(d_n)}\Big).$$

Obviously, the natural bijection between parabolic subgroups $P$ of $G$ and 
subsets of $\Delta_0$ or better $\Delta$ is given explicitly
 as follows: The parabolic subgroup $P=P_I$ corresponds to the 
subset $I(P)$ of $\Delta=\Big\{1,2,\cdots, r-1\Big\}$ such that 
$\Delta-I(P)=\{r_1,r_2,\cdots,r_n\}$. That is, 
$$I(P)=\Big\{1,2,\cdots, r_1-1\Big\}\bigcup\Big\{r_1+1,r_1+2,\cdots,r_2-1
\Big\}\bigcup\cdots
\bigcup\Big\{r_{n-1}+1, r_{n-1}+2,\cdots, r_n-1\Big\}.$$ 
For simplicity, from now on, 
we will use $I$ to indicate both the partition and its corresponding 
subset $I(P_I)$. For example, the subset $I=\emptyset$ corresponds
to the partition $\empty=(1,1,\cdots,1)$ of $r=1+1+\cdots+1$ 
(since $\Delta-\emptyset=\{1,2,\cdots,r-1\}$) whose 
corresponding parabolic subgroup is simply $P_0$ the minimal parabolic 
subgroup.

Consider now $\emptyset\subset I\subset J$. From 
$\Delta-J\subset\Delta-I\subset \Delta$, we get 
$P_0=P_\emptyset\subset P_I\subset P_J$. Set
$\Delta-J=(s_1,s_2,\cdots,s_m)$ corresponding to the partition
$r=f_1+f_2+\cdots+f_m$. Then
$f_1=d_1+d_2+\cdots+d_{k_1}, f_2=d_{k_1+1}+d_{k_1+2}+\cdots+d_{k_2},\cdots,
f_m=d_{k_{m-1}+1}+d_{k_{m-1}+2}+\cdots+d_{k_m}$ and
 $\frak a_J$ is a subspace of $\frak a_0$ defined by
$$\frak a_J:=\Big\{(H_1^{(f_1)},H_2^{(f_2)},\cdots,H_m^{(f_m)})
\in\frak a_0\Big\}.$$
Clearly $\frak a_J$ is in $\frak a_I$ since $\frak a_I$ is a 
subspace of $\frak a_0$ defined by
$$\frak a_I:=\Big\{\Big(A_1^{(d_1)},A_2^{(d_2)},\cdots,A_n^{(d_n)}\Big)
\in\frak a_0\Big\}$$
and we have the map $$\begin{matrix}
\frak a_J&\to&\frak a_I\\ 
\Big(H_1^{(f_1)},H_2^{(f_2)},\cdots,H_m^{(f_m)}\Big)&\mapsto&
\Big(H_1^{(d_1)},H_1^{(d_2)},\cdots,H_1^{(d_{k_1})},\cdots,H_m^{(d_{k_m})}
\Big).
\end{matrix}$$
Moreover, with the help of Iwasawa decomposition,
the corresponding logarithmic map may be simply viewed as the map 
$H_J:A_{J+}\to\frak a_J$ defined by $$\mathrm{diag}
\Big(b_1I_{f_1},b_2I_{f_2},\cdots,b_mI_{f_m}\Big)
\mapsto \Big((\log b_1)^{(f_1)},(\log b_2)^{(f_2)},\cdots,
(\log b_m)^{(f_m)}\Big).$$

To go further, for $i\in\Delta$, 
set $$\varpi_i:=\varpi_i^\Delta=\frac{r-i}{r}\sum_{k=1}^ik\cdot 
\alpha_k+\frac {i}{r}\sum_{k=i+1}^{r-1}(r-k)\cdot\alpha_k$$ be the 
corresponding fundamental weight. 
$\varpi_i$ is the $\mathbb Q$-linear form on $\frak a_\emptyset=\frak a_0$ 
given by
$$\begin{aligned}\varpi_i(H)=&H_1+H_2+\cdots+H_i
-\frac{i}{r}\Big(H_1+H_2+\cdots+H_r\Big)\\
=&H_1+H_2+\cdots+H_i\end{aligned}$$ 
for all $H=\Big(H_1,H_2,\cdots,H_r\Big)\in\frak a_0,$ 
so that $\widehat \Delta_I=\Big\{\varpi_i:i\in\Delta-I\Big\}$ 
since $$\Delta_I=
\Big\{\alpha_i\Big|_{\frak a_I}:i\in\Delta-I\Big\}.$$
 
Note that $\frak a_I^*=\langle\Delta_I\rangle_{\mathbb R}$ and
elements of $\widehat \Delta_I$ are non-degenerate positive linear 
combination of elements in 
$\Delta_I$, so  both $\Delta_I$ and $\widehat \Delta_I$ 
generate $\frak a_I^*$. Consequently, we have 
$$\frak a_I^*=\Big\langle\varpi_i:i\in\Delta-I\Big\rangle_{\mathbb R}
\supset 
\Big\langle\varpi_j:j\in\Delta-J\Big\rangle_{\mathbb R}=\frak a_J^*.$$ 
Moreover, by definition, the space $\frak a_I^J$ is the  
subspace of $\frak a_I$ annihilated by $\frak a_J^*$. Therefore,
$$\begin{aligned}\frak a_I^J=&\Big\{H\in\frak a_I:\alpha(H)=0\ 
\forall\alpha\in\frak a_J^*\Big\}\\
=&\Big\{H\in\frak a_I:\varpi_j(H)=0\ \forall j\in\Delta-J\Big\}\\
=&\Big\{(H_1^{(d_1)},H_2^{(d_2)},\cdots,H_n^{(d_n)})\in\frak a_0:
d_1H_1+d_2H_2+\cdots+ d_{k_1}H_{k_1}=0,\\
&\hskip 5.0cm d_{k_1+1}H_{k_1+1}+d_{k_1+2}
H_{k_1+2}+\cdots+ d_{k_2}H_{k_2}=0,\cdots\Big\}.\end{aligned}$$
Note that clearly $\dim \frak a_I^J+\dim \frak a_J=\dim\frak a_I$ and 
$\frak a_I^J\cap \frak a_J=\emptyset$, so we have the natural decomposition 
$$\frak a_I=\frak a_I^J\oplus\frak a_J.$$

Furthermore, note that 
$$\begin{aligned}I=&\bigcup_{i=1}^n\Big\{r_{j-1}+1,r_{j-1}+2,\cdots, 
r_{j}-1\Big\}\\
=&\Big\{1,2,\cdots,r-1\Big\}-\Big\{r_1,r_2,\cdots,r_{n-1}\Big\}.\end{aligned}$$
Set, for $k=1,2,\cdots, d_j-1$,
$$\varpi_{r_{j-1}+k}^I=\frac{d_j-k}{d_j}\sum_{l=1}^kl\cdot 
\alpha_{r_{j-1}+l}+\frac{k}{d_j}\sum_{l=k+1}^{d_j-1}\Big(d_j-l\Big)\cdot
\alpha_{r_{j-1}+l},$$
motivated by the construction of fundamental weights above. Then if, as usual, 
we define the coroot $\alpha_i^\vee\in\frak a_0$ by 
$$(\alpha_i^\vee)_k=\begin{cases}
1& \mathrm{if}\ $k=i$,\\
-1&\mathrm{if}\ $k=i+1$,\\
0&\mathrm{otherwise,}\end{cases}$$ 
then $\Big\{\alpha_i^\vee:i\in I\Big\}$ is a basis 
of $\frak a_\emptyset^I$ of $\frak a_\emptyset$ and easily one checks that 
$\varpi_i(\alpha_j^\vee)=\delta_{ij}$ for all $i,j\in I$. Consequently, 
$\Big\{\varpi_i^I\Big|_{\frak a_\emptyset^I}: i\in I\Big\}$ 
is its dual basis in $(\frak a_\emptyset^I)^*$.
In particular, we also get $$\widehat\Delta_I^J=\Big\{\varpi_k^J:k\in J-I=
(\Delta-I)-(\Delta-J)\Big\}.$$
 
We end this preperation with the following remarks about the notation. 
With respct to the decomposition
$$\Big(\frak a_\emptyset^J\Big)^*=\Big(\frak a_\emptyset^I\Big)^*
\oplus\Big(\frak a_I^J\Big)^*,$$ we 
have that the space $\Big(\frak a_\emptyset^J\Big)^*$ is generated by 
$$\Delta_\emptyset^J
=\Big\{\alpha_l:l\in J-\emptyset=J\Big\}\subset \Delta_\emptyset$$ 
the root system, 
similarly $\Big(\frak a_\emptyset^I\Big)^*$ is generated by 
$$\Delta_\emptyset^I
=\Big\{\alpha_k:k\in I-\emptyset=I\Big\}\subset \Delta_\emptyset$$
while $\Big(\frak a_I^J\Big)^*$ is generated by $\Delta_I^J$ 
which is not a subset of 
$\Delta_\emptyset$ but we have $\Delta_I^J\subset \Delta_I$ consists of the 
restriction of $\alpha_k$ to $\frak a_I^J$ with 
$k\in J-I=(\Delta-I)-(\Delta-J)$. That is 
$$\Delta_I^J=\Big\{\alpha_k\Big|_{\frak a_I^J}:k\in 
J-I=(\Delta-I)-(\Delta-J)\Big\}.$$ In particular, if $J=\Delta$ we have 
$$\frak a_\emptyset^\Delta=\frak a_\emptyset,\ \frak a_I^\Delta=\frak a_I,
\qquad\mathrm{and}\qquad \Delta_\emptyset^\Delta=\Delta_\emptyset,\ 
\Delta_I^\Delta=\Delta_I.$$ 

\subsection{An Algebraic Truncations}

For simplicity, we in this subsection work only with 
the field of rationals $\mathbb Q$ and use mixed
languages of adeles and lattices. Also by the discussion in Chapter 1, 
without loss of generality, we assume that $\mathbb Z$-lattices are of 
volume one. Accordingly, set $G=SL(r,\mathbb R)$. 

Denote by $\mathcal P$ the collection of all parabolic subgroups of $G$ 
both standard and non-standard, and by $\mathcal P_0$ the collection of all
standard parabolic subgroups. 
For a rank $r$ lattice $\Lambda$ of volume one,  denote the sublattices 
filtration associated to  $P\in \mathcal P$ by $$0=\Lambda_0\subset
\Lambda_1\subset\Lambda_2\subset\cdots\subset\Lambda_{|P|}=\Lambda.$$ 
Assume that $P$ corresponds to the partition $I=(d_1,d_2,\cdots,d_{n=:|P|})$. 
Consequently, we have $$\mathrm{rk}(\Lambda_i)=r_i:=d_1+d_2+\cdots+d_i,\qquad 
\mathrm{for}\ i=1,2,\cdots,|P|.$$

\noindent
{\bf {\large Definition}.} (1) Let $p,q:[0,r]\to {\mathbb R}$ be two polygons 
such that $p(0)=q(0)=p(r)=q(r)=0$. Then we say $q$ is {\it strongly
bigger than} $p$ {\it with respect to} $P$ and denote it by 
$q\triangleright_Pp$, if 
$$\begin{aligned}&\frac{q(r_1)-p(r_1)}{r_1}>\frac{q(r_2)-p(r_2)}{r_2}\\
>&\cdots\\
>&\frac{q(r_i)-p(r_i)}{r_i}>\frac{q(r_{i+1})-p(r_{i+1})}{r_{i+1}}\\
>&\cdots\\
>&\frac{q(r_{|P|-1})-p(r_{|P|-1})}{r_{|P|-1}}>
\frac{q(r_{|P|})-p(r_|{P|})}{r_{|P|}}=0.\end{aligned}$$ 

\noindent
{\bf Remark.} While this definition is motivated by Larfforgue [Laf], 
it is quite different. Recall that in Lafforgue's definition, $p>_Pq$ 
 means that $q(r_i)-p(r_i)>0$ for all $i=1,2,\cdots,|P|-1$. So our conditions
are much stronger than that of Lafforgue. The reason for making such a  
definition will become clear in the following discussion.
\vskip 0.30cm
For our own convenience, we call a polygon 
$p:[0,r]\to\mathbb R$ as a {\it normalized} one if $p(0)=p(r)=0$. From now 
on, all polygons are assumed to be normalized. 
\vskip 0.30cm
\noindent
{\bf {\large Definition}.} (2) 
Let  $p:[0,r]\to\mathbb R$ be a (normalized) polygon. 
Define the associated (real) character $T=T(p)$ of 
$M_0$ by the condition that $$\alpha_i(T)=\Big[p(i)-p(i-1)\Big]-
\Big[p(i+1)-p(i)\Big]$$ 
for all $i=1,2,\cdots, r-1.$
\vskip 0.30cm
\noindent
{\bf {\large Lemma 1}.} {\it With the same notation as above, $$T(p)=\Big(p(1), p(2)-p(1),
\cdots, p(i)-p(i-1),\cdots,p(r-1)-p(r-2),-p(r-1)\Big)$$.}

\noindent
{\it Proof.} Let $T=T(p)=\Big(t_1(p),t_2(p),\cdots,t_r(p)\Big)
=\Big(t_1,t_2,\cdots,t_r\Big)\in\frak a_0$. 
Note that $\alpha_i(T)=t_i-t_{i+1}$. Hence
$$\begin{aligned}t_i-t_{i+1}=&[p(i)-p(i-1)]-[p(i+1)-p(i)]\\
=&\Delta p(i)-\Delta p(i+1),\qquad i=1,2,\cdots, r-1\end{aligned}$$
where  $\Delta p(i):=p(i)-p(i-1)$. Consequently,
$$t_1-t_{i+1}=\Delta p(1)
-\Delta p(i+1),\qquad i=1,2,\cdots, r-1.$$ In particular, 
$$t_1-t_r=\Delta p(1)-\Delta p(r)=p(1)-p(r-1),$$ and 
$$t_{i+1}=t_1+\Delta p(i+1)-p(1),\qquad i=1,2,\cdots, r-1,$$
since $$\Delta p(1)=p(1)-p(0)=p(1)\qquad\mathrm{and}\qquad 
\Delta(r)=p(r)-p(r-1)=-p(r-1).$$  
But $$\begin{aligned}0=&\sum_{i=1}^{r} t_i\\
=&\sum_{i=0}^{r-1}\Big(t_1-p(1)+\Delta p(i+1)\Big)\\
=&(r-1)\Big(t_1-p(1)\Big)+\Big(p(0)-p(r)\Big)\\
=&(r-1)\Big(t_1-p(1)\Big).\end{aligned}$$ Therefore 
$t_1=p(1)$ and hence $t_{i+1}=\Delta p(i+1).$ This completes the proof.
\vskip 0.30cm
Now take $g=g(\Lambda)\in G(\mathbb A)$ such that the corresponding 
lattice is simply $\Lambda^g$. Then for $P\in\mathcal P$, we have the 
associated filtration
$$0=\Lambda_0^{g,P}\subset \Lambda_1^{g,P}\subset\cdots\subset 
\Lambda_{|P|}^{g,P}=\Lambda^g.$$
Recall that then we may define the associated polygon 
$p_P^g=p_P^{\Lambda^g}:[0,r]\to\mathbb R$ by the conditions that

\noindent
(1) $p_P^g(0)=p_P^g(r)=0$;

\noindent
(2) $p_P^g$ is affine on $[r_i,r_{i+1}]$, $i=1,2,\cdots,|P|-1$; and

\noindent
(3) $p_P^g(r_i)=\mathrm{deg}\Big(\Lambda_i^{g,P}\Big)-r_i\cdot
\frac{\mathrm{deg}\Big(\Lambda^{g}\Big)}{r},\ i=1,2,\cdots,|P|-1$.
\vskip 0.30cm
\noindent
Note that the volume of $\Lambda$ is assumed to be one, therefore 
(3) is equivalent to 

\noindent
(3)$'$ $p_P^g(r_i)=\mathrm{deg}\Big(\Lambda_i^{g,P}\Big),\ i=1,2,\cdots,|P|-1$.
\vskip 0.30cm
Now we can see the advantage of partially using adelic language:  
the values of $p_P^g$ may be written down more precisely.
Indeed, using Langlands decompositon $g=n\cdot m\cdot a(g)\cdot k$ 
with $n\in N_P(\mathbb A), m\in M_P(\mathbb A)^1, a\in A_+$ and $k\in K$ 
the naturally associated maximal compact subgroup 
$\prod_pSL(\mathcal O_{{\mathbb Q}_p})\times SO(r)$. Write 
$a=a(g)=\mathrm{diag}\Big(a_1I_{d_1},a_2I_{d_2},\cdots,a_{|P|}I_{d_{|P|}}
\Big)$ where $r=d_1+d_2+\cdots+d_{|P|}$ is the partition corresponding to $P$.
Then it is a standard fact that
$$\mathrm{deg}\Big(\Lambda_i^{g,P}\Big)=-\log\Big(\prod_{j=1}^i a_j^{d_j}\Big)
=-\sum_{j=1}^i d_j\log a_j,\qquad i=1,\cdots,|P|.$$

Set now ${\bf 1}(p_P^*\triangleright_Pp)$ to be the characteristic function 
of the subset
of $g$'s such that $p_P^g\triangleright_Pp$. Then we see that
$${\bf 1}(p_P^g\triangleright_Pp)=1$$ if and only if
$$\begin{aligned}&\frac{p_P^g(r_1)-p(r_1)}{r_1}>\frac{p_P^g(r_2)-p(r_2)}{r_2}\\
>&\cdots>\frac{p_P^g(r_i)-p(r_i)}{r_i}>\frac{p_P^g(r_{i+1})-
p(r_{i+1})}{r_{i+1}}\\
>&\cdots >\frac{p_P^g(r_{|P|-1})-p(r_{|P|-1})}
{r_{|P|-1}}>\frac{p_P^g(r_{|P|})-p(r_|{P|})}{r_{|P|}}=0.\end{aligned}$$ 
if and only if
$$\begin{aligned}&\frac{\mathrm{deg}\Big(\Lambda_1^{g,P}\Big)-p(r_1)}{r_1}>
\frac{\mathrm{deg}\Big(\Lambda_2^{g,P}\Big)-p(r_2)}{r_2}\\
>&\cdots>
\frac{\mathrm{deg}\Big(\Lambda_i^{g,P}\Big)-p(r_i)}{r_i}>
\frac{\mathrm{deg}\Big(\Lambda_{i+1}^{g,P}\Big)-p(r_{i+1})}{r_{i+1}}\\
>&\cdots \frac{\mathrm{deg}\Big(\Lambda_{|P|-1}^{g,P}\Big)-p(r_{|P|-1})}
{r_{|P|-1}}>\frac{\mathrm{deg}\Big(\Lambda_{|P|}^{g,P}\Big)-p(r_|{P|})}
{r_{|P|}}=0.\end{aligned}$$ That is to say, we have established the following
\vskip 0.30cm
\noindent
{\bf {\large Lemma 2}.} {\it With the same notation as above,
$${\bf 1}(p_P^g\triangleright_Pp)=1$$ if and only if
$$\begin{aligned}\frac{-d_1\log a_1-p(r_1)}{r_1}
>&\frac{-d_1\log a_1-d_2
\log a_2-p(r_2)}{r_2}\\
>&\cdots>\frac{-d_1\log a_1-d_2\log a_2-\cdots-
d_i\log a_i-p(r_i)}{r_i}\\
>&\frac{-d_1\log a_1-d_2\log a_2-\cdots-d_{i+1}
\log a_{i+1}-p(r_{i+1})}{r_{i+1}}\\
>&\cdots >\frac{-d_1\log a_1-d_2\log a_2
-\cdots-d_{|P|-1}\log a_{|P|-1}-p(r_{|P|-1})}{r_{|P|-1}}\\
>&\frac{-d_1\log a_1-d_2\log a_2-\cdots-d_{|P|}\log a_{|P|}
-p(r_|{P|})}{r_{|P|}}=0.\end{aligned}$$}
\eject
\section{A Bridge between Algebraic and Analytic Truncations}

In this subsection, we expose a beautiful intrinsic relation
between algebraic and analytic truncations.  
To begin with, let us study the function $\tau_P(-H_0(g)-T(p)).$ 
When does it take the value one?

Recall that for $P=P_I$, $\tau_P$ is the characteristic function of the 
following subset of $\frak a_0=\frak a_0^I+\frak a_I$
defined to be $$\frak a_0^I\oplus
\Big\{H\in\frak a_I:\alpha_i(H)>0\quad\forall i\in\Delta-I\Big\}.$$
That is to say, if $H\in\frak a_0$ belongs to  $$\frak a_0^I
\oplus\Big\{H\in\frak a_i:\alpha_I(H)>0\quad i=r_1,r_2,\cdots,r_{|P|-1}
\Big\},$$ then $\tau_P(H)=1$. 

Recall also that $H\in\frak a_I$ if and only if 
$H=\Big(H_1^{(d_1)},H_2^{(d_2)},\cdots, H_{|P|}^{(d_{|P|})}\Big)$ and 
$\sum_{i=1}^{|P|} d_i H_i=0$. Therefore, the conditions that 
$$H\in\frak a_I\quad\mathrm{such\ that}\qquad\alpha_i(H)>0\quad 
\mathrm{for\ all}\ i=r_1,r_2,\cdots,r_{|P|-1}$$
is equivalent to the conditions that 
$$H=\Big(H_1^{(d_1)},H_2^{(d_2)},
\cdots, H_{|P|}^{(d_{|P|})}\Big)\qquad\mathrm{where}\qquad
\sum_{i=1}^{|P|} d_i H_i=0$$ and $$H_1-H_2>0, H_2-H_3>0,\cdots, 
H_{|P|}-H_{|P|-1}>0.$$
 
Furthermore, by definition, 
$$-H_0(g)-T(p)=\Big(-t_1-\log b_1,-t_2-\log b_2,\cdots,-t_r-\log b_r\Big)$$ 
where $$a(g)=\mathrm{diag}\Big(b_1,b_2,
\cdots,b_r\Big)=\mathrm{diag}\Big(a_1^{(d_1)},a_2^{(d_2)},\cdots,
a_{|P|}^{(d_{|P|})}\Big)$$ and $$T(p)=\Big(t_1,t_2,\cdots,t_r\Big)=
\Big(p(1), p(2)-p(1),\cdots,p(r-1)-p(r-2),-p(r-1)\Big).$$ 

On the other hand, for $H=\Big(H_1,H_2,\cdots,H_r\Big)\in \frak a_0$, 
the projections of $H$ to $\frak a_I$  (resp. to  $\frak a_0^I$) is given by
$$\begin{aligned}\bigg(\Big(\frac{H_1+H_2+\cdots+H_{r_1}}{d_1}\Big)^{(d_1)},
&\Big(\frac{H_{r_1+1}+H_{r_1+2}+\cdots+H_{r_2}}{d_2}\Big)^{(d_2)},\\
&\cdots,\Big(\frac{H_{r_{|P|-1}+1}+H_{r_{|P|-1}+2}+\cdots+
H_{r_{|P|}}}{d_{|P|}}\Big)^{(d_{|P|})}\bigg)\in\frak a_I\end{aligned}$$
(resp.
$$\begin{aligned}&\bigg(H_1-\frac{H_1+H_2+\cdots+H_{r_1}}{d_1},
H_2-\frac{H_1+H_2+\cdots+H_{r_1}}{d_1},\cdots,H_{r_1}-\frac{H_1+H_2+
\cdots+H_{r_1}}{d_1},\\
&\cdots,
H_{r_{|P|-1}+1}-\frac{H_{r_{|P|-1}+1}+H_{r_{|P|-1}+2}+\cdots+
H_{r_{|P|}}}{d_{|P|}},
H_{r_{|P|-1}+2}-\frac{H_{r_{|P|-1}+1}+H_{r_{|P|-1}+2}+\cdots+
H_{r_{|P|}}}{d_{|P|}},\\
&\qquad\cdots,H_{r_{|P|}}-\frac{H_{r_{|P|-1}+1}+H_{r_{|P|-1}+2}+\cdots+
H_{r_{|P|}}}{d_{|P|}}\bigg)\\
=&\Big(H_1,H_2,\cdots,H_r\Big)-\bigg(\Big(\frac{H_1+H_2+\cdots+
H_{r_1}}{d_1}\Big)^{(d_1)},
\Big(\frac{H_{r_1+1}+H_{r_1+2}+\cdots+H_{r_2}}{d_2}\Big)^{(d_2)},\\
&\hskip 5.0cm\cdots,
\Big(\frac{H_{r_{|P|-1}+1}+H_{r_{|P|-1}+2}+\cdots+H_{r_{|P|}}}{d_{|P|}}
\Big)^{(d_{|P|})}\bigg)\in\frak a_0^I.\end{aligned}$$
Therefore,  
$$\tau_P\Big(-H_0(g)-T(p)\Big)=1$$ if and only if the coordinates $H_i$ of
 $-H_0(g)-T(p):=H=\Big(H_1,H_2,\cdots,H_r\Big)$  satisfy the 
conditions that
  $$\frac{H_1+H_2+\cdots+H_{r_1}}{d_1}>\frac{H_{r_1+1}+H_{r_1+2}+
\cdots+H_{r_2}}{d_2}>
  \cdots>\frac{H_{r_{|P|-1}+1}+H_{r_{|P|-1}+2}+\cdots+H_{r_{|P|}}}
{d_{|P|}}.$$ Clearly, this latest group 
  of inequalities, after replacing $H_i$'s with the precise coordinates 
of $-H_0(g)-T(p)$,  are then changed to the following group conditions:
 $$\begin{aligned}&\frac{\Big(-\log b_1-p(1)\Big)+\Big(-\log b_2-p(2)+p(1)\Big)
+\cdots +
\Big(-\log b_{r_1}-p(r_1)+p(r_1-1)\Big)}{r_1}>\\
&\frac{\Big(-\log b_{r_1+1}-p(r_1+1)+p(r_1)\Big)+\Big(-\log b_{r_1+2}-p(r_1+2)+p(r_1+1)\Big)+
\cdots+\Big(-\log b_{r_2}-p(r_2)+p(r_2-1)\Big)}{r_2-r_1}\\
 &>\cdots>\\
&\frac{\Big(-\log b_{r_{i-1}+1}-p(r_{i-1}+1)+p(r_i)\Big)+\Big(-
\log b_{r_{i-1}+2}-p(r_{i-1}+2)+p(r_{i-1}+1)\Big)+
 \cdots+\Big(-\log b_{r_i}-p(r_i)+p(r_i-1)\Big)}{r_i-r_{i-1}}\\
 &>\cdots>\\
&\frac{\Big(-\log b_{r_{|P|-1}+1}-p(r_{|P|-1}+1)+p(r_{|P|})\Big)+\Big(-
\log b_{r_{|P|-1}+2}-p(r_{|P|-1}+2)+p(r_{|P|-1}+1)\Big)
 +\cdots+\Big(-\log b_{r_{|P|}}-p(r_{|P|})+p(r_{|P|}-1)\Big)}
{r_{|P|}-r_{|P|-1}}.\end{aligned}$$ 
 Or, after the obvious simplification, with $$B_1:=-\log b_1, B_2
=-\log b_2,\cdots,B_r=-\log b_r,\ r_0:=0,$$
 $$\begin{aligned}&\frac{B_1+B_2+\cdots+B_{r_1}-p(r_1)+p(0)}{r_1-r_0}\\
&>\frac{B_{r_1+1}+B_{r_1+2}+
 \cdots+B_{r_2}-p(r_2)+p(r_1)}{r_2-r_1}\\
 &>\cdots\\
&>\frac{B_{r_{i-1}+1}+B_{r_{i-1}+2}+\cdots+B_{r_i}-p(r_i)
+p(r_{i-1})}{r_i-r_{i-1}}\\
&>\frac{B_{r_{i}+1}+B_{r_{i}+2}+\cdots+B_{r_{i+1}}-p(r_{i+1})+
p(r_{i})}{r_{i+1}-r_{i}}\\
 &>\cdots\\
&>\frac{B_{r_{|P|-1}+1}+B_{r_{|P|-1}+2}+\cdots+B_{r_{|P|}}-p(r_{|P|})
+p(r_{|P|-1})}{r_{|P|}-r_{|P|-1}}.\end{aligned}$$
 Thus note that $$\Big(b_1,b_2,\cdots, b_r\Big)=a(g)=
\Big(a_1^{(r_1-r_0)},a_2^{(r_2-r_1)},
\cdots, a_{|P|}^{(r_{|P|}-r_{|P|-1})}\Big),$$
and set $A_i=-\log a_i$, we arrive at the following
\vskip 0.30cm
\noindent
{\bf {\large Sublemma}}. {\it With the same notation as above,} 
$$\tau_P\Big(-H_0(g)-T(p)\Big)=1$$ 
{\it if and only if the following conditions
 are satisfied}:

\noindent 
(1) $\displaystyle{\sum_{i=1}^{|P|}(r_i-r_{i-1})A_i=0};$ and

\noindent 
(2) $\displaystyle{\frac{(r_i-r_{i-1})A_i-p(r_i)+p(r_{i-1})}{r_i-r_{i-1}}>
\frac{(r_{i+1}-r_{i})A_{i+1}
 -p(r_{i+1})+p(r_{i})}{r_{i+1}-r_{i}},}$ for $1\leq i\leq |P|-1$.
\vskip 0.30cm
Now we are ready to expose an intrinsic relation 
between algebraic and analytic truncations:
\vskip 0.30cm
\noindent
{\bf {\large Bridge Lemma}}.  
{\it For a fixed convex normalized polygon 
$p:[0,r]\to\mathbb R$, and $g\in SL_r(\mathbb A)$, with respect to
 any parabolic subgroup $P$, we have
 $$\tau_P\Big(-H_0(g)-T(p)\Big)=\bold 1\Big(p_P^g\triangleright_Pp\Big).$$}
 
\noindent
{\it Proof.} This is based on the following very easy trick: If $a,c>0$,
 then $$\frac{b}{a}>\frac{d}{c} \Leftrightarrow
\frac{b}{a}>\frac{b+d}{a+c} \Leftrightarrow\frac{b+d}{a+c}>\frac{d}{c} 
\Leftrightarrow\frac{b}{a}>\frac{b+d}{a+c}>\frac{d}{c}.$$
Consequently,
$$\frac{b_1}{a_1}>\frac{b_2}{a_2}>\frac{b_3}{a_3}>\cdots>\frac{b_s}{a_s}$$
if and only if
$$\frac{b_1}{a_1}>\frac{b_1+b_2}{a_1+a_2}>\frac{b_3}{a_3}>\cdots>
\frac{b_s}{a_s}$$
if and only if 
$$\frac{b_1}{a_1}>\frac{b_1+b_2}{a_1+a_2}>\frac{b_1+b_2+b_3}{a_1+a_2+a_3}
>\frac{b_4}{a_4}>\cdots>\frac{b_s}{a_s}$$
$$\cdots\cdots\cdots$$
if and only if 
$$\begin{aligned}&\frac{b_1}{a_1}>\frac{b_1+b_2}{a_1+a_2}
>\frac{b_1+b_2+b_3}{a_1+a_2+a_3}
>\cdots\\
&>\frac{b_1+b_2+\cdots+b_i}{a_1+a_2+\cdots+a_i}
>\cdots
>\frac{b_1+b_2+\cdots+b_{s-1}}{a_1+a_2+\cdots+a_{s-1}}\\
&>\frac{b_s}{a_s}\end{aligned}$$
if and only if 
$$\begin{aligned}&\frac{b_1}{a_1}
>\frac{b_1+b_2}{a_1+a_2}\\
&>\frac{b_1+b_2+b_3}{a_1+a_2+a_3}\\
&>\cdots\\
&>\frac{b_1+b_2+\cdots+b_i}{a_1+a_2+\cdots+a_i}\\
&>\cdots\\
&>\frac{b_1+b_2+\cdots+b_{s-1}}{a_1+a_2+\cdots+a_{s-1}}\\
&>\frac{b_1+b_2+\cdots+b_{s-1}+b_s}{a_1+a_2+\cdots+a_{s-1}+a_s}.\end{aligned}$$

Therefore, if we start with $$\tau_P\Big(H_0(g)-T(p)\Big)=1,$$ using the Sublemma  
above, we see that it is equivalent to
$$\begin{aligned}&\frac{(r_1-r_{0})A_1-p(r_1)+p(r_{0})}{r_1-r_{0}}\\
&>\frac{(r_{2}-r_{1})A_{2}-p(r_{2})+p(r_{1})}{r_{2}-r_{1}}\\
&>\cdots\\
&>\frac{(r_i-r_{i-1})A_i-p(r_i)+p(r_{i-1})}{r_i-r_{i-1}}\\
&>\frac{(r_{i+1}-r_{i})A_{i+1}-p(r_{i+1})+p(r_{i})}{r_{i+1}-r_{i}}\\
&>\cdots\\
&>\frac{(r_{|P|}-r_{|P|-1})A_{|P|}-p(r_{|P|})+p(r_{|P|-1})}
{r_{|P|}-r_{|P|-1}}.\end{aligned}$$
So by applying the above argument on inequalities, we further conclude that
this latest group of inequalities are equivalent to
$$\begin{aligned}&\frac{r_1A_1-p(r_1)}{r_1}>\frac{r_1A_1+(r_{2}-r_{1})
A_{2}-p(r_{2})}{r_{2}}\\
&>\cdots\cdots\\
&> \frac{r_1A_1+(r_2-r_1)A_2+\cdots+(r_i-r_{i-1})A_i-p(r_i)}{r_i}\\
&>\cdots\cdots\\
&>\frac{r_1A_1+(r_2-r_1)A_2+\cdots+(r_{|P|-1}-r_{|P|-2})
A_{|P|-1}-p(r_{|P|-1})}{r_{|P|-1}}\\
&>\frac{r_1A_1+(r_2-r_1)A_2+\cdots+(r_{|P|-1}-r_{|P|-2})A_{|P|-1}+
(r_{|P|}-r_{|P|-1})A_{|P|}-p(r_{|P|})}{r_{|P|}},\end{aligned}$$
since $r_0=0$. 
Now, using the fact that 
$$\mathrm{deg}\Lambda_i^g=-\log\prod_{j=1}^i a_j^{(d_j)}
=r_1A_1+(r_2-r_1)A_2+\cdots+(r_i-r_{i-1})A_i,$$ we conclude that, by 
the Lemma 2 in the previous subsection, 
this latest group of inequalities is equivalent to the condition that
$$\bold 1(p_P^g\triangleright_Pp)=1,$$ provided that
the last term $$\frac{r_1A_1+(r_2-r_1)A_2+
\cdots+(r_{|P|-1}-r_{|P|-2})A_{|P|-1}+(r_{|P|}-r_{|P|-1})A_{|P|}-
p(r_{|P|})}{r_{|P|}}=0.$$ 
But note that $p(r_{|P|})=0$ since $r_{|P|}=r$ and that the volumes of the 
lattices are fixed to be one, 
so $$\begin{aligned}&0=\mathrm{deg}
\Lambda^g=\mathrm{deg}\Lambda_{|P|}^g\\
&=-\log\prod_{j=1}^{|P|} a_j^{(d_j)}
=r_1A_1+(r_2-r_1)A_2+\cdots+(r_{|P|}-r_{|P|-1})A_{|P|}.\end{aligned}$$ 
This completes the proof.

\noindent
{\bf Remark.} It is well-known that in Arthur's analytic truncation, 
what is used is the function $\widehat\tau$ associated to positive cone.
In this sense, to everybody's surprise, it is $\tau$, associated to
positive chamber, not $\widehat\tau$ that
is appeared as a part of the bridge. However, as we will see later, 
such a subtle difference creates new room for non-abelian parts 
of the theory.
\eject
\section{Global Relation between 
Analytic Truncation and Geometric Truncation}

In this subsection, we give a global relation between 
analytic truncation and geometric truncation. 

Let $\Lambda=\Lambda^g$ be a rank $r$ lattice associated to 
$g\in\mathrm{GL}_r(\mathbb A)$ and  $P\in\mathcal P$ a parabolic 
subgroup.  Denote the sublattices 
filtration associated to  $P$ by $$0=\Lambda_0\subset
\Lambda_1\subset\Lambda_2\subset\cdots\subset\Lambda_{|P|}=\Lambda.$$ 
Assume that $P$ corresponds to the partition $I=(d_1,d_2,\cdots,d_{n=:|P|})$. 
Consequently, we have $$\mathrm{rk}(\Lambda_i)=r_i:=d_1+d_2+\cdots+d_i,\qquad 
\mathrm{for}\ i=1,2,\cdots,|P|.$$ 
Let $p,q:[0,r]\to {\mathbb R}$ be two polygons 
such that $p(0)=q(0)=p(r)=q(r)=0$. Then following Lafforgue, 
we say $q$ is {\it bigger than} $p$ {\it with respect to} $P$ and denote it by 
$q>_Pp$, if $q(r_i)-p(r_i)>0$ for all $i=1,\cdots,|P|-1.$ 
Introduce also  the characteristic function $\bold 1(\overline p^*\leq p)$ by
$$\bold 1(\overline p^g\leq p)
=\begin{cases} 1,&\text{if $\overline p^g\leq p$;}\\
0,&\text{otherwise}.\end{cases}$$ 
Then we have the following
\vskip 0.30cm
\noindent
{\bf {\large Fundamental Relation.}} 
{\it For a fixed convex polygon $p:[0,r]\to 
{\mathbb R}$ such that $p(0)=p(r)=0$, the following relations
$$\bold 1(\overline p^g\leq p)=\sum_{P:\, \text{standard\, 
parabolic}}(-1)^{|P|-1}\sum_{\delta\in
P({\mathbb Z})\backslash G({\mathbb Z})}\bold 1(p_P^{\delta g}>_Pp)
\qquad\forall g\in G.$$}

\noindent
{\bf Remarks.} (1) This result and its proof below are 
motivated by a similar result of Lafforgue
for vector bundles over function fields. 

\noindent
(2) The right hand side may be naturally decomposite into 
two parts according to whether $P=G$ or not.
In such  away,  the right hand side becomes
$$\bold 1_G-\sum_{P:\, \mathrm{proper\,standard\, parabolic}}(-1)^{|P|-1}
\cdots.$$ This then exposes two aspects of our geometric truncation:
First of all, if a lattice is not stable, then there will be parabolic
subgroups which take the responsibility; Secondly, each parabolic subgroup
has its fix role -- Essentially, they should be counted only once 
each time. In other words, if more are substracted, then we need to 
add one fewer back to make sure the whole process is not overdone.
\vskip 0.30cm
\noindent 
{\it Proof.} (Following Lafforgue) Note that all parabolic 
subgroups of $G$ may be obtained from taking $\delta$-conjugates of standard 
parabolic subgroups with $\delta\in P(F)\backslash G(F)$, and that 
$\overline p_{P}^{\delta g}=\overline p_{\delta P\delta^{-1}}^g$. Therefore, 
it suffices to prove the following
\vskip 0.30cm
\noindent
{\bf {\large Fundamental Relation.}}$'$ 
$$\bold 1(\overline p^g\leq p)=
\sum_{P:\,\mathrm{parabolic}}(-1)^{|P|-1}\bold 1(p_P^{g}>_Pp).$$
\vskip 0.30cm
To establish this equality, we consider two cases depending on
\begin{description}
\item[whether] (a) $\overline p^g\leq p$;

\item[or] (b) $\overline p^g$ is not bounded from above by $p$.
\end{description}
Assume that (a) holds. Then the LHS is simply 1. On the other hand,  
for a proper parabolic subgroup $Q\not=G$, by definition, $p_Q\leq p$. 
Hence $\bold 1(p_Q^{g}>_Qp)=0$.
Consequently, the right hand side degenerates 
to the one consisting of a single term involving $G$ only, which, 
by definition, it is simply 1 as well. We are done.
\vskip 0.30cm 
Next, consider (b), for which $\overline p^g$ is 
not bounded from above by $p$. This is a bit complicated.
To start with, we set
$\overline\Lambda_*^{g,P}$ to be the unique filtration of $\Lambda^g$ 
such that the associated polygon $\overline p_P^g$  is 
{\bf maximal} among $p_Q^g$ where $Q$ runs over all parabolic subgroups 
contained in $P$. That is
 $$\overline p_P^g:=\max\Big\{p_Q^g:\,Q\subset P,\ \mathrm{parabolic\ 
subgroup}\Big\}.$$
Denote the associated parabolic subgroup by $\overline Q^g_P$ and call it 
$P$-{\it canonical}. Clearly, $G$-canonical is simply canonical.

Moreover, since $\overline p^g\not\leq p$, the LHS=0. So we should prove that
the RHS=0 as well. For this, regroup the parabolic subgroups $P$ 
appeared in the RHS according to the refined canonical parabolic subgroup 
$\overline Q_P^g$.
Then clearly, it suffices to establish the following
\vskip 0.30cm
\noindent
{\bf {\large Proposition.}} {\it With the same notation as above, if 
$\overline p^g\not\leq p$, then for any fixed parabolic subgroup $Q$,
$$\sum_{P:\overline Q_P^g=Q}(-1)^{|P|}\bold 1(p_P^p>_Pp)=0.$$}

\noindent
{\it Proof.} We break this into the following 3 steps.

\noindent
{\bf Step A}. Fix a non-negative real number $\mu$. For the lattice 
$\Lambda=\Lambda^g$ and a fixed parabolic subgroup $P$ of $G$, set
\vskip 0.30cm
\noindent
\begin{description}
\item[(1)] $^\mu\overline Q^g$ to be the unique parabolic subgroup of 
$G$ such that if
$^\mu\overline\Lambda^g_*$ is the induced filtration of $\Lambda^g$, then
$$\Big\{^\mu\overline\Lambda^g_*\Big\}\,=\,\Big\{\overline\Lambda_j^g:\ \mu\Big
(\overline\Lambda_j^g\Big/\overline\Lambda_{j-1}^g\Big)>
\mu\Big(\overline\Lambda_{j+1}^g\Big/\overline\Lambda_{j}^g\Big)+\mu\Big\},$$
where $\overline\Lambda_*^g$ is the canonical  filtration associated with 
$\Lambda^g$;

\item[(2)] $^\mu\overline Q_P^g$ to be the unique parabolic subgroup of 
$G$ such that if
$^\mu\overline\Lambda^{g,P}_*$ is the induced filtration of $\Lambda^g$, then
$$\Big\{^\mu\Lambda^{g,P}_*\Big\}\,=\,\Big\{\Lambda_j^{g,P}\Big\}\cup
\Big\{\overline\Lambda_j^{g,P}:\ \mu\Big(\overline\Lambda_j^{g,P}\Big/
\overline\Lambda_{j-1}^{g,P}\Big)>
\mu\Big(\overline\Lambda_{j+1}^{g,P}\Big/\overline\Lambda_{j}^{g,P}\Big)+\mu
\Big\},$$
where $\Lambda_j^{g,P}$ is the filtration of $\Lambda^g$ induced from the
 parabolic subgroup $P$, and $\overline\Lambda_*^{g,P}$ is the
$P$-canonical  filtration associated with $\Lambda^g$.
\end{description}

In particular, from this construction, we have:
\vskip 0.30cm
\noindent
(i) For $\overline\Lambda_*^{g,P}$, the so-called $P$-canonical filtration,
the Harder-Narasimhan conditions hold, i.e., $\overline\Lambda_j^{g,P}\Big/\overline\Lambda_{j-1}^{g,P}$ are semi-stable, and
$$\mu\Big(\overline\Lambda_j^{g,P}\Big/\overline\Lambda_{j-1}^{g,P}\Big)>
\mu\Big(\overline\Lambda_{j+1}^{g,P}\Big/\overline\Lambda_{j}^{g,P}\Big),$$
except for those $j$ appeared as that associated to $P$ itself;
\vskip 0.30cm
\noindent
(ii) While the canonical filtration  $\overline\Lambda_*^g$ of $\Lambda^g$ 
satisfying
$$\mu\Big(\overline\Lambda_j^{g}\Big/\overline\Lambda_{j-1}^{g}\Big)>
\mu\Big(\overline\Lambda_{j+1}^{g}\Big/\overline\Lambda_{j}^{g}\Big),
\quad\forall j=1,\cdots,|\overline Q^g|-1,$$
as a partial filtration, $\Big\{^\mu\overline\Lambda^g_*\Big\}$ 
is the one obtained from this canonical one by
selecting only the part where much stronger conditions 
$$\mu\Big(\overline\Lambda_j^g\Big/\overline\Lambda_{j-1}^g\Big)\,>\,
\mu\Big(\overline\Lambda_{j+1}^g\Big/\overline\Lambda_{j}^g\Big)+\mu$$ are 
satisfied. Consequently, 
$$\overline Q^g\subset\,^\mu\overline Q^g;$$

\noindent
(iii) The filtration $\{^\mu\Lambda^{g,P}_*\}$ contains not only
that part of the $P$-canonical filtration
 where stronger conditions $$\mu\Big(\overline\Lambda_j^{g,P}\Big/
\overline\Lambda_{j-1}^{g,P}\Big)\,>\,
\mu\Big(\overline\Lambda_{j+1}^{g,P}\Big/\overline\Lambda_{j}^{g,P}\Big)+\mu$$ 
are satisfied, it contains also the
full filtration of $\Lambda^g$ associated to $P$. So, 
$$\overline Q^g_P\subset\, ^\mu\overline
Q^g_P,\qquad\mathrm{while}\qquad ^\mu\overline Q^g_P\subset P.$$
\vskip 0.30cm
\noindent
{\bf Step B}. For $g\in G(\mathbb A)$, set
$$\begin{aligned}I(Q_1):=&\Big\{\,\mathrm {rk}\Lambda_i^{g,Q_1}:\ 
\ i=1,2,\cdots,|Q_1|-1\,\Big\};\\
I(Q_2):=&\Big\{\,\mathrm {rk}\Lambda_j^{g,Q_2}:\ j=1,2,\cdots,|Q_2|-1\,
\Big\};\\
J_0(Q_1):=&\Big\{\,\mathrm {rk}\Lambda_k^{g,Q_1}:\ 
\mu\Big(\overline\Lambda_k^{g,Q_1}\Big/
\overline\Lambda_{k-1}^{g,Q_1}\Big)\,>\,
\mu\Big(\overline\Lambda_{k+1}^{g,Q_1}\Big/\overline\Lambda_{k}^{g,Q_1}\Big)
\,\Big\};\\
J_\mu(Q_1):=&\Big\{\,\mathrm {rk}\Lambda_l^{g,Q_1}:\ 
\mu\Big(\overline\Lambda_l^{g,Q_1}\Big/
\overline\Lambda_{l-1}^{g,Q_1}\Big)\,>\,
\mu\Big(\overline\Lambda_{l+1}^{g,Q_1}\Big/\overline\Lambda_{l}^{g,Q_1}\Big)
+\mu\,\Big\};\\
K_p(Q_2):=&\Big\{\,\mathrm {rk}\Lambda_m^{g,Q_2}:\ 
\Big(p_{Q_2}^g-p\Big)\Big(\mathrm{rk}\Lambda_m^{g,Q_2}\Big)>0\,\Big\}.
\end{aligned}$$

\noindent
{\bf {\large Lemma.}} {\it With the same notation as above, the following
correspondence
$$\Big\{\,P\in\mathcal P:\ \overline Q_P^g=Q_1,\ \,^\mu\overline Q_P^g=Q_2,\
p_P^g>_Pp\,\Big\}\longrightarrow$$
$$\Big\{\,\Sigma:\ \Sigma\subset I(Q_2),\ J_\mu(Q_1)\subset I(Q_2), 
\Big(I(Q_1)-J_0(Q_1)\Big)\cup 
\Big(I(Q_2)-J_\mu(Q_1)\Big)\subset \Sigma\subset K_p(Q_2)\,\Big\}$$ defined by
$$P\mapsto I(P):=\Big\{\,\mathrm {rk}\Lambda_i^{g,P}:\ i=1,2,\cdots,|P|-1\,
\Big\}$$ is a well-defined bijection.}
\vskip 0.30cm
\noindent
{\it Proof of the Lemma.} (I) {\bf Well-defined}.
For this, let us look at the conditions for $P$ in 
$$\Big\{\,P\in\mathcal P:\ \overline Q_P^g=Q_1,\ \,^\mu\overline Q_P^g=Q_2,\
p_P^g>_Pp\,\Big\}.$$

\noindent 
(i) $\overline Q_P^g=Q_1$ gives the filtration
$$\{0\}=\Lambda_{0}^{g,Q_1}\subset \Lambda_{1}^{g,Q_1}\subset
\Lambda_{2}^{g,Q_1}\subset\cdots\subset \Lambda_{|Q_1|-1}^{g,Q_1}\subset
\Lambda_{|Q_1|}^{g,Q_1}=\Lambda^{g}$$ with
$p_{Q_1}^g$ be the maximal among all $\Big\{ p_Q^g:Q\subset P\Big\}$.
Hence,

\noindent
(a) $Q_1\subset P$; and

\noindent
(b) Except possibly at indices corresponding to the partion coming from $P$,
$$\mu\Big(\overline\Lambda_k^{g,Q_1}\Big/
\overline\Lambda_{k-1}^{g,Q_1}\Big)\,>\,
\mu\Big(\overline\Lambda_{k+1}^{g,Q_1}\Big/\overline\Lambda_{k}^{g,Q_1}\Big),$$
by the maximal property, or more clearly, the Harder-Narasimhan property.
Therefore, we have
$$I(Q_1)-J_0(Q_1)\subset I(P).$$

\noindent
(ii) $\,^\mu\overline Q_P^g=Q_2$ gives the filtration
$$\{0\}=\Lambda_{0}^{g,Q_2}\subset \Lambda_{1}^{g,Q_2}\subset
\Lambda_{2}^{g,Q_2}\subset\cdots\subset \Lambda_{|Q_2|-1}^{g,Q_2}\subset
\Lambda_{|Q_2|}^{g,Q_2}=\Lambda^{g},$$
which by definition coincides with
$$\Big\{\,\Lambda_{j}^{g,P}:\ j=1,2,\cdots,|P|-1\,\Big\}\cup
\Big\{\,\Lambda_{k}^{g,Q_1}:\ \mu\Big(\overline\Lambda_k^{g,Q_1}\Big/
\overline\Lambda_{k-1}^{g,Q_1}\Big)\,>\,
\mu\Big(\overline\Lambda_{k+1}^{g,Q_1}\Big/\overline\Lambda_{k}^{g,Q_1}\Big)+\mu\,\Big\}.$$ By taking the rank for each lattices, we get the relation
$I(Q_2)=I(P)\cup J_\mu(Q_1)$. Therefore,

\noindent
(a)  $Q_2\subset P$ or equivalently $I(P)\subset I(Q_2)$; and

\noindent
(b) $I(Q_2)-J_\mu(Q_1)\subset I(P)$.

\noindent
(iii) $p_P^g>_Pp$ is, by definition, equivalent to the following 
group of inequalities
$$\Big(p_{P}^g-p\Big)\Big(\mathrm{rk}\Lambda_j^{g,P}\Big)>0.$$
Consequently, $I(P)\subset K_p(Q_2)$ since $I(P)\subset I(Q_2)$.

Therefore, $P\mapsto I(P)$ is well-defined.

\noindent
(II) {\bf Injectivity}. This is clear since with a fixed partition $I$, 
there exists
only one unique standard parabolic subgroup $P_I$ such that $I(P_I)=\Delta-I$.

\noindent
(III) {\bf Surjectivity}. For a subset $\Sigma$ of $\Big\{1,2,\cdots,r\Big\}$ 
such that
$$\Sigma\subset I(Q_2),\ J_\mu(Q_1)\subset I(Q_2), 
\Big(I(Q_1)-J_0(Q_1)\Big)\cup 
\Big(I(Q_2)-J_\mu(Q_1)\Big)\subset \Sigma\subset K_p(Q_2),$$
clearly, there exists a parabolic subgroup $P=P_\Sigma$. So we need to check
whether $\overline Q_P^g=Q_1, \,^\mu\overline Q_P^g=Q_2$ and  $p_P^g>_Pp$.

It is clear that $p_P^g>_Pp$ since $\Sigma\subset K_p(Q_2)$. On the other hand,
$\overline Q_P^g=Q_1, \,^\mu\overline Q_P^g=Q_2$ are direct consequences
of the definition of P-canonical filtration and the conditions that
$\Big(I(Q_1)-J_0(Q_1)\Big)\subset\Sigma$ and  
$\Big(I(Q_2)-J_\mu(Q_1)\Big)\subset \Sigma$ as explained in the proof above
of the well-defined statement. This completes the proof of the Lemma.
\vskip 0.30cm
\noindent
{\bf Step C.} With the Lemma, we have
$$\begin{aligned}
&\Big|\sum_{P:\overline Q_P^g=Q_1,\,^\mu\overline Q_P^g=Q_2}(-1)^{|P|}
\bold 1(p_P>_Pp)\Big|\\
=&\Big|\sum_{\Sigma:\Big(I(Q_1)-J_0(Q_1)\Big)\cup 
\Big(I(Q_2)-J_\mu(Q_1)\Big)\subset \Sigma\subset 
K_p(Q_2)}(-1)^{\#(\sigma)}\Big|\\
=&\begin{cases}
0,&\text{if $\Big(I(Q_1)-J_0(Q_1)\Big)\cup 
\Big(I(Q_2)-J_\mu(Q_1)\Big)\not=K_p(Q_2)$;}\\
1,&\text{if $\Big(I(Q_1)-J_0(Q_1)\Big)\cup 
\Big(I(Q_2)-J_\mu(Q_1)\Big)=K_p(Q_2)$.}\end{cases}\end{aligned}$$
Surely, we want to show that the value  1 cannot be taken for $g$ such that
$\overline p^g$ is not bounded by $p$.
To this end, as the final push to establish our Fundamental Relation, 
set $\mu=0$ and $Q_1=Q_2=Q$ in the above discussion so that
$$\begin{aligned}
\sum_{P:\overline Q_P^g=Q_1,\,^\mu\overline Q_P^g=Q_2}(-1)^{|P|}
\bold 1(p_P>_Pp)=&\sum_{P:\overline Q_P^g=Q=\,^0\overline Q_P^g}
(-1)^{|P|}\bold 1(p_P>_Pp)\\
=&\sum_{P:\overline Q_P^g=Q}(-1)^{|P|}\bold 1(p_P>_Pp)\end{aligned}$$ 
is the summation in the Proposition. Then,
by definition and the convexity property of $p$, 
$\Big(I(Q_1)-J_0(Q_1)\Big)\cup 
\Big(I(Q_2)-J_\mu(Q_1)\Big)=K_p(Q_2)$ implies that 
$I(Q_1)=J_0(Q_1)$, or better $K_p^g(Q_2)=\emptyset.$ 
Consequently, this then shows that 
$\overline p^g\leq p$, a contradiction. This completes the proof of the 
Proposition and hence the Fundamental Relation.

\chapter{Non-Abelian L-Functions}
\section{Moduli Spaces as Integration Domains}

Let $F$ be a number field with $\Delta_F$ the absolute of its discriminant. 
Denote by $\mathbb A$ its ring of adeles. Fix a positive integer
$r\in {\mathbb Z}_{>0}$ and a convex polygon $p:[0,r]\to{\mathbb R}$. 
Consider the moduli space 
$${\mathcal M}_{F,r}^{\leq p}\Big[\Delta_F^{r\over 2}\Big]:=
\Big\{g\in GL_r(F)\backslash GL_r({\mathbb A}):
\mathrm{deg} (g)=-{r\over 2}\log\Delta_F,\ \overline p^g\leq p\Big\}.$$ 
Denote by $d\mu$ the induced Tamagawa measures on 
${\mathcal M}_{F,r}^{\leq p}\Big[\Delta_F^{r\over 2}\Big]$. 
For example,  in the case when $p=0$ and $F=\mathbb Q$, this moduli space 
coincides with the one defined using rank $r$ semi-stable lattices of volume 
one introduced in Chapter 1.

More generally, for any standard parabolic subgroup $P$ of $GL_r$, we 
introduce the moduli spaces 
$${\mathcal M}_{F,r}^{P;\leq p}\Big[\Delta_F^{r\over 2}\Big]
:=\Big\{g\in P(F)\backslash GL_r({\mathbb A}):
\mathrm{deg} (g)=-{r\over 2}\log\Delta_F,\ \overline p_P^g\leq p,
\ \overline p_P^g\geq -p\Big\}.$$ By the discussion in Chapter 2, all these 
moduli spaces 
${\mathcal M}_{F,r}^{P;\leq p}\Big[\Delta_F^{r\over 2}\Big]$ are  compact, 
a key property which plays a central role in our definition of non-abelian 
$L$-functions below.

\section{Choice of Eisenstein Series: First Approach to Non-Abelian 
$L$-Function}

To faciliate our ensuing discussion, we start with  some preperations. We 
follow [MW] closely. So from now on until the end of this chapter, 
a new system of notations will be used for our own convenience.

Fix a connected reduction group $G$ defined over $F$, denote by $Z_G$ its 
center. Fix a minimal parabolic subgroup $P_0$ of $G$. Then $P_0=M_0U_0$, 
where as usual we fix once and for all the Levi $M_0$ and  the unipotent
 radical $U_0$. A parabolic subgroup $P$ of $G$ is called {\it standard} if 
$P\supset P_0$. For such groups write $P=MU$ with $M_0\subset M$ the 
{\it standard Levi} and $U$ the unipotent radical. Denote by 
$\mathrm {Rat}(M)$ the group of rational characters of $M$, i.e, the morphism 
$M\to {\mathbb G}_m$  where ${\mathbb G}_m=\mathrm{GL}_1$ denotes the 
multiplicative group. Set 
$$\frak a_M^*:=\mathrm {Rat}(M)\otimes_{\mathbb Z}{\mathbb C},\qquad \frak a_M
:=\mathrm{Hom}_{\mathbb Z}\Big(\mathrm {Rat}(M),{\mathbb C}\Big),$$ and 
$$\mathrm {Re}\frak 
a_M^*:=\mathrm {Rat}(M)\otimes_{\mathbb Z}{\mathbb R},\qquad 
\mathrm{Re}\frak a_M
:=\mathrm{Hom}_{\mathbb Z}\Big(\mathrm {Rat}(M),{\mathbb R}\Big).$$ 
For any $\chi\in 
\mathrm {Rat}(M)$, we obtain a (real) character $|\chi|:M({\mathbb A})\to 
{\mathbb R}^*$ defined by $m=(m_v)\mapsto m^{|\chi|}:=\prod_{v\in S}
|m_v|_v^{\chi_v}$ with $|\cdot|_v$ the $v$-absolute values. Set then 
$M({\mathbb A})^1:=\cap_{\chi\in \mathrm {Rat}(M)}\mathrm{Ker}|\chi|$, 
which is 
a normal subgroup of $M({\mathbb A})$. Set $X_M$ to be the group of complex 
characters which are trivial on $M({\mathbb A})^1$. Denote by 
$H_M:=\log_M:M({\mathbb A})\to \frak a_M$ the map such that 
for all $\chi\in \mathrm {Rat}(M)\subset \frak a_M^*,$
$\displaystyle{\langle\chi,\log_M(m)\rangle:=\log(m^{|\chi|})}$. 
Clearly, $$M({\mathbb A})^1
=\mathrm{Ker}(\log_M);\qquad \log_M\Big(M({\mathbb A})
\Big/M({\mathbb A})^1\Big)\simeq 
\mathrm{Re}\frak a_M.$$ Hence in particular there is a natural isomorphism 
$\kappa:\frak a_M^*\simeq X_M.$
Set $$\mathrm{Re}X_M:=\kappa \Big(\mathrm{Re}\frak a_M^*\Big),
\qquad \mathrm{Im}X_M:=\kappa \Big(i\cdot \mathrm{Re}\frak a_M^*\Big).$$ 
Moreover define our working space 
$X_M^G$ to be the subgroup of $X_M$ consisting of complex characters of 
$M({\mathbb A})\Big/M({\mathbb A})^1$ which are trivial on 
$Z_{G({\mathbb A})}$.
\vskip 0.30cm
Fix a maximal compact subgroup ${\mathbb K}$ such that for all standard 
parabolic subgroups $P=MU$ as above, 
$P({\mathbb A})\cap{\mathbb K}=\Big(M({\mathbb A})
\cap{\mathbb K}\Big)\cdot\Big(U({\mathbb A})\cap{\mathbb K}\Big).$ 
Hence we get the Langlands decomposition 
$G({\mathbb A})=M({\mathbb A})\cdot U({\mathbb A})\cdot {\mathbb K}$. 
Denote by 
$m_P:G({\mathbb A})\to M({\mathbb A})/M({\mathbb A})^1$ the map 
$g=m\cdot n\cdot
 k\mapsto M({\mathbb A})^1\cdot m$ where $g\in G({\mathbb A}), m\in 
M({\mathbb A}), 
n\in U({\mathbb A})$ and 
$k\in {\mathbb K}$.
\vskip 0.30cm 
Fix Haar measures on $M_0({\mathbb A}), U_0({\mathbb A}), {\mathbb K}$ 
respectively 
such that 

\noindent
(1) induced measure on $M(F)$ is the counting measure and the volume of
 the induced measure on $M(F)\backslash M({\mathbb A})^1$ is 1. 
(It is  a fundamental fact that $M(F)\backslash M({\mathbb A})^1$ is 
compact.)

\noindent
(2) induced measure on $U_0(F)$ is the counting measure and the volume 
of $U_0(F)\backslash U_0({\mathbb A})$ is 1. (Being unipotent 
radical, $U_0(F)\backslash U_0({\mathbb A})$ is compact.)

\noindent
(3) the volume of ${\mathbb K}$ is 1.

Such measures then also induce Haar measures via $\log_M$ to $\frak a_{M_0}, 
\frak a_{M_0}^*$, etc. Furthermore, if we denote by $\rho_0$ the half of the 
sum of the positive roots of  the maximal split torus $T_0$ of the central 
$Z_{M_0}$
of $M_0$, then $$f\mapsto \int_{M_0({\mathbb A})\cdot U_0({\mathbb A})\cdot 
{\mathbb K}}f(mnk)\,dk\,dn\,m^{-2\rho_0}dm$$ defined for continuous functions 
with compact supports on $G({\mathbb A})$  defines a Haar measure $dg$ on 
$G({\mathbb A})$. This in turn gives measures on $M({\mathbb A}), 
U({\mathbb A})$ 
and hence on $\frak a_{M}, \frak a_{M}^*$, $P({\mathbb A})$, etc, for all 
parabolic subgroups $P$. In particular, one checks that the following 
{\it compactibility condition} holds
$$\int_{M_0({\mathbb A})\cdot U_0({\mathbb A})\cdot {\mathbb K}}f(mnk)\,
dk\,dn\,
m^{-2\rho_0}dm=\int_{M({\mathbb A})\cdot U({\mathbb A})\cdot 
{\mathbb K}}f(mnk)\,dk\,
dn\,m^{-2\rho_P}dm$$ for all continuous functions $f$ with compact supports 
on $G({\mathbb A})$, where $\rho_P$ denotes the half of the sum of the 
positive 
roots of  the maximal split torus $T_P$ of the central $Z_{M}$ of $M$. 
For later use, 
denote also by $\Delta_P$ the set of positive roots determined by $(P,T_P)$ 
and $\Delta_0=\Delta_{P_0}$.
\vskip 0.30cm
Fix an isomorphism $T_0\simeq {\mathbb G}_m^R$. Embed ${\mathbb R}_+^*$ 
by the map 
$t\mapsto (1;t)$. That is, at all finite places, the local components are 1, 
and at all infinite places, the local components are all $t$.  
Then we obtain a natural injection $({\mathbb R}_+^*)^R
\hookrightarrow T_0({\mathbb A})$ which splits. Denote by 
$A_{M_0({\mathbb A})}$ 
the unique connected subgroup of $T_0({\mathbb A})$ which projects onto 
$({\mathbb R}_+^*)^R$. More generally, for a standard parabolic subgroup 
$P=MU$, set $A_{M({\mathbb A})}:=
A_{M_0({\mathbb A})}\cap Z_{M({\mathbb A})}$ where as used above $Z_*$ denotes 
the center of the group $*$. Clearly, $M({\mathbb A})=A_{M({\mathbb A})}\cdot 
M({\mathbb A})^1$. For later use, set also 
$A_{M({\mathbb A})}^G:=\{a\in A_{M({\mathbb A})}:\log_Ga=0\}.$ 
Then $A_{M({\mathbb A})}=A_{G({\mathbb A})}\oplus
A_{M({\mathbb A})}^G.$
\vskip 0.30cm
Note that ${\mathbb K}$, $M(F)\backslash M({\mathbb A})^1$ and $U(F)\backslash 
U({\mathbb A})$ are all compact, thus with the Langlands decomposition 
$G({\mathbb A})=U({\mathbb A})M({\mathbb A}){\mathbb K}$ in mind, the 
reduction theory 
for $G(F)\backslash G({\mathbb A})$ or more generally 
$P(F)\backslash G({\mathbb A})$ is reduced to that for 
$A_{M({\mathbb A})}$ since $\Big(Z_G(F)\cap Z_{G({\mathbb A})}\Big)
\backslash \Big(Z_{G({\mathbb A})}\cap G({\mathbb A})^1\Big)$ 
is compact as well. As such for $t_0\in M_0({\mathbb A})$ set
$$A_{M_0({\mathbb A})}(t_0):=
\Big\{a\in A_{M_0({\mathbb A})}:a^\alpha>t_0^\alpha\,
\forall\alpha\in\Delta_0\Big\}.$$ Then, for a fixed compact subset 
$\omega\subset P_0({\mathbb A})$, we have the corresponding {\it Siegel set} 
$$S(\omega;t_0):=\Big\{p\cdot a\cdot k:p\in \omega, a\in A_{M_0({\mathbb A})}
(t_0),k\in {\mathbb K}\Big\}.$$ In 
particular, for big enough $\omega$ and small 
enough $t_0$, i.e, $t_0^\alpha$ is very close to 0 for all 
$\alpha\in\Delta_0$, the classical reduction theory may be restated as 
$G({\mathbb A})=G(F)\cdot S(\omega;t_0)$. More generally, set 
$$A_{M_0({\mathbb A})}^P(t_0):=
\Big\{a\in A_{M_0({\mathbb A})}:a^\alpha>t_0^\alpha\,
\forall\alpha\in\Delta_0^P\Big\},$$ and $$S^P(\omega;t_0):=
\Big\{p\cdot a\cdot k:
p\in \omega, a\in A_{M_0({\mathbb A})}^P(t_0),k\in {\mathbb K}\Big\}.$$  
Then similarly as above
for big enough $\omega$ and small enough $t_0$, $G({\mathbb A})=P(F)\cdot 
S^P(\omega;t_0)$. (Here 
$\Delta_0^P$ denotes the set of positive roots for $(P_0\cap M,T_0)$.)
\vskip 0.30cm
Fix an embedding $i_G:G\hookrightarrow SL_n$ sending $g$ to $(g_{ij})$. 
Introducing a hight function on $G({\mathbb A})$ by setting 
$\|g\|:=\prod_{v\in S}\mathrm{sup}\Big\{|g_{ij}|_v:\forall i,j\Big\}$. It is 
well-known 
that up to $O(1)$, hight functions are unique. This implies that the following 
growth conditions do not depend on the height function we choose.

A function $f:G({\mathbb A})\to {\mathbb C}$ is said to have {\it moderate 
growth} if  there exist 
$c,r\in {\mathbb R}$ such that $\Big|f(g)\Big|\leq c\cdot \Big\|g\Big\|^r$ 
for all $g\in G({\mathbb A})$. Similarly, for a standarde parabolic 
subgroup $P=MU$, 
a function $f:U({\mathbb A})M(F)\backslash G({\mathbb A})\to {\mathbb C}$ 
is said 
to have {\it moderate growth} if there exist $c,r\in {\mathbb R}$ and
$\lambda\in \mathrm{Re}X_{M_0}$ such that for any 
$a\in A_{M({\mathbb A})},k\in 
{\mathbb K}, 
m\in M({\mathbb A})^1\cap S^P(\omega;t_0)$,
$$\Big|f(amk)\Big|\leq c\cdot \Big\|a\Big\|^r\cdot m_{P_0}(m)^\lambda.$$

Also a function $f:G({\mathbb A})\to {\mathbb C}$ is said to be smooth 
if for any 
$g=g_f\cdot g_\infty \in G({\mathbb A}_f)\times G({\mathbb A}_\infty)$, 
there exist 
open neighborhoods $V_*$ of $g_*$ in $G({\mathbb A})$ and a 
$C^\infty$-function 
$f':V_\infty\to {\mathbb C}$ such that $f(g_f'\cdot g_\infty')=
f'(g_\infty')$ for 
all $g_f'\in V_f$ and $g_\infty'\in V_\infty$.
\vskip 0.30cm
By contrast, a function $f: S(\omega;t_0)\to {\mathbb C}$ is said to 
be {\it rapidly 
decreasing} if there exists $r>0$ and for all $\lambda\in \mathrm{Re}X_{M_0}$ 
there exists $c>0$ such that for $a\in A_{M({\mathbb A})}, g\in 
G({\mathbb A})^1\cap 
S(\omega;t_0)$, $\Big|\phi(ag)\Big|\leq c
\cdot\Big\|a\Big\|\cdot m_{P_0}(g)^\lambda$. And a 
function $f:G(F)\backslash G({\mathbb A})\to {\mathbb C}$ is said to 
be {\it rapidly decreasing} if $f|_{S(\omega;t_0)}$ is so.
\vskip 0.30cm
By definition,  a function 
$\phi:U({\mathbb A})M(F)\backslash G({\mathbb A})\to {\mathbb C}$ is 
called {\it automorphic} (of level $P$) if

\noindent
(i) $\phi$ has moderate growth;

\noindent
(ii) $\phi$ is smooth;

\noindent
(iii) $\phi$ is ${\mathbb K}$-finite, i.e, the ${\mathbb C}$-span of all 
$\phi(k_1\cdot *\cdot k_2)$ parametrized by $(k_1,k_2)\in {\mathbb K}\times 
{\mathbb K}$ is finite dimensional; and 

\noindent
(iv) $\phi$ is $\frak z$-finite, i.e, the ${\mathbb C}$-span of all 
$\delta(X)\phi$ parametrized by all $X\in \frak z$ is finite dimensional. Here 
$\frak z$ denotes the center of the universal enveloping algebra 
$\frak u:=\frak U(\mathrm{Lie}G({\mathbb A}_\infty))$ of the Lie algebra of  
$G({\mathbb A}_\infty)$ and $\delta(X)$
denotes the derivative of $\phi$ along $X$.

For such a function $\phi$, set $\phi_k:M(F)\backslash M({\mathbb A})\to 
{\mathbb C}$ by $m\mapsto m^{-\rho_P}\phi(mk)$ for all $k\in {\mathbb K}$. 
Then one checks that $\phi_k$ is an automorphic form in the usual sense. Set
 $A\Big(U({\mathbb A})M(F)\backslash G({\mathbb A})\Big)$ be the space of 
automorphic 
forms on $U({\mathbb A})M(F)\backslash G({\mathbb A})$.
\vskip 0.30cm
For a measurable locally $L^1$-function 
$f:U(F)\backslash G({\mathbb A})\to 
{\mathbb C}$ define its {\it constant term} along with the standard parabolic 
subgroup 
$P=UM$ to be the function $f_P:U({\mathbb A})\backslash G({\mathbb A})
\to {\mathbb C}$ given by
$g\mapsto\int_{U(F)\backslash G({\mathbb A})}f(ng)\,dn.$ 
Then an automorphic form $\phi\in A\Big(U({\mathbb A})M(F)\backslash 
G({\mathbb A})\Big)$ is 
called a {\it cusp form} if for any standard parabolic subgroup $P'$ properly 
contained in $P$, $\phi_{P'}\equiv 0$. Denote by
 $A_0\Big(U({\mathbb A})M(F)\backslash G({\mathbb A})\Big)$  
the space of cusp forms on $U({\mathbb A})M(F)\backslash G({\mathbb A})$. 
One checks easily that

\noindent
(i) all cusp forms are rapidly decreasing; and hence

\noindent
(ii) there is a natural pairing
 $$\Big\langle\cdot,\cdot\Big\rangle:A_0(U({\mathbb A})M(F)\backslash 
G({\mathbb A}))\times 
A\Big(U({\mathbb A})M(F)\backslash G({\mathbb A})\Big)\to 
{\mathbb C}$$ defined by
$\Big\langle \psi,\phi\Big\rangle:=\int_{Z_{M({\mathbb A})}U({\mathbb A})
M(F)\backslash 
G({\mathbb A})}\psi(g)\overline{\phi(g)}\,dg.$

Moreover, for a (complex) character $\xi:Z_{M({\mathbb A})}\to 
{\mathbb C}^*$ of 
$Z_{M({\mathbb A})}$ set 
$$\begin{aligned}A&\Big(U({\mathbb A})M(F)\backslash G({\mathbb A})
\Big)_\xi\\
:=&\Big\{\phi\in A\Big(U({\mathbb A})M(F)\backslash G({\mathbb A})\Big):
\phi(zg)=z^{\rho_P}\cdot\xi(z)\cdot\phi(g),\forall 
z\in Z_{M({\mathbb A})}, g\in G({\mathbb A})\Big\}\end{aligned}$$ and 
$$A_0\Big(U({\mathbb A})M(F)\backslash G({\mathbb A})\Big)_\xi:=
A_0\Big(U({\mathbb A})M(F)\backslash G({\mathbb A})\Big)\bigcap 
A\Big(U({\mathbb A})M(F)\backslash G({\mathbb A})\Big)_\xi.$$
 
More globally, set 
$$A\Big(U({\mathbb A})M(F)\backslash G({\mathbb A})\Big)_Z:=
\sum_{\xi\in \mathrm{Hom}\Big(Z_{M({\mathbb A})},{\mathbb C}^*\Big)}
A\Big(U({\mathbb A})M(F)\backslash 
G({\mathbb A})\Big)_\xi$$ and
$$A_0\Big(U({\mathbb A})M(F)\backslash G({\mathbb A})\Big)_Z:=
\sum_{\xi\in \mathrm{Hom}\Big(Z_{M({\mathbb A})},{\mathbb C}^*\Big)}
A_0\Big(U({\mathbb A})M(F)\backslash G({\mathbb A})\Big)_\xi.$$ 
One checks that the natural morphism $${\mathbb C}
\Big[\mathrm{Re}\frak a_M\Big]\bigotimes 
A\Big(U({\mathbb A})M(F)\backslash G({\mathbb A})\Big)_Z
\to A\Big(U({\mathbb A})M(F)\backslash G({\mathbb A})\Big)$$ defined by 
$(Q,\phi)\mapsto \Big(g\mapsto Q(\log_M(m_P(g))\Big)\cdot \phi(g)$ is an 
isomorphism, using the special structure of $A_{M({\mathbb A})}$-finite 
functions and the Fourier analysis over the compact space $A_{M({\mathbb A})}
\backslash Z_{M({\mathbb A})}$. Consequently, we also obtain a natural 
isomorphism 
$${\mathbb C}\Big[\mathrm{Re}\frak a_M\Big]\bigotimes 
A_0\Big(U({\mathbb A})M(F)\backslash G({\mathbb A})\Big)_Z\to 
A_0\Big(U({\mathbb A})M(F)\backslash G({\mathbb A})\Big).$$
In such a way, we may trace back where the automorphic forms at hands come.

Set also $\Pi_0\Big(M({\mathbb A})\Big)_\xi$ be isomorphism classes of 
irreducible representations of $M({\mathbb A})$ occuring in the space 
$A_0\Big(M(F)\backslash M({\mathbb A})\Big)_\xi$, and
$$\Pi_0\Big(M({\mathbb A})\Big):=
\bigcup_{\xi\in \mathrm{Hom}(Z_{M({\mathbb A})},{\mathbb C}^*)}
\Pi_0\Big(M({\mathbb A})\Big)_\xi.$$ (More precisely, we should 
use $M({\mathbb A}_f)\times \Big(M({\mathbb A})\cap {\mathbb K},
\mathrm{Lie}(M({\mathbb A}_\infty))\otimes_{\mathbb R}{\mathbb C}\Big)$ 
instead of 
$M({\mathbb A})$.) For any 
$\pi\in \Pi_0\Big(M({\mathbb A})\Big)_\xi$ set 
$A_0\Big(M(F)\backslash M({\mathbb A})\Big)_\pi$ 
to 
be the isotypic component of type $\pi$ of $A_0\Big(M(F)\backslash 
M({\mathbb A})\Big)_\xi$, i.e, the set of these
cusp forms of $M({\mathbb A})$ which generate the 
semi-simple isotypic $M({\mathbb A}_f)\times 
\Big(M({\mathbb A})\cap {\mathbb K},
\mathrm{Lie}(M({\mathbb A}_\infty))\otimes_{\mathbb R}{\mathbb C}\Big)$-module 
of type $\pi$. Set $$\begin{aligned}
A_0&\Big(U({\mathbb A})M(F)\backslash G({\mathbb A})\Big)_\pi\\
:=&\Big\{\phi\in A_0\Big(U({\mathbb A})M(F)\backslash G({\mathbb A})\Big):
\phi_k\in A_0\Big(M(F)\backslash M({\mathbb A})\Big)_\pi,\,\forall k\in 
{\mathbb K}\Big\}.\end{aligned}$$ Clearly
$$A_0\Big(U({\mathbb A})M(F)\backslash G({\mathbb A})\Big)_\xi
=\bigoplus_{\pi\in
\Pi_0\big(M({\mathbb A})\big)_\xi} A_0\Big(U({\mathbb A})M(F)\backslash 
G({\mathbb A})\Big)_\pi.$$

More generally, let $V\subset A\Big(M(F)\backslash M({\mathbb A})\Big)$ 
be an irreducible 
$M({\mathbb A}_f)\times \Big(M({\mathbb A})\cap {\mathbb K},\mathrm{Lie}
(M({\mathbb A}_\infty))
\otimes_{\mathbb R}{\mathbb C}\Big)$-module with $\pi_0$ the induced 
representation of 
$M({\mathbb A}_f)\times \Big(M({\mathbb A})\cap {\mathbb K},\mathrm{Lie}
(M({\mathbb A}_\infty))
\otimes_{\mathbb R}{\mathbb C}\Big)$. Then we call $\pi_0$ an {\it automorphic 
representation} of $M({\mathbb A})$. Denote by $A\Big(M(F)\backslash 
M({\mathbb A}\Big)_{\pi_0}$  the isotypic 
subquotient module of type $\pi_0$ of 
$A\Big(M(F)\backslash M({\mathbb A}\Big)$. One checks that
$$V\otimes \mathrm{Hom}_{M({\mathbb A}_f)\times (M({\mathbb A})\cap
 {\mathbb K},
\mathrm{Lie}(M({\mathbb A}_\infty))\otimes_{\mathbb R}{\mathbb C})}
\Big(V,A\big(M(F)\backslash M({\mathbb A})\big)\Big)\simeq 
A\Big(M(F)\backslash M({\mathbb A})\Big)_{\pi_0}.$$ 
Set $$\begin{aligned}A&\Big(U({\mathbb A})M(F)\backslash 
G({\mathbb A})\Big)_{\pi_0}\\
:=&\Big\{\phi\in 
A\Big(U({\mathbb A})M(F)\backslash G({\mathbb A})\Big):\phi_k\in A(M(F)
\backslash M({\mathbb A}))_{\pi_0},\forall k\in {\mathbb K}\Big\}.
\end{aligned}$$
Moreover if 
$A\Big(M(F)\backslash M({\mathbb A})\Big)_{\pi_0}\subset
A_0\Big(M(F)\backslash M({\mathbb A})\Big)$, we call $\pi_0$ a {\it 
cuspidal representation}.
\vskip 0.30cm
Two automorphic representations $\pi$ and $\pi_0$ of $M({\mathbb A})$ are said 
to be {\it equivalent} if there exists $\lambda\in X_M^G$ such that $\pi\simeq 
\pi_0\otimes\lambda$. This, in practice, simply means that the following
relation holds: $A\Big(M(F)\backslash M({\mathbb A})\Big)_\pi
=\lambda\cdot A\Big(M(F)\backslash M({\mathbb A})\Big)_{\pi_0}.$ 
That is for any 
$\phi_\pi\in A\Big(M(F)\backslash M({\mathbb A})\Big)_\pi$ there exists a 
$\phi_{\pi_0}\in A\Big(M(F)\backslash M({\mathbb A})\Big)_{\pi_0}$ 
such that 
$\phi_\pi(m)=m^\lambda\cdot \phi_{\pi_0}(m)$. Consequently, 
$$A\Big(U({\mathbb A})M(F)\backslash G({\mathbb A})\Big)_\pi=
\Big(\lambda\circ m_P\Big)\cdot 
A\Big(U({\mathbb A})M(F)\backslash G({\mathbb A})\Big)_{\pi_0}.$$ 
Denote by 
$\varpi:=[\pi_0]$ the equivalence class of $\pi_0$. Then $\varpi$ is an 
$X_M^G$-principal homogeneous space, hence admits a natural complex structure. 
Usually we call $(M,\varpi)$ a {\it cuspidal datum} 
of $G$ if $\pi_0$ is cuspidal. 
Also for $\pi\in\varpi$, 
set $$\mathrm{Re}\pi:=\mathrm{Re}\chi_\pi=|\chi_\pi|\in 
\mathrm{Re}X_M,\qquad \mathrm{Im}\pi:=\pi\otimes(-\mathrm{Re}\pi),$$ where 
$\chi_\pi$ is the central character of $\pi$.
\vskip 0.30cm
Now fix an irreducible automorphic representation $\pi$ of $M({\mathbb A})$ and
for an automorphic form
$\phi\in A\Big(U({\mathbb A})M(F)\backslash G({\mathbb A})\Big)_\pi$, define the  associated {\it Eisenstein series} 
$E(\phi,\pi):G(F)\backslash G({\mathbb A})\to 
{\mathbb C}$ by $$E(\phi,\pi)(g):=\sum_{\delta\in P(F)\backslash G(F)}\phi
(\delta g).$$ Then one checks that there is an open cone ${\mathcal C}\subset 
\mathrm{Re}X_M^G$ such that if $\mathrm{Re}\pi\in {\mathcal C}$, 
the Eisenstein series $E(\lambda\cdot \phi,
\pi\otimes\lambda)(g)$ converges uniformly for $g$ in a compact subset of 
$G({\mathbb A})$ and $\lambda$ in an open neighborhood of 0 in $X_M^G$. For 
example, if $\varpi=[\pi]$ is cuspidal, we may even take ${\mathcal C}$ to be 
the cone $\Big\{\lambda\in \mathrm{Re}X_M^G:\langle\lambda-\rho_P,
\alpha^\vee\rangle>0,\forall\alpha\in \Delta_P^G\Big\}$.
As a direct consequence, then $E(\phi,\pi)\in A\Big(G(F)\backslash 
G({\mathbb A})\Big)$. That is, it is an automorphic form of level $G$.
\vskip 0.30cm
As noticed above, being an automorphic form, $E(\phi,\pi)$ is of moderate 
growth. However, in general it is not integrable over 
$Z_{G({\mathbb A})}G(F)\backslash 
G({\mathbb A})$. To remedy this, classically, as initiated in the so-called 
Rankin-Selberg method, analytic truncation is used: From Fourier analysis, 
we understand 
that the probelmatic terms are the so-called constant terms, which are of 
slow growth only. So by cutting off these constant terms suitably, 
the reminding one is expacted to be rapidly decreasing, and hence integrable.
(Yet, in general, it is very difficult to make such an analytic truncation 
intrinsically related with arithmetic properties of number fields. 
See however, the Rankin-Selberg method [Bu], [Z] and the Arthur-Selberg 
trace formula [Ar1-4].) 

Furthermore, note that Eisenstein series themselves are quite intrinsic 
arithmetical invariants.
Thus it is natural for us on one hand to keep Eisenstein series unchanged 
while on the other to find new moduli spaces, which themselves are 
intrinsically parametrized certain 
modular objects, and over which Eisenstein series are integrable.

This then naturally leads to our non-abelian $L$-functions: 
As said, we are going to view Eisenstein
 series as something intrinsically defined. 
In contrast, using a geo-arithmetical truncation 
for the space $G(F)\backslash G({\mathbb A})$,  we can make 
the integrations of Eisenstein series (over the newly obtained compact
moduli spaces) well-defined.

More precisely, let us return to the group $G=GL_r$. Recall that we 
obtain the moduli space ${\mathcal M}_{F,r}^{\leq p}
\Big[\Delta_F^{r\over 2}\Big]$ 
and hence also a well-defined integration $$L_{F,r}^{\leq p}(\phi,\pi):=
\int_{{\mathcal M}_{F,r}^{\leq p}[\Delta_F^{r\over 2}]}E(\phi,\pi)(g)\,dg,
\qquad \mathrm {Re}\pi\in {\mathcal C}.$$

\section{New Non-Abelian $L$-Functions}

In general, however, we do not know whether the above  
defined integrations satisfy nice properties such as meromorphic continuation 
and functional equations etc... It is to remedy this that 
we make a further choice of automorphic forms.

For $G=GL_r$, fix a 
standard parabolic subgroup $P_I=U_IM_I$ corresponding to the partition 
$I=(r_1,\cdots,r_{|P|})$ of $r$ with $M_I$ the standard Levi and $U_I$ the 
unipotent radical. Then for a fixed irreducible automorphic representation 
$\pi$ of $M_I({\mathbb A})$, choose $$\begin{aligned}\phi\in 
A&\Big(U_I({\mathbb A})M_I(F)\backslash 
G({\mathbb A})\Big)_\pi\cap 
L^2\Big(U_I({\mathbb A})M_I(F)\backslash G({\mathbb A})\Big)\\
:=&A^2\Big(U_I({\mathbb A})M_I(F)
\backslash G({\mathbb A})\Big)_\pi,\end{aligned}$$ where 
$L^2\Big(U_I({\mathbb A})M_I(F)\backslash 
G({\mathbb A})\Big)$ denotes the space of $L^2$ functions on the quotient
space  $Z_{G({\mathbb A})}
U_I({\mathbb A})M_I(F)\backslash G({\mathbb A})$. Then we have the 
associated Eisenstein 
series  $E(\phi,\pi)\in A\Big(G(F)\backslash G({\mathbb A})\Big)$.
\vskip 0.30cm
\noindent
{\bf {\large Main Definition A.}} {\it A rank $r$ 
non-abelian $L$-function $L_{F,r}^{\leq p}(\phi,\pi)$ for the number field $F$
associated to an
$L^2$-automorphic form $\phi\in A^2(U_I({\mathbb A})M_I(F)\backslash 
G({\mathbb A}))_\pi$ is defined to be the integration
$$L_{F,r}^{\leq p}(\phi,\pi):=
\int_{{\mathcal M}_{F,r}^{\leq p}[\Delta_F^{r\over 2}]}E(\phi,\pi)(g)\,dg,
\qquad \mathrm {Re}\,\pi\in {\mathcal C}.$$}

More generally, for any standard parabolic subgroup $P_J=U_JM_J\supset P_I$ 
(so that the partition $J$ is a refinement of $I$), we have the corresponding 
relative Eisenstein series $$E_I^J(\phi,\pi)(g):=\sum_{\delta\in P_I(F)
\backslash P_J(F)}\phi(\delta g),\qquad\forall g\in P_J(F)\backslash 
G({\mathbb A}).$$ 
It is well-known that there is  an open cone ${\mathcal C}_I^J$
in $\mathrm {Re}X_{M_I}^{P_J}$ such that for $\mathrm {Re}\pi\in 
{\mathcal C}_I^J$, 
$E_I^J(\phi,\pi)\in A\Big(P_J(F)\backslash G({\mathbb A})\Big).$ 
Here $X_{M_I}^{P_J}$ 
is defined similarly as $X_M^G$ with $G$ replaced by $P_J$. Then we have 
a well-defined {\it relative non-abelian $L$-function}
$$L_{F,r}^{P_J;\leq p}(\phi,\pi):=
\int_{{\mathcal M}_{F,r}^{P_J;\leq p}[\Delta_F^{r\over 2}]}E_I^J
(\phi,\pi)(g)\,dg,
\qquad \mathrm {Re}\pi\in {\mathcal C}_I^J.$$ 
 \vskip 0.30cm
\noindent
{\bf Remarks.} (1) Here when defining non-abelian $L$-functions we assume that 
$\phi$ comes from a single irreducible automorphic representations. But this 
restriction is rather artifical and can be removed easily. However, to present
our constructions and results in a very neat way, we in the sequel will keep 
using it. 

\noindent
(2) The discussion for non-abelian 
$L$-functions holds  for the just defined relative
 non-abelian $L$-functions as well. So from now on, we will leave such a 
 modification to the reader while concentrate ourselves only on non-abelian 
$L$-functions.
\section{Meromorphic Extension and Functional Equations}

With the same notation as above, set $\varpi=[\pi]$. For $w\in W$ the Weyl 
group of $G$, fix once and for all representative $w\in G(F)$ of $w$. Set 
$M':=wMw^{-1}$ and denote the associated parabolic subgroup by $P'=U'M'$. 
$W$ acts naturally on  automorphic representations, from which we obtain 
an equivalence
classes $w{\varpi}$ of automorphic representations of $M'({\mathbb A})$. As 
usual, define the associated {\it intertwining operator} $M(w,\pi)$ by 
$$\Big(M(w,\pi)\phi\Big)(g):=
\int_{U'(F)\cap wU(F)w^{-1}\backslash U'({\mathbb A})}
\phi(w^{-1}n'g)\,dn',\qquad
\forall g\in G({\mathbb A}).$$ One checks that if $\langle \mathrm {Re}\pi,
\alpha^\vee\rangle\gg 0,\forall \alpha\in\Delta_P^G$,

\noindent
(i) for a fixed $\phi$, $M(w,\pi)\phi$ depends only on the double coset 
$M'(F)wM(F)$. So $M(w,\pi)\phi$ is well-defined for $w\in W$;

\noindent
(ii) the above integral converges absolutely and uniformly for $g$ varying 
in a compact subset of $G({\mathbb A})$;

\noindent
(iii) $M(w,\pi)\phi\in A\Big(U'({\mathbb A})M'(F)\backslash 
G({\mathbb A})\Big)_{w\pi}$;
and if $\phi$ is $L^2$, which from now on we always assume, so is  
$M(w,\pi)\phi$.
\vskip 0.30cm
\noindent
{\bf {\large Basic Facts of Non-Abelian $L$-Functions.}} 
{\it With the same notation above,} 

\noindent
(I) ({\bf Meromorphic Continuation}) {\it $L_{F,r}^{\leq p}(\phi,\pi)$ for 
$\mathrm{Re}\pi\in {\mathcal C}$ is well-defined and admits a unique 
meromorphic continuation to the whole space $\varpi$;}

\noindent
(II) ({\bf Functional Equations}) {\it As meromorphic functions on $\varpi$,
$$L_{F,r}^{\leq p}\Big(\phi,\pi\Big)=
L_{F,r}^{\leq p}\Big(M(w,\pi)\phi,w\pi\Big),\qquad \forall w\in W.$$}
 
\noindent
{\it Proof.} This is a direct consequence of the fundamental 
results of 
Langlands on Eisenstein series and spectrum decompositions. [See e.g, [Ar1], 
[La1,2] and [MW1]).

Indeed, if $\phi$ is cuspidal, by definition, (I) is a direct 
consequence of Prop. II.15, Thm. IV.1.8 of [MW] and (II) is a direct 
consequence of Thm. IV.1.10 of [MW].

More generally, if $\phi$ is only $L^2$, then by Langlands' theory of 
Eisenstein series and spectral decomposition, $\phi$ may be obtained as 
the residue of
relative Eisenstein series coming from cuspidal forms, since $\phi$ is $L^2$ 
automorphic. As such then (I) and (II) are direct consequences of the proof 
of VI.2.1(i) at p.264 of [MW].
\section{Holomorphicity and Singularities}

Let $\pi\in \varpi$ and $\alpha\in\Delta_M^G$. Define the function 
$h:\varpi\to {\mathbb C}$ by $\pi\otimes\lambda\mapsto \langle \lambda,
\alpha^\vee\rangle,\forall\lambda\in X_M^G\simeq \frak a_M^G$. Here as usual, 
$\alpha^\vee$ denotes the coroot associated to $\alpha$. Set 
$H:=\Big\{\pi'\in\varpi:h(\pi')=0\Big\}$ 
and call it a root hyperplane. Clearly the 
function $h$ is determined by $H$, hence we also denote  $h$ by $h_H$. 
Note also that root hyperplanes depend on the base point $\pi$ we choose.

Let $D$ be a set of root hyperplanes. Then 

\noindent
(i) the singularities of a meromorphic function $f$ on $\varpi$ is said to 
lie on $D$ if for all $\pi\in\varpi$, there exist $n_\pi:D\to 
{\mathbb Z}_{\geq 0}$ zero almost everywhere such that $\pi'\mapsto 
\Big(\Pi_{H\in D}h_H(\pi')^{n_\pi(H)}\Big)\cdot f(\pi')$ is holomorphic at
$\pi'$;

\noindent
(ii) the singularities of $f$ are said to be without multiplicity at $\pi$ if 
$n_\pi\in\{0,1\}$;

\noindent
(iii) $D$ is said to be locally finite,  if for any compact subset 
$B\subset\varpi$, $\Big\{H\in D:H\cap B\not=\emptyset\Big\}$ is finite.
\vskip 0.30cm
\noindent
{\bf {\large Basic Facts of Non-Abelian $L$-Functions}.} 
{\it With the same notation 
above,} 

\noindent
(III) ({\bf Holomorphicity}) (i) {\it When $\mathrm{Re}\pi\in {\mathcal C}$, 
$L_{F,r}^{\leq p}(\phi,\pi)$ is holomorphic;}

\noindent
(ii) {\it $L_{F,r}^{\leq p}(\phi,\pi)$ is holomorphic at $\pi$ where 
$\mathrm{Re}\pi=0$;}

\noindent
(IV) ({\bf Singularities}) {\it Assume further that $\phi$ is a cusp form. 
Then}

\noindent
(i) {\it There is a locally finite set of root hyperplanes 
$D$ such that the singularities of $L_{F,r}^{\leq p}(\phi,\pi)$ are lied on 
$D$;}

\noindent
(ii) {\it The singularities of $L_{F,r}^{\leq p}(\phi,\pi)$ are without 
multiplicities at $\pi$ if $\langle \mathrm{Re}\pi,\alpha^\vee\rangle\geq 0,
\forall \alpha\in\Delta_M^G$;}

\noindent 
(iii) {\it There are only finitely many of singular hyperplanes of 
$L_{F,r}^{\leq p}(\phi,\pi)$ which intersect $\Big\{\pi\in\varpi:
\langle\mathrm{Re}\pi,\alpha^\vee\rangle\geq 0,\forall\alpha\in\Delta_M
\Big\}$.}
\vskip 0.30cm
\noindent
{\it Proof.} As above, this is a direct consequence of the fundamental 
results of 
Langlands on Eisenstein series and spectrum decompositions. [See e.g, [Ar1], 
[La1,2] and [MW1]).
Indeed, if $\phi$ is a cusp form, (III.i) is a direct consequence of 
Lemma IV.1.7 of [MW], while  (III.ii)  and (IV) are direct consequence of
Prop. IV.1.11 of [MW].

In general when $\phi$ is only $L^2$ automorphic, then we have to use the 
theory of Langlands to realize $\phi$ as the residue of relative Eisenstein 
series defined using cusp forms. (See e.g., item (5) at p.198 and the second 
half part of p.232 of [MW].) 

As such, (III) and (IV) are direct consequence of the definition of 
residue datum and the compactibility between residue and Eisenstein series 
as stated
 for example under item (3) at p.263 of [MW]. 
\vskip 0.30cm
\noindent 
{\bf Remark.}  Since $G=GL_r$, one can write down the functional equations
concretely, and give a much more refined result about the singularities of
the non-abelian $L$-functions, for instance, with the 
use of [MW2] about the  residue of Eisenstein 
series. We discuss this elsewhere.

\chapter{Arthur's Analytic Truncation}
\section{Langlands' Combinatorial Lemma and Arthur's Partition}
\subsection{Height Functions}

We shall recall some basic properties of height functions associated to 
rational representations of $G$.

Let $V$ be a vector space defined over $F$. Suppose that 
$\Big\{v_1,v_2,\cdots,v_n\Big\}$ is an $F$-basis of $V(F)$. If 
$\xi_v\in V(F_v)$ and 
$$\xi_v=\sum_i\xi_v^i\,v_i,\qquad \xi_v^i\in F_v,$$
define $$\|\xi_v\|_v:=\begin{cases} \max_i|\xi_v^i|_v, & \text{if $v$ 
is finite}\\ \Big(\sum_i|\xi_{\mathbb R}^i|^2\Big)^{1\over 2},& 
\text{if $v=\mathbb R$}\\
\sum_i|\xi_{\mathbb C}^i|^2,& \text{if $v=\mathbb C$.}\end{cases}$$
An element $\xi=\prod_v\xi_v$ in $V(\mathbb A)$ is said to be {\it primitive} 
if $\|\xi_v\|_v=1$ for almost all $v$, in which case we set
$$\|\xi\|:=\prod_v\|\xi_v\|_v.$$ $\|\cdot\|$ is called the 
{\it height function} associated to the basis $\Big\{v_1,v_2,\cdots,v_n\Big\}$.
\vskip 0.30cm
Suppose that $\Lambda:G\to GL(V)$ is a homomorphism defined over $F$. 
Let $K_\Lambda$ be the group of elements $k\in K$ such that 
$\|\Lambda(k)v\|=\|v\|$ for any primitive $v\in V(\mathbb A)$. It is 
possible to choose the basis $\Big\{v_1,v_2,\cdots,v_n\Big\}$ such that 

\noindent
(i) $K_\Lambda$ is of finite index in $K$; and 

\noindent
(ii) for each $a\in A_0$, the operator $\Lambda(a)$ is diagonal.

\noindent
We shall always assume that for a given $\Lambda$, the basis has been chosen 
to satisfy these two conditions. From our basis on $V(F)$, we obtain a basis 
for the vector space of endomorphisms of $V(F)$. We have

\noindent
(a) {\it every element in $G(\mathbb A)$ is primitive with respect to the 
corresponding height function}; and

\noindent
(b) {\it for every primitive $v\in V(\mathbb A)$ and every $x\in G(\mathbb A)$,
$$\|\Lambda(x) v\|\leq \|\Lambda(x)\|\cdot\|v\|.$$}

If $t>0$, define $$G_t:=\Big\{x\in G(\mathbb A):\|\Lambda(x)\|\leq t\Big\}.$$
Suppose that $\Lambda$ has the further property that

\noindent
(c) {\it $G_t$ is compact for every $t$}.

\noindent
It is known that

\noindent
(c$'$) {\it there are constants $C$ and $N$ such that for any $t$, 
the volume of 
$G_t$ (with respect to our Haar measure) is bounded by $C t^N$.}

\noindent
For the rest of this paper we shall simply assume that {\it some $\Lambda$, 
satisfying this additional property, has been fixed}, and we shall write 
$\|x\|$ for $\|\Lambda(x)\|$.
\vskip 0.30cm
This \lq norm' function on $G(\mathbb A)$ satisfies the following properties:

\noindent
(d) $$\|x\|\geq 1;$$
$$\|k_1xk_2\|=\|x\|;$$
$$\|xy\|\leq \|x\|\cdot\|y\|;$$
and $$\|x^{-1}\|\leq C\|x\|^N$$
for constants $C$ and $N$, elements $x,y\in G(\mathbb A)$, and $k_1,k_2$ 
belonging to a subgroup of finite index in $K$.

Once $\|\cdot\|$ has been fixed, we shall want to consider different 
rational representations $\Lambda$ of $G$. In particular, suppose that 
{\it the highest weight of $\Lambda$ is $\lambda$}, for some element 
$\lambda$ in $\frak a_0^*$. Then

\noindent
(e) {\it there are constants $c_1$ and $c_2$ such that 
$$c_1e^{-\lambda(H_0(x))}\leq\|\Lambda(x)^{-1}v\|\leq 
c_2e^{-\lambda(H_0(x))},$$ for all $x\in G(\mathbb A)$.}

The point here is that uniformly on $x\in G(\mathbb A)$, 
$\|\Lambda(x)^{-1}v\|$is approximately of growth $e^{-\lambda(H_0(x))}$. 
By varying the linear functional $\lambda$, we can then show that

\noindent
(e$'$) {\it for any Euclidean norm $\|\cdot\|$ on $\frak a_0$, we can choose a 
constant $c$ so that $$\|H_0(x)\|\leq c\Big(1+\log\|x\|\Big),\qquad 
x\in G(\mathbb A).$$} 

Finally, we state the following non-trivial property whose proof may be 
found in [MW]:

\noindent
(f) {\it there is a constant $c$ such that for any $\delta\in G(F)$ and 
$x\in \frak s$, the Siegel domain,  $$\|x\|\leq c\|\delta x\|.$$}

\subsection{Partial Truncation and First Estimations}

If $P_1\subset P_2$ are two (standard) parabolic subgroups of $G$, 
following Arthur [Ar2], let $$\tau_{P_1}^{P_2}:=\tau_1^2\qquad\mathrm{and}
\qquad \widehat\tau_{P_1}^{P_2}:=\widehat\tau_1^2$$ be the characteristic 
functions on $\frak a_0$ of
$$\Big\{H\in \frak a_0:\alpha(H)>0,\ \alpha\in\Delta_1^2\Big\}$$ and
$$\Big\{H\in \frak a_0:\varpi(H)>0,\ \varpi\in\widehat\Delta_1^2\Big\}.$$
We shall denote $\tau_P^G$ and $\widehat\tau_P^G$ simply by $\tau_P$ and 
$\widehat\tau_P$.
\vskip 0.30cm
\noindent
{\bf Basic Estimation.} (Arthur)  Suppose that we are given a parabolic 
subgroup $P$, and a Euclidean norm $\|\cdot\|$ on $\frak a_P$. Then 
{\it there are constants $c$ and $N$ such that for all $x\in G(\mathbb A)^1$ 
and $X\in\frak a_P$,
$$\sum_{\delta\in P(F)\backslash G(F)}\widehat\tau_P\Big(H(\delta x)-X\Big)
\leq c\Big(\|x\|e^{\|X\|}\Big)^N.$$} Moreover, {\it the sum is finite}.

\noindent
{\it Proof.} The idea is to show that for fixed $x$ and $X$, the condition
that $\hat\tau_P\Big(H(\delta x)-X\Big)\not=0$ for $\delta\in P(F)$
forms a compact subset. Hence by the discreteness of $G(F)$, the sum above 
is a finite one. For this purpose, we go as follows (please pay special 
attention to inequalities (1), (2) and (3) below):

First, we prove the following

\noindent
{\bf Lemma.} {\it there is a constant $c$ such that 
$$\varpi\Big(H_0(\delta x)\Big)\leq c\Big(1+\log\|x\|\Big)\eqno(1)$$ 
for all $\varpi\in\widehat\Delta_0,\ x\in G(\mathbb A)^1$ and 
$\delta\in G(F)$.}

\noindent
{\it Proof.} Suppose that $\varpi\in\widehat\Delta_0$. Let $\Lambda$ be 
a rational 
representation of $G$ on the vector space $V$, with highest weight 
$d\varpi,\ d>0.$ Choose a height function relative to a basis on $V(F)$ 
as above. We can assume that {\it the basis contains a highest weight 
vector $v$}. 

According to the Bruhat decomposition, any element $\delta\in G(F)$ can 
be written in the form $\pi\, w_sn$ for $\pi\in P_0(F),\, n\in N_0(F)$ 
and $s\in \Omega$. Hence $H_0(\delta)=0???$. It follows from (e) above that,
 after suitable normalization on $v$ if necessary,
$$\|\Lambda(\delta)^{-1}v\|\geq 1\qquad\forall \delta\in G(F).$$

On the other hand, there are constants $c_1$ and $N_1$ such that for any
$x\in G(\mathbb A)^1$,
$$\begin{aligned}
\|\Lambda(\delta)^{-1}v\|=&\|\Lambda(x)\Lambda(\delta x)^{-1}v\|
 \leq \|\Lambda(x)\|\cdot \|\Lambda(\delta\, x)^{-1}v\|\qquad 
\mathrm{(by\ (d))}\\
=&\|x\|\cdot \|\Lambda(\delta\,x)^{-1}v\|\leq c_1\|x\|
e^{-d\varpi\Big(H_0(\delta x)\Big)}\qquad(\mathrm{by\ (e)\ again}).\end{aligned}$$ 
This combines the established relation 
$$\|\Lambda(\delta)^{-1}v\|\geq 1$$ completes the proof of the Lemma.

For each $x$, let $\Gamma(x)$ be a fixed set of representatives of 
$P(F)\backslash G(F)$ in $G(F)$ such that for any $\delta\in\Gamma(x)$, 
$\delta x$ belongs to $\omega\frak s A(\mathbb R)^0K$, where $\omega$ is a 
fixed compact subset of $N(\mathbb A)$ and $\frak s$ is a fixed Siegel set 
in $M(\mathbb A)^1$.
By reduction theory, or better, from the definition of Siegel domain $\frak s$,
we conclude that

\noindent
{\it there is a compact subset $\omega_0$ of $N_0(\mathbb A)M_0(\mathbb A)^1$ 
and a point $T_0$ in $\frak a_0$ such that for any $x$, and any 
$\delta\in \Gamma(x)$, $\delta x$ belongs to $\omega_0A_0(\mathbb R)^0K$, 
and in addition,
$$\alpha\Big(H_0(\delta x)\Big)\geq\alpha(T_0),\eqno(2)$$ for every 
$\alpha \in \Delta_0^P$.}

Next, note that what we are interested in are those $\delta$ such that
 $\widehat\tau_P\Big(H_0(\delta x)-X\Big)=1$, that is, such that
$$\varpi\Big(H_0(\delta x)\Big)>\varpi(X)\eqno(3)$$ for every 
$\varpi \in\widehat\Delta_P.$

Clearly, {\it the set of points $H_0(\delta x)$ in $\frak a_0^G$ which 
satisfy (1), (2) and (3) above is compact.}
Moreover, from our discussion, it follows that for 
$x\in G(\mathbb A)^1$,
and $\delta\in \Gamma(x)$, if $\widehat\tau_P\Big(H(\delta x)-X\Big)=1$, then
$\|\delta x\|$ is bounded by a constant multiple of a power of 
$\|x\|\cdot e^{\|X\|}$.

On the other hand, by (d),  $$\|\delta\|\leq\|\delta x\|\cdot\|x^{-1}\|\leq 
c\|\delta x\|\cdot\|x\|^N,$$ for some $c$ and $N$. Hence, $\|\delta\|$ too 
is bounded by a constant multiple of a power of $\|x\|\cdot e^{\|X\|}$.
Because $G(F)$ is a discrete subgroup of $G(\mathbb A)^1$, the Basic Estimation
follows from (c$'$) on the volume of $G_t$ above. This completes the proof.
\vskip 0.30cm
As a direct consequence,  we have the following

\noindent
{\bf Corollary.} Suppose that $T\in\frak a_0$ and $N\geq 0$. Then
{\it there exist constants $c'$ and $N'$ such that for any function $\phi$ 
on $P(F)\backslash G(\mathbb A)^1$, and $x,y\in G(\mathbb A)^1$,
$$\sum_{\delta\in P(F)\backslash G(F)}\Big|\phi(\delta x)\Big|
\cdot\widehat\tau_P
\Big(H(\delta x)-H(y)-X\Big)\eqno(4)$$ is bounded by $$c'\|x\|^{N'}\cdot
\|y\|^{N'}\cdot
\sup_{u\in G(\mathbb A)^1}\Big(|\phi(u)|\cdot \|u\|^{-N}\Big).$$}

\noindent
{\it Proof.} The expression (4) is bounded by the product of
$$\sup_{u\in G(\mathbb A)^1}\Big(|\phi(u)|\cdot \|u\|^{-N}\Big)$$ and 
$$\sum_{\delta\in P(F)\backslash G(F)}\|\delta x\|^N\cdot \widehat\tau_P
\Big(H(\delta x)-H(y)-X\Big).$$
We have shown in the above proof that when $\widehat\tau_P
\Big(H(\delta x)-H(y)-X\Big)$ is equal to 1, $\|\delta x\|$ is bounded
 by a constant multiple of a power of $\|x\|e^{\|H_P(y)+T\|}.$ The 
corollary therefore follows from the Basic
Estimation itself.
\vskip 0.30cm
In particular, we obtain the following

\noindent
{\bf Basic Fact I.} {\it If $\phi$ is a slow growth function, then so is 
$\Lambda^T\phi$.}

\subsection{\bf Langlands' Combinatorial Lemma}

In this section, following Arthur, we give a proof of what is usually 
called Langlands' combinatorial lemma. As an application, we give a 
refomulation of classical reduction theory which we call Arthur's partition 
for the total space.

If $P_1\subset P_2$, following Arthur [Ar3], set
$$\sigma_1^2(H):=\sigma_{P_1}^{P_2}:=\sum_{P_3:P_2\supset P_2}
(-1)^{\mathrm{dim}(A_3/A_2)}\tau_1^3(H)\cdot\hat\tau_3(H),$$ for 
$H\in\frak a_0$.

Before exposing detailed properties for $\sigma_1^2$,
we make the following comment: 
While the summation is for all (standard) parabolic subgroups $P_3$ 
which contains $P_2$, each term is with respect to $P_1$ in the sense 
that it is $\tau_1^3$, not $\tau_2^3$ which is used here. 
\vskip 0.30cm
\noindent
{\bf Lemma 1.} If $P_1\subset P_2$, $\sigma_1^2$ {\it is a characteristic 
function of the subset of $H\in\frak a_1$} such that

\noindent
(i) $\alpha(H)>0$ for all $\alpha\in\Delta_1^2$;

\noindent
(ii) $\sigma(H)\leq 0$ for all $\sigma\in\Delta_1\backslash \Delta_1^2$; and

\noindent
(iii) $\varpi(H)>0$ for all $\varpi\in\hat\Delta_2$.

Maybe it is better to recall that we have the following decomposition
$$\frak a_1=\frak a_1^2\oplus\frak a_2,$$ and that $\frak a_2$ may be 
identified with the subspace of $\frak a_1$ annihilated by $\Delta_1^2$. 
Consequently (i) says that for each $H=H_1^2+H_2$ according the above 
decomposition of $\frak a_0$, $H_1^2$ is in the 
positive chamber (of $\frak a_1^2$); (iii) says that $H_2$ is in the 
positive cone of $\frak a_2$; while (ii) says that $H=H_1^2+H_2$ is 
negative with respect to the roots of $\Delta_1$ outside $\Delta_1^2$.

With this understanding, the proof becomes a bit easy. But before that, 
let us state the following important consequence: 
\vskip 0.30cm
\noindent
{\bf Langlands' Combinatorial Lemma.} (I) If $Q\subset P$ are parabolic 
subgroups, then for all $H\in\frak a_0$,
$$\sum_{R:Q\subset R\subset P}(-1)^{\mathrm{dim} (A_R/A_P)}\tau_Q^R(H)
\hat\tau_R^P(H)=\delta_{QP}.$$

\noindent
{\it Proof.} Clearly, we may assume that $P=G$. Thus it suffices to 
show the following relation:
$$\sum_{P_3:P_3\supset P_1}(-1)^{\mathrm{dim} (A_1/A_3)}\tau_1^3(H)
\hat\tau_3(H)=\delta_{QP}.$$ But this is a direct consequence of the 
previous lemma. Indeed, in the lemma,  take the special case where 
$P_1=P_2\not=G.$ Then conditions (i) and (ii)
contradict to each other. That is to say, $\sigma_1^2(H)$ in this case 
is a characteristic function of the empty set, which should the constant 
function 0.
\vskip 0.30cm
Now let us come back to give a proof of the lemma itself.

\noindent
{\it Proof of the Lemma.} Fix $H\in\frak a_1$. Consider the subset of those 
$\varpi$ in $\hat\Delta_2$ for which $\varpi(H)>0$. While depending on 
$H$, this subset is of the form $\hat \Delta_R$, for a unique parabolic 
subgroup $R\supset P_2$. Accordingly,
$$\sigma_1^2(H)=\sum_{P_3:P_3\supset R}(-1)^{\mathrm{dim}(A_3/A_2)}
\tau_1^3(H).$$

Suppose now that $\tau_1^3(H)=1$ for a given $P_3: P_3\supset R$. Then 
$\tau_1^3(H)=1$ for all smaller $P_3$. It follows from the Ground 0 
Relation that the above sum vanishes unless the original $P_3$ equals $R$. 
Thus, for a fixed $R: R\supset P_2$, 
$$\Big|\sum_{P_3:P_3\supset R}(-1)^{\mathrm{dim}(A_3/A_2}\tau_1^3(X)\Big|$$ 
is the characteristic function of
$$\Big\{X\in\frak a_1:\alpha(X)>0, \alpha\in\Delta_1^R; \beta(X)\leq 0,\ 
\beta\in\Delta_1\backslash\Delta_1^R\Big\}.$$
Indeed, the first condition on $\alpha$ is rather clear while the condition 
for $\beta$ comes from the statement about Grand 0.

Now there are two cases. When $R=P_2$, clearly $\sigma_1^2(H)=1$ if and only 
if $H$ satisfies conditions (i), (ii) and (iii) of the Lemma, as required.

When $R\not= P_2$, i.e., $R$ is strictly larger than $P_2$, we mush show that 
$\sigma_1^2(H)=0$. Suppose not. Then $H$ belongs to the subset in the 
Lemma by the above discussion. In particular, the projection of 
$H$ onto $\frak a_1^R$ lies in the positive chamber as 
$\Delta_1^R\subset\Delta_1^2$, which is contained in the set of positive 
linear combinations of roots in $\Delta_1^R$. Thus $\varpi(H)>$ for all 
$\varpi\in\hat\Delta_1^R$. Also, by the defintion of $R$, $\varpi(H)>0$ 
for all $\varpi\in\hat\Delta_R$. From this, we shall show the following

\noindent
{\bf Claim.} {\it if $H=\sum_{\alpha\in\Delta_1}c_\alpha\alpha^\vee,$ 
each $c_\alpha$ is positive.}

\noindent
{\it Proof of the Claim.} Suppose that 
$\varpi\in\hat \Delta_R\subset\hat\Delta_1$, and that 
$\alpha_\varpi$ is the element in $\Delta_1$ which is paired with 
$\varpi$. Then $c_{\alpha_\varpi}=\varpi(H)$ is {\it positive}. Therefore 
the projection of $$H_R:=\sum_{\varpi\in\hat\Delta_R}c_{\alpha_\varpi}
\alpha_\varpi^\vee$$ onto $\frak a_1^R$ is in the negative chamber, so 
that if $\nu\in\hat\Delta_1^R$, $\nu(H_R)$ is negative.
Hence, if $\alpha_\nu$ is the root in $\Delta_1^R$ corresponding to 
$\nu$, $c_{\alpha_\nu}=\nu(H)-\nu(H_R)$ is {\it positive}.
Thus each $c_\alpha$ is positive. This completes the proof of the claim.

Consequently, $\varpi(H)$ is positive 
for each $\varpi\in\hat\Delta_1$. But $\hat\Delta_2\subset\hat\Delta_1$ 
which shows that $\varpi(H)>0$ for all $\varpi\in\hat\Delta_2$. So by the 
definition of $R$ again, $R=P_2$ so we have a contradiction. This completes the proof of the Lemma.
\vskip 0.30cm
We end this section with a direct consequence of Langlands' Combinatorial 
Lemma. Before stating it, let us make the following preperation:

Suppose that $Q\subset P$ are parabolic subgroups. Fix a vector
$\Lambda\in\frak a_0^*$. Let $$\varepsilon_Q^P(\Lambda):=
(-1)^{\#\{\alpha\in\Delta_Q^P:\Lambda(\alpha^\vee)\leq 0\}},$$ 
and let $$\phi_Q^P(\Lambda,H),\qquad H\in\frak a_0,$$ be the characteristic 
function of the set  $$\Big\{
H\in\frak a_0: 
\begin{matrix}
\varpi(H)>0, &\mathrm{if}\ \Lambda(\alpha^\vee)\leq 0\\ 
\varpi(H)\leq 0,&\mathrm{if}\ \Lambda(\alpha^\vee)>0
\end{matrix},
\forall\alpha\in\Delta_Q^P
\Big\}.$$

\noindent
{\bf Lemma 2.} With the same notation as above,
$$\sum_{R:Q\subset R\subset P}\varepsilon_Q^R(\Lambda)\cdot
\phi_Q^R(\Lambda,H)\cdot\tau_R^P(H)=\begin{cases} 0, &\text{if $\Lambda(\alpha^\vee)\leq 0,\ 
\exists \alpha\in\Delta_Q^P$}\\
1,&\text{otherwise}\end{cases}.$$

\noindent
{\it Proof.} By Langlands' combinatorial lemma, if $R\not=P$, then
$$\sum_{P_1:R\subset P_1\subset P}(-1)^{\mathrm{dim}{A_1/A_P}}\tau_R^1(H)
\hat\tau_1^P(H)=0$$ for all $H$. Therefore
$$\sum_{R:Q\subset R\subset P}\varepsilon_Q^R(\Lambda)\phi_Q^R(\Lambda,H)
\tau_R^P(H)$$ is the difference between $$\varepsilon_Q^P(\Lambda)
\phi_Q^P(\Lambda,H)$$ and
$$\sum_{R,P_1:Q\subset R\subset P_1\subsetneq P}\varepsilon_Q^R(\Lambda)
\cdot\phi_Q^R(\Lambda,H)\cdot\tau_R^1(H)\cdot
(-1)^{\mathrm{dim}(A_1/A_P)}\cdot\hat\tau_1^P(H).\eqno(1)$$
We shall prove the lemma by induction on $\mathrm{dim}(A_0/A_P)$.

Define $\Delta_Q^\Lambda$ to be the set of roots $\alpha\in\Delta_Q^P$ 
such that $\alpha(\Lambda)>0$. Associated to $\Delta_A^\Lambda$ we have a 
parabolic subgroup $P_\Lambda$ with $Q\subset P_\Lambda\subset P$. By our 
induction assumption, the sum over $R$ in (1) vanishes unless $P_1\subset
 P_\Lambda$, in which case it equals $(-1)^{\mathrm{dim}(A_1/A_P)}
\hat\tau_1^P(H)$. Thus, (1) equals
$$\begin{cases}
\sum_{P_1:Q\subset P_1\subset P_\Lambda}(-1)^{\mathrm{dim}(A_1/A_P)}
\hat\tau_1^P(H),&\text{if $P_\Lambda\not=P$,}\\
\sum_{P_1:Q\subset P_1\subset P_\Lambda}(-1)^{\mathrm{dim}(A_1/A_P)}
\hat\tau_1^P(H)-1,&\text{if $P_\Lambda=P$}\end{cases}.$$ Thus we only 
need to show that 
$$\varepsilon_Q^P(\Lambda)\phi_Q^P(\Lambda,H)=\sum_{P_1:Q\subset 
P_1\subset P_\Lambda}(-1)^{\mathrm{dim}(A_1/A_P)}\hat\tau_1^P(H).$$ 
This is a consequence of Grand 0 Relation and the definition.
\vskip 0.30cm
As a direct consequence, we have the following

\noindent
\vskip 0.30cm
\noindent
{\bf Langlands' Combinatorial Lemma.} (II) If $Q\subset P$ are parabolic 
subgroups, then for all $H\in\frak a_0$,
$$\sum_{R:Q\subset R\subset P}(-1)^{\mathrm{dim} (A_Q/A_R)}\widehat\tau_Q^R(H)
\tau_R^P(H)=\delta_{QP}.$$

\noindent
{\it Proof.}(Arthur) {\it Method 1.} In Lemma 2, let $\Lambda\in -(\frak a_0^*)^+$, then
$$\varepsilon_Q^R(\Lambda):=
(-1)^{\#\{\alpha\in\Delta_Q^R:\Lambda(\alpha^\vee)\leq 0\}}=
(-1)^{\#\Delta_Q^R}=(-1)^{\mathrm{dim} (A_Q/A_R)}.$$ 
Moreover, since, by definition, 
$$\phi_Q^R(\Lambda,H),\qquad H\in\frak a_0,$$ is nothing but 
the characteristic function of the set  $$\Big\{
H\in\frak a_0: 
\begin{matrix}
\varpi(H)>0, &\mathrm{if}\ \Lambda(\alpha^\vee)\leq 0\\ 
\varpi(H)\leq 0,&\mathrm{if}\ \Lambda(\alpha^\vee)>0
\end{matrix},
\forall\alpha\in\Delta_Q^R
\Big\}.$$ So, with $\Lambda\in -(\frak a_0^*)^+$,
$$\phi_Q^R(\Lambda,H)=\widehat\tau_Q^R(H).$$ Therefore,
$$\sum_{R:Q\subset R\subset P}(-1)^{\mathrm{dim} (A_Q/A_R)}\widehat\tau_Q^R(H)
\tau_R^P(H)=\sum_{R:Q\subset R\subset P}\varepsilon_Q^R(\Lambda)\cdot
\phi_Q^R(\Lambda,H)\cdot\tau_R^P(H)=\delta_{QP},$$ as required.

\noindent
{\it Method 2.} It is enough to assume that $P=G$ and $Q$ is proper in $G$. 
One then can deduce (II) form (I) by evaluating the expression
$$\sum_{R,P:Q\subset R\subset P}(-1)^{\mathrm{dim}(A_R/A_P)}
\widehat\tau_Q^R(H)\tau_R^P(H)\widehat\tau_QP(H)$$ as two different iterated 
sums. For if one takes $P$ to index the inner sum, and assumes inductively 
that (II) holds whenever $G$ is replaced by a proper Levi subgroup, one 
finds that the expression equals the sum $\widehat\tau_Q(H)$ with the 
left hand side of (II). On the other hand, by taking the inner sum to be 
over $R$, one sees from (I) that the expression reduces simply to 
$\widehat\tau_Q(H)$. It follows that the left hand side of (II) vanishes, 
as required.

\subsection{Langlands-Arthur's Partition: Reduction Theory}

Our aim here is to derive Langlands-Arthur's partition of 
$G(F)\backslash G(\mathbb A)$ into disjoint subsets, one 
for each (standard) parabolic subgroup. 
This partition is similar to a construction of Langlands, in which 
disjoint subsets of $G(F)\backslash G(\mathbb A)$ are associated to 
maximal parabolic subgroups. More generally, we shall partition 
$N(\mathbb A)M(F)\backslash G(\mathbb A)$, where $P=NM$ is a parabolic 
subgroup. Essentially, this is a restatement of the following basic 
lemma in reduction theory:

Suppose that $\omega$ is a compact subset of $N_0(\mathbb A)M_0(\mathbb A)^1$
 and that $T_0\in-\frak a_0^+$. For any parabolic subgroup $P_1$, let 
$\frak s^{P_1}(T_0,\omega)$ be the set of
$$pak,\qquad p\in\omega,\ a\in A_0(\mathbb R)^0,\ k\in K,$$
such that $\alpha\Big(H_0(a)-T_0\Big)$ is positive for each 
$\alpha\in\Delta_0^1$. 

\noindent
{\bf {\large Basic Fact.}} (Siegel Sets) {\it We can fix $\omega$ and $T_0$ so that
for any $P_1$, $G(\mathbb A)=P_1(F)\frak s^{P_1}(T_0,\omega).$}

The proof is based on the following

\noindent
{\bf {\large Basic Lemma.}} (Reduction Theory) Let $T$ in $\frak a_0^+$ be any 
suitably regular point. Suppose that $P_1\subset P$ are parabolic subgroups, 
and that
$x$ and $\delta x$ belong to $\frak s^{P_1}(T_0,\omega)$ for points 
$x\in G(\mathbb A)$ and $\delta\in P(F)$. Then {\it if 
$\alpha\Big(H_0(a)-T_0\Big)>0$
 for all $\alpha\in\Delta_0^P\backslash \Delta_0^{P_1}$,  $\delta$ 
belongs to $P_1(F)$.}

For the proof, please see [MW].
\vskip 0.30cm
Suppose that $P_1$ is given. Let $\frak s^{P_1}(T_0,T,\omega)$ be the set 
of $x$ in $\frak s^{P_1}(T_0,\omega)$ such that 
$\varpi\Big(H_0(x)-T\Big)\leq 0$ for each $\varpi\in\hat\Delta_0^1$. 
Let $$F^{P_1}(x,T):=F^1(x,T)$$ be the characteristic function of the 
set of $x\in G(\mathbb A)$ such that $\delta x$ belongs to 
$\frak s^{P_1}(T_0,T,\omega)$ for some $\delta\in P_1(F)$.

As such, $F^1(x,T)$ is left $A_1(\mathbb R)^0N_1(\mathbb A)M_1(F)$ 
invariant, and can be regarded as the characteristic function of the 
projection of $\frak s^{P_1}(T_0,T,\omega)$ onto 
$A_1(\mathbb R)^0N_1(\mathbb A)M_1(F)\backslash G(\mathbb A)$, 
a compact subset of $A_1(\mathbb R)^0N_1(\mathbb A)M_1(F)\backslash 
G(\mathbb A)$.

For example,  $F(x,T):=F^G(x,T)$ admits the following more direct description
which will play a key role in our study of Arthur's periods:

If $P_1\subset P_2$ are (standard) parabolic subgroups, we write 
$A_1^\infty:=A_{P_1}^\infty$ for $A_{P_1}(\mathbb A)^0$, the identity 
component of $A_{P_1}(\mathbb R)$, and $$A_{1,2}^\infty:=
A_{P_1,P_2}^\infty:=A_{P_1}\cap M_{P_2}(\mathbb A)^1.$$ Then 
$H_{P_1}$ maps $A_{1,2}^\infty$ isomorphically onto $\frak a_1^2$, 
the orthogonal complement of $\frak a_2$ in $\frak a_1$. If $T_0$ and 
$T$ are points in $\frak a_0$, set $A_{1,2}^\infty(T_0,T)$ equal to the set
$$\Big\{a\in A_{1,2}^\infty:\alpha\Big(H_1(a)-T\Big)>0,\ \alpha\in\Delta_1^2;
\
\varpi\Big(H_1(a)-T\Big)<0,\,\varpi\in\hat\Delta_1^2\Big\},$$ where
$$\Delta_1^2:=\Delta_{P_1\cap M_2}\qquad\mathrm{and}\qquad\hat\Delta_1^2
:=\hat\Delta_{P_1\cap M_2}.$$ Now fix $T_0$ so that $-T_0$ is suitably regular.
Then 

\noindent
{\it $F(x,T)$ is the characteristic function of the compact subset of 
$G(F)\backslash G(\mathbb A)^1$ obtained by projecting 
$$N_0(\mathbb A)\cdot M_0(\mathbb A)^1\cdot A_{P_0,G}^\infty(T_0,T)\cdot K$$ 
onto $G(F)\backslash G(\mathbb A)^1.$}

(The reader should know that the compactness comes from that for 
$A_{1,2}^\infty(T_0,T)$: being positive for roots and negative for weights, 
it is compact, since $N_0(F)\backslash N_0(\mathbb A)$ and
 $M(F)\backslash M(\mathbb A)^1$ are compact.)

\noindent
{\bf {\large Arthur's Partition.}} {\it Fix $P$ and let $T$ be any suitably 
point in $T_0+\frak a_0^+$. Then $$\sum_{P_1:P_0\subset P_1\subset P}
\sum_{\delta\in P_1(F)\backslash G(F)}F^1(\delta x)\cdot\tau_1^P
\Big(H_0(\delta x)-T\Big)=1\qquad\forall x\in G(\mathbb A).$$}

\noindent
{\it Proof.} We will use Lemma 2 in the previous section.

\noindent
{\bf Step 1:  The given sum is at least one.} Fix $x\in G(\mathbb A)$. Choose 
$\delta\in P(F)$ such that $\delta x$ belongs to $\frak s^P(T_0,\omega)$.
 Apply Lemma 2 with $Q=P_0,\Lambda\in(\frak a_0^*)^+$ and $H=H_0(\delta x)-T.$ 
Then we see that the right hand side (of Lemma 2) takes the value 1, 
since $\Lambda\in (\frak a_0^*)^+$ implies that the condition 
$\Lambda(\alpha^\vee)>0$ holds for all $\alpha\in\Delta_0^P$. Thus,
 on the left hand (of Lemma 2), there exists at least one term, say, 
$R=P_1$ so that $$\varepsilon_{P_0}^{P_1}(\Lambda)\cdot
\phi_{P_0}^{P_1}\Big(\Lambda,H_0(\delta x)-T\Big)\cdot
\tau_{P_1}^P\Big(H_0(\delta x)-T\Big)$$ takes values 1.
Since $\Lambda\in(\frak a_0^*)^+$, so $\varepsilon_{P_0}^{P_1}(\Lambda)=
(-1)^0=1$.
Thus, both factors $$\phi_{P_0}^{P_1}\Big(\Lambda,H_0(\delta x)-T\Big)\qquad
\mathrm{and}\qquad\tau_{P_1}^P\Big(H_0(\delta x)-T\Big)$$ take values 1. By 
definition of $\phi_{P_0}^{P_1}$, note that $\Lambda(\alpha^\vee)>0$ 
holds, we have $$\varpi\Big(H_0(\delta x)-T\Big)\leq 0\qquad \forall 
\varpi\in\hat\Delta_0^{P_1}.\eqno(1)$$ Similarly, by the definition of 
$\tau_0^1$, we conclude that $$\alpha\Big(H_0(\delta x)-T\Big)> 0\qquad 
\forall \alpha\in\Delta_1^{P}.\eqno(2)$$
Thus by definition, from (1) and (2), we have
$$F^1(\delta x)\cdot\tau_1^P\Big(H_0(\delta x)-T\Big)=1,$$ as required.

\noindent
{\bf Step 2: The given sum is at most one.} Suppose that there are elements 
$\delta_1,\delta_2\in G(F)$, and parabolic subgroups $P_1$ and $P_2$ 
contained in $P$ such that
$$F^1(\delta_1 x)\cdot\tau_1^P\Big(H_0(\delta_1 x)-T\Big)=F^2(\delta_2 x)
\cdot\tau_2^P\Big(H_0(\delta_2 x)-T\Big)=1.$$
From the reduction theory, after left translating $\delta_1$ by an element 
in $P_i(F)$ if necessary, we may assume that
$$\delta_i x\in\frak s^{P_i}(T_0,T,\omega),\qquad i=1,2.$$
Moreover, by definition (see also the dicussion in A), the projection of 
$H(\delta_i x)-T$ onto $\frak a_0^P$ can be written as 
$$-\sum_{\alpha\in\Delta_0^i}c_\alpha\alpha^\vee+
\sum_{\varpi\in\hat\Delta_i^P}c_\varpi\varpi^\vee,$$ where each 
$c_\alpha$ and $c_\varpi$ is positive. It follows that 
$\alpha\Big(H_0(\delta_i x)-T\Big)>0$ for every 
$\alpha\in\Delta_0^P\backslash\Delta_0^i$. In particular, since 
$T$ lies in $T_0+\frak a_0^+$, $\delta_i x$ belongs to $
\frak s^P(T_0,\omega)$. The reduction theoretic result cited as the 
Basic Lemma above implies that $\delta_2\delta_1^{-1}$ belongs to $P_1(F)$ and 
$\delta_1\delta_2^{-1}$ belongs to $P_2(F)$. In other words, 
$\delta_2=\xi\delta_1$ for some element $\xi\in P_1(F)\cap P_2(F)$. 
Let $Q=P_1\cap P_2$. Then 
$H_0(\delta_1 x)-T$ and $H_0(\delta_2 x)-T$ project onto the same point, 
say $H_Q^P$ on $\frak a_Q^P$. If $R$ equals either $P_1$ or $P_2$, we have
$\varpi(H_Q^P)\leq 0$ for $\varpi\in\hat\Delta_Q^R$ and $\alpha(H_Q^P)>0$ 
for $\alpha\in\Delta_R^P$. Applying Lemma 2  with $\Lambda\in(\frak a_0^*)^+$ 
again, as discussed in A, by definition, we see that there is exactly one $R$,
 with $Q\subset R\subset P$, for which these inequalities hold. Therefore 
$P_1=P_2$ and $\delta_1$ and $\delta_2$ belong to the same $P_1(F)$ coset in
 $G(F)$. This proves that the given sum is at most one. We are done.
\vskip 0.30cm
As Arthur points out  (see e.g., [Ar6]), this  partition can be restated 
geometrically in terms of the subsets
$$G_P(T):=\Big\{x\in G(F)\backslash G(\mathbb A): F^P(x,T)=1,\ 
\tau_P\Big(H_P(x)-T\Big)=1\Big\}$$ of 
$P(F)\backslash G(\mathbb A)$. Indeed, Arthur's partition says that for any 
$P$, the projection of 
$P(F)\backslash G(\mathbb A)$ onto $G(F)\backslash G(\mathbb A)$ maps 
$G_P(T)$ {\it injectively}
onto a subset $\overline G_P(T)$ of $G(F)\backslash G(\mathbb A)$, and that
$G(F)\backslash G(\mathbb A)$ is a disjoint union over $P$ of the sets 
$\overline G_P(T)$. Otherwise said,
$G(F)\backslash G(\mathbb A)^1$  has a partition parametrized by the 
set of standard parabolic subgroups, which separates the problem of  
noncompactness from the topologically
complexity of $G(F)\backslash G(\mathbb A)^1$. The subset corresponding 
to $P=P_0$ is topologically
simple but highly noncompact. The subset corresponding to a group 
$P\not\in\{P_0,G\}$ is mixed, being a product of a compact set of
intermediate complexity with a simple set of intermediate degree of 
noncompactness. The partition of $G(F)\backslash G(\mathbb A)^1$ is, 
incidentally, closely related to the compactification of this space defined 
by Borel and Serre.

\section{Arthur's Analytic Truncation}
\subsection{Definition}
Fix a suitably regular point $T\in\frak a_0^+$. If $\phi$ is a continuous 
function on $G(F)\backslash G(\mathbb A)^1$, define Arthur's analytic 
trunction $\Big(\Lambda^T\phi\Big)(x)$ to be the function
$$\Big(\Lambda^T\phi\Big)(x):=
\sum_P(-1)^{\mathrm{dim} (A/Z)}\sum_{\delta\in P(F)\backslash G(F)}
\phi_P(\delta x)\cdot\hat\tau_P\Big(H(\delta x)-T\Big),$$ where
$$\phi_P(x):=\int_{N(F)\backslash N(\mathbb A)}\phi(nx)\,dn$$ denotes the 
constant term of $\phi$ along $P$, and the sum is over all (standard) 
parabolic subgroups.
 
The main purpose for introducing analytic truncation is to give a natural way
to construct integrable functions: even from the example of $GL_2$, we 
know that automorphic forms are generally not integrable over the total 
fundamental domain $G(F)\backslash G(\mathbb A)^1$ mainly due to the fact 
that in the Fourier expansions of such functions, the so-called constant 
terms are only of moderate growth (hence not integrable). Thus in order 
to naturally obtain  integrable functions, we should truncate the original 
function along the cuspidal regions by removing constant terms. Thus, 
naturally, in Arthur's more general definition,
the summation may be divided into different levels according to the ranks 
of $P$. Say, at the grand level zero where $P=G$, the term we obtain is 
nothing but the function $\phi$ itself; and at the level one, where $P$ 
is maximal, then what we obtain is the product of the constant term $\phi_P$ 
of $\phi$ along $P$, together with a minus sign. That is to say, we have to
 substract from $\phi$ the constant term $\phi_P$ of $\phi$ along with 
these maximal parabolic subgroups. However the term substracted is not 
over the whole space, rather, it is only over the cuspidal region 
corresponding to maximal parabolics; with this being done, then we go to the 
next level, and so on. Simply put, Arthur's analytic truncation is a 
well-designed divice in which constant terms are tackled in such a way 
that  different levels of parabolic subgroups are suitably counted at 
the corresponding cuspidal region so that the whole truncation will not 
be overdone while there will be no parabolic subgroups left untackled.

Note that all parabolic subgroups of $G$ can be obtained from standard 
parabolic subgroups by taking conjugations with elements from 
$P(F)\backslash G(F)$. So we have:

\noindent
(a) $\displaystyle{\Big(\Lambda^T\phi\Big)(x)=\sum_P(-1)^{\mathrm{dim} (A/Z)}
\phi_P(x)
\cdot\hat\tau_P\Big(H(x)-T\Big),}$ {\it where the sum is over all, both
standard and non-standard, parabolic subgroups};
 
\noindent
(b) {\it If $\phi$ is a cusp form, then $\Lambda^T\phi=\phi$};

This is because by definition, all constant terms along proper $P: P\not=G$ 
are zero.
Moreover, it is a direct consequence of the Basic Estimation for partial 
truncation, (see e.g. the Corollary there), we have

\noindent
(c) {\it If $\phi$ is of moderate growth} in the sense that there exist some 
constants $C,N$ such that $$\Big|\phi(x)\Big|\leq c\|x\|^N$$ for all 
$x\in G(\mathbb A)$, {\it then so is $\Lambda^T\phi$}.

\subsection{Idempotence}
Next property for Arhur's analytic truncation is the following
\vskip 0.30cm
\noindent
{\bf {\large Idempotence}}. $\Lambda^T\circ\Lambda^T=\Lambda^T$.

Indeed this is a direct consequence of the following much stronger

\noindent
{\bf Lemma.} {\it Fix $P_1$. Then for $\phi\in C\Big(G(F)\backslash 
G(\mathbb A)\Big)$,
$$\int_{N_1(F)\backslash N_1(\mathbb A)}\Lambda^T\phi(n_1x)\,dn_1=0$$ 
unless $\varpi\Big(H_0(x)-T\Big)<0$ for each $\varpi\in\hat \Delta_1$.}

(We say that this is much stronger because for the Idempotance of 
$\Lambda^T$, by definition it is suffices to show that 
$$\Big(\Lambda^T\phi\Big)_{P_1}=0$$ 
unless there exists a certain $\varpi\in\hat \Delta_1$ such that 
$\varpi\Big(H_0(x)-T\Big)<0$. But the lemma says that
$$\Big(\Lambda^T\phi\Big)_{P_1}=0$$ 
unless for all $\varpi\in\hat \Delta_1$, $\varpi\Big(H_0(x)-T\Big)<0$).

\noindent
{\it Proof.} {\bf Step 1.} For any $P$, let $$\Omega(\frak a_0;P):=
\Big\{s\in\Omega:s^{-1}\alpha>0,\ \forall\alpha\in\Delta_0^P\Big\}.$$
Applying the Bruhat decomposition to $P(F)\backslash G(F)$, we find that
$$\begin{aligned}
&\int_{N_1(F)\backslash N_1(\mathbb A)}\Lambda^T\phi(n_1x)\,dn_1\\
=&
\int_{N_1(F)\backslash N_1(\mathbb A)}
\sum_P(-1)^{\mathrm{dim} (A/Z)}
\sum_{\delta\in P(F)\backslash G(F)}\phi_P(\delta n_1x)\cdot\hat\tau_P
\Big(H(\delta n_1x)-T\Big)\,dn_1\\
=&\sum_P(-1)^{\mathrm{dim} (A/Z)}\int_{N_1(F)\backslash N_1(\mathbb A)}
\sum_{\delta\in P(F)\backslash G(F)}\int_{N(F)\backslash N(\mathbb A)}\phi
(n\delta n_1x)\,dn\cdot\hat\tau_P\Big(H(\delta n_1x)-T\Big)
\,dn_1\\
=&\sum_P(-1)^{\mathrm{dim} (A/Z)}\sum_{\delta\in P(F)\backslash G(F)}
\int_{N(F)\backslash N(\mathbb A)}\Big(\int_{N_1(F)\backslash N_1(\mathbb A)}
\phi
(n\delta n_1x)\cdot\hat\tau_P\Big(H(\delta n_1x)-T\Big)
\,dn_1\Big)\,dn\\
=&\sum_P(-1)^{\mathrm{dim} (A/Z)}\sum_{s\in\Omega(\frak a_0;P)}\int_{N(F)
\backslash N(\mathbb A)}dn\\
&\qquad\times\Big(
\int_{N_1(F)\backslash N_1(\mathbb A)}\sum_{\nu\in w_s^{-1}N_0(F)w_s\cap 
N_0(F)\backslash N_0(F)}\phi
(n w_s\nu n_1x)\cdot\hat\tau_P\Big(H(w_s\nu n_1x)-T\Big)
\,dn_1\Big).\end{aligned}$$
Now let us analyse the inner integration $$
\int_{N_1(F)\backslash N_1(\mathbb A)}\sum_{\nu\in w_s^{-1}N_0(F)w_s\cap 
N_0(F)\backslash N_0(F)}\phi
(n w_s\nu n_1x)\cdot\hat\tau_P\Big(H(w_s\nu n_1x)-T\Big)
\,dn_1.\eqno(1)$$ 
Note that $N_1\supset N_0$ and we may have a decomposition
$$N_1=N_0\cdot N_0^1$$ for a suitable nilpotent subgroup $N_0^1$.
 (Try use the group $GL_n$ to understand this decomposition.)
Thus $N_1(F)\backslash N_1(\mathbb A)=N_0(F)\backslash N_0^1(F)N_1(\mathbb A).$
Consequently, (1) is equal to
$$\int_{w_s^{-1}N_0(F)w_s\cap N_0(F)\backslash N_0^1(F)N_1(\mathbb A)}
\phi(n w_s n_1x)\cdot\hat\tau_P\Big(H(w_s n_1x)-T\Big)
\,dn_1.\eqno(2)$$
Decompose $w_s^{-1}N_0(F)w_s\cap N_0(F)\backslash N_0^1(F)N_1(\mathbb A)$ as
$$\begin{aligned}&\Big(w_s^{-1}N_0(F)w_s\cap N_0(F)\backslash
w_s^{-1}N_0(\mathbb A)w_s\cap N_0^1(F)N_1(\mathbb A)\Big)
\times \Big(w_s^{-1}N_0(\mathbb A)w_s\cap N_0^1(F)N_1(\mathbb A)
\backslash N_0^1(F)N_1(\mathbb A)\Big)\\
=&\Big(w_s^{-1}N_0(F)w_s\cap N_0(F)\backslash
w_s^{-1}N_0(\mathbb A)w_s\cap N_1(\mathbb A)\Big)
\times \Big(w_s^{-1}N_0(\mathbb A)w_s\cap N_0^1(F)N_1(\mathbb A)
\backslash N_0^1(F)N_1(\mathbb A)\Big).\end{aligned}$$ 
This induces a decomposition of 
the measure $dn_1$ as $dn_*dn^*.$ Then write
$$w_sn_*n^*=w_sn_*w_s^{-1}w_sn^*=\tilde n_*w_sn^*,$$ and finally combine 
the integral over $\tilde n_*$ with the integral over $n$ in $N(F)
\backslash N(\mathbb A)$. Because $s$ lies in $\Omega(\frak a_0;P)$, 
$N_0\cap w_sN_1w_s^{-1}\cap M$ {\it is the nilpotent radical of a standard 
parabolic subgroup of $M$.} It follows that $$\Big(N_0\cap w_sN_1w_s^{-1}
\cap M\Big)N=N_s$$ is the unipotent radical of a {\it uniquely determined}
 parabolic subgroup $P_s$ of $G$, which is contained in $P$. That is to 
say, we have shown the following

\noindent
{\bf Sublemma 1.} {\it With the same notation as above,
$$\begin{aligned}
&\int_{N_1(F)\backslash N_1(\mathbb A)}\Lambda^T\phi(n_1x)\,dn_1\\
=&\sum_P(-1)^{\mathrm{dim} (A/Z)}\sum_{s\in\Omega(\frak a_0;P)}
\int_{w_s^{-1}N_0(\mathbb A)w_s\cap N_0^1(F)N_1(\mathbb A)\backslash N_0^1(F)
N_1(\mathbb A)}dn^*\\
&\qquad\Big(
\int_{N_s(F)\backslash N_s(\mathbb A)}\phi
(n w_s n^*x)\cdot\hat\tau_P\Big(H(w_s n^*x)-T\Big)
\,dn\Big).\end{aligned}$$}

\noindent
{\bf Step 2.} We shall change the order of summation, and consider the 
set of $P$ which give rise to a fixed $P_s$. For a fixed $s\in \Omega$, 
decompose $\Delta_0=S^1\sqcup S_1\sqcup S_0$ into disjoint unions with $S_1$ 
(resp. $S_1$) being the set of $\alpha\in\Delta_0$ such that $s^{-1}\alpha$ 
is a positive root which is {\it
orthogonal} (resp. {\it not orthogonal}) to $\frak a_1$. Then by definition

\noindent
(1) If $P_s$ is one of the groups that appear in the above formula, 
and $\Delta_0^s$ is a subset of $S^1$;

\noindent
(2) Those $P$ which give rise to a fixed $P_s$ are exactly the groups for 
which $\Delta_0^P$ is the union of $\Delta_0^s$ and a subset $S$ of $S_1$.

Thus, for fixed $s$ with $\Delta_0^s\subset S^1$, we will obtain an 
alternating sum over subsets $S\subset S_1$ of the corresponding functions 
$\hat\tau_P$. Consequently, we may apply the Basic Relation of Grand 0 
to conclude the follows:

Let $\chi_s$ be the characteristic function of the set 
 $$\Big\{H\in\frak a_0:
\begin{matrix} \varpi_\alpha(H)>0,&\mathrm{if}\ 
\alpha\in\Delta_0\backslash(\Delta_0^s\cup S_1)\\ 
\varpi_\alpha(H)\leq 0& 
\mathrm{if}\ \alpha\in S_1
\end{matrix}
\Big\}.$$ Here 
$\varpi_\alpha$ is the element in $\hat\Delta_0$ corresponding to $\alpha$. 

\noindent
{\bf Sublemma 2.} {\it With the same notation as above,
$$\begin{aligned}
&\int_{N_1(F)\backslash N_1(\mathbb A)}\Lambda^T\phi(n_1x)\,dn_1\\
=&\sum_{s\in\Omega}\sum_{\Delta_0^s\subset S^1}
\int_{w_s^{-1}N_0(\mathbb A)w_s\cap N_0^1(F)N_1(\mathbb A)\backslash N_0^1(F)
N_1(\mathbb A)}dn^*\\
\qquad\qquad&\times\Big(\int_{N_s(F)\backslash N_s(\mathbb A)}\phi
(n w_s n^*x)\cdot(-1)^{\#(\Delta_0\backslash(S^1\cup S_1))}\chi_s
\Big(H(w_s n^*x)-T\Big)\,dn\Big).\end{aligned}$$}

\noindent
{\bf Step 3.} Suppose that for some $s$, $\chi_s\Big(H(w_s n^*x)-T\Big)
\not=0$. Clearly, by definition, if $$H(w_s n^*x)-T=\sum_{\alpha\in \Delta_0} 
t_\alpha\alpha^\vee,\qquad t_\alpha\in\mathbb R,$$ then $t_\alpha$ is
 positive for $\alpha\in \Delta_0\backslash(\Delta_0^s\cup S_1)$, and 
is not positive for $\alpha\in S_1$.

Now if $\varpi\in\hat\Delta_1$,
$$\begin{aligned}&\varpi\Big(s^{-1}\Big(H(w_s n^*x)-T\Big)\Big)
=\sum_{\alpha\in \Delta_0} t_\alpha\cdot \varpi(s^{-1}\alpha^\vee)\\
=&\sum_{\alpha\in \Delta_0\backslash S^1} t_\alpha\cdot 
\varpi(s^{-1}\alpha^\vee)+\sum_{\alpha\in  S^1} t_\alpha\cdot 
\varpi(s^{-1}\alpha^\vee)
=\sum_{\alpha\in \Delta_0\backslash S^1} t_\alpha\cdot 
\varpi(s^{-1}\alpha^\vee)\end{aligned}$$ 
since $S^1$ is defined to be the set of 
$\alpha\in\Delta_0$ such that $s^{-1}\alpha$ is orthogonal to $\frak a_1$. 

Moreover, by definition $$\sum_{\alpha\in \Delta_0\backslash S^1} 
t_\alpha\cdot \varpi(s^{-1}\alpha^\vee)
=\sum_{\alpha\in  S_1} t_\alpha\cdot \varpi(s^{-1}\alpha^\vee)+
\sum_{\alpha\in  S_0} t_\alpha\cdot \varpi(s^{-1}\alpha^\vee).$$ 
In the first case, if $\alpha\in S_1$, then $t_\alpha\leq 0$ while
 $\varpi(s^{-1}\alpha^\vee)>0$ since $s^{-1}\alpha$ is by defintion 
a positive root. Consequently, 
{\it the first sum is always less than or equal to 0}.
 In the second case, if $\alpha\in S_0$, then $\alpha\in \Delta_0
\backslash(\Delta_0^s\cup S_1)$, so $t_\alpha>0$ as stated above, 
while $s^{-1}\alpha$ is not a positive root, hence a negative of a 
positive root, which implies that $\varpi(s^{-1}\alpha^\vee)\leq 0$. 
That is to say, {\it the second sum is also always less than or equal 
to 0}. All in all, we have shown that following

\noindent
{\bf Claim.} {\it If $\varpi\in\hat\Delta_1$, then 
$$\varpi\Big(s^{-1}\Big(H(w_s n^*x)-T\Big)\Big)\leq 0.$$}

On the other hand, since $w_sn^*x=w_sn^*w_s^{-1}\cdot w_sx$. Hence
$$s^{-1}H_0(w_sn^*x)=s^{-1}H_0(w_s\nu w_s^{-1})+s^{-1}H_0(w_sx)
=s^{-1}H_0(w_s\nu w_s^{-1})+H_0(x)$$ for suitable element $\nu\in N_0
(\mathbb A)$. Therefore
$$H_0(x)-T=s^{-1}\Big(H(w_s n^*x)-T\Big)-s^{-1}H_0(w_s\nu w_s^{-1})-
\Big(T-s^{-1}T\Big).$$ If $\varpi\in\hat \Delta_0$, it is well-known that
$\varpi\Big(s^{-1}H_0(w_s\nu w_s^{-1})\Big)$ is non-negative. Moreover,
since $T$ is suitably regular, 
$\varpi\Big(T-s^{-1}T\Big)$ is strictly positive. 
Therefore, for any $\varpi\in\hat\Delta_1\subset\hat\Delta_0$, 
$\varpi\Big(H_0(x)-T\Big)$ is negative, as required.

\subsection{Self-Adjointness}

In this section, we show that $\Lambda^T$ is a self-adjoint operator. 
However,  we must make sure that the intagrations involved are all 
well-defined. With this in mind, we have the following:

\noindent
{\bf {\large Self-Adjointness.}} {\it Suppose that $\phi_1$ and $\phi_2$ are continuous 
functions on $G(F)\backslash G(\mathbb A)^1$. Assume that $\phi$ is of
 moderate growth, and that $\phi_2$ is rapidly decreasing,} in the sense 
that for any $N$, the function $\|x\|^N\cdot|\phi_2(x)|$ is bounded on 
any Siegel set. {\it Then
$$\Big(\Lambda^T\phi_1,\phi_2\Big)=\Big(\phi_1,\Lambda^T\phi_2\Big).$$}

\noindent
{\it Proof.} Since $\phi_2$ is rapidly decreasing, 
the inner product $\Big(\Lambda^T\phi_1,\phi_2\Big)$ is 
defined by an absolutely convergent integral. It equals
$$\begin{aligned}
&\int_{G(F)\backslash G(\mathbb A)^1}\sum_P(-1)^{\mathrm{dim}(A/Z)}
\sum_{\delta\in P(F)\backslash G(F)}\int_{N(F)\backslash N(\mathbb A)}
\phi_1(n\delta x)\hat\tau_P\Big(H(\delta x)-T\Big)\overline{\phi_2(x)}dn\,dx\\
=&\sum_P(-1)^{\mathrm{dim}(A/Z)}\int_{N(F)\backslash N(\mathbb A)}
\int_{P(F)\backslash G(\mathbb A)^1}\phi_1(n x)\overline{\phi_2(x)}
\hat\tau_P\Big(H(x)-T\Big)dx\,dn\\
&\hskip 5.0cm\mathrm{(by\ unfolding\ trick)}\\
=&\sum_P(-1)^{\mathrm{dim}(A/Z)}\int_{N(F)\backslash N(\mathbb A)}
\int_{P(F)\backslash G(\mathbb A)^1}\phi_1(x)\overline{\phi_2(n x)}
\hat\tau_P\Big(H(x)-T\Big)dx\,dn\\
&\hskip 5.0cm\mathrm{ (by\ invariance\ of\ Haar\ measure)}.\end{aligned}$$ This 
last expression reduces to $\Big(\phi_1,\Lambda^T\phi_2\Big)$ as required.

\subsection{$\Lambda^T\phi$ is Rapidly Decreasing}

We would like to show that under suitable conditions, $\Lambda^T\phi(x)$ 
is rapidly decreasing at infinity. 
 
\noindent
{\bf Step A.  Application of Arthur's Partition}. 
To start with, suppose that $\phi$ is a continuous function on 
$G(F)\backslash G(\mathbb A)^1$. Also for any two parabolic subgroups 
$P_1\subset P_2$, set
$$\phi_{P_1,P_2}(x):=\sum_{P:P_1\subset P\subset P_2}(-1)^{\mathrm{dim} 
(A/Z)}\phi_P(y).$$

\noindent
{\bf Lemma 1.} With the same notation as above,
$$\Lambda^T\phi(x)=\sum_{P_1,P_2:P_0\subset P_1\subset P_2}\sum_{\delta\in P_1
(F)\backslash G(F)}F^1(\delta x,T)\cdot\sigma_1^2\Big(H_0(\delta x)-T\Big)\cdot
\phi_{P_1,P_2}(\delta x).$$

\noindent
{\it Proof.} (Arthur) By definition,
$$\begin{aligned}\Lambda^T\phi(x)=&\sum_P(-1)^{\mathrm{dim}(A_P/A_G)}
\sum_{\delta\in P(F)\backslash G(F)}\phi_P(\delta x)
\widehat\tau\Big(H_P(\delta x)-T\Big)\\
=&\sum_P(-1)^{\mathrm{dim}(A_P/A_G)}\sum_{\delta\in P(F)\backslash G(F)}
\Big(\sum_{P_1:P_1\subset P}
F^{P_1}\Big(\delta_1\delta x,T\Big)\cdot\tau_{P_1}^P
\Big(H_{P_1}(\delta_1\delta x)-T\Big)\Big)\\
&\hskip 5.0cm\cdot\widehat\tau\Big(H_P(\delta x)-T\Big)\cdot \phi_P(\delta x)
\end{aligned}$$ 
by using Arthur's partition. We then write $$\phi_P(\delta x)=
\phi_P(\delta_1\delta x)\qquad\mathrm{and}\qquad
\widehat\tau\Big(H_P(\delta x)-T\Big)=
\widehat\tau\Big(H_P(\delta_1\delta x)-T\Big),$$ since both functions are 
left $P(F)$-invariant. Combining the double sum over $\delta$ and 
$\delta_1$ into a single sum over $\delta\in P_1(F)\backslash G(F)$, we
 write $\Lambda^T\phi(s)$ as the sum over pairs $P_1\subset P$ of the 
 product of $(-1)^{\mathrm{dim}(A_P/A_G)}$
with $$\sum_{\delta\in P_1(F)\backslash G(F)}
F^{P_1}\Big(\delta x,T\Big)\cdot\tau_{P_1}^P\Big(H_{P_1}(\delta x)-T\Big)
\widehat\tau\Big(H_P(\delta x)-T\Big)\cdot \phi_P(\delta x).$$
 
Next, consider the product $$\tau_{P_1}^P\Big(H_{P_1}(\delta x)-T\Big)
\widehat\tau_P\Big(H_P(\delta x)-T\Big)
=\tau_{P_1}^P(H_1)\widehat\tau_P(H_1)$$ for the vector 
$H_1=H_{P_1}(\delta x)-T_{P_1}$ in $\frak a_{P_1}$. (We have written
$T_{P_1}$ for the projection of $T$ onto $\frak a_{P_1}$.) 

\noindent
{\bf Claim.} For fixed groups $P_1:P_1\subset P$,
$$\tau_{P_1}^P(H_1)\widehat\tau_P(H_1)=\sum_{P_2,Q:P\subset P_2\subset Q}
(-1)^{\mathrm{dim}(A_{P_2}/A_Q)}
\tau_{P_1}^Q(H_1)\widehat\tau_Q(H_1).$$

\noindent{\it proof of the Claim.} Indeed, for a given pair of parabolic 
subgroups 
$P\subset Q$, the set of $P_2$ with $P\subset P_2\subset Q$ is bijective 
with the 
collection of subsets $\Delta_P^{P_2}$ of $\Delta_P^Q$. Since 
$$(-1)^{\mathrm{dim}(A_{P_2}/A_Q)}=(-1)^{|\Delta_P^Q|-|\Delta_P^{P_2}|},$$ 
the claim follows from the BAsic Lemma of Grand Zero.
We can therefore write
$$\tau_{P_1}^P(H_1)\widehat\tau_P(H_1)=\sum_{P_2:P_2\supset P}
\sigma_1^2(H_1),$$
where we recall that $$\sigma_1^2(H_1):=\sigma_{P_1}^{P_2}(H_1):=
\sum_{Q:Q\supset P_2}(-1)^{\mathrm{dim}(A_{P_2}/A_Q)}
\tau_{P_1}^Q(H_1)\widehat\tau_Q(H_1).$$ This completes the proof.

\noindent
{\bf Step B.} {\bf Compactness}.
For the moment, fix $\delta$ and $x$. Regard $\delta$ as an 
element in $G(F)$ which we are free to left multiply by an element in 
$P_1(F)$. We can therefore assume, as in the proof of idempotence of
 $\Lambda^T$ that $$\delta x=n^*n_*mak$$ where $k\in K, n^*, n_*$ and 
$m$ belong to a fixed compact subsets of $N_1(\mathbb A), N^2(\mathbb A)$ 
and $M_1(\mathbb A)^1$ respectively, and $a$ is an element in 
$A_1(\mathbb R)^0$ with $\sigma_1^2\Big(H_0(a)-T\Big)\not=0$. Therefore
$$\phi_{P_1,P_2}(\delta x)=\phi_{P_1,P_2}(n^*n_*mak)=\phi_{P_1,P_2}(n_*mak)
=\phi_{P_1,P_2}(aa^{-1}n_*mak)=\phi_{P_1,P_2}(ac),$$ where $c$ belongs 
to a fixed compact subset of $G(\mathbb A)^1$ which depends only on $G$.

\noindent
{\bf Step C.} {\bf Harish-Chandra's Trick} In this step, our aim is to 
give an estimation of the function $\phi_{P_1,P_2}$. Roughly speaking, 
such a function is of rapidly decreasing since it is obtained by removing 
constant terms, due to the fact that the structure of the nilpotent is 
rather simple -- infinitesmially, it is simply certain copies of 
$\mathbb R$. Hence we can use a step by step trace of corresponding Fourier 
expansions to arrive at our conclusion. More precisely, it goes as follows.
(However, for simplicity, we here only treat the case when $F=\mathbb Q$
following Arthur. The reader may consult pp.30-34 of [MW] 
for general number fields.)

\noindent
(C.1) {\bf Structure of $N$.} If $\alpha\in\Delta_1^2$, let 
$P_\alpha: P_1\subset P_\alpha\subset P_2$, be the parabolic 
subgroup such that $\Delta_1^\alpha=\Delta_{P_1}^{P_\alpha}$ is 
the complement of $\alpha$ in $\Delta_1^2$.

For each $\alpha$, let $\Big\{Y_{\alpha,1},\cdots, Y_{\alpha,n_\alpha}\Big\}$
 be a basis of $\frak n_\alpha^1(F)$, the Lie algebra of $N_\alpha^2(F)$. 
We shall assume that the basis is compatible with the action of $A_1$ so 
that each $Y_{\alpha,i}$ is a root vector corresponding to the root 
$\beta_{\alpha,i}$ of $(M_2\cap P_1,A_1)$. We shall also assume that if 
$i\leq j$, the height of $\beta_{\alpha,i}$ is not less that the height 
of $\beta_{\alpha,i}$.

Define $\frak a_{\alpha,j}$, $0\leq j\leq n_\alpha$, to be the direct sum 
of $\Big\{Y_{\alpha,1},\cdots, Y_{\alpha,j}\Big\}$ within the Lie algebra 
of $N_2$, and let $N_{\alpha,j}=\exp\frak n_{\alpha,j}$. Then $N_{\alpha,j}$ 
is a normal subgroup of $N_1$ which is defined over $F$.

Now for any subgroup $V$ of $N_1$ defined over $F$, let $\pi(V)$ be the 
operator which sends $\phi$ to $$\int_{V(F)\backslash V(\mathbb A)}
\phi(nx)\,dn,\qquad x\in G(\mathbb A).$$

\noindent
{\bf Sublemma.} With the same notation as above, {\it $\phi_{P_1,P_2}$ is the 
transform of $\phi$ by the product over $\alpha\in\Delta_1^2$ of the operators
$$\pi(N_2)-\pi(N_\alpha)=\sum_{i=1}^{n_\alpha}\pi(N_{\alpha,i-1})
-\pi(N_{\alpha,i}).$$}

\noindent
(C.2) {\bf Application of Fourier Analysis.} 
If $K_0$ is an open compact subgroup of $G(\mathbb A_f)$, $G(F)\backslash 
G(\mathbb A)^1/K_0$ is a differentiable manifold. Assume from now in addition 
that $\phi$ is a function on this space which is {\it differentiable of 
sufficiently high order}. Suppose that $I$ is a collection of indices
$$\Big\{i_\alpha:\alpha\in\Delta_1^2, 1\leq i_\alpha\leq n_\alpha\Big\}.$$
Then $$N_I^-:=\prod_\alpha N_{\alpha,i_\alpha-1}\qquad\mathrm{and}\qquad
N_I:=\prod_\alpha N_{\alpha,i_\alpha}$$ are normal subgroups of $N_1$. 
Let $\frak n^I$ be the span of $\{Y_{\alpha,i_\alpha}\}$ and let 
$\frak n^I(F)'$ be the set of elements $$\xi:=\sum_\alpha r_\alpha 
Y_{\alpha,i_\alpha},\qquad r_\alpha\in F^*.$$ Then if $n$ is any positive 
integer, $$\xi^n:=\prod_\alpha(r_\alpha^n)$$ is a non-zero real number. 
By the Fourier inversion formula for the group $\mathbb A/F$, 
$\phi_{P_1,P_2}(y)$ is the sum over all $I$ of 
$$\sum_{\xi\in \frak n^I(F)'}\int_{\frak n^I(F)\backslash 
\frak n^I(\mathbb A)}dX\cdot
\int_{N_I^-(F)\backslash N^-(\mathbb A)}du\cdot \phi\Big(u e^X y\Big)\cdot
\psi\Big(\langle X,\xi\rangle\Big).$$ Here $e^*$ denote the exponential 
map while $\psi$ denotes an additive character over $\mathbb A$ used as a 
reference for the Fourier transform, and $\langle\cdot,\cdot\rangle$ is the 
inner product defined by our basis on $\frak n^I$. Clearly, if $n$ is a 
positive integer, then
$$Y_I^n:=\prod_\alpha\Big(-\sqrt{-n}Y_{\alpha,i_\alpha}\Big)^n$$ can be 
regarded as an element in $\mathcal U\Big(\frak g(\mathbb R)^1
\otimes\mathbb C\Big)$. In this way, we have shown the following

\noindent
{\bf Lemma 2.} With the same notation as above, $$\begin{aligned}
&\phi_{P_1,P_2}(y)\\
=&
\sum_I\sum_{\xi\in \frak n^I(F)'}(\xi^n)^{-1}\int_{\frak n^I(F)
\backslash \frak n^I(\mathbb A)}dX\cdot
\int_{N_I^-(F)\backslash N^-(\mathbb A)} R_y(\mathrm {Ad}
(y^{-1})Y_I^n\Big)\cdot\phi\Big(u e^X y\Big)\cdot
\psi\Big(\langle X,\xi\rangle\Big)\,du.\end{aligned}$$

Now set $$y=\delta x=ac$$ as above. Since $\sigma_1^2\Big(H_0(a)-T\Big)
\not=0$, $a$ belongs to a fixed Siegel set in $M_2(\mathbb A)$. It 
follows that the integrand above, as a function of $X$, is invariant 
by an open compact subgroup of $\frak n^I(\mathbb A_f)$ which is 
independent of $a$ and $c$. Consequently, the involved {\it integration 
vanishes unless $\xi$ belongs to a fixed lattice, $L^I(K_0)$, in 
$\frak n^I(\mathbb R)$.}

But for $n$ sufficiently large $$\sum_{\xi\in\frak n^I(F)'
\cap L^I(K_0)}|\xi^n|^{-1}$$ is finite for all $I$. Let $c_n(K_0)$ be
 the supremum over all $I$ of these numbers. Then 
$\Big|\phi_{P_1,P_2}(ac)\Big|$ is bounded by
$$c_n(K_0)\sum_I\int_{N_I(F)\backslash N_I(\mathbb A)}
\Big|\bigg(R\Big(\mathrm {Ad}(c)^{-1}\mathrm {Ad}(a)^{-1}Y_I^n\Big)\phi\bigg)
\Big(u ac\Big)\Big|du.$$ Let $\beta_I=\sum_\alpha\beta_{\alpha,i_\alpha}$.
 Then $\beta_I$ is a positive sum of roots in $\Delta_1^2$. For any $n$,
$$\mathrm {Ad}(a)^{-1}Y_I^n=e^{-n\beta_I(H_0(a))Y_I^n}=e^{-n\beta_I
(H_0(\delta x))Y_I^n}.$$ 

We can {\it choose a finite set of elements} $\{X_i\}$ in 
$\mathcal U\Big(\frak g(\mathbb R)^1\otimes\mathbb C\Big)$, depending 
only on $n$ and $K_0$, such that for any $P_1,P_2, I$ and $c$, such 
that $$c_n(K_0)\mathrm {Ad}(c)^{-1}Y_I^n$$ is a linear combination of 
$\{X_i\}$. 

Since $c$ lies in a compact set, we may assume that each of the coefficients
has absolute value less than 1. We have thus far shown the following

\noindent
{\bf Lemma 3.} With the same notation as above,
{\it $\Big|\Lambda^T\phi(s)\Big|$ is bounded by the sum over all $P_1,P_2$ 
and $\delta\in P_1(F)\backslash G(F)$ of the product of
$$F^1(\delta x,T)\sigma_1^2\Big(H_0(\delta x)-T\Big)$$ with $$\sum_I\sum_i
\int_{N_I(F)\backslash N_I(\mathbb A)}\Big|R(X_i)\phi(u\delta x)\Big|
du\cdot e^{-n\beta_I(H_0(\delta x))}.\eqno(1)$$}

\noindent
{\bf Step D.} Finally we are ready to tackle the problem of rapidly decreasing
property of $\Lambda^T\phi$. However as a bit more general version 
is ready to be stated without any additional work, we give the following
\vskip 0.30cm
\noindent
{\bf Theorem.} (Arthur) Let $\frak s$ be a Siegel set in $G(\mathbb A)^1$.
 For any pair of positive integers $N'$ and $N$, and any open compact 
subgroup $K_0$ of $G(\mathbb A_f)$, we can choose a finite subset
 $\{X_i\}$ of $\mathcal U\Big(\frak g(\mathbb R)^1\otimes\mathbb C\Big)$, 
and a positive integer $r$ which satisfy the following property: Suppose 
that $(S,d\sigma)$ is a measure space and that $\phi(\sigma,x)$ is a 
measurable function from $S$ to 
$C^r\Big(G(F)\backslash G(\mathbb A)^1/K_0\Big)$. Then
{\it for any $x\in\frak s$,
$$\int_S\Big|\Lambda^T\phi(\sigma,x)\Big|\,d\sigma$$
is bounded by $$\sum_i\sup_{y\in G(\mathbb A)}\Big(\int_S\Big|R(X_i)\phi
(\sigma,y)\Big|d\sigma\cdot\|y\|^{-N}\Big)\cdot \|x\|^{-N'}.$$}

\noindent
{\it Proof.} Substitute $\phi(\sigma,\cdot)$ for $\phi$ in (1) and integrate 
over $\sigma$. The result is $$\sum_I\sum_i
\int_{N_I(F)\backslash N_I(\mathbb A)}\int_S\Big|R(X_i)\phi(u\delta x)
\Big|d\sigma\,
du\cdot e^{-n\beta_I(H_0(\delta x))}.\eqno(2)$$ If $\delta x=ac$, with 
$a$ and $c$ as above, $$\|\delta x\|\leq\|a\|\cdot\|c\|.$$ We are assuming 
that $\sigma_1^2\Big(H_0(a)-T\Big)\not=0$. With this condition, we next 
want to give a natural bound for $\|a\|$. For this we need the following

\noindent
{\bf Sublemma.} Fix $T\in\frak a_0^+$, let 
$H\in\frak a_1^G=\frak a_1^2\oplus\frak a_2^G$ with corresponding decomposition
$H=H_1^2+H_2$. {\it If $\sigma_1^2\Big(H-T\Big)\not=0$, then}

\noindent
(i) {\it $\alpha(H_1^2)$ is positive for each $\alpha\in\Delta_1^2$;} and

\noindent
(ii) {\it $\|H\|\leq c\Big(1+\|H_1^2\|\Big)$ for any Euclidean norm 
$\|\ \|$ on $\frak a_0$ and some constant $c$.}

\noindent
{\it Proof.} The first condition follows directly from  
the characterization of $\sigma_1^2$ and the fact that $T\in\frak a_0^+$.
To prove the second one, note that the value at $H_2$ of any root in 
$\Delta_2$ equals $\alpha(H_2)$ for some root 
$\alpha\in\Delta\backslash\Delta_1^2$. But
$$\alpha(H_2)=\alpha(H-T)-\alpha(H_1^2)+\alpha(T)<-\alpha(H_1^2)+\alpha(T),$$ 
bythe characterization of $\sigma_1^2$ again. Since 
$$\varpi(H_2)=\varpi(H)>0\qquad\forall\varpi\in\widehat\Delta_2,$$ $H_2$ 
belongs to a compact set. In fact the norm of $H_2$ is bounded by a 
constant multiple of $1+\|H_1^2\|$ as required. This completes the proof 
of the Sublemma.

Since $\beta_I$ is a positive sum of roots in $\Delta_1^2$ we conclude from 
\S1 that $\|a\|$ is bounded by a fixed power of $$e^{\beta_I(H_0(a))}=
e^{\beta_I(H_0(\delta x))}.$$ It follows that for any positive integres 
$N$ and $N_1$ we may choose $n$ so that (2) is bounded by a constant multiple
 of $$\sum_i
\sup_{y\in G(\mathbb A)^1}\Big(\int_S\Big|R(X_i)\phi(u\delta y)\Big|d\sigma
\cdot\|y\|^{-N'}\Big)\|\delta x\|^{-N_1}.$$
It is well-known that there is a constant $c_1$ such that for any 
$\gamma\in G(F)$ and $c\in\frak s$, $$\|\gamma x\|^{-N_1}\leq 
c_1\|x\|^{-N_1}.$$

\noindent
{\bf Step E}. As such, the only thing left to estimate is
$$\sum_{\delta\in P_1(F)\backslash G(F)}F^1(\delta x,T)\sigma_1^2
\Big(H_0(\delta x)-T\Big).$$ The summand is the characteristic function, 
evaluated at $\delta x$, of a certain subset of $$\Big\{y\in G(\mathbb A)^1:
\varpi\Big(H_0(y)-T\Big)>0,\varpi\in\hat\Delta_1\Big\}.$$
The sum is bounded by $$\sum_{\delta\in P_1(F)\backslash G(F)}\hat\tau_1
\Big(H_0(\delta x)-T\Big).$$
It follows from the Basic Estimation 
that we can find constants $C_2$ and $N_2$ such that
 for all $P_1$ this last expression is bounded by $C_2\|x\|^{N_2}$. Set 
$N_1=N'+N_2$. $N_1$ dictates our choice of $n$, from which we obtain the 
differential operators $\{X_i\}$. The theorem follows with any $r$ greater 
than all the degree of the operators $X_i$.

\chapter {Arthur's Period and Abelian Part of Non-Abelian L-Functions}
\section {Arthur's Periods}
\vskip 0.30cm
Let $G$ be a reductive group defined over a number field $F$. Fix a 
minimal $F$-parabolic subgroup $P_0$ with minimal $F$-Levi component $M_0$.
 Let $P$ be a 
parabolic subgroup of $G$ with $M$ the Levi component and $N$ the unipotent
 radical. Fix $T\in \mathrm{Re}(\frak a_0)$, and assume that $T$ is 
sufficiently regular in the sense that $\langle\alpha,T\rangle$ are 
sufficiently large for all $\alpha\in\Delta_0$.
Then, following Arthur, for any locally $L^1$ function 
$\phi:G(F)\backslash 
G(\mathbb A)\to \mathbb C$,  we define a locally $L^1$ function 
$\Lambda^T\phi$, the {\it analytic truncation} of $\phi$ 
with parameter $T$, on $G(F)\backslash G(\mathbb A)$ by
$$\Lambda^T\phi(g):=\sum_{P=MN; P\supset P_0}
(-1)^{\mathrm{rk}(G)-\mathrm{rk}(M)}\sum_{\gamma\in P(F)\backslash G(F)}
\widehat \tau_P\Big(\log_M(m_P(\gamma g))-T_M\Big)\cdot \phi_P(\gamma g),$$
where as usual, $T_M$ denote the projection of $T$ onto 
$\mathrm{Re}(\frak a_M)$ with respect to the decomposition
$\frak a_0=\frak a_0^M\oplus\frak a_M$, and $\phi_P$ denotes the 
constant term of $\phi$ with respect to $P$, i.e., $$\phi_P(g):=
\int_{N(F)\backslash N(\mathbb A)}\phi(gn)\,dn.$$
\vskip 0.30cm
Fundamental properties of Arthur's analytic truncation, as proved in Ch. 4 
may be summarized as follows:

\noindent
{\bf {\Large Theorem.}} (Arthur) (1) {\it Let $\phi:G(F)\backslash 
G(\mathbb A)\to\mathbb C$ be a locally $L^1$ function. Then we have 
$$\Lambda^T\Lambda^T\phi(g)=\Lambda^T\phi(g)$$ for almost all $g$. 
If $\phi$ is also locally bounded,  then the above
is true for all $g$};

\noindent
(2) {\it Let $\phi_1,\,\phi_2$ be two locally $L^1$ functions on 
$G(F)\backslash G(\mathbb A)$. Suppose that 
$\phi_1$ is of moderate growth and $\phi_2$ is  rapidly decreasing. Then
$$\int_{Z_{G(\mathbb A)}G(F)\backslash G(\mathbb A)}\overline{\Lambda^T
\phi_1(g)}\cdot \phi_2(g)\,dg
=\int_{Z_{G(\mathbb A)}G(F)\backslash G(\mathbb A)}\overline{\phi_1(g)}\cdot
\Lambda^T\phi_2(g)\,dg;$$}

\noindent
(3) {\it Let $K_f$ be an open compact subgroup of $G(\mathbb A_f)$, and 
$r, r'$ are 
two positive real numbers. Then there exists a finite subset 
$\Big\{X_i:i=1,2,\cdots,N\Big\}\subset\mathcal U$, the universal enveloping 
algebra of $\frak g_\infty$, such that the 
following is satisfied: Let $\phi$ be a smooth function on 
$G(F)\backslash G(\mathbb A)$, right invariant under $K_f$ and let 
$a\in A_{G(\mathbb A)},\ g\in G(\mathbb A)^1\cap S$. Then 
$$\Big|\Lambda^T\phi(ag)\Big|\leq \|g\|^{-r}\sum_{i=1}^N
\sup\Big\{|\delta(X_i)\phi(ag')|\,\|g'\|^{-r'}:g'\in G(\mathbb A)^1\Big\},$$
where $S$ is a  Siegel domain with respect to $G(F)\backslash G(\mathbb A)$.}

\noindent
{\it Sketch of the Proof.}  Simply put, 
(1) says that $\Lambda^T$ is idempotent; 
(2) says that $\Lambda^T$ is self-adjoint; and
(3) says that if $\phi$ is of moderate growth, then 
$\Lambda^T$ is rapidly decreasing. 
Thus in particular, (1) and (2) implies 
that $\Lambda^T$ is an orthogonal projection while (3) implies that the 
integration $\int_{G(F)\backslash G(\mathbb A)}\Lambda^T\phi(g)\,dg$ is 
well-defined for any automorphic form $\phi$. 

As fora proof, the reader may find all details form the previous chapter. 
Roughly, (1) is based on the following 
{\it miracle}, which is much stronger than what is actually need: 
{\it For a fix parabolic subgroup $P_1$, unless $\varpi(H_0(x)-T)<0$ for all 
$\varpi\in\widehat \Delta_1$, the constant term of $\Lambda^T\phi$ along 
$P_1$ is identically zero.} 

\noindent
(2) is merely a manipulation of the definition.  

\noindent
(3) is a bit involved. To explain it, for $P_1\subset P_2$ two
parabolic subgroups, set
$$\phi_{12}(g):=\sum_{P:P_1\subset P\subset P_2}(-1)^{d(P)-d(G)}\phi_P(g)$$  
be the alternating sum of constant terms along parabolic subgroups $P$ 
between $P_1$ and $P_2$. Then,

\noindent
(i) with the help of Langlands' Combinatorial Lemma, or better, Arthur's
partition, $\Lambda^T\phi$ may be 
rewritten as $$\sum_{P_1\subset P_2}\sum_{\gamma\in P_1(F)
\backslash G(F)} \Big(F_G^1(\gamma g, T)\cdot\sigma_1^2(H(\gamma g)-T)\Big)
\cdot \phi_{12}(\gamma g)$$ where the definitions of the characteristic 
functions $F_G^1$ and $\sigma_1^2$ is the same as in 4.1.4 and 4.1.3. 
Indeed, this 
step may be understood as a refined version of classical reduction theory in 
which the combinatorial lemma plays a key role;

\noindent
(ii) The characteristic function $\sum_{\gamma\in P_1(F)\backslash G(F)}
\Big(F_G^1(\gamma g, T)\cdot\sigma_1^2(H(\gamma g)-T)\Big)$ may be 
naturally bounded by the summation $$\sum_{\gamma\in P_1(F)\backslash 
G(F)}\widehat\tau_1\Big(H(\gamma g)-T\Big)$$ which is known to be a finite 
summation and 
bounded from above by $c\|x\|^{C}$ for certain constants $c>0$ and $C$; 
(While  simple, this latest upper bound, the Basic Estimation, is indeed very 
basic to Arthur's analytic truncation.) And 

\noindent
(iii) Natural upper bound for $\phi_{12}$: For any $N>0$, there exists a 
constant $c>0$ such that for all $\gamma\in G(F), g\in G(\mathbb A)^1$, 
if $F_G^1(\gamma g, T)\cdot\sigma_1^2(H(\gamma g)-T)=1$, then 
$$\Big|\phi_{12}(\gamma g)\Big|\leq C\|g\|^{-N}.$$ This is basically 
a direct generalization of the fact that
after cutting off the constant term in the Fourier expansion for 
functions over $\mathbb R$, the 
remainning part is rapid decreasing. However,  the situation at hands is 
much more complicated: We are not just working over a one dimensional 
affine space, but the Lie of the unipotent radicals. Hence Arthur uses 
a trick of Harish-Chandra, albert in a gereralized version to establish 
the above inequality.
\vskip 0.30cm
With this, now we are ready to define {it Arthur's period} 
associated to an automorphic form $\phi$ on $G(F)
\backslash G(\mathbb A)$ by the integral $$A(\phi;T):=
\int_{G(F)\backslash G(\mathbb A)}\Lambda^T\phi(g)\,dg.$$ From the 
above definition, for sufficiently regular $T\in\frak a_0$, 
$A(\phi;T)$ is well-defined, since being an automorphic 
form, $\phi$ is of moderate growth, so $\Lambda^T$ is rapidly 
decreasing, and hence integrable.
 
\section{Arthur's Periods as Integrals over Truncated Fundamental Domains}
  
Next we want to express Arthur's period $A(\phi;T)$ as an integration of 
the original automorphic form $\phi$ over a compact subset. 
This goes as follows:

To start with, note that for Arthur's analytic truncation 
$\Lambda^T$, we have $\Lambda^T\circ \Lambda^T=\Lambda^T$. Hence,
$$A(\phi;T)=\int_{Z_{G({\mathbb A})}G(F)\backslash G({\mathbb A})}\Lambda^T
 \phi\ d\mu(g)=\int_{Z_{G({\mathbb A})}G(F)\backslash G({\mathbb A})}
\Lambda^T\Big(\Lambda^T \phi \Big)(g)\ d\mu(g).$$ Moreover, by the 
self-adjoint property, we have, for the constant function $\bold 1$ on 
$G(\mathbb A)$, 
$$\begin{aligned}&\int_{Z_{G({\mathbb A})}G(F)
\backslash G({\mathbb A})}{\bf 1}(g)\cdot 
\Lambda^T\Big(\Lambda^T \phi \Big)(g)\ d\mu(g)\\
=&\int_{Z_{G({\mathbb A})}G(F)\backslash G({\mathbb A})}
\Big(\Lambda^T {\bf 1}\Big)(g)\cdot
\Big(\Lambda^T \phi \Big)(g)\ d\mu(g)\\
=&\int_{Z_{G({\mathbb A})}G(F)\backslash G({\mathbb A})}
\Lambda^T\Big(\Lambda^T {\bf 1}\Big)(g)\cdot
\phi(g)\ d\mu(g),\end{aligned}$$ 
since $\Lambda^T\phi$ and $\Lambda^T{\bf 1}$ are 
rapidly decresing. Therefore, using $\Lambda^T\circ \Lambda^T=\Lambda^T$ again,
we arrive at
$$A(\phi;T)=\int_{Z_{G({\mathbb A})}G(F)\backslash G({\mathbb A})}
\Lambda^T{\bold 1}(g)\cdot \phi(g)\ d\mu(g).\eqno(1)$$ 

To go further, let us give a much more detailed study of Authur's analytic 
truncation for the constant function ${\bold 1}$.
Introduce the truncated subset
$\Sigma(T):=\Big(Z_{G({\mathbb A})}G(F)\backslash G({\mathbb A})\Big)_T$ of 
the space $G(F)\backslash G(\mathbb A)^1$. That is,
$$\Sigma(T):=\Big(Z_{G({\mathbb A})}G(F)\backslash G({\mathbb A})\Big)_T:=
\Big\{g\in Z_{G({\mathbb A})}G(F)\backslash G({\mathbb A}):
\Lambda^T{\bf 1}(g)\not=0\Big\}.$$ We claim that
$\Sigma(T)$ or the same 
$\Big(Z_{G({\mathbb A})}G(F)\backslash G({\mathbb A})\Big)_T$, is compact. In fact, much stronger result is correct. Namely, we have the following

\noindent
{\bf {\bf large Basic Fact.}}  (Arthur) With the same notation as above, 
$$\Lambda^T{\bold 1}(x)=F(x,T).$$ That is to say,
{\it $\Lambda^1{\bold 1}$ is the characteristic function of the compact subset
$\Sigma(T)$ of $G(F)\backslash G(\mathbb A)^1$ obtained by projecting 
$$N_0(\mathbb A)\cdot M_0(\mathbb A)^1\cdot A_{P_0,G}^\infty(T_0,T)\cdot K$$ 
onto $G(F)\backslash G(\mathbb A)^1.$} 

\noindent
{\bf Remark.} The discussion in the first version of this mauscript
is a bit different: at that time, I was not so sure on whether 
$\Lambda^T{\bold 1}$ is a characteristic function.
I thank Arthur for informing me 
that the above result is proved in [Ar5].
\vskip 0.30cm
Before proving this fundamental result, let me continue our discussion on 
Arthur's period. From relation (1),  $$
\int_{Z_{G({\mathbb A})}G(F)\backslash G({\mathbb A})}\Lambda^T\phi(g)\, 
d\mu(g)=
\int_{Z_{G({\mathbb A})}G(F)\backslash G({\mathbb A})}
\Lambda^T{\bold 1}(g)\cdot\phi(g)\, d\mu(g),$$ which, by the Basic Fact above,
is nothing but
$$\int_{\Sigma(T)}\phi(g)\, d\mu(g).$$ 
That is to say, we have obtained the following very beautiful relation:

\noindent
{\bf {\large Modified Geometric Truncation
=Analytic Truncation}:}  {\it For a sufficiently regular $T\in\frak a_0$, 
and an automorphic form $\phi$ on $G(F)\backslash G(\mathbb A)$,
$$\int_{\Sigma(T)}\phi(g)\, 
d\mu(g)=\int_{G(F)\backslash G({\mathbb A})^1}\Lambda^T\phi(g)\, d\mu(g).$$} 

As such, accordingly, we may have two different
ways to study Arthur's periods. Say, when taking $\phi$ to be Eisenstein series
associated with an $L^2$ automorphic forms:  
\vskip 0.30cm
\noindent
(a) ({\it Analytic Evaluation for RHS}) 
The right hand side can be evaluated using certain analytic methods. For 
example, 
in the case for Eisensetin series associated to cusp forms, as done in the next
chapter, the RHS is evaluated in terms of integrations
of cusp forms twisted by intertwining 
operators, which normally are close related with abelian $L$-functions;
\vskip 0.30cm
\noindent
(b) ({\it Geometric Evaluation for LHS}) 
The left hand, by taking the residue (for suitable Eisenstein series), 
can be eveluated using geometric methods totally independently. For example,
in the so-called rank one case, essentially
the LHS is the volume of the geometrically truncated domain, i.e., 
the difference
between the volume of the total fundamental domain, which can be 
written in terms of some special values of abelian $L$-functions,
 thanks to the fundamental work of Siegel [S], and a certain combination of
volumes of cuspidal regions corresponding to proper parabolic subgroups, 
which are supposed to be easily calculated thanks to the simple structure 
of cuspidal regions.@
\vskip 0.30cm
We end this section by offering a proof of the Basic Fact. This is based on
Arthur's partition for $G(F)\backslash G(\mathbb A)$ and the inversion 
formula, 
which itself is a direct consequence of Langlands' Combinatorial lemma.

Recall that for $T\in\frak a_0$, level $P$-{\it Arthur's analytic truncation} 
of $\phi$ is defined by the formula
$$\Lambda^{T,P}\phi(g):=\sum_{R:R\subset P}(-1)^{d(R)-d(P)}
\sum_{\delta\in R(F)\backslash P(F)}\phi_R(\delta g)\cdot
\widehat\tau_R^P\Big(H(\delta g)-T\Big).$$ Then $\Lambda^T$ stands simply
for $\Lambda^{T,G}$. Thus
\vskip 0.30cm
\noindent
{\bf (i)} $\displaystyle{\Lambda^{T}\phi(g):=\sum_{R}(-1)^{d(R)-d(G)}
\sum_{\delta\in R(F)\backslash G(F)}\phi_R(\delta g)\cdot
\widehat\tau_R\Big(H(\delta g)-T\Big)};$ and

\noindent
{\bf (ii)} ({\bf Inversion Formula}) {\it For a $G(F)$-invariant 
function $\phi$,} 
$$\phi(g)=\sum_P\sum_{\delta\in P(F)\backslash G(F)}\Lambda^{T,P}
\phi(\delta g)\cdot \tau_P\Big(H(\delta g)-T\Big).$$ 
\vskip 0.30cm
\noindent
{\it Proof.} Indeed, recall that we have the following
\vskip 0.30cm
\noindent
{\bf Langlands' Combinatorial Lemma.} {\it If $Q\subset P$ are parabolic 
subgroups, 
then for all $H\in \frak a_0$, we have
$$\sum_{R:Q\subset R\subset P}(-1)^{d(R)-d(P)}\tau_Q^R(H)\widehat\tau_R^P(H)=
\delta_{QP},$$ and 
$$\sum_{R:Q\subset R\subset P}(-1)^{d(Q)-d(R)}\widehat\tau_Q^R(H)\tau_R^P(H)=
\delta_{QP}.$$}

Thus in particular, the RHS of (ii) is simply
$$\begin{aligned}&\sum_P\sum_{\delta\in P(F)\backslash G(F)}\bigg(
\sum_{R:R\subset P}(-1)^{d(R)-d(P)}\sum_{\gamma\in R(F)\backslash P(F)}
\phi_R(\gamma \delta g)\cdot\widehat\tau_R^P\Big(H(\gamma\delta g)-T\Big)\bigg)
\cdot \tau_P\Big(H(\delta g)-T\Big)\\
=&\sum_{R,P:R\subset P}\sum_{\delta\in P(F)\backslash G(F)}
\sum_{\gamma\in R(F)\backslash P(F)}\phi_R(\gamma\delta g)(-1)^{d(R)-d(P)}
\widehat\tau_R^P\Big(H(\gamma\delta g)-T\Big)\cdot 
\tau_P\Big(H(\delta g)-T\Big)\\
=&\sum_{R,P:R\subset P}\sum_{\delta\in R(F)\backslash G(F)}
\phi_R(\delta g)(-1)^{-d(R)+d(P)}\widehat\tau_R^P\Big(H(\delta g)-T\Big)
\cdot \tau_P\Big(H(\delta g)-T\Big)\\
=&\sum_{R}\sum_{\delta\in R(F)\backslash G(F)}
\phi_R(\delta g)\sum_{P:R\subset P}
(-1)^{-d(R)+d(P)}\widehat\tau_R^P\Big(H(\delta g)-T\Big)
\cdot \tau_P\Big(H(\delta g)-T\Big)\\
=&\sum_{R}\sum_{\delta\in R(F)\backslash G(F)}
\phi_R(\delta g)\sum_{P:R\subset P}
(-1)^{d(P)-d(R)}\widehat\tau_R^P\Big(H(\delta g)-T\Big)
\cdot\tau_P\Big(H(\delta g)-T\Big)\\
=&\sum_{R}\sum_{\delta\in R(F)\backslash G(F)}
\phi_R(\delta g)\cdot \delta_{RG}\end{aligned}$$ 
by the second relation in Langlands' combinatorial lemma. 
But this latest quantity is clearly $\phi(g)$. This completes the proof.
\vskip 0.30cm
With all this done, now we are ready to give the following:

\noindent
{\it Proof of Basic Relation.} (Arthur) By Arthus's partition for 
$G(F)\backslash G(\mathbb A)$, we have 
$$\sum_P\sum_{\delta\in P(F)\backslash G(F)}
F^P\Big(\delta x,T\Big)\cdot\tau_P\Big(H_P(\delta x)-T\Big)=1$$ where
$\tau_P$ is the characteristic function of 
$$\Big\{H\in \frak a_0:\alpha(H)>0,
\ \alpha\in\Delta_P\Big\}$$ and $$F^P\Big(nmk,T\Big)=F^{M_P}\Big(m,T\Big),
\qquad n\in N_P(\mathbb A), m\in M_P(\mathbb A), k\in K.$$
On the other hand, by the inversion formula, applying to the constant 
function ${\bold 1}$, we have
$$\sum_P\sum_{\delta\in P(F)\backslash G(F)}
\Big(\Lambda^{T,P}{\bold 1}\Big)(\delta x)\cdot\tau_P
\Big(H_P(\delta x)-T\Big)=1,$$
where $\Lambda^{T,P}$ is the partial analytic truncation operator 
defined above.With this, the desired result is immediately obtained by 
induction.

\section{Geometrically Oriented Analytic Truncation}

Motivated by Arthur's analytic truncation [Ar3] recalled above,
 Lafforgue's corresponding work on function fields [Laf] and the 
discussion in Chapter 2,  for  a normalized polygon $p:[0,r]\to\mathbb R$,
we introduce a geometrically oriented  analytic truncation as follows:

For any locally $L^1$ function $\phi:G(F)\backslash G(\mathbb A)\to 
\mathbb C$,  we define a locally $L^1$ function $\Lambda^T\phi$, the 
{\it modified analytic truncation} of $\phi$ associated with the polygon
$p$ on $G(F)\backslash G(\mathbb A)$ by
$$\Lambda_p\phi(g):=
\sum_P(-1)^{\mathrm{dim}(A_P/Z_G)}\sum_{\delta\in 
P(F)\backslash G(F)}\phi_P(\delta g)\cdot{\bold 1}
(p_P^{\delta g}>_Pp).$$
 
Moreover, recall that we also have the relation 
$${\bold 1}\Big(p_P^{g}\triangleright_Pp\Big)
=\tau_P\Big(-\log_M(m_P(g))-T(p)\Big).$$ 
Formally, we can also introduce the following modified
analytic truncation
$$\begin{aligned}\Lambda_{T(p)}\phi(g):=
&\sum_{P=MN; P\supset P_0}(-1)^{\mathrm{rk}(G)-
\mathrm{rk}(M)}\sum_{\gamma\in P(F)\backslash G(F)}
\tau_P\Big(-\log_M(m_P(\gamma g))-T(p)_M\Big)\cdot \phi_P(\gamma g)\\
=&
\sum_P(-1)^{\mathrm{dim}(A_P/Z_G)}\sum_{\delta\in 
P(F)\backslash G(F)}\phi_P(\delta g)\cdot{\bold 1}
(p_P^{\delta g}\triangleright_Pp),\end{aligned}$$
where as usual, $T(p)_M$ denote the projection of $T(p)$ onto 
$\mathrm{Re}(\frak a_M)$ with respect to the decomposition
$\frak a_0=\frak a_0^M\oplus\frak a_M$, and $\phi_P$ denotes the 
constant term of $\phi$ with respect to $P$, i.e., $$\phi_P(g):=
\int_{N(F)\backslash N(\mathbb A)}\phi(gn)\,dn.$$

That is to say, in these new trunctions, instead of the 
characteristic function $\widehat \tau_P$ associated with positive cones 
as used in the original Arthur analytic truncation, we use the 
characteristic functions ${\bold 1}(p_P^{*}>_Pp)$ and
a much more restricted characteristic function ${\bold 1}
(p_P^{*}\triangleright_Pp)=\tau_P\Big(-H_P(*)-T(p)_M\Big)$ 
associated with positive chambers. 
In the follows, we will show that while
the first new geo-analytic trunction beautifully makes extra 
rooms for housing what we are going to call the essential non-abelian 
part of non-abelian zeta functions, the second one plays no role as 
with such a restricted truncation the resulting functions are generally not
of rapidly decreasing (and hence may even not be well-defined).
 
\section{Abelian and Essential Non-Abelain Parts}

To see this, naturally, motivated by the above discussion on 
Arthur's analytic truncation, at the beginning, we may try to compare
our new truncations with that of Arthur.
Say, we may expect that the following nice properties hold 
for the geo-analytic truncation $\Lambda_p$ as well. That is to say, 
we expect that:

\noindent
(1) $\Lambda_p\circ\Lambda_p=\Lambda_p$;

\noindent
(2) $\Lambda_p$ is self-adjoint;

\noindent
(3) $\Lambda_p\phi$ for automorphic forms $\phi$ are all rapidly decreasing.

However, the life is not that simple: For Arthur's analytic truncation, 
similar properties are established with lucky and great care: Say to show (1), 
one needs a miracle, while to get (3), we need a certain refined reduction 
theory. 

Accordingly, we will not try to pursue this direction any further, 
as our final aim is quite different.

Instead, we are trying to explore the differnce between the geometrically 
oriented integration and the analytically oriented integration. 
As the discussions are parallel,
we here only give the details for $\Lambda_p$, while leave the case
for $\Lambda_T$ to the reader. (See however Ch.6.)

First of all, for a fixed normalized polygon $p$, we do 
not really have the property that $\Lambda_p\circ\Lambda_p=\Lambda_p$.
Therefore, even if the integration 
$$\int_{Z_{G({\mathbb A})}G(F)\backslash G({\mathbb A})}\Lambda_p
 \phi(g)\, d\mu(g)$$ is well-defined, 
we cannot say that such an integration is equal to the integration
$$\int_{Z_{G({\mathbb A})}G(F)\backslash G({\mathbb A})}\Lambda_p{\bf 1}(g)
\cdot\phi(g)\, d\mu(g);$$

Secondly, we need to show that $\Lambda_p{\bf 1}$ is a characteristic function.
For this, we are extremely lucky:
By the Fundamental Relation established in section 2.4,
$$\bold 1(\overline p^g\leq p)=\sum_{P:\, \text{standard\, 
parabolic}}(-1)^{|P|-1}\sum_{\delta\in
P({\mathbb Z})\backslash G({\mathbb Z})}\bold 1(p_P^{\delta g}>_Pp)
\qquad\forall g\in G.$$
Therefore, 
$$\bold 1(\overline p^g\leq p)=\Lambda_p{\bold 1}(g)\qquad\forall g\in G.$$
Accordingly, set 
$\Big(Z_{G({\mathbb A})}G(F)\backslash G({\mathbb A})\Big)_p$ 
be the subset in the fundamental domain obtained by applying our new 
geo-analytic truncation to the constant function ${\bold 1}$. That is to say, 
$$\Big(Z_{G({\mathbb A})}G(F)\backslash G({\mathbb A})\Big)_p:=\Big\{g\in 
Z_{G({\mathbb A})}G(F)\backslash G({\mathbb A}):\Lambda_p{\bf 1}(g)=1\Big\}.$$
Then,  by the Main Theorem in Ch.2, we have the folowing
\vskip 0.30cm
{\bf {\large Basic Fact}}. {\it There is a natural identification} 
$$\Big(Z_{G({\mathbb A})}G(F)\backslash G({\mathbb A})\Big)_p
\simeq \mathcal M_{F,r}^{\leq p}\Big[\Delta_F^{\frac{r}{2}}\Big].$$ 

\noindent
In particular, it is  compact. Consequently, the integration 
$$\int_{\Big(Z_{G({\mathbb A})}G(F)\backslash G({\mathbb A})\Big)_p}
 \phi\, d\mu(g)$$ is well-defined, which is simply our non-abelian 
$L$-function if $\phi$ is replaced by Eisenstein series associated with 
$L^2$-automorphic forms.

As such, obviously, there is a discrepency
between $$\int_{\Big(Z_{G({\mathbb A})}G(F)\backslash G({\mathbb A})\Big)_p}
 \phi(g)\, d\mu(g)\eqno(*)$$ and 
$$\int_{Z_{G({\mathbb A})}G(F)\backslash G({\mathbb A})}
\Lambda_p \phi(g)\, d\mu(g),\eqno(**)$$ even if the integration (**) is 
well-defined.
 Such a discrepency results what we 
call {\it essential non-abelian part} of non-abelin $L$ functions. 
 
Indeed, as said above, by definition, if we take $\phi$ to be 
Eisenstein series associated to 
$L^2$-automorphic forms, then
(*) represents our general non-abelian $L$-function. On the other hand,
(**) is a kind of modified period. Moreover, as we will see in the following 
chapter, (**) may be evaluated by Rankin-Selberg method, or better, 
the regularized
integration in the sense of Jacquet-Lapid-Rogawski as a further development of
Zagier's regularization. 
It is for this reason that we will call 
the integration (**) as an {\it abelian} one, 
while call the difference between (*) and (**) as an 
{\it essential non-abelian} one. In particular, in the case when
$\phi$ is an Eisenstein series associated with an
$L^2$-automorphic form, we will call the corresponding integration 
(**) as the {\it abelian part} of the non-abelian $L$-function, 
while call the difference between (*) and (**) as the 
{\it essential non-abelian part} of the non-abelian $L$-function.

\chapter{Rankin-Selberg Method: Evaluation of Abelian Part}

In this chapter, we mainly concentrate on abelian part of our 
non-abelian $L$-functions,
using a generalized Rankin-Selberg method, so as to give a precise expression 
of these abelian part. However, our success is quite limited: We only can 
write these abelian parts down for
non-abelian $L$-functions associated with cuspidal automorphic forms. 

In general, there are two types of difficulties in studying abelian parts: 
the convergence problems and the precise expression of constant terms of 
Eisenstein series. For the first one, what we have done here is to indicate 
how a parallel discussion can be carried out for our more geometrically 
oriented truncation, following the work of Jacquet-Lapid-Rogawski
on periods of automrphic forms; 
while for the second one,  our limited achievement 
is a precise formula for the abelian parts of non-abelian $L$-functions 
associated with cusp forms. Clearly, this chapter is
motivated by the work of Arthur [Ar], follows
Jacquet-Lapid-Rogawski [JLR] and our RIMS notes on: Analytic Truncation
and Rankin-Selberg Methods versus Geometric Truncation and Non-abelian 
Zeta Functions.

\section{Abelian Part of Non-Abelian $L$-Functions: Fomal Calculation}

\subsection{From Geometric Truncation to Periods}

Let $G=GL_r$ be the general linear group defined over a number field $F$. 
Fix a minimal ($F$-)parabolic subgroup $P_0$ to be the Borel subgroup 
consisting of upper-triangle matrices. Let $p:[0,r]\to\mathbb R$ be a 
normalized convex polygon, i.e., a convex polygon such that $p(0)=p(r)=0$.

Recal that for any locally $L^1$-function 
$\phi:G(F)\backslash G(\mathbb A)\to \mathbb C$, the geometrically
oriented analytic truncation $\Lambda_p \phi$ is an 
$L^1$-function on $G(F)\backslash G(\mathbb A)$
defined by $$\Big(\Lambda_p\phi\Big)(g):=
\sum_P(-1)^{d(P)-d(G)}\sum_{\delta\in P(F)\backslash G(F)}
\phi_P(\delta g)\cdot\bold 1\Big(p_P^{\delta g}>_Pp\Big),$$ 
where $P$ runs over all (standard) parabolic subgroups of $G$, and $d(P)$ 
denotes the rank of $P$. Thus in particular, formally, we have
$$\begin{aligned}
&\int_{G(F)\backslash G(\mathbb A)^1}\Big(\Lambda_p\phi\Big)(g)\,dg\\
=&
\int_{G(F)\backslash G(\mathbb A)^1}\sum_P(-1)^{d(P)-d(G)}
\sum_{\delta\in P(F)\backslash G(F)}
\phi_P(\delta g)\cdot\bold 1\Big(p_P^{\delta g}>_Pp\Big)\,dg\\
=&\sum_P(-1)^{d(P)-d(G)}\int_{G(F)\backslash G(\mathbb A)^1}
\sum_{\delta\in P(F)\backslash G(F)}
\phi_P(\delta g)\cdot\bold 1\Big(p_P^{\delta g}>_Pp\Big)\,dg\\
=&\sum_P(-1)^{d(P)-d(G)}\int_{P(F)\backslash G(\mathbb A)^1}
\phi_P(g)\cdot\bold 1\Big(p_P^{g}>_Pp\Big)\,dg.\end{aligned}$$

\subsection{Regularization}

Surely, if all the integrations involved are absolutely convergent, 
this formal calculation works without any further work. 
However, for most of functions $\phi$, we cannot expect such a cheap way out. 
It is to remedy this, in literatures, in particular, in the work of Arthur, 
along with the line of the classical Ranking-Selberg methods, introduced is
an analytic truncation, for which
Langlands' combinatorial lemma is used so as to get an inversion formula 
for the function $\phi$ in terms of (a combination of) 
its partial analytic truncations. (See the discussion below.) 
Consequently, by the rapid decreasing property for the associated 
partial truncation, then the above discussion is justified. 
However, here we will not try to pursue this. Instead, we go 
further with our aim to understand the partial integration 
$$\int_{P(F)\backslash G(\mathbb A)^1}
\phi_P(g)\cdot\bold 1\Big(p_P^{g}>_Pp\Big)\,dg.$$ This then shows that we 
should develop a kind of regularization, say $\int_{P(F)\backslash G(\mathbb A)^1}^*$ for the above integration so as to get a well-defined
$$\int_{P(F)\backslash G(\mathbb A)^1}^*
\phi_P(g)\cdot\bold 1\Big(p_P^{g}>_Pp\Big)\,dg.$$
 
Clearly, one should not expect such a regularization works for all kinds of
$L^1$-function. But still, it is reasonable to expect that such a 
regularization works when $\phi$ is an automorphic form. As such, the first 
property we except is the following
\vskip 0.30cm
\noindent
{\bf {\large Basic Equality}.} (I) {\it For an automorphic form $\phi$ on 
$G(F)\backslash G(\mathbb A)$, we have
$$\int_{G(F)\backslash G(\mathbb A)^1}^*\Big(\Lambda_p\phi\Big)(g)\,dg
=\sum_P(-1)^{d(P)-d(G)}\int_{P(F)\backslash G(\mathbb A)^1}^*
\phi_P(g)\cdot\bold 1\Big(p_P^{g}>_Pp\Big)\,dg.$$}

Thus, from now on, assume in addition that $\phi$ is an {automorphic form}. 
Then, by [MW], as a function on 
$P(F)N(\mathbb A)\backslash G(\mathbb A)$, the 
constant term $\phi_P$ of $\phi$ along $P$ may be written in the  form
$$\phi_P(namk)=\sum_{j=1}^s\phi_j(m,k)\cdot\alpha_j\Big(H_P(a),k\Big)
e^{\langle\lambda_j+\rho_P,H_P(a)\rangle},$$ for
$n\in N(\mathbb A),\ a\in A_P,\ m\in M(\mathbb A)^1)$ and $k\in K$, where 
for all $j$,

\noindent
(a) $\phi_j(m,k)$ is an automorphic form on $M(F)\backslash M(\mathbb A)^1
\times K$; and

\noindent
(b) $\lambda_j\in\frak a_P^*$ and $\alpha_j(X,k)$ is a continuous family of 
polynomial on $\frak a_P$.

\noindent
Thus formally, we may further carry on the above calculation as follows:
$$\begin{aligned}&\int_{P(F)\backslash G(\mathbb A)^1}
\phi_P(g)\cdot\bold 1\Big(p_P^{g}>_Pp\Big)\,dg\\
=&\int_K\int_{M(F)\backslash M(\mathbb A)^1}\int_{\frak a_P}\Big(\sum_{j=1}^s
\phi_j(m,k)\cdot\alpha_j\Big(X,k\Big)e^{\langle\lambda_j+\rho_P,X\rangle}\Big)
e^{-\langle\rho_P,X\rangle}\cdot\bold 1\Big(p_P^{e^X}>_Pp\Big)\,dX\,dm\,dk\\
=&\sum_{j=1}^s
\int_K\int_{M(F)\backslash M(\mathbb A)^1}\phi_j(m,k)\,dm\,dk
\cdot\int_{\frak a_P}\alpha_j\Big(X,k\Big)
e^{\langle\lambda_j,X\rangle}\cdot\bold 1\Big(p_P^{e^X}>_Pp\Big)\,dX.
\end{aligned}$$

With this, then we need a new type of regularization
$\int_{\frak a_P}^{\#}$ such that we get the following
\vskip 0.30cm
\noindent
{\bf {\large Basic Equality}.} (II) {\it With the same notation as above, 
we have
$$\begin{aligned}
&\int_{G(F)\backslash G(\mathbb A)^1}^*\Big(\Lambda_p\phi\Big)(g)\,dg\\
=&\sum_P(-1)^{d(P)-d(G)}\sum_{j=1}^s
\int_K\Big(\int_{M(F)\backslash M(\mathbb A)^1}^*\phi_j(m,k)\,dm\Big)\,dk
\cdot\int_{\frak a_P}^{\#}
\alpha_j\Big(X,k\Big)e^{\langle\lambda_j,X\rangle}\cdot
\bold 1\Big(p_P^{e^X}>_Pp\Big)\,dX.\end{aligned}$$}
\vskip 0.30cm
\noindent
{\bf Remarks.} (1) Clearly, it is enough to build the whole regularization
from that $\int_{\frak a_P}^{\#}$ using the parabolic induction: after all,
for minimal parabolic subgroup $P_0$, $M_0$ is discrete;

\noindent
(2) Roughly speaking, modulo  subtle difference such as here we use
geometric truncation $\bold 1$ instead of the analytic truncation $\widehat\tau$, 
the above process defining regularizations $\int^*$
on $P(F)\backslash G(\mathbb A)^1$ and $\int^{\#}$ on $\frak a_P$ 
is the technique heart of Jacquet-Lapid-Rogawski's  paper [JLR].

\subsection{Berstein's Principle}

Instead of giving a detailed discussion on such regularizations for our 
geometric truncation, in the follows, we will first give an application of the 
above discussion and then show how JLR's discusion works for 
$\bold 1(p_P^*>Pp)$ as well -- With the help of our detailed discussion
on geometric truncation $\bold 1$ in Chapter 2, 
we conclude that $\bold 1(p_P^*>Pp)$ are also chapacteristic functions
of  cones. Consequently, we can then be confident that the final formula 
for the abelian part of non-abelian $L$-functions (associated to cusp forms) 
given below should work correctly.

However, as it turns out, for this purpose, there is still one key point which
is missing from the above discussion: Bernstein's Principle, as discovered by
Jacquet-Lapid-Rogawski ([JLR]).

Let $P=MN$ be a parabolic subgroup and let $\sigma$ be an automorphic 
subrepresentation of $L^2\Big(M(F)\backslash M(\mathbb A)^1\Big)$. Let 
$\mathcal A_P(G)_\sigma$ be the subspace of functions 
$\phi\in \mathcal A_P(G)$, the space of level $P$-automorphic forms, 
such that $\phi$ is left-invariant under 
$A_P$ and for all $k\in K$, the function $m\mapsto \phi(mk)$ belongs 
to the space of $\sigma$. For $\phi\in \mathcal A_P(G)_\sigma$ and
 $\lambda\in\frak a_P^*$, we write $E(g,\phi,\lambda)$ for the
{\it Eisenstein series} which is given in its domian of absolute 
convergence by the infinite series
$$E(g,\phi,\lambda):=\sum_{\gamma\in P(F)\backslash G(F)}
\phi(\gamma g)e^{\langle \lambda+\rho_P,H_P(\gamma g)\rangle}.$$

Then by a work of JLR, based on a result of Bernstein,
we  have the following
\vskip 0.30cm
\noindent
{\bf {\large Bernstein's Principle}}. 
{\it Let $P=MN$ be a proper parabolic sungroup and let $\sigma$ be 
an irredcible cuspidal representation in 
$L^2\Big(M(F)\backslash M(\mathbb A)^1\Big)$. Let $E(g,\phi,\lambda)$ 
be an Eisenstein series associated to $\phi\in \mathcal A_P(G)_\sigma$. Then
$$\int_{P(F)\backslash G(\mathbb A)^1}^*E(g,\phi,\lambda)\,dg =0$$ 
for generic $\lambda$.}

As will be seen later, the main reason for such a result to be hold is that 

\noindent
(a) with a natural regularization, we may further assume that
{\it the map $$\phi\mapsto \int_{P(F)\backslash G(\mathbb A)^1}^*
E(g,\phi,\lambda)\,dg$$ defines a 
$G(\mathbb A_f)^1$-invariant functional on 
$\mathrm{Ind}_P^G(\sigma\otimes e^\lambda)$ for nice $\lambda$} motivated by a result of JLR [JLR] whose analogue will be discussed below; and

\noindent
(b) there does not exist any such invariant functional for generic values of 
$\lambda$, by a result of Bernstein. 

\subsection{Constant Terms of Eisenstein Series}

As above, let $P=MN$ be a parabolic subgroup and let 
$\sigma$ be an automorphic 
subrepresentation of $L^2\Big(M(F)\backslash M(\mathbb A)^1\Big)$. Let 
$\phi\in \mathcal A_P(G)_\sigma$, 
 $\lambda\in\frak a_P^*$, and $E(g,\phi,\lambda)$ be the associated
Eisenstein series. 

Furthermore, let $N_{G(F)}(A_0)$ be the normalizer of $A_0$ in $G(F)$ and let 
$$\Omega:=N_{G(F)}(A_0)\Big/C_{G(F)}(A_0)$$ be the {\it Weyl group} of $G$. 
Recall that a parabolic subgroup $Q$ is {\it associated} to $P$, written as
$Q\sim P$ if $M_Q$ 
is conjugate to $M_P$ under $\Omega$. If $Q$ is associated to $P$, let 
$\Omega(P,Q)$ be the set of maps $A_P\to A_Q$ obtained by restriction 
of elements $w\in \Omega$ such that $wM_Pw^{-1}=M_Q$.

Suppose that $Q$ is associated to $P$ and let $w$ be an element of 
$\Omega(P,Q)$ with representative $\tilde w$ in $N_{G(F)}(A_0)$. 
Recall that the {\it standard intertwining operator} is defined by 
$$\Big(M(w,\lambda)\phi\Big)(g):=
e^{-\langle w\lambda+\rho_Q,H_Q(g)\rangle}\int_{N_w(\mathbb A)
\backslash N_Q(\mathbb A)}\phi(\tilde w^{-1}ng)e^{\langle \lambda+
\rho_P,H_P(\tilde w^{-1}ng)\rangle}dn$$ where $N_w:=N_Q\cap \tilde
 w N\tilde w^{-1}$. The operator $M(w,\lambda)$ depends on $w$ but 
not on the choice of representative $\tilde w$.
\vskip 0.30cm
Assume that $\sigma$ is {\bf cuspidal}. For $Q'\subset Q$, denote by 
$E^Q(g,\psi,\lambda)$  the Eisenstein 
series induced from $M_{Q'}$ to $M_Q$:
$$E^Q(g,\psi,\lambda):=\sum_{\gamma\in Q'(F)\backslash Q(F)}
\psi(\gamma g)e^{\langle \lambda+\rho_{Q'},H_{Q'}(\gamma g)\rangle}.$$

Assume that in particular, $\phi$ is {\it cuspidal}.
 Then the constant term $E_Q(g,\phi,\lambda)$ of $E(g,\phi,\lambda)$ 
relative to a parabolic subgroup $Q$ has a 
simple expression. More precisely,
\vskip 0.30cm
\noindent
(a) If $Q$ does not contain any associated of $P$, 
then $E_Q(g,\phi,\lambda)$ is identically zero. 

\noindent
(b) If $Q$ properly contains an association of $P$, then ([A2], or/and [MW]) 
$$E_Q(g,\phi,\lambda)=\sum_{Q'}\sum_{w\in \Omega(P,Q'),w^{-1}\alpha>0
\forall\alpha\in\Delta_{Q'}^Q} E^Q\Big(g, M(w,\lambda)\phi,w\lambda\Big)$$ 
where the sum is over the standard parabolic subgroups $Q'\subset Q$ 
associated to $P$. Thus in particular, 

\noindent
(b)$'$ If $Q$ is associated to $P$, then
$$E_Q(g,\phi,\lambda)=\sum_{w\in\Omega(P,Q)}\Big(M(w,\lambda)\phi\Big)(g)
\cdot e^{\langle w\lambda+\rho_Q,H_Q(g)\rangle}.$$

\subsection{Abelian Part of Non-Abelian $L$-Functions}

Thus, when $\phi\in \mathcal A_P(G)_\sigma$ is {\it cuspidal},
the abelian part of the non-abelian zeta function associated to $\phi$ may be 
given using the following calculation:
$$\begin{aligned}&L^{\mathrm{ab}}_{\leq p}(\phi;\lambda)\\
:=&\int_{G(F)\backslash G(\mathbb A)^1}^*
\Big(\Lambda_pE(g,\phi,\lambda)\Big)\,dg\\
=&\int_{G(F)\backslash G(\mathbb A)^1}^*\sum_Q(-1)^{d(Q)-d(G)}
\sum_{\delta\in Q(F)\backslash G(F)}E_Q(\delta g,\phi,\lambda)\cdot
\bold 1\Big(p_Q^{\delta g}>_Qp\Big)\,dg\\
&\hskip 3.0cm\mathrm{(by\ the\ definition\ of\ geometric\ truncation)}\\
=&\sum_Q(-1)^{d(Q)-d(G)}\int_{G(F)\backslash G(\mathbb A)^1}^*
\sum_{\delta\in Q(F)\backslash G(F)}E_Q(\delta g,\phi,\lambda)\cdot
\bold 1\Big(p_Q^{\delta g}>_Qp\Big)\,dg\\
&\hskip 3.0cm\mathrm{(by\ the\ desired\ property\ of\ regularization)}\\
=&\sum_Q(-1)^{d(Q)-d(G)}\int_{Q(F)\backslash G(\mathbb A)^1}^*
E_Q(g,\phi,\lambda)\cdot\bold 1\Big(p_Q^g>_Qp\Big)\,dg\\
&\hskip 3.0cm\mathrm{(by\ the\ desired\ property\ of\ regularization)}\\
=&\sum_{Q:\exists Q'\ s.t. Q\supset Q'\sim P}(-1)^{d(Q)-d(G)}
\int_{Q(F)\backslash G(\mathbb A)^1}^*\sum_{Q'}\sum_{w\in \Omega(P,Q'): 
w^{-1}\alpha>0,\forall \alpha\in\Delta_{Q'}^Q}\\
&\hskip 3.0cm
E^Q\Big(g,M(w,\lambda)\phi,w\lambda\Big)\cdot\bold 1\Big(p_Q^g>_Qp\Big)\,dg\\
&\hskip 3.0cm\mathrm{(by\ the\ formula\ for\ constant\ 
terms\ of\ Eisenstein\ series)}\\
=&\sum_{Q: Q\sim P}(-1)^{d(Q)-d(G)}
\int_{Q(F)\backslash G(\mathbb A)^1}^*
E_Q\Big(g,\phi,\lambda\Big)\cdot\bold 1\Big(p_Q^g>_Qp\Big)\,dg\\
&\hskip 3.0cm\mathrm{(by\ Bernstein's\ Principal)}\\
=&\sum_{Q: Q\sim P}(-1)^{d(Q)-d(G)}
\int_{Q(F)\backslash G(\mathbb A)^1}^*\sum_{w\in\Omega(P,Q)}
M(w,\lambda)\phi(g)\cdot e^{\langle w\lambda+\rho_Q,H_Q(g)\rangle}
\cdot\bold 1\Big(p_Q^g>_Qp\Big)\,dg\\
&\hskip 3.0cm\mathrm{(by\ the\ formula\ for\ constant\ terms\ 
of\ Eisenstein\ series\ again)}\\
=&\sum_{Q: Q\sim P}(-1)^{d(Q)-d(G)}\sum_{w\in\Omega(P,Q)}
\int_K
\Big(\int_{M_Q(F)\backslash M_Q(\mathbb A)^1}^*
M(w,\lambda)\phi(mk)dm\Big)\,dk\\
&\hskip 5.0cm
\cdot\int_{\frak a_P}^{\#} e^{\langle w\lambda,X\rangle}
\cdot\bold 1\Big(p_Q^{e^X}>_Qp\Big)\,dX\\
&\hskip 3.0cm\mathrm{(by\ the\ desired\ property\ of\ regularization)}\\
=&\sum_{Q: Q\sim P}(-1)^{d(Q)-d(G)}\sum_{w\in\Omega(P,Q)}
\int_K
\Big(\int_{M_Q(F)\backslash M_Q(\mathbb A)^1}
M(w,\lambda)\phi(mk)dm\Big)\,dk\\
&\hskip 3.0cm
\cdot\int_{\frak a_P}^{\#} e^{\langle w\lambda,X\rangle}
\cdot\bold 1\Big(p_Q^{X}>_Qp\Big)\,dX\\
&\hskip 3.0cm\mathrm{(by\ the\ fact\ that}\ \phi\ \mathrm{is\ cuspidal\ 
and\ hence\ integrable)}.\end{aligned}$$
For simplicity, here we write
$\bold 1\Big(p_Q^{X}>_Qp\Big)$ for  
$\bold 1\Big(p_Q^{e^X}>_Qp\Big)$.
In such a way, we arrive at the following
\vskip 0.30cm
\noindent
{\bf {\large Abelian Part of Non-Abelian $L$-Functions.}} {\it For
$\phi\in \mathcal A_P(G)_\sigma$ a  cusp form,
the abelian part of the non-abelian zeta function associated to $\phi$ 
is given by
$$\begin{aligned}&L^{\mathrm{ab}}_{\leq p}(\phi;\lambda)\\
=&\sum_{Q: Q\sim P}(-1)^{d(Q)-d(G)}
\sum_{w\in\Omega(P,Q)}\int_K
\Big(\int_{M_Q(F)\backslash M_Q(\mathbb A)^1}
M(w,\lambda)\phi(mk)dm\Big)\,dk\\
&\hskip 5.0cm
\cdot\int_{\frak a_P}^{\#} e^{\langle w\lambda,X\rangle}
\cdot\bold 1\Big(p_Q^{X}>_Qp\Big)\,dX.\end{aligned}$$}
 
In the sequel, we will not try to pursue further the rigorous
arguments in establishing the claim above. Instead, we will 
give an example in showing how such regularizations can be worked out
following [JLR].

\section{Regularized Integration: Integration over Cones}

We here give a review  of Jacquet-Lapid-Rogawski's 
beautiful yet elementary treatment of integrations over cones.

Let $V$ be a real finite-dimensional vector space of dimension $r$. Let $V^*$ 
be the space of complex linear forms on $V$. Denote by $S(V^*)$ the symmetric 
algebra of $V^*$, which may be regarded as the space of polynomial 
functions on $V$ as well. 
By definition, an {\it exponential polynomial} function on $V$ is a function 
of the form $$f(x)=\sum_{i=1}^r e^{\langle\lambda_i,x\rangle}P_i(x)$$ where 
the $\lambda_i$ are distinct elements of $V^*$ and the $P_i(x)$ are non-zero 
elements of $S(V^*)$. One checks easily that 
\vskip 0.30cm
\noindent
(1) {\it The $\lambda_i$ are uniquely determined and called 
the exponents of $f$.}
\vskip 0.30cm
\noindent
{\bf Remark.} As being seen in the previous section,
the reason for introducing these functions is that
constant terms of an automorphic form along with parabolic subgroups are
 of such types.
\vskip 0.30cm
By a {\it cone} in $V$ we shall mean a closed subset of the form
$$\mathcal C:=\Big\{x\in V:\langle \mu_i,x\rangle\geq 0\ \forall i\Big\}$$ 
where $\{\mu_i\}$ is a
 basis of $\mathrm{Re}(V^*)$, the space of {\it real} linear forms on $V$. Let 
$\big\{e_j\big\}$ be the dual basis of $V$. We shall say that 

\noindent
$\lambda\in V^*$ is {\it negative with respect to} $\mathcal C$ if 
$\mathrm {Re}\,\langle\lambda,e_j\rangle<0$ for each $j=1,\cdots,r$; and that

\noindent
$\lambda$ is {\it non-degenerate with respect to} $\mathcal C$ if 
$\langle\lambda,e_j\rangle\not=0$ for each $j=1,\cdots,r$.
\vskip 0.30cm
\noindent
Easily we have

\noindent
(2) {\it The function $f(s)=\sum_{i=1}^r e^{\langle\lambda_i,x\rangle}P_i(x)$ 
is integrable over $\mathcal C$ if and only if $\lambda_i$ is negative with
 respect to $\mathcal C$ for all $i$.}
\vskip 0.30cm
To define the regularized integral over $\mathcal C$, we study the integral 
$$I_{\mathcal C}\Big(f;\lambda\Big):=\int_{\mathcal C} f(x)
e^{\langle\lambda,x\rangle}dx.$$
The integral converges absolutely for $\lambda$ in the open set 
$$\Big\{\lambda\in V^*:\mathrm{Re}\langle\lambda_i-\lambda,e_j\rangle<0\ 
\forall 1\leq i\leq n, 1\leq j\leq r\Big\}.$$ Moreover, it may be 
analytically continued as follows:
 
First,  we have for $\mathrm{Re}\,(\lambda)>0$ and 
for any polynomial $P$ in one variable:
$$\int_0^\infty e^{-\lambda x}P(x)dx=\sum_{m\geq 0}\frac{\Big(D^mP\Big)(0)}
{\lambda^{m+1}}.$$ 

More generally, fix an index $k$ and let $\mathcal C_k$ be the intersection of 
$\mathcal C$ and the hyperplane $V_k:=\Big\{x:\langle \mu_k,x\rangle=0\Big\}$. 
Then in the integral we can write $x=\langle\mu_k,x\rangle e_k+y$ with 
$y\in V_k$. Thus, for a suitable choice of Haar measures:
$$I_{\mathcal C}\Big(f;\lambda\Big)=\sum_{i=1}^n\sum_{m\geq 0}
\frac{1}{\Big(\langle\lambda-
\lambda_i,e_k\rangle\Big)^{m+1}}\int_{\mathcal C_k}e^{-\langle
\lambda-\lambda_i,y\rangle}\Big(D_{e_k}^m\Big)P_i(y)\,dy.$$
This formula gives the analytic continuation of 
$I_{\mathcal C}\Big(f;\lambda\Big)$ 
to the tube domain defined by $\mathrm {Re}\,\langle\lambda_i-\lambda,
e_k\rangle>0$ for $j\not=k$, $1\leq i\leq n$ with hyperplane singularities. 
Clearly, the singular hyperplanes are given by
$$H_{k,i}:=\Big\{\lambda:\langle\lambda,e_k\rangle=\langle\lambda_i,e_k\rangle
\Big\},
\qquad 1\leq i\leq n.$$ Since this is true for all $k=1,2,\dots,r$, the 
function $I_{\mathcal C}\Big(f;\lambda\Big)$ has analytic continuation to 
$V^*$ by Hartogs' lemma with (only) hyperplane singularities along $H_{i,k}$, 
$1\leq k,\leq r, 1\leq i\leq r$. Moreover, from such an induction, one checkes
\vskip 0.30cm
\noindent
(3) {\it Suppose that $f$ is absolute integrable over $\mathcal C$. Then 
$I_{\mathcal C}\Big(f;\lambda\Big)$ is holomorphic at 0 and 
$I_{\mathcal C}\Big(f;0\Big)=\int_{\mathcal C}f(x)\,dx$}; and 
\vskip 0.30cm
\noindent
(4) {\it The function $I_{\mathcal C}\Big(f;\lambda\Big)$ is holomorphic at 0  
if and only if for all $i$, $\lambda_i$ is non-degenerate with respect to 
$\mathcal C$, i.e.,$\langle\lambda_i,e_k\rangle\not=0$ for all pairs $(i,k)$.}
\vskip 0.30cm
Now we are ready to define the regularized integral. Denote the characteristic 
function of a set $Y$ by $\tau^Y$. For $\lambda\in V^*$ such that 
$\lambda_i-\lambda$ is negative with respect to $\mathcal C$ for all 
$i=1,2,\cdots,n$, set
$$\begin{aligned}\widehat F\Big(\lambda;\mathcal C,T\Big):=&
\int_Vf(x)\cdot \tau^{\mathcal C}(x-T)\cdot 
e^{-\langle \lambda,x\rangle}dx\\
=&\sum_{i=1}^n\int_VP_i(x)\tau^{\mathcal C}(x-T)e^{-\langle \lambda-\lambda_i,
x\rangle}dx\\
=&\sum_{i=1}^ne^{-\langle \lambda_i-\lambda,T\rangle}I_{\mathcal C}\Big(
P_i(*+T)e^{\langle \lambda_i,*\rangle};\lambda\Big).\end{aligned}$$
The integrals are absolutely convergent and 
$\widehat F\Big(\lambda;\mathcal C,T\Big)$ 
extends to a meromorphic function on $V^*$. 
\vskip 0.30cm
\noindent
{\bf {\large Definition}}. (JLR) The function 
$f(x)\cdot \tau^{\mathcal C}(x-T)$ 
is called \#-{\it integrable} if $\widehat F\Big(\lambda;\mathcal C,T\Big)$ 
is holomorphic at $\lambda=0$. 
In this case, set $$\int_V^{\#}f(x)\cdot 
\tau^{\mathcal C}(x-T)dx:=\widehat F\Big(0;\mathcal C,T\Big).$$ 

By the basic property (4) above, the \#-integral exists if and only if each 
exponent $\lambda_j$ is non-degenerate with respect to $\mathcal C$.
\vskip 0.30cm
To go further, suppose that $V=W_1\oplus W_2$ is a decomposition of 
$V$ as a direct sum, and 
let $\mathcal C_j$ be a cone in $W_j$. Write $T=T_1+T_2$ and $x=x_1+x_2$ 
relative to this decomposition. If the exponents $\lambda_j$ of $f$ are 
non-degenerate with respect to $\mathcal C=\mathcal C_1+\mathcal C_2$, then 
the function
$$w_2\mapsto \int_{W_1}^{\#}f(w_1+w_2)\tau_1^{\mathcal C_1}(w_1-T_1)\,dw_1$$ 
is defined and is an exponential polynomial. And it follows by analytic 
continuation that 

\noindent
(5) $\displaystyle{\int_{V}^{\#}f(x)\tau^{\mathcal C}(x-T)\,dx}$
$$=\int_{W_2}^{\#}\bigg(\int_{W_1}^{\#}f(w_1+w_2)
\tau_1^{\mathcal C_1}(w_1-T_1)\,dw_1\bigg)\cdot \tau_2^{\mathcal C_2}
(w_2-T_2)\,dw_2;$$ 
{\it and as a function of $T$, it is an exponential polynomial 
with the same exponents as $f$.}
\vskip 0.30cm
More generally, let $g(x)$ be a compactly supported function on $W_1$ and 
consider  functions of the form 
$g(w_1-t_1)\tau_{\mathcal C}(w_2-T_2)$, which we call functions of 
{\it type (C)}. It is clear that the integral 
$$\widehat F(\lambda):=\int_Vf(x)\cdot
g(w_1-t_1)\tau_{\mathcal C}(w_2-T_2)\cdot e^{-\langle\lambda,x\rangle}dx$$ 
converges 
absolutely for an open set of $\lambda$ whose restriction to $W_2$ is 
negative with respect to $\mathcal C_2$. Furthermore, $\widehat F(\lambda)$ 
has a meromorphic continuation to $V^*$. Following JLR, the function $f(x)
g(w_1-t_1)\tau_{\mathcal C}(w_2-T_2)$ is called \#-{\it integrable} if 
$\widehat F(\lambda)$ is holomorphic at $\lambda=0$. If so we define
$$\int_V^{\#}f(x)\cdot
g(w_1-t_1)\tau_{\mathcal C}(w_2-T_2)\,dx:=\widehat F(0).$$
Note that as a function of $w_2$, $\int_{W_1}f(w_1+w_2)g(w_1-T_1)\,dw_1$ 
is an exponential polynomial on $W_2$, hence
$$\begin{aligned}&\int_V^{\#}f(x)\,
g(w_1-t_1)\tau_{\mathcal C}(w_2-T_2)\,dx\\
=&\int_{W_2}^{\#}\bigg(
\int_{W_1}f(w_1+w_2)g(w_1-T_1)\,dw_1\bigg)\tau_2^{\mathcal C}
(w_2-T_2)\,dw_2.\end{aligned}$$

Moreover, if $V=W_{i1}\oplus W_{i2}$ are decompositions of $V$ with
 $\mathcal C_{i2}$ a cone in $W_{i2}$, and $g_i$ are compactly 
supported functions on $W_{i1}$. Set $G_i(w_1+w_2)=g_1(w_1)
\tau_{i2}^{\mathcal C}(w_2)$ for $w_j\in W_{ij}$.
Then for cones $\mathcal C$ and $\mathcal C'$ in $V$ such that 
$\mathcal C_{i2},\mathcal C\subset\mathcal C'$ for all $i$, and that
 $\tau^{\mathcal C}(x)=\sum_{i=1}^n a_iG_i(x)$ for some constants $a_i$,
we have
\vskip 0.30cm
\noindent
(6) {\it If each of the integrals $\int_V^{\#}f(x)\cdot G_i(x-T)\,dx$ 
exists for all $i$, then $f(x)\tau^{\mathcal C}(x-T)$ is \#-integrable and
$$\int_V^{\#}f(x)\tau^{\mathcal C}(x-T)\,dx=\sum_i a_i\int_V^{\#}f(x)G_i(x-T)\,
dx.$$}

We end this discussion by the following explicit formulas:
\vskip 0.30cm
\noindent
(7) {\it $\displaystyle{
\int_V^{\#}e^{\langle\lambda,x\rangle}\tau^{\mathcal C}(x-T)\,dx\,=\,(-1)^r
\mathrm{Vol}\Big(e_1,e_2,\cdots,e_r\Big)\cdot 
\frac{e^{\langle\lambda,T\rangle}}{\prod_{j=1}^r
\langle\lambda,e_j\rangle},}$
where the cone $\mathcal C$ is given by 
$\Big\{\sum_{j=1}^r a_je_j:a_j\geq 0\Big\}$ and 
$\mathrm{Vol}\Big(e_1,e_2,\cdots,e_r\Big)$ is the volume of the associated 
parallelepiped 
$\Big\{\sum_{j=1}^r a_je_j:a_j\in[0,1]\Big\}$.}
 
In particular, if $V=\mathbb R$, we have 
$\displaystyle{\int_T^\infty e^{\lambda t}dt=
-\frac{e^{\lambda T}}{\lambda}}$ and hence for any cone $\mathcal C\subset 
\mathbb R$, 
$$\int_V^{\#}e^{\langle\lambda,x\rangle}
\Big(1-\tau^{\mathcal C}\big(x-T\big)\Big)\,dx
=-\int_V^{\#}e^{\langle\lambda,x\rangle}\tau^{\mathcal C}(x-T)\,dx$$ since 
$1-\tau^{\mathcal C}$ is the characteristic function of the cone $-\mathcal C$.

\section{Positive Cones and Positive Chambers}

\subsection{Modified Truncation}

For $T\in\frak a_0$, define  {\it analytic truncations} of $\phi$ by the 
formula
$$\Lambda^{T,P}\phi(g):=\sum_{R:R\subset P}(-1)^{d(R)-d(P)}
\sum_{\delta\in R(F)\backslash P(F)}\phi_R(\delta g)\cdot
\widehat\tau_R^P\Big(H(\delta g)-T\Big),$$ and
$$\Lambda_{T,P}\phi(g):=\sum_{R:R\subset P}(-1)^{d(R)-d(P)}
\sum_{\delta\in R(F)\backslash P(F)}\phi_R(\delta g)\cdot
\tau_R^P\Big(H(\delta g)-T\Big).$$ 
As they stand, the first one is simply Arthur's truncation, while
the second one is an obvious imitation in which positive chambers, instead of
positive cones, are used.
Write $\Lambda^T$ and $\Lambda_T$ for $\Lambda^{T,G}$ and $\Lambda_{T,G}$
respecitvely. Then
\vskip 0.30cm
\noindent
(i) $\displaystyle{\Lambda^{T}\phi(g):=\sum_{R}(-1)^{d(R)-d(G)}
\sum_{\delta\in R(F)\backslash G(F)}\phi_R(\delta g)\cdot
\widehat\tau_R\Big(H(\delta g)-T\Big)},$ while
$$\Lambda_{T}\phi(g):=\sum_{R}(-1)^{d(R)-d(G)}
\sum_{\delta\in R(F)\backslash G(F)}\phi_R(\delta g)\cdot
\tau_R\Big(H(\delta g)-T\Big);$$

\noindent
(ii) ({\bf Inversion Formula}) {\it For a $G(F)$-invariant function $\phi$,} 
$$\begin{aligned}
\phi(g)=&\sum_P\sum_{\delta\in P(F)\backslash G(F)}\Lambda^{T,P}
\phi(\delta g)\cdot \tau_P\Big(H(\delta g)-T\Big)\\
=&\sum_P\sum_{\delta\in P(F)\backslash G(F)}\Lambda_{T,P}
\phi(\delta g)\cdot\widehat \tau_P\Big(H(\delta g)-T\Big).\end{aligned}$$ 
\vskip 0.30cm
\noindent 
{\it Proof.} Let us recall the following
\vskip 0.30cm
\noindent
{\bf Langlands' Combinatorial Lemma.} {\it If $Q\subset P$ are parabolic 
subgroups, 
then for all $H\in \frak a_0$, we have
$$\sum_{R:Q\subset R\subset P}(-1)^{d(R)-d(P)}\tau_Q^R(H)\widehat\tau_R^P(H)=
\delta_{QP},$$ and 
$$\sum_{R:Q\subset R\subset P}(-1)^{d(Q)-d(R)}\widehat\tau_Q^R(H)\tau_R^P(H)=
\delta_{QP}.$$}
 
We start with the second one for positive cone.
Thus in particular, the RHS of (ii) is simply
$$\begin{aligned}&\sum_P\sum_{\delta\in P(F)\backslash G(F)}\bigg(
\sum_{R:R\subset P}(-1)^{d(R)-d(P)}\sum_{\gamma\in R(F)\backslash P(F)}
\phi_R(\gamma \delta g)\cdot\tau_R^P\Big(H(\gamma\delta g)-T\Big)\bigg)
\cdot\widehat \tau_P\Big(H(\delta g)-T\Big)\\
=&\sum_{R,P:R\subset P}\sum_{\delta\in P(F)\backslash G(F)}
\sum_{\gamma\in R(F)\backslash P(F)}\phi_R(\gamma\delta g)(-1)^{d(R)-d(P)}
\tau_R^P\Big(H(\gamma\delta g)-T\Big)\cdot\widehat \tau_P\Big(H(\delta g)-T\Big)\\
=&\sum_{R,P:R\subset P}\sum_{\delta\in R(F)\backslash G(F)}
\phi_R(\delta g)(-1)^{-d(R)+d(P)}\tau_R^P\Big(H(\delta g)-T\Big)
\cdot\widehat \tau_P\Big(H(\delta g)-T\Big)\\
=&\sum_{R}\sum_{\delta\in R(F)\backslash G(F)}
\phi_R(\delta g)\sum_{P:R\subset P}
(-1)^{d(G)-d(R)+d(P)-d(G)}\tau_R^P\Big(H(\delta g)-T\Big)
\cdot\widehat \tau_P\Big(H(\delta g)-T\Big)\\
=&\sum_{R}\sum_{\delta\in R(F)\backslash G(F)}(-1)^{d(R)-d(G)}
\phi_R(\delta g)\sum_{P:R\subset P}
(-1)^{d(P)-d(G)}\tau_R^P\Big(H(\delta g)-T\Big)
\cdot\widehat \tau_P\Big(H(\delta g)-T\Big)\\
=&\sum_{R}\sum_{\delta\in R(F)\backslash G(F)}(-1)^{d(R)-d(G)}
\phi_R(\delta g)\cdot \delta_{RG}\end{aligned}$$ 
by the first equality in Langlands' combinatorial lemma. 
But this latest quantity is clearly $\phi(g)$. 

The proof of the first equality is similar: here the second relation in 
Langlands' Combinatorial Lemma is needed. We leave this to the reader.
This completes the proof.
\vskip 0.30cm
Despite the fact the at formally level the truncations $\Lambda^{T,P}$
and $\Lambda_{T,P}$ are quite parallel, they bear essential differences:
In the case for Arthur's truncation $\Lambda^{T,P}$, we have the following 

\noindent
{\bf Basic Fact.} (Arthur) {\it Assume that $T\in \mathcal C_P$, and that
$\phi$ is an automorphic form on $N(\mathbb A)M(F)\backslash G(\mathbb A)$. 
Then $\Lambda_{T,P}\phi(m)$ is rapidly decreasing 
on $M(F)\backslash M(\mathbb A)^1$. In particular, 
the integration  $\int_{M(F)\backslash M(\mathbb A)^1}\Lambda_{T,P}\phi(m)\,dm$
is well-defined.}

\noindent
By contrast, for the new truncation $\Lambda_{T,P}$,  as can be checked 
with $GL_3$, a similar statement does not hold. That is, we have the following 
\vskip 0.30cm
\noindent
{\bf {\large False Assumption.}} {\it For $T\in -\mathcal C_P$,   
the function defined by (a possibly regilarized integration)
$g\mapsto \int_{M(F)\backslash M(\mathbb A)^1}\Lambda_{T,P}\phi(mg)\,dm$ 
behaves
as if the function $m\mapsto \Lambda_{T,P}\phi(mg)$ is rapidly decreasing 
on $M(F)\backslash M(\mathbb A)^1$.}
\subsection{JLR Period of an Automorphic Form}

Now let us follow [JLR] to apply the above discussion on integration over 
cones to the theory of periods of automorphic forms. 

Let $\tau_k(X)$ be a function of type (C) on 
$\frak a_P$ that depends continuously on $k\in K$, i.e., we assume that there 
is a decomposition $\frak a_P=W_1\oplus W_2$ such that $\widehat\tau_k$ 
has the form
$$g_k(w_1-T_1)\tau_{\mathcal C_{2k}}(w_2-T_2)$$ where the compactly 
supported function $g_k$ varies continuously in the $L^1$-norm and linear 
inequalities defining the cone $\mathcal C_{2k}$ vary continuously.

Also with the application to  automorphic forms in mind,
let $f$ be a function on $M(F)N(\mathbb A)\backslash G(\mathbb A)$ 
of the form
$$f(namk)=\sum_{j=1}^s\phi_j(m,k)\cdot \alpha_j\Big(H_P(a),k\Big)
e^{\langle\lambda_j+\rho_P,H_P(a)\rangle}$$ 
for $n\in N(\mathbb A), a\in A_P, m\in M(\mathbb A)^1$ and $k\in K$, 
where, for all $j$,

\noindent
(a) $\phi_j(m,k)$ is absolutely integrable on $M(F)\backslash 
M(\mathbb A)^1\times K$, so that $\int_{M(F)\backslash 
M(\mathbb A)^1\times K}\phi_j(m,k)dm\,dk$ is well-defined; and

\noindent
(b)  $\lambda_j\in\frak a_P^*$ and $\alpha_j(X,k)$ is a continuous family of 
polynomials on $\frak a_P$ such that, for all $k\in K$, 
$\alpha_J(X,k)e^{\langle\lambda_j,X\rangle}\tau_k(X)$ is \#-integrable;

As such, following [JLR], define the \#-{\it integral}
$$\int_{P(F)\backslash G(\mathbb A)}^{\#}f(g)\cdot 
\tau_k\Big(H_P(g)\Big)\,dg:=\sum_{j=1}^s\int_K\bigg(\int_{M(F)
\backslash M(\mathbb A)^1}\phi_j(mk)\,
dm\bigg)\cdot \bigg(\int_{\frak a_P}^{\#}\alpha_j(X,k)
\tau_k(X)\,dX\bigg)\,dk.$$ 

\noindent
{\bf Remark.} The reader should notice that in this definition there is 
no integration involving $\int_{N(F)\backslash N(\mathbb A)}dn:$ In practice,
$f$ is supposed to be the constant term along $P$, so is 
$N(\mathbb A)$-invariant. Thus, this integration may be omitted since we 
normalize our Haar measure so that the total volume of 
$N(F)\backslash N(\mathbb A)$ is 1.
\vskip 0.30cm
As usual, write $\mathcal E_P(f)$ for the set of 
distinct exponents $\{\lambda_j\}$ of $f$. As can be easily seen,
this set is uniquely determined 
by $f$, but the functions $\phi_j$ and $\alpha_j$ are not. So we need to show 
that the \#-integral is independent of the choices of these functions. 
Clearly, the function
$$X\mapsto \int_{M(F)\backslash M(\mathbb A)^1} f(e^Xmk)\,dm$$ is an 
exponential polynomial on $\frak a_P$ and  the above \#-integral is equal to
$$\int_K\int_{\frak a_P}^{\#}\bigg(\int_{M(F)\backslash M(\mathbb A)^1}
f(e^Xmk)\,dm\bigg) e^{-\langle\rho_P,X\rangle}\tau_k(X)\,dX\,dk.$$ 
This 
shows in particular that the original \#-integral is independent of the 
decomposition of $f$ as used in the definition. Moreover, 
if each of the exponents
 $\lambda_j$ is negative with respect to $\mathcal C_{2k}$ for all
 $k\in K$, then  the ordinary integral $$\int_{P(F)\backslash 
G(\mathbb A)^1}f(g)\cdot\tau_k(H_P(g))\,dg$$ is absolutely convergent
 and equals to the \#-integral by Basic Property (3) of integration over cones.
That is to say, we have the following
\vskip 0.30cm
\noindent
{\bf {\large Lemma 1}.} (1) {\it The above \#-integral 
$\int_{P(F)\backslash G(\mathbb A)}^{\#}f(g)\cdot \tau_k\Big(H_P(g)\Big)\,dg$
is well-defined;} and

\noindent
(2) {\it  If each of the exponents
 $\lambda_j$ of $f$ is negative with respect to $\mathcal C_{2k}$ for all
 $k\in K$, then  the ordinary integral $$\int_{P(F)\backslash 
G(\mathbb A)^1}f(g)\cdot\tau_k\Big(H_P(g)\Big)\,dg$$ is absolutely convergent
 and equals to the \#-integral.}
\vskip 0.30cm
Fix a sufficiently regular element $T\in \frak a_0^+$. Then the above 
construction applies to Arthur's analytic truncation
$\Lambda^{T,P}\Psi(g)$ and the {\bf characteristic 
function} $\tau_P\Big(H_P(g)-T\Big)$ where $\Psi\in \mathcal A_P(G)$,
the spaec of level $P$ automorphic forms. 
Indeed, with
$$\Psi(namk)=\sum_{j=1}^sQ_j\Big(H_P(a)\Big)\cdot \psi_j(amk)$$ 
for $n\in N(\mathbb A), a\in A_P, m\in M(\mathbb A)^1$ and $k\in K$,
where $Q_j$ are polynomials and $\psi_j\in \mathcal A_P(G)$ satisfies 
$$\psi_j(ag)=e^{\langle\lambda_j+\rho_P,H_P(a)\rangle}\psi_j(g)$$ for 
some exponent $\lambda_j\in\frak a_P^*$ and all $a\in A_P$, the integration
$\int_{M(F)\backslash M(\mathbb A)^1\times K}\Lambda^{T,P}\psi_j(mk)dm\,dk$
is well-defined, since by the Basic Fact about Arthur's truncation,  
the function $m\mapsto \Lambda^{T,P}\psi_j(mk)$ is rapidly 
decreasing on
$M(F)\backslash M(\mathbb A)^1\times K$. Since $\tau_P$ is the 
characteristic function of the cone spanned by the coweights  
${\widehat \Delta}_P^\vee$, we see that
$$\int_{P(F)\backslash G(\mathbb A)^1}^{\#}\Lambda^{T,P}\Psi(g)
\cdot  \tau_P\Big(H_P(g)-T\Big)\,dg$$ exists if and only if
$$\langle\lambda_j,\varpi^\vee\rangle\not=0 \qquad \forall 
\varpi\in  {\widehat \Delta}_P^\vee\quad\mathrm{and}\quad 
\lambda_j\in\mathcal E_P(\Psi).$$ The same is true for $\Lambda^{T,P}\Psi(g)
\cdot \tau_P(H_P(gx)-T)$ for $x\in G(\mathbb A)$.
Indeed, for $x\in G(\mathbb A)$, let $K(g)\in K$ be an element such that 
$gK(g)^{-1}\in P_0(\mathbb A)$, then $H_P(gx)=H_P(g)+H_P(K(g)x)$ and hence 
$\tau_P\Big(X-T-H_P(K(g)x)\Big)$ is the characteristic function of a 
cone depending continuously on $g$.
Set then $$I^{P,T}_G(\Psi):=\int_{P(F)\backslash G(\mathbb A)^1}^{\#}
\Lambda^{T,P}\Psi(g)\cdot\tau_P\Big(H_P(g)-T\Big)\,dg.$$
\vskip 0.30cm
For any automorphic form $\phi\in \mathcal A(G)$, write 
$\mathcal E_P(\phi)$ for the set of exponents 
$\mathcal E_P(\phi_P)$. Set 
$$\mathcal A(G)^*:=\Big\{\phi\in \mathcal A(G):
\langle\lambda,\varpi^\vee\rangle\not=0 \qquad \forall 
\varpi^\vee\in {\widehat\Delta}_P^\vee,\ \lambda\in\mathcal E_P(\phi),\ 
P\not=G\Big\}.$$ If $\phi\in \mathcal A(G)^*$, then $I_{P,T}^{G}(\phi_P)$ 
exists for all $P$; and define JLR's {\it regularized period} by 
$$\int_{G(F)\backslash G(\mathbb A)}^*\phi(g)\,dg
=\sum_{P}I^{P,T}_G(\phi_P):=I^{T}_G(\phi).$$ 
\vskip 0.30cm
In the follows, let us expose some of the basic properties of JLR's 
regularized periods. For this, we make the following preperation.
As above, for $x\in G(\mathbb A_f)$, let $\rho(x)$ denote 
right translation by $x$: 
$\rho(x)\phi(g)=\phi(gx)$. The space $\mathcal A(G)$ is stable 
under right translation by $G(\mathbb A_f)$. Furthermore, 
$\rho(x)\phi$ has the same set of exponents as $\phi$. 
Consequently, the space $\mathcal A(G)_*$ is 
invariant under right translation by $G(\mathbb A_f)$.
Indeed, for $k\in K$, write the Iwasawa decomposition of $kx$ as 
$kx=n'a'm'K(kx)$ and write $\phi_P$ in the form 
$$\phi_P\Big(namk\Big)
=\sum_jQ_j\Big(H_P(a)\Big)\cdot e^{\langle \lambda_j+\rho_P,
H(a)\rangle}\phi_j(mk).$$ Since $amkx=n^* aa'mm'K(kx)$ for some 
$n^*\in N(\mathbb A)$, we have
$$\phi_P\Big(namk\cdot x\Big)=\sum_jQ_j\Big(H_P(a)+H_P(a')\Big)\cdot 
e^{\langle \lambda_j+\rho_P,H(a)+H(a')\rangle}\phi_j\Big(mm'K(kx)\Big).$$ 
This shows that $\mathcal E_P\Big(\rho(x)\phi\Big)=\mathcal E_P(\phi).$
\vskip 0.30cm
\noindent
{\bf {\large Basic Facts.}} ([JLR]) 
{\it With the same notation as above,} 

\noindent
(1) {\it $I^{T}_G$ 
defines a $G(\mathbb A_f)^1$-invariant linear functional 
on $\mathcal A(G)_*$;}

\noindent
(2) {\it $I^{T}_G(\phi):=
\int_{G(F)\backslash G(\mathbb A)^1}^*\phi(g)\,dg$ is 
independent of the choice of $T$};

\noindent
(3) {\it If $\phi\in\mathcal A(G)$ is integrable over 
$G(F)\backslash G(\mathbb A)^1$, then $\phi\in\mathcal A(G)_*$ and 
$$\int_{G(F)\backslash G(\mathbb A)^1}^*\phi(g)\,dg=
\int_{G(F)\backslash G(\mathbb A)^1}\phi(g)\,dg.$$}

\noindent
{\bf Remark.} Recall that we have the inversion formula
$$\phi(g)=\sum_P\sum_{\delta\in P(F)\backslash G(F)}\Lambda^{T,P}
\phi(\delta g)\cdot\tau_P\Big(H(\delta g)-T\Big).$$
Thus, formally,
$$\begin{aligned}&\int_{G(F)\backslash G(\mathbb A)^1}\phi(g)\,dg\\
=&\int_{G(F)\backslash G(\mathbb A)^1}\sum_P\sum_{\delta\in P(F)
\backslash G(F)}\Lambda^{T,P}\phi(\delta g)\cdot\tau_P
\Big(H(\delta g)-T\Big)\,dg\\
=&\sum_P
\int_{G(F)\backslash G(\mathbb A)^1}\sum_{\delta\in P(F)\backslash G(F)}
\Lambda^{T,P}\phi(\delta g)\cdot\tau_P\Big(H(\delta g)-T\Big)\,dg\\
=&\sum_P\int_{P(F)\backslash G(\mathbb A)^1}
\Lambda^{T,P}\phi(g)\cdot\tau_P\Big(H(g)-T\Big)\,dg\\
=&\sum_P\int_{P(F)\backslash G(\mathbb A)^1}^*
\Lambda^{T,P}\phi_P(g)\cdot\tau_P\Big(H(g)-T\Big)\,dg.\end{aligned}$$  
In this sense, in the Basic Facts, (2) is natural, while (3) gives a 
justification of the above calculation, since by definition
$$\begin{aligned}
I_G(\phi):=&\int_{G(F)\backslash G(\mathbb A)^1}^*\phi(g)\,dg\\
=&\sum_P I_P^T(\phi_P)\\
:=&\sum_P\int_{P(F)\backslash G(\mathbb A)}^{\#}
\Lambda^{T,P}\phi_P(g)\cdot\tau_P\Big(H(g)-T\Big)\,dg.\end{aligned}$$
\vskip 0.30cm
\noindent
{\it Proof.} ([JLR]) (1) This breaks into several steps.
\vskip 0.30cm
\noindent
{\bf Step A}. We start with the following natural:

\noindent
{\bf Claim 1} {\it For $f$ as above, we have
$$\int_{P(F)\backslash G(\mathbb A)^1}^{\#}f(g)\tau_k(gx)\,dg=
\int_{P(F)\backslash G(\mathbb A)^1}^{\#}f(g)\tau_k(g)\,dg,$$ for 
all $x\in G(\mathbb A_f)^1$.}
 
This is easy: Set $f_\mu(g):=e^{\langle \mu, H_P(g)\rangle}
f(g)$ for $\mu\in\frak a_P^*$. If $\langle \mathrm{Re}\,(\mu),
\alpha^\vee\rangle\ll 0$ for all $\alpha\in\Delta_P$,
then $$\int_{P(F)\backslash G(\mathbb A)^1}f(gx)\tau_k(gx)\,dg=
\int_{P(F)\backslash G(\mathbb A)^1}f_\mu(g)\tau_k(g)\,dg,$$
by the invariance of Haar measure, based on the fact that
 both sides are 
absolutely convengent. Since both sides have the same meromorphic 
continuation, evaluate it at $\mu=0$ gives the claim.
\vskip 0.30cm
\noindent
{\bf Step B}. Now fix $x\in G(\mathbb A_f)^1$ and set 
$F_P(g):=\Big(\Lambda^{T,P}\rho(x^{-1})\phi\Big)(gx).$ Then 
from the previous claim, we have
$$\int_{P(F)\backslash G(\mathbb A)^1}^{\#}
\Big(\Lambda^{T,P}\rho(x^{-1})\phi\Big)(gx)\cdot 
\tau_P\Big(H_P(g)-T\Big)\,dg=
\int_{P(F)\backslash G(\mathbb A)^1}^{\#}F_P(g)\cdot 
\tau_P\Big(H_P(g)-T\Big)\,dg.$$
Hence, we mush show the following

\noindent
{\bf Claim 2.} With the same notation as above,
$$\int_{G(F)\backslash G(\mathbb A)^1}^*\phi(g)\,dg=\sum_P
\int_{P(F)\backslash G(\mathbb A)^1}^{\#}F_P(g)\cdot
\tau_P(H_P(gx)-T)\,dg.
\eqno(*)$$

For this purpose, we want to write $F_P$ in terms of certain auxiliary function
$\Gamma_R^P$ defined by
$$\Gamma_Q^P(H,X):=\sum_{R:Q\subset R\subset P}(-1)^{d(R)-d(P)}
\tau_Q^R(H)\widehat\tau_R^P(H-X),$$
and $\nabla_Q^P$ by $$\nabla_Q^P(H,X):=
\sum_{R:Q\subset R\subset P}(-1)^{d(Q)-d(R)}
\tau_Q^R(H-X)\widehat\tau_R^P(H).$$
Dually, we also introduce $\widehat \Gamma_R^P$ defined by
$$\widehat\Gamma_Q^P(H,X):=\sum_{R:Q\subset R\subset P}(-1)^{d(P)-d(R)}
\widehat\tau_Q^R(H)\tau_R^P(H-X),$$
and also $\widehat\nabla_Q^P$ by $$\widehat\nabla_Q^P(H,X):=
\sum_{R:Q\subset R\subset P}(-1)^{d(Q)-d(R)}
\widehat\tau_Q^R(H-X)\tau_R^P(H).$$

\noindent
{\bf Remarks.} (1) The simplest case is when $\frak a_P^G$ is one-dimensional. 
Suppose for example that $G=SL(3)$ amd $P=P_1$ is a maximal parabolic 
subgroup. Then $Q$ is summed over the set $\{P_1,G\}$. Taking $X$ to be a 
fixed point in positive chamber in $\frak a_P^G$, we see that 
$H\mapsto \Gamma_P^G(H,X)$ is the difference of characteristic functions 
of two open half lines, and is hence the characteristic function of the 
bounded half open interval. As Arthur illustrates, with thism one may 
even get a non-trivial better picture by considering the case with $P=P_0$, 
the minimal parabolic subgroup of $SL(3)$. (See e.g. [Ar6].)

\noindent
(2) Please note that our $\widehat\Gamma$ is different from that 
in [JLR], which is supposed to be $\nabla$ in our system of notations.

Then we have the following
\vskip 0.30cm
\noindent
{\bf {\large Lemma 2.}} {\it With the same notation as above,}

\noindent
(i) $\displaystyle{\Gamma_Q^P(H-X,-X)=(-1)^{d(P)-d(Q)}
\nabla_Q^P(H,X)}$;

\noindent
(ii) $\displaystyle{\tau_S^P(H-X)=\sum_{R:S\subset R\subset P}
\Gamma_S^R(H-X,-X)\cdot\tau_R^P(H);}$

\noindent
(iii) $\displaystyle{\widehat\tau_Q^P(H-X)=
\sum_{R:Q\subset R\subset P}
(-1)^{d(R)-d(P)}
\widehat\tau_Q^R(H)\cdot\Gamma_R^P(H,X).}$

Dually, we also have
\vskip 0.30cm
\noindent
{\bf {\large Lemma $\hat 2$.}} {\it With the same notation as above,}

\noindent
(i) $\displaystyle{\widehat\Gamma_Q^P(H-X,-X)=(-1)^{d(P)-d(Q)}
\widehat\nabla_Q^P(H,X)}$;

\noindent
(ii) $\displaystyle{\widehat\tau_S^P(H-X)=\sum_{R:S\subset R\subset P}
\widehat\Gamma_S^R(H-X,-X)\cdot\widehat\tau_R^P(H);}$

\noindent
(iii) $\displaystyle{\tau_Q^P(H-X)=\sum_{R:Q\subset R\subset P}
(-1)^{d(R)-d(P)}
\tau_Q^R(H)\cdot\widehat\Gamma_R^P(H,X).}$
\vskip 0.30cm
\noindent
{\it Proof of the Lemmas.} The proof of Lemma 2 and Lemma $\hat 2$ are similar.
We here give one for the later.

\noindent
(i) By definition,
$$\begin{aligned}\widehat\Gamma_Q^P(H-X,-X):=&
\sum_{R:Q\subset R\subset P}(-1)^{d(P)-d(R)}
\widehat\tau_Q^R(H-X)\tau_R^P\Big(H-X-(-X)\Big)\\
=&\sum_{R:Q\subset R\subset P}
(-1)^{d(P)-d(R)}\widehat\tau_Q^R(H-X)\tau_R^P(H)\\
=&(-1)^{d(P)-d(Q)}\sum_{R:Q\subset R\subset P}
(-1)^{d(Q)-d(R)}\widehat\tau_Q^R(H-X)\tau_R^P(H)\\
=&(-1)^{d(P)-d(Q)}\widehat\nabla_Q^P(H,X);\end{aligned}$$

\noindent
(ii) By definition,
$$\begin{aligned}&\sum_{R:S\subset R\subset P}
\widehat\Gamma_S^R\Big(H-X,-X\Big)\cdot\widehat\tau_R^P(H)\\
=&\sum_{R:S\subset R\subset P}(-1)^{d(R)-d(S)}\widehat\nabla_S^R(H,X)
\cdot\widehat\tau_R^P(H)\\
=&\sum_{R:S\subset R\subset P}(-1)^{d(R)-d(S)}
\sum_{T:S\subset T\subset R}(-1)^{d(S)-d(T)}
\widehat\tau_S^T(H-X)\tau_T^R(H)
\cdot\widehat\tau_R^P(H)\\
=&\sum_{R,T:S\subset T\subset R\subset P}(-1)^{d(R)-d(T)}
\widehat\tau_S^T(H-X)\tau_T^R(H)
\cdot\widehat\tau_R^P(H)\\
=&\sum_{T:S\subset T}\widehat\tau_S^T(H-X)\sum_{R:
T\subset R\subset P}(-1)^{d(R)-d(T)}
\tau_T^R(H)\cdot\widehat\tau_R^P(H)\\
=&\sum_{T:S\subset T}\widehat\tau_S^T(H-X)\delta_{TP}\end{aligned}$$ 
by Langlands' Combinatorial Lemma. Therefore,
$$\sum_{R:S\subset R\subset P}\widehat\Gamma_S^R\Big(H-X,-X\Big)
\cdot\widehat\tau_R^P(H)=\widehat\tau_S^P(H-X)$$ as required;

\noindent
(iii) By definition,
$$\begin{aligned}
&\sum_{R:Q\subset R\subset P}
(-1)^{d(R)-d(P)}
\tau_Q^R(H)\cdot\widehat\Gamma_R^P(H,X)\\
=&\sum_{R:Q\subset R\subset P}
(-1)^{d(R)-d(P)}\tau_Q^R(H)
\sum_{S:R\subset S\subset P}(-1)^{d(P)-d(S)}
\widehat\tau_R^S(H)\tau_S^P(H-X)\\
=&\sum_{R,S:Q\subset R\subset S\subset P}
(-1)^{d(R)-d(S)}\tau_Q^R(H)
\widehat\tau_R^S(H)\tau_S^P(H-X)\\
=&\sum_{S:S\subset P}\Big(\sum_{R:Q\subset R\subset S}
(-1)^{d(R)-d(S)}\tau_Q^R(H)
\widehat\tau_R^S(H)\Big)\tau_S^P(H-X)\\
=&\sum_{S:S\subset P}\Delta_{QS}\tau_S^P(H-X)\end{aligned}$$
by Langlands' Combinatorial Lemma. Therefore,
$$\sum_{R:Q\subset R\subset P}
(-1)^{d(R)-d(P)}
\tau_Q^R(H)\cdot\widehat\Gamma_R^P(H,X)=\tau_Q^P(H-X)$$ as required. 
This completes the proof.
\vskip 0.30cm
To continue the proof of Claim 2, as before, for $g\in G(\mathbb A)$, let $K(g)\in K$ be an element 
such that $gK(g)^{-1}\in P_0(\mathbb A)$.
\vskip 0.30cm
\noindent
{\bf Subclaim 1.} $\displaystyle{F_P(g)=\sum_{S:S\subset P}
\sum_{\eta\in S(F)\backslash P(F)}
\Lambda^{T,S}\phi_S(\eta g)\cdot\Gamma_S^P
\Big(H_S(\eta g)-T,-H_S(K(\eta g)x)\Big)}$
$$=\sum_{S:S\subset P}
\sum_{\eta\in S(F)\backslash P(F)}
\Lambda_{T,S}\phi_S(\eta g)\cdot\widehat\Gamma_S^P
\Big(H_S(\eta g)-T,-H_S(K(\eta g)x)\Big).$$
\vskip 0.30cm
The proof of the first is similar as that for the second. So here we only give 
a proof for the later one. Indeed, by definition, 
$H(\delta g x)-T=H(\delta g)-T+H(K(\delta g)x)$ and
$$\begin{aligned}F_P(g)=&\sum_{R:R\subset P}(-1)^{d(R)-d(P)}
\sum_{\delta\in R(F)\backslash P(F)}\phi_R(\delta g)
\cdot\tau_R^P\Big(H(\delta gx)-T\Big)\\
=&\sum_{R:R\subset P}(-1)^{d(R)-d(P)}\sum_{\delta\in R(F)\backslash P(F)}
\phi_R(\delta g)\cdot\tau_R^P\Big(H(\delta g)-T-
\big(-(H(K(\delta g)x))\big)\Big)\\
=&\sum_{R:R\subset P}(-1)^{d(R)-d(P)}\sum_{\delta\in R(F)\backslash P(F)}
\phi_R(\delta g)\cdot\sum_{S:R\subset S\subset P}(-1)^{d(S)-d(P)}\tau_R^S(H)
\widehat\Gamma_S^P(H,X)\end{aligned}$$ with $H=H(\delta g)-T,\ 
X=-H(K(\delta g)x).$
Here, we have used the inverse relation (iii) above, i.e.,
$$\tau_Q^P(H-X)=\sum_{R:Q\subset R\subset P}(-1)^{d(R)-d(P)}\tau_Q^R(H)
\widehat\Gamma_R^P(H,X).$$
Therefore,
$$\begin{aligned}&F_P(g)\\
=&\sum_{S:S\subset P}\sum_{\eta\in S(F)\backslash P(F)}
\sum_{R:R\subset S}\sum_{\delta\in R(F)\backslash S(F)}(-1)^{d(R)-d(S)}
\phi_R(\delta\eta g)\cdot\tau_R^S(H)\widehat\Gamma_S^P
\Big(H_S(\eta g)-T,-H_S(K(\eta g)x)\Big)\\
=&\sum_{S:S\subset P}\sum_{\eta\in S(F)\backslash P(F)}\Lambda_{T,S}
\phi_S(\eta g)\cdot \widehat\Gamma_S^P\Big(H_S(\eta g)-T,-H_S(K(\eta g)x)\Big).
\end{aligned}$$ 
This completes of the proof of the subclaim.
\vskip 0.30cm
\noindent
{\bf Step C}. Thus we have
$$\begin{aligned}&\int_{P(F)\backslash G(\mathbb A)^1}^{\#}F_P(g)\cdot
\tau_P\Big(H_P(gx)-T\Big)\,dg\\
=&\int_{P(F)\backslash G(\mathbb A)^1}^{\#}
\sum_{S:S\subset P}\sum_{\eta\in S(F)\backslash P(F)}\Lambda^{T,S}
\phi_S(\eta g)\cdot \Gamma_S^P\Big(H_S(\eta g)-T,-H_S(K(\eta g)x)\Big)
\cdot\tau_P\Big(H_P(gx)-T\Big)\,dg.\end{aligned}\eqno(**)$$

\noindent 
{\bf Subclaim 2.} {\it With our assumption effective, the sum over $S$ can 
be taken outside the integral.}
\vskip 0.30cm
\noindent
{\it Proof of the Subclaim.} (a) First, note that the functions 
$\Gamma_S^P(Z,W)$ depend only on 
the projections of $Z$ and $W$ on $\frak a_S^P$.
Then we note the following
\vskip 0.30cm
\noindent
{\bf Lemma 3.} (Arthur) {\it For each $X$ in a  fixed compact subset of 
$\frak a_P/\frak a_G$, the support of the function 
$H\mapsto \Gamma_P^G(H,X),\ H\in
\frak a_P/\frak a_G$ is contained in a fixed compact subset, which is 
independent of $X$.}

\noindent
{\it Proof.} (Arthur)  Step 1. {\it An auxiliary function which is compact-supported}:
If $Q\supset P$, set $\widehat\tau^{Q/P}(H)$ equal to the characteristic 
function of $$\Big\{H:\varpi(H)>0,\varpi\in\widehat\Delta_P\backslash
\widehat\Delta_Q\Big\}.$$ Define the the function $\widetilde\Gamma_Q^G(H,X)$ 
inductively by the condition that for all $Q\supset P_0$.
$$\widehat\tau_Q(H-X)=\sum_{R:R\supset Q}(-1)^{d(R)-d(G)}
\widehat\tau^{R/Q}(H)\widetilde\Gamma_R^G(H,X).$$ Then
$$\widetilde\Gamma_Q(H,X)=\sum_{R:R\supset Q}(-1)^{d(R)-d(G)}
\widehat\tau^{R/Q}(H)\widehat\tau(H-X).$$
Directly from te definition, we see that, modulo sign,
$\widetilde \Gamma_Q^G(\cdot,X)$ is just the characteristic function in 
$\frak a_Q/\frak a_G$ of a parallelepiped with opposite vertices 0 and $X$. 
In particular, $\widetilde \Gamma_Q^G(H,X)$ is compactly supported as a 
function of $H\in \frak a_Q/\frak a_G$.

\noindent
Step 2. {\it Seperate $X$ from $X$}. By definition,
$$\begin{aligned}\Gamma_P^G(H,X)=&
\sum_{Q:Q\subset P}(-1)^{d(Q)-d(G)}\tau_P^Q(H)\widehat
\tau_Q(H-X)\\
=&\sum_{R:R\subset P}\tilde\Gamma_R^G(H,X)\sum_{Q:P\subset Q\subset R}
(-1)^{d(Q)-d(R)}\tau_P^Q(H)\widehat\tau^{R/Q}(H).\end{aligned}$$

\noindent
Step 3. {\it An induction on $\mathrm{dim}(G)$}.
By Langlands' combinatorial lemma, if $P\not=G$,
$$\sum_{Q:P\subset Q\subset R}
(-1)^{d(Q)-d(R)}\tau_P^Q(H)\widehat\tau_Q(H)=0.$$
Thus, if $R=G\not=P$, $$\sum_{Q:P\subset Q\subset R}
(-1)^{d(Q)-d(R)}\tau_P^Q(H)\widehat\tau^{R/Q}(H)=\sum_{Q:P\subset Q\subset R}
(-1)^{d(Q)-d(R)}\tau_P^Q(H)\widehat\tau_Q(H)=0.$$
Thereore, in the above expression of $\Gamma_P^G(H,X)$, the outer sum may 
be taken over only those $R$ which are not equal to $G$. For a given 
$R\not=G$, and 
$H=H_P^R+H_R^G\in\frak a_P^G$ according to the decomposition
$\frak a_P^G=\frak a_P^R\oplus\frak a_R^G$, clearly
$\widetilde\Gamma_R^G(H,X)=\widetilde\Gamma_R^G(H_R^G,X).$ Moreover, by 
definition, we may rewrite
$$\widehat\tau^{R/Q}(H)=\widehat\tau^{R/Q}\Big(H_P^R-L(H_R^G)\Big),$$ for 
a certain linear map $H_R^G\mapsto L(H_R^G)$ from $\frak a_R^G$ to 
$\frak a_P^R$
 which is independent of $Q$. If the summand corresponding to $R$ does not 
vanish, $H_R^G$ will lie in a fixed compact set. So, therefore, will 
$L(H_R^G)$. Applying the induction assumption to the group $M_R$, we see 
that $H_P^R$ must lie in a fixed compact subset of $\frak a_P^R$. It follows 
that $H$ is contained in a fixed compact subset of $\frak a_P^G$. This 
completes the proof.

\noindent
{\bf Remark.} A simple example with $GL_3$ will show that this subclaim is 
not really true when $\widehat\Gamma$ is used.

Back to our proof. Clearly, a similar argument works for 
$\Gamma_S^P$ as well. Consequently, for $W=-H_S(K(\eta g)x)$ belonging to 
a fixed compact subset of $\frak a_S^P$, there exists a compact subset 
$\mathcal Y\subset\frak a_S^P$ such that the function $g\mapsto 
\Gamma_S^P\Big(H_S(\eta g)-T,-H_S(K(\eta g)x)\Big)$ is supported 
inside a subset of elements $g$ for which the projection of $H_S(g)$ onto 
$\frak a_S^P$ lies in  a compact set depending only on $x$. 

Therefore, being a finite sum, $$\sum_{\eta\in S(F)\backslash P(F)}
\Lambda^{T,S}\phi_S(\eta g)\cdot\Gamma_S^P
\Big(H_S(\eta g)-T,-H_S(K(\eta g)
\cdot x)\Big)$$  is integrable over $M(F)\backslash 
M(\mathbb A)^1$.
\vskip 0.30cm
\noindent 
(b) Now let $g=n\cdot e^X\cdot m\cdot k$ be an Iwasawa decomposition of 
$g$ relative to $P$ with $X\in \frak a_P$. There exist polynomials $Q_j$ 
on $\frak a_P$, automorphic forms $\psi_j\in \mathcal A_P(G)$, and exponents 
$\lambda_j$ in the $\mathcal E_P(\phi)$ such that for all $S\subset P$,
$$\Lambda^{T,S}\phi_S(g)=\sum_j Q_j(X)\cdot e^{\langle \lambda_j+\rho_P,
X\rangle}\cdot\Lambda^{T,S}\psi_j(mk).$$

Since $\Gamma_S^P(Z,W)$ depends only on the projections 
of $Z$ and $W$ onto 
$\frak a_S^P$, we may write the right hand side of (**) 
as the integral over $k\in K$ and sum over 
$j$ of
$$\begin{aligned}&\int_{\frak a_P}^{\#}\int_{M(F)\backslash M(\mathbb A)^1}
\sum_{S:S\subset P}\sum_{\eta\in S(F)\backslash P(F)}Q_j(X)\cdot 
e^{\langle\lambda_j,X\rangle}\cdot\Lambda^{T,S}\psi_j(\eta mk)\\
&\qquad\Gamma_S^P\Big(H_S(\eta m)-T,-H_S(K(\eta mk)x)\Big)\,dm
\cdot\tau_P\Big(X+H(kx)-T\Big)\,dX\\
=&\int_{\frak a_P}^{\#}\sum_jQ_j(X)\cdot e^{\langle\lambda_j,X\rangle}
\cdot\tau_P\Big(X+H(kx)-T\Big)\,dX\\
&\int_{M(F)\backslash M(\mathbb A)^1}\sum_{S:S\subset P}
\sum_{\eta\in S(F)\backslash P(F)}\Lambda^{T,S}\psi_j(\eta mk)\cdot
\Gamma_S^P\Big(H_S(\eta m)-T,-H_S(K(\eta mk)x)\Big)\,dm.\end{aligned}$$
Now,  by our basic fact for Arthur's truncation, 
each term in the sum over $S$  is seperately integrable over 
$M(F)\backslash M(\mathbb A)^1$, we may take the sum over $S$ outside 
the integral as claimed.
\vskip 0.30cm
\noindent 
{\bf Step D}. We also have the following natural

\noindent{\bf Subclaim 3.} {\it With the same notation as above,}
 $$\begin{aligned}
\int_{P(F)\backslash G(\mathbb A)^1}^{\#}&
\sum_{\eta\in S(F)\backslash P(F)}\Lambda^{T,S}\phi_S(\eta g)\cdot
\Gamma_S^P\Big(H_S(\eta g)-T,-H_S(K(\eta g)x)\Big)
\cdot\tau_P(H(gx)-T)\,dg\\
=&\int_{S(F)\backslash G(\mathbb A)^1}^{\#}
\Lambda^{T,S}\phi_S(g)\cdot
\Gamma_S^P(H_S(g)-T,-H_S(K(g)x))\cdot\tau_P(H(gx)-T)\,dg.
\end{aligned}$$

\noindent
{\it Proof of the Subclaim.}  (a) LHS is equal to the integral 
over $k\in K$ and sum over $j$ of
$$\begin{aligned}
\int_{\frak a_P}^{\#}&Q_j(X)\cdot e^{\langle\lambda_j,X\rangle}
\cdot\tau_P\Big(X+H(kx)-T\Big)\,dX\\
&\cdot\int_{M(F)\backslash M(\mathbb A)^1}\sum_{\eta\in S(F)\backslash P(F)}
\Lambda^{T,S}\psi_j(\eta mk)\cdot
\Gamma_S^P\Big(H_S(\eta m)-T,-H_S(K(\eta mk)x)\Big)\,dm\end{aligned}$$ 
which can be written as
$$\begin{aligned}
\int_{\frak a_P}^{\#}&Q_j(X)\cdot e^{\langle\lambda_j,X\rangle}\cdot
\tau_P\Big(X+H(kx)-T\Big)\,dX\\
&\cdot\int_{S_M(F)\backslash M(\mathbb A)^1}\Lambda^{T,S}\psi_j(mk)\cdot
\Gamma_S^P\Big(H_S(m)-T,-H_S(K(mk)x)\Big)\,dm\end{aligned}$$ 
where $S_M:=S\cap M$.

\noindent 
(b) Expressing the integral over $M_S(F)\backslash M(\mathbb A)^1$ using 
the Iwasawa decomposition $m=n'e^{X'}m'k'$ of $M(\mathbb A)^1$ relative 
to $S_M$ gives
$$\begin{aligned}
\int_{\frak a_P}^{\#}&Q_j(X)\cdot e^{\langle\lambda_j,X\rangle}\cdot
\tau_P\Big(X+H(kx)-T\Big)\,dX\\
&\cdot\int_{K_M}\int_{\frak a_S^P}^{\#}
\int_{M_S(F)\backslash M_S(\mathbb A)^1}e^{-\langle \rho_S^P,X'\rangle}
\Lambda^{T,S}\psi_j(e^{X'}m'k'k)\cdot
\Gamma_S^P\Big(X'-T,-H_S(k'kx)\Big)\,dm'\,dX'\,dk',\end{aligned}$$ 
where $K_M=K\cap M(\mathbb A)^1$.

\noindent
(c) Since we are integrating over $K$, we may drop the integration 
over $K_M$.

\noindent
(d) Moreover, each function $\Lambda^{T,S}\psi_j(e^{X'}m'k)$ has an analogous
decomposition with respect to $S$. Therefore for fixed $k$,
the function $$(m',X')\mapsto \sum_j\Lambda^{T,S}\psi_j(e^{X'}m'k)\cdot
\Gamma_S^P\Big(X'-T,-H_S(kx)\Big)$$ is a sum of terms each of which is the 
product of 

\noindent
(i) a function of $m'$ which  is absolutely integrable over 
$M_S(F)\backslash M_S(\mathbb A)^1$ by our basic assumption; and 

\noindent
(ii) a function of $X'$ which 
itself is equal to an exponential polynomials, times the compactly 
supported function $\Gamma_S^P\Big(X'-T,-H_S(kx)\Big)$.
\vskip 0.30cm
Moreover, by basic property (5) of integration over cones, 
we may combine the integrals over $\frak a_S^P$ and 
$\frak a_P$ to a $\#$-integral over $\frak a_S$ and we find that
the LHS of the claim 
is equal to the integral over $k\in K$ and $m'\in M_S(F)\backslash 
M_S(\mathbb A)^1$ of $$\int_{\frak a_S}^{\#} e^{-\langle\rho_P,X\rangle}
\Lambda^{T,S}\phi_S(e^{X}m'k)\cdot
\Gamma_S^P\Big(X-T,-H_S(kx)\Big)\cdot\tau_P\Big(X+H(kx)-T\Big)\,dX$$ 
this equals the RHS of the claim. We are done.
\vskip 0.30cm
\noindent
{\bf Step E}. In Subclaim 3, summing the RHS over all $S$ and $P$ such 
that $S\subset P$ we see that
$$\sum_P\int_{P(F)\backslash G(\mathbb A)^1}^{\#}F_P(g)\cdot 
\tau_P\Big(H(gx)-T\Big)\,dg$$ is equal to the sum 
over parabolic subgroups $S$ of
$$\sum_{P:P\supset S}\int_{S(F)\backslash G(\mathbb A)^1}^{\#}
\Lambda^{T,S}\phi_S(g)\cdot\Gamma_S^P\Big(H_S(g)-T,-H_S(K(g)\cdot x)\Big)
\cdot \tau_P\Big(H(gx)-T\Big)\,dg.$$  
Thus the proof of (*), i.e., Claim 2, will be completed 
if we prove the following

\noindent
{\bf Claim 2}$'$. $\displaystyle{\sum_{P:P\supset S}
\int_{S(F)\backslash G(\mathbb A)^1}^{\#}
\Lambda^{T,S}\phi_S(g)\cdot\Gamma_S^P\Big(H_S(g)-T,-H_S(K(g)\cdot x)\Big)
\cdot \tau_P(H(gx)-T)\,dg}$
$$=\int_{S(F)\backslash G(\mathbb A)^1}^{\#}
\Lambda^{T,S}\phi_S(g)\cdot \tau_S\Big(H(g)-T\Big)\,dg.$$
 
This is based on the  relation (ii) in the Lemma of Step B: 
$$\tau_S(Y-X)=\sum_{P\supset S}\Gamma_S^P(Y-X,-X)
\cdot\tau_P(Y).$$
Indeed,  applying this claim  to $Y=H_S(gx)-T$ and $X=H_S(gx)-H_S(g)$
$=H_S(K(g)x)$ 
gives $$\sum_{P:P\supset S}\Gamma_S^P\Big(H_S(g)-T,-H_S(K(h)x)\Big)
\cdot\tau_P\Big(H(gx)-T\big)=\tau_S\Big(H(g)-Y\Big).$$
Thus we need to show the

\noindent 
{\bf Subclaim 4.} {\it The summation over $P$ in the RHS of the claim can be 
taken inside the \#-integral.}

\noindent
Proof. Let $g=ne^Xmk$ with $X\in\frak a_S$ be the Iwasawa decomposition of 
$g\in G(\mathbb A)$ relative to $S$. Then 
$$\Gamma_S^P\Big(H_S(g)-T,-H_S(K(h)x)\Big)
\cdot\tau_P\Big(H(gx)-T\Big)=
\Gamma_S^P\Big(X-T,-H_S(kx)\Big)\cdot\tau_P\Big(X+H(kx)-T\Big).$$ Since
the subset $\Big\{H_S(kx):k\in K\Big\}$ of $\frak a_S$ (for $x$ fixed) 
is compact,  there exists a compact subset 
$\mathcal Y\subset\frak a_S^P$ such that the function $g\mapsto 
\Gamma_S^P\Big(H_S(\eta g)-T,-H_S(K(\eta g)x)\Big)$ 
is supported inside a subset 
of elements $g$ for which the projection of $H_S(g)$ onto $\frak a_S^P$ lies 
in  a compact set depending only on $x$.
The cone defining $\tau_P$ is the positive Weyl chamber in $\frak a_P$ 
which is contained in the positive cone of $\frak a_S$
Thus we may apply basic property (6) of integration over cones 
to take the sum over $P$ inside the integral. 
This completes the proof of assertion (1) of the Theorem.
\vskip 0.30cm
\noindent
The proof of (2) is nearly identical. 
First, suppose that $T'\in\frak a_0$ is regular. We have 
$$\Lambda^{T+T',P}\phi(g)=\sum_{Q:Q\subset P}\sum_{\delta\in 
Q(F)\backslash P(F)}\Lambda^{T,Q}\phi(\delta g)\cdot\Gamma_Q^P
\Big(H(\delta g)-T,T'\Big).$$
Indeed, using the formula (iii) proved in the Lemma of Step B:
$$\widehat\tau_R^P\Big(H(\delta g)-T-T'\Big)=\sum_{Q:R\subset Q
\subset P}(-1)^{d(Q)-d(P)}\widehat\tau_R^Q\Big(H(\delta g)-T\Big)
\Gamma_Q^P\Big(H(\delta g)-T,T'\Big)$$ we have
$$\begin{aligned}&\Lambda^{T+T',P}\phi(g)\\
=&\sum_{R:R\subset P}(-1)^{d(R)-d(P)}
\sum_{\delta\in R(F)\backslash P(F)}
\phi_R(\delta g)\cdot\widehat\tau_R^P\Big(H(\delta g)-T-T'\Big)\\
=&\sum_{R, Q:R\subset Q\subset P}\sum_{\delta\in Q(F)\backslash P(F)}
(-1)^{d(R)-d(Q)}\bigg(\sum_{\gamma\in R(F)\backslash Q(F)}\phi_R(\delta g)
\widehat\tau_R^Q\Big(H(\gamma \delta g)-T\Big)\bigg)\cdot
\Gamma_Q^P\Big(H(\delta g)-T,T'\Big)\\
=&\sum_{Q:Q\subset P}\sum_{\delta\in Q(F)\backslash P(F)}\bigg(
\sum_{R:R\subset Q}(-1)^{d(R)-d(Q)}
\sum_{\gamma\in R(F)\backslash Q(F)}\phi_R(\gamma\delta g)
\widehat\tau_R^Q\Big(H(\gamma \delta g)-T\Big)\bigg)\cdot\Gamma_Q^P
\Big(H(\delta g)-T,T'\Big)\\
=&\sum_{Q:Q\subset P}\sum_{\delta\in Q(F)\backslash P(F)}\Lambda^{T,Q}
\phi(\delta g)\cdot\Gamma_Q^P\Big(H(\delta g)-T,T'\Big).\end{aligned}$$
Therefore, $$\begin{aligned}&
\sum_P\int_{P(F)\backslash G(\mathbb A)^1}^{\#}
\Lambda^{T+T',P}\phi(g)\cdot\tau_P\Big(H(g)-T-T'\Big)\,dg\\
=&\sum_P\int_{P(F)\backslash G(\mathbb A)^1}^{\#}\sum_{Q:Q\subset P}
\sum_{\delta\in Q(F)\backslash P(F)}\Lambda^{T,Q}\phi(\delta g)\cdot
\Gamma_Q^P\Big(H(\delta g)-T,T'\Big)\cdot
\tau_P\Big(H(g)-T-T'\Big)\,dg\\
=&\sum_{Q,P:Q\subset P}\int_{Q(F)\backslash G(\mathbb A)^1}^{\#}
\Lambda^{T,Q}\phi(g)\cdot\Gamma_Q^P\Big(H(g)-T,T'\Big)\cdot
\tau_P\Big(H(g)-T-T'\Big)\,dg\\
=&\sum_{Q}\int_{Q(F)\backslash G(\mathbb A)^1}^{\#}
\Lambda^{T,Q}\phi(g)\sum_{P:P\supset Q}\Gamma_Q^P
\Big(H(g)-T,T'\Big)\cdot
\tau_P\Big(H(g)-T-T'\Big)\,dg,\end{aligned}$$ 
where the second equality is justified 
in the same way as in Step D and the third 
equality is justified as in the Subclaim of Step E above.
 
Thus, using the relation (ii) of the Lemma in Step B: $$\sum_{P:P\supset Q}
\Gamma_Q^P\Big(H(\delta g)-T,T'\Big)\cdot
\tau_P\Big(H(g)-T-T'\Big)=\tau_Q\Big(H(g)-T\Big),$$  we
obtain
$$\sum_{Q}\int_{Q(F)\backslash G(\mathbb A)^1}^{\#}
\Lambda^{T,Q}\phi(g)\tau_Q\Big(H(g)-T\Big)\,dg=
\int_{G(F)\backslash G(\mathbb A)^1}^*\phi(g)\,dg$$ as required.
\vskip 0.30cm
\noindent
We now prove (3)  By the inversion formula
$$\phi(g)=\sum_{P\subset G}\sum_{\delta\in P(F)\backslash G(F)}
\Lambda^{T,P}\phi(\delta g)\cdot\tau_P\Big(H(\delta g)-T\Big),$$ 
it suffices to check that  if $\phi\in\mathcal A(G)^*$ is integrable over 
$G(F)\backslash G(\mathbb A)^1$, then the integral
$$\int_{P(F)\backslash G(\mathbb A)^1}^{\#}
\Lambda^{T,P}\phi_P(g)\tau_P(H(g)-T)\,dg$$  is well-defined and is 
equal to $I_P^{G,T}(\phi_P).$ 
Expand $\phi_P$ as usual, as shown in [MW], top of p.50, 
for all $j$, there exist
a parabolic subgroup $Q\subset P$ and a cuspidal exponent $\mu$ of 
$\phi_Q^{cusp}$ such that the exponent $\lambda_j$ of $\phi_P$
is equal to the restriction of $\mu$ to $\frak a_P$ relative to the 
decomposition $\frak a_Q=\frak a_Q^P\oplus\frak a_P$.
 
Moreover, according to Lemma I.4.11 of [MW], p.75,
if $\phi$ is $L^2$ on $G(F)\backslash G(\mathbb A)^1$, then the exponent 
$\mu$ can be written in the form $\sum_{\alpha\in\Delta_P}x_\alpha\alpha$ 
with $x_\alpha<0$. A nearly identical argument shows that this remains 
true if $\phi$ is assumed to be integrable over $G(F)\backslash 
G(\mathbb A)^1$. This says that $\mu$ is negative with respect to
the cone $\Big\{X\in\frak a_P:\tau_P(X)=1\Big\}$. Therefore, the integral 
$I_P^{G,T}(\phi_P)$ is absolutely convergent and coincides with 
the ordinary integral over
$P(F)\backslash G(\mathbb A)^1$ as required.

\subsection{Arthur's Period}

Recall that for $P=MN$  a parabolic subgroup and $f\in\mathcal A_P(G)$,
we may generalize the construction of the previous section so as to get the 
regularized integral
$$\int_{P(F)\backslash G(\mathbb A)^1}f(g)\cdot
\tau\Big(H_P(g)-T\Big)\,dg.$$
where $\tau$ is a function of type (C) on $\frak a_P$. Indeed,
if $$f\Big(namk\Big)=\sum_{j=1}^l\phi_j(mk)\cdot\alpha_j\Big(H(a)\Big)
e^{\langle\lambda_j+
\rho_P,H(a)\rangle}$$ for $n\in N(\mathbb A),\,a\in A_P,\,
m\in M(\mathbb A)^1$ and $k\in K$, where for all $j$, $\alpha_j(X)$ 
is a polynomial, and $\phi_j(g)$ is an automorphic form in 
$\mathcal A_P(G)$ such that $\phi_j(ag)=\phi_j(g)$ for $a\in A_P$,
then we can define 
$$\int_{P(F)\backslash G(\mathbb A)^1}^*f(g)\cdot
\tau_P\Big(H(g)-T\Big)\,dg$$ by
$$\sum_{j=1}^l\int_K\int_{M(F)\backslash M(\mathbb A)^1}^*
\phi_j(mk)\,dm\,dk\Bigg(\int_{\frak a_P}^{\#}\alpha_j(X)
e^{\langle\lambda_j,X\rangle}\tau_P(X-T)\,dX\bigg).$$
This is well-defined provided that the following two conditions are satisfied:

\noindent
(i$^*$) $\langle\mu,\varpi^\vee\rangle\not=0 \qquad \forall
Q\subset P,\ \varpi^\vee\in ({\widehat \Delta}^\vee)_Q^P,\ 
\mu\in\mathcal E_Q(\phi)$; and

\noindent
(ii$^*$) $\langle\lambda,\alpha^\vee\rangle\not=0 \qquad \forall
 \alpha\in \Delta_P,\ \lambda\in\mathcal 
E_P(\phi).$

Let $\mathcal A(G)^{**}$ be the space of $\phi\in\mathcal A(G)$ such 
that (i$^*$) and hence also (ii$^*$) are satisfied for any $P$. As a direct 
consequence of the Theorem in the previous subsection, we have the following
\vskip 0.30cm
\noindent
{\bf {\Large Basic Fact.}} ([JLR]) {\it For $\phi\in\mathcal A(G)^{**}$, 
$$\int_{G(F)\backslash 
G(\mathbb A)^1)}\Lambda^T\phi(g)\,dg$$ is equal
to $$\sum_{P}(-1)^{d(P)-d(G)}\int_{P(F)\backslash 
G(\mathbb A)^1}^* \phi_P(g)\cdot\widehat\tau_P\Big(H(g)-T\Big)\,dg,$$
where $P$ runs through all standard parabolic subgroups.}

\noindent
{\bf Remark.} By definition, 
$$\Lambda^T\phi(g)=\sum_{P:P\subset G}(-1)^{d(P)-d(G)}\phi_P(g)
\cdot\widehat\tau_P(H(g)-T).$$ Thus formally,
$$\begin{aligned}&\int_{G(F)\backslash 
G(\mathbb A)^1)}\Lambda^T\phi(g)\,dg\\
=&\int_{G(F)\backslash 
G(\mathbb A)^1)}\sum_{P}(-1)^{d(P)-d(G)}\sum_{\delta\in P(F)
\backslash G(F)}\phi_P(\delta g)\cdot\widehat
\tau_P\Big(H(\delta g)-T\Big)\,dg\\
=&\sum_{P}(-1)^{d(P)-d(G)}\int_{G(F)\backslash 
G(\mathbb A)^1)}\sum_{\delta\in P(F)\backslash G(F)}\phi_P(\delta g)
\cdot\widehat\tau_P\Big(H(\delta g)-T\Big)\,dg\\
=&\sum_{P}(-1)^{d(P)-d(G)}\int_{P(F)\backslash 
G(\mathbb A)^1)}\phi_P(g)\cdot\widehat\tau_P(H(g)-T)\,dg.\end{aligned}$$
In this sense, the present theorem serves the sole purpose of
giving a justification of this formal calculation.

\noindent
{\it Proof}. ([JLR]) 
Step 1. By induction on the rank, we may assume that the theorem
 holds for the Levi subgroup $M$ of a {\it proper} parabolic subgroup 
$P$ of $G$. We will show below the
\vskip 0.30cm
\noindent
{\bf Claim 1:} {\it This induction hypothesis implies that, for any proper 
$P$, $$\int_{P(F)\backslash  
G(\mathbb A)^1}^{\#} \Lambda^{T,P}\phi_P(g)\cdot
\tau_P\Big(H(g)-T\Big)\,dg\eqno(*)$$ is
 equal to $$\sum_{R:R\subset P}(-1)^{d(R)-d(P)}\int_{R(F)
\backslash G(\mathbb A)^1}^* \phi_R(g)\cdot
\widehat\tau_R^P\Big(H(g)-T\Big)\cdot\tau_P\Big(H(g)-T\Big)\,dg.\eqno(**)$$} 

Assuming this, we may sum over $P$ to write
$$\int_{G(F)\backslash G(\mathbb A)^1}^*\phi(g)\,dg-
\int_{G(F)\backslash G(\mathbb A)^1}\Lambda^T\phi(g)\,dg$$ as
 
$$\sum_{P\not=G, P\supset R}(-1)^{d(R)-d(P)}\int_{R(F)\backslash 
G(\mathbb A)^1}^* \phi_R(g)\cdot \widehat\tau_R^P(H(g)-T)\cdot
\tau_P(H(g)-T)\,dg\eqno(*3)$$
\vskip 0.30cm
\noindent
{\bf  Claim 2.} {\it The summation can be taken inside the integral
on the RHS of (*3).}

If so, then
$$\begin{aligned}&\int_{G(F)\backslash G(\mathbb A)^1}^*\phi(g)\,dg-
\int_{G(F)\backslash G(\mathbb A)^1}\Lambda^T\phi(g)\,dg\\
=&\sum_R \int_{R(F)\backslash 
G(\mathbb A)^1}^* \phi_R(g)\sum_{P:P\not=G, P\supset R}
(-1)^{d(R)-d(P)+d(G)-d(G)}\widehat\tau_R^P\Big(H(g)-T\Big)\cdot
\tau_P\Big(H(g)-T\Big)\,dg\\
=&\sum_R \int_{R(F)\backslash 
G(\mathbb A)^1}^*(-1)^{d(R)-d(G)} \phi_R(g)\sum_{P:G\not=P\supset R}
(-1)^{d(G)
-d(P)}\widehat\tau_R^P\Big(H(g)-T\Big)\cdot\tau_P\Big(H(g)-T\Big)\,dg\\
=&\sum_R (-1)^{d(R)-d(G)}\int_{R(F)\backslash 
G(\mathbb A)^1}^* \phi_R(g)\cdot(-1)(-1)^{d(G)
-d(G)}\widehat\tau_R^G(H(g)-T)\tau_G(H(g)-T)\,dg,\end{aligned}$$ since
 for $R\not=G$, Langlands' Combinatorial Lemma gives
$$\sum_{P:R\subset P\not=G}(-1)^{d(R)-d(P)}\widehat\tau_R^P\Big(H(g)-T\Big)
\cdot\tau_P\Big(H(g)-T\Big)
=-(-1)^{d(R)-d(G)}\widehat\tau_R\Big(H(g)-T\Big).$$
Thus, $$\begin{aligned}&\int_{G(F)\backslash 
G(\mathbb A)^1)}\Lambda^T\phi(g)\,dg\\
=&\int_{G(F)\backslash G(\mathbb A)^1)}^*\phi(g)\,dg
+\sum_{R:R\not=G}(-1)^{d(R)-d(G)}
\int_{R(F)\backslash G(\mathbb A)^1)}^*\phi_R(g)\widehat\tau_R(H(g)-T)\,dg\\
=&\sum_{P:P\not=G}(-1)^{d(P)-d(G)}
\int_{P(F)\backslash G(\mathbb A)^1)}^*\phi_P(g)\widehat\tau_P(H(g)-T)\,dg,
\end{aligned}$$ as required.

\noindent
Step 2.  {\it Proof of Claim 1.} Write the constant term $\phi_P$ 
as a sum $$\phi_P(namk)=\sum_j Q_j\Big(H(a)\Big)\psi_j(amk)$$ for 
$n\in N(\mathbb A), a\in A_P, m\in M(\mathbb A)^1$ and $k\in K$, 
where $Q_j$ are polynomials and $\psi_j\in \mathcal A_P(G)$ satisfies
$$\psi_j(ag)=e^{\langle \lambda_j+\rho_P,H(a)\rangle}\psi_j(g)$$ for 
some exponents $\lambda_j\in\frak a_P^*$ for all $a\in A_P$. Then (*) 
is equal to
$$\sum_j\Big(\int_K\int_{M(F)\backslash M(\mathbb A)^1}\Lambda^{T,M}
\psi_j(mk)dm\,dk\Big)\Big(\int_{\frak a_P}^{\#}
Q_j(X)e^{\langle \lambda_j,X\rangle}\tau_P(X-T)\,dX\Big).$$ Using our
 induction hypothesis, we may further write it as the sum over $j$ and 
$R\subset P$ of $(-1)^{d(R)-d(P)}$ times
$$\int_K\int_{R_M(F)\backslash M(\mathbb A)^1}^*(\psi_j)_{M_R}(mk)
\cdot\widehat\tau_R^P\Big(H(m)-T\Big)dm\,dk\eqno(*4)$$ times
$$\int_{\frak a_P}^{\#}
Q_j(X)e^{\langle \lambda_j,X\rangle}\tau_P(X-T)\,dX,$$ where 
$R_M=R\cap M$. Choose an analogous decomposition  for the constant 
term $(\psi_j)_{R_M}$:
$$\Big(\psi_j\Big)_{R_M}(namk)=\sum_l P_{j,l}\Big(H(a)\Big)\cdot
\psi_{j,l}\Big(amk\Big)$$ for 
$n\in \Big(N_R\cap M\Big)(\mathbb A), 
a\in A_R\cap M(\mathbb A)^1, m\in M_R(\mathbb A)^1$ and $k\in K$, 
where $P_{j,l}$ are polynomials and $\psi_{j,l}$  
$$\psi_{j,l}(ag)=e^{\langle \lambda_{j,l}+\rho_R^P,H(a)\rangle}
\psi_{j,l}(g)$$ for some exponent $\lambda_{j,l}\in
(\frak a_R^P)^*$.  Then we may write (*4) as a sum over $l$ of
$$\int_K\int_{M_R(F)\backslash M_R(\mathbb A)^1}^*
(\psi_{j,l})(mk)dm\,
dk\int_{\frak a_R^P}^{\#}
P_{j,l}(X)e^{\langle \lambda_{j,l},X\rangle}\widehat\tau_R^P(X-T)\,dX.$$
By basic property (5) for integration over cones, 
we may combine the \#-integrals over $\frak a_R^P$ and 
$\frak a_P$ into a single \#-integral over $\frak a_R$ and we see that
$$\bigg(\int_K\int_{R_M(F)\backslash M(\mathbb A)^1}^*(\psi_j)_{M_R}(mk)
\cdot\widehat\tau_R^P\Big(H(m)-T\Big)dm\,dk\bigg)
\bigg(\int_{\frak a_P}^{\#}
Q_j(X)e^{\langle \lambda_j,X\rangle}\tau_P(X-T)\,dX\bigg)$$ is equal
 to the sum over $l$ of
$$\int_K\int_{M_R(F)\backslash M_R(\mathbb A)^1}^*\psi_{j,l}(mk)dm\,dk
\int_{\frak a_R}^{\#}
P_{j,l}(X)Q_j(X)e^{\langle \mu_{j,l},X\rangle}\widehat\tau_R^P\Big(H(g)-T\Big)
\tau_P\Big(X-T\Big)\,dX,$$
where $\mu_{j,l}=\lambda_j+\lambda_{j,l}$ and this equals
$$\int_{R(F)\backslash G(\mathbb A)^1}^*Q_j\Big(H_P(g)\Big)(\psi_j)_{R_M}(g)
\cdot\widehat\tau_R^P\Big(H(g)-T\Big)\cdot\tau_P\Big(H(g)-T\Big)\,dg.$$
Summing over $j$ gives
$$\int_{R(F)\backslash G(\mathbb A)^1}^*\phi_R(g)
\cdot\widehat\tau_R^P\Big(H(g)-T\Big)\tau_P\Big(H(g)-T\Big)\,dg.$$
This shows that (*) is equal to
$$\sum_{R\subset P}(-1)^{d(R)-d(P)}\int_{R(F)\backslash G(\mathbb A)^1}^*
\phi_R(g)\widehat\tau_R^P(H(g)-T)\tau_P(H(g)-T)\,dg,$$ as required.
\vskip 0.30cm
\noindent
Step 3. {\it Proof of Claim II.} This is an easy consequence of 
Basic Property (6) of integration over cones:
Consider the three cones $$\widehat{\mathcal C}_R:=\{X\in\frak a_R:
\widehat\tau_R(X)=1\},
\mathcal C_P:=\{X\in\frak a_P:\tau_P(X)=1\}, \widehat{\mathcal C}_R^P:=
\{X\in\frak a_R^P:\widehat\tau_R^P(X)=1\}.$$

\noindent
(a) The product $\mathcal C_P\times \widehat{\mathcal C}_R^P$ is contained in
$\widehat{\mathcal C}_R$. Indeed, $\widehat{\mathcal C}_R$ is the positive span
 of the coroots $\{\alpha^\vee:\alpha\in \Delta_R\}$, 
$\mathcal C_P$ is the positive span of the coweights in 
$\widehat\Delta_P^\vee$ and $\widehat {\mathcal C}_R^P$ is the positive 
span of the coroots $\{\alpha^\vee:\alpha\in\Delta_R^P\}$, so the 
assertion follows from the fact that all coweights in 
$\widehat\Delta_P^\vee$ are non-negative linear combinations of 
coroots in $\Delta_P^\vee$.

\noindent
(b) Now let $\lambda\in \mathcal E_R(\phi)$ and for $P$ containing $R$, 
write $\lambda=\lambda_R^P+\lambda_P$ relative
to the decomposition $\frak a_R=\frak a_R^P\oplus \frak a_R.$ By our
 hypothesis, $\langle\lambda_P,\varpi^\vee\rangle\not=0$ for all 
$\varpi^\vee\in\Delta_R^P$, and hence $\lambda_P$ is non-negative 
with respect to $\widehat{\mathcal C}_R^P$. Since 
$\mathcal C_P\times \widehat {\mathcal C}_R^P\subset \widehat {\mathcal C}_R$ 
for all $P$, we may apply Basic Property (6) for integration over cones
to conclude that the summation can be taken inside the integral
on the RHS of (*3). This completes the proof of Claim 2 and hence the theorem
as well.

\subsection{Bernstein Principle}

As such, with our basic assumption on modified analytic truncations, 
similarly, we have the following
\vskip 0.30cm
\noindent
{\bf {\large Bernstein's Principle}}$'$. ([JLR]) {\it Let $P=MN$ be a proper 
parabolic subgroup and let $\sigma$ be an irredcible cuspidal representation 
in $L^2\Big(M(F)\backslash M(\mathbb A)^1\Big)$. Let $E(g,\phi,\lambda)$ 
be an Eisenstein series associated to $\phi\in \mathcal A_P(G)_\sigma$. Then
$$I_G\Big(E(g,\phi,\lambda)\Big)=0$$ for all $\lambda$ such that
$E(g,\phi,\lambda)$ and $I_G\Big(E(g,\phi,\lambda)\Big)$ are defined.}

\noindent
{\it Proof.}  This is because 

\noindent
(a) the map $\phi\mapsto I_G\Big(E(g,\phi,\lambda)\Big)$ defines a 
$G(\mathbb A_f)^1$-invariant functional on 
$\mathrm{Ind}_P^G(\sigma\otimes e^\lambda)$ for nice $\lambda$ by JLR's result
above; and

\noindent
(b) there does not exist any such invariant functional for generic values of 
$\lambda$, by a result of Bernstein [B]. 
\vskip 0.30cm
We believe that this principle is the main reason behind special yet very 
important relations among multiple and ordinary zeta functions in literature.
\subsection{Eisenstein Periods}
 
As a direct consequence of Bernstein principle, with the explicit formula 
for constant terms, we have the following close formula for Eisenstein periods
associated with cusp forms:
\vskip 0.30cm
\noindent
{\bf {\large Basic Fact.}} ([JLR])  
{\it Let $P=MN$ be a parabolic sungroup and 
let $\sigma$ be an irredcible cuspidal representation in 
$L^2\Big(M(F)\backslash M(\mathbb A)^1\Big)$. Let $E(g,\phi,\lambda)$ 
be an Eisenstein series associated to $\phi\in \mathcal A_P(G)_\sigma$. Then
$$\begin{aligned}
&\int_{G(F)\backslash G(\mathbb A)}\Lambda^TE(g,\phi,\lambda)\,dg
=\mathrm{Vol}\Big(\Big\{\sum_{\alpha\in\Delta_0}
a_\alpha\alpha^\vee:a_\alpha\in[0,1)\Big\}\Big)\\
&\hskip 4.0cm\cdot\sum_{w\in\Omega}\frac{e^{\langle w\lambda-\rho,T\rangle}}
{\prod_{\alpha\in\Delta_0}\langle w\lambda-\rho,\alpha^\vee\rangle}
\cdot\int_{M(F)\backslash M(\mathbb A)^1\times K}\Big(M(w,\lambda)
\phi\Big)(mk)\,dm\,dk.\end{aligned}$$}

\noindent
{\it Proof.} (a) For $\lambda\in\frak a_P^*$ generic, by the
Basic Facts in the previous subsection,
$$\int_{G(F)\backslash G(\mathbb A)^1}\Lambda_TE(g,\phi,\lambda)=
\sum_Q(-1)^{d(Q)-d(G)}\int_{Q(F)\backslash G(\mathbb A)^1}^*
E_Q(g,\phi,\lambda)\cdot\widehat\tau_Q\Big(H(g)-T\Big)\,dg.$$

\noindent
(b) Since $\phi$ is cuspidal, $E_Q(g,\phi,\lambda)$ vanishes identically 
unless $Q$ contains as association of $P$.

\noindent
(c) As a function of $m$, $E_Q(namk,\phi,\lambda)$ is a sum of Eisenstein 
series induced from parabolic subgroups associated to $P$. Therefore, 
by Bernstein's principal and the definition,
$$\int_{Q(F)\backslash G(\mathbb A)^1}^*E_Q(g,\phi,\lambda)\cdot\widehat\tau_Q
\Big(H(g)-T\Big)\,dg$$ vanishes identically unless $Q$ is associated to $P$.

\noindent
(d) If $Q$ is associated to $P$, then $E_Q$ is a sum of cusp forms on 
$M_Q(\mathbb A)$. In this case, $Q=P$ and the only non-zero term is 
$$(-1)^{d(P)-d(G)}\sum_{w\in\Omega}\Big(\int_{\frak a_P}^*
e^{\langle w\lambda-\rho_P\rangle}\cdot\tau_P(H-T)\,dH\Big)\cdot
\int_{M(F)\backslash M(\mathbb A)^1\times K}M(w,\lambda)\phi(mk)\,dm\,dk.$$ 

\noindent
(e) Finally, by Basic Property (7) for the integration over cones, we see
$$\begin{aligned}
&(-1)^{d(P)-d(G)}\int_{\frak a_P}^{\#}e^{\langle w\lambda-\rho_P,H\rangle}
\widehat \tau_P(H-T)\,dH\\
=&\mathrm{Vol}\Big(\Big\{\sum_{\alpha\in\Delta_P}
a_\alpha\alpha^\vee:a_\alpha\in[0,1)\Big\}\Big)
\cdot\frac{e^{\langle w\lambda-\rho_P,T\rangle}}
{\prod_{\alpha\in\Delta_P}\langle w\lambda-\rho_P,\alpha^\vee\rangle}.
\end{aligned}$$
This completes the proof.
\vskip 0.30cm
This final close formula in turn may be used to understand Bernstein
principle. Indeed, by the inversion formula,
$$E(g,\phi,\lambda)=\sum_{Q}\sum_{\delta\in Q(F)\backslash G(F)}
\Lambda^{T,Q}E(\delta g,\phi,\lambda)\cdot\tau_Q\Big(H(\delta g)-T\Big).$$ 
Therefore,
$$\begin{aligned}
0=&\int_{G(F)\backslash G(\mathbb A)^1}^*E(g,\phi,\lambda)\,dg\\
=&
\sum_Q\int_{Q(F)\backslash G(\mathbb A)^1}^*\Lambda^{T,Q}
E(g,\phi,\lambda)\cdot \tau_Q\Big(H(g)-T\Big)\,dg.\end{aligned}$$ 
Thus by definition,
$$\begin{aligned}
&-\int_{G(F)\backslash G(\mathbb A)^1}^*\Lambda^TE(g,\phi,\lambda)\,dg\\
=&\sum_{Q:Q\subsetneqq G}\int_{Q(F)\backslash G(\mathbb A)^1}^*
\Lambda^{T,Q}E(g,\phi,\lambda)\cdot \tau_Q\Big(H(g)-T\Big)\,dg.\end{aligned}$$
That is to say,
$$\begin{aligned}
&\sum_{Q:Q\subsetneqq G}\int_{Q(F)\backslash G(\mathbb A)^1}^*\Lambda^{T,Q}
E(g,\phi,\lambda)\cdot \tau_Q\Big(H(g)-T\Big)\,dg\\
=&-\mathrm{Vol}\Big(\Big\{\sum_{\alpha\in\Delta_0}
a_\alpha\alpha^\vee:a_\alpha\in[0,1)\Big\}\Big)\\
&\qquad\cdot\sum_{w\in\Omega}\frac{e^{\langle w\lambda-\rho,T\rangle}}
{\prod_{\alpha\in\Delta_0}\langle w\lambda-\rho,\alpha^\vee\rangle}
\cdot\int_{M(F)\backslash M(\mathbb A)^1\times K}\Big(M(w,\lambda)
\phi\Big)(mk)\,dm\,dk\end{aligned}\eqno(*)$$
It is our understanding that this formula is the main reason behind 
special yet very important relations among multiple and ordinary zeta 
functions in literature.

\section{Abelian Part of Non-Abelian L-Functions: Open Problems}

\subsection{Cone Corresponding to ${\bold 1}\Big(p_P^g>_Pp\Big)$}
In this quite elementary section, we show that in fact the subset in 
$\frak a_P$ defined by
the characteristic function ${\bold 1}\Big(p_P^g>_Pp\Big)$ is also a cone, 
in the sense of JLR used in defining  regularized integration of cones. 
Consequently, JLR's
regularized integration over cones applies to our function 
${\bold 1}\Big(p_P^g>_Pp\Big)$ as well.
 
To start with, as before, fix a standard parabolic subgroup 
$P:P\supset P_0$ corresponding to the partition $r=d_1+d_2+\cdots+d_P$. 
Set $$H_P(g):=H=\Big(H_1^{(d_1)}, H_2^{(d_2)},\cdots, 
H_{|P|}^{(d_{|P|})}\Big)\in\frak a_P,$$ 
where $d_j=r_j-r_{j-1}$ for $j=1,\cdots, |P|$ and $r_0=0$.
By definition,
$$p_P^g>_Pp\Leftrightarrow p_P^g(r_i)>p(r_i)\qquad i=1,2,\cdots,|P|-1.$$
Or equivalently,
$$d_1H_1+d_2H_2+\cdots+d_iH_i>p(r_i).$$ 

Clearly,
$$\begin{aligned}&d_1H_1+d_2H_2+\cdots+d_iH_i\\
=&r_1H_1+(r_2-r_1)H_2+\cdots+(r_i-r_{i-1})H_i\\
=&r_1(H_1-H_2)+r_2(H_2-H_3)+\cdots+r_{i-1}(H_{i-1}-H_i)+r_iH_i\\
=&r_1\alpha_{r_1}(H)+r_2\alpha_{r_2}(H)+\cdots+r_{i-1}\alpha_{i-1}(H)
+r_iH_i\\
=&\Big(r_1\alpha_{r_1}+r_2\alpha_{r_2}+\cdots+r_{i-1}\alpha_{i-1}\Big)(H)
+r_iH_i.\end{aligned}$$ Here as before, we use 
$$\Delta:=\Big\{\alpha_i=e_i-e_{i+1}:i=1,2,\cdots,|P|-1\Big\}.$$
Next, we want to write $r_iH_i$ as a linear combination of $\alpha_{r_j}$ 
acting on $H$. Thus set
$$r_iH_i=:\sum a_{ji}\alpha_{r_j}(H).$$ 
Then $$
\begin{aligned}r_iH_i=&\sum_{j=1}^{|P|-1}a_{ji}(H_j-H_{j+1})
=\sum_{j=1}^{|P|-1}a_{ji}H_j-\sum_{j=2}^{|P|}a_{j-1\,i}H_j\\
=&a_{1i}H_1+\sum_{j=2}^{|P|-1}\Big(a_{ji}-a_{j-1\,i}\Big)H_j
-a_{|P|-1\,i}H_{|P|}.\end{aligned}$$ 
But $$d_1H_1+d_2H_2+\cdots+d_{|P|}H_{|P|}=0,$$ and $d_1=r_1$,
so $$-H_{|P|}=\frac{1}{d_{|P|}}\Big(\sum_{j=2}^{|P|-1}d_j H_j+r_1H_1\Big).$$
Thus we have
$$r_iH_i=a_{1i}H_1+\sum_{j=2}^{|P|-1}(a_{ji}-a_{j-1\,i})H_j+
\frac{1}{d_{|P|}}\Big(\sum_{j=2}^{|P|-1}d_jH_j+r_1H_1\Big)\cdot a_{|P|-1\,i}.$$
In other words,
$$r_iH_i=\Big(a_{1i}+\frac{r_1}{d_{|P|}}a_{|P|-1\,i}\Big)H_1
+\sum_{j=2}^{|P|-1}\Big(a_{ji}-a_{j-1\,i}+\frac{d_j}{d_{|P|}}\cdot 
a_{|P|-1\,i}\Big)H_j.$$

As an example, take $i=1$. Then we have
$$\begin{aligned}a_{11}+&\frac{r_1}{d_{|P|}}a_{|P|-1\,1}=r_1,\\
a_{j1}-a_{j-1\,1}+&\frac{d_j}{d_{|P|}}\cdot a_{|P|-1\,1}=0,\qquad j=2,3,
\cdots,|P|-1.\end{aligned}$$
That is, we have
$$\begin{aligned}
a_{|P|-1\,1}-a_{|P|-2\,1}+&\frac{r_{|P|-1}-r_{|P|-2}}{d_{|P|}}\cdot 
a_{|P|-1\,1}=0,\\
a_{|P|-2\,1}-a_{|P|-3\,1}+&\frac{r_{|P|-2}-r_{|P|-3}}{d_{|P|}}\cdot 
a_{|P|-1\,1}=0,\\
\cdots\cdots&\cdots\\
a_{21}-a_{11}+&\frac{r_2-r_1}{d_{|P|}}\cdot a_{|P|-1\,1}=0,\\
a_{11}+&\frac{r_1}{d_{|P|}}a_{|P|-1\,1}=r_1.
\end{aligned}$$ 
Therefore, $$\begin{aligned}a_{|P|-1\,1}+
&\frac{r_{|P|-1}}{d_{|P|}}a_{|P|-1\,1}=r_1,\\
a_{|P|-1\,1}-a_{|P|-j\,1}+&\frac{r_{|P|-1}-r_{|P|-j}}{d_{|P|}}\cdot 
a_{|P|-1\,1}=0, \qquad j=2,3,\cdots, |P|-1.\end{aligned}$$
That is to say,
$$\begin{aligned}
a_{|P|-1\,1}=&\frac{r_1\cdot (r_{|P|}-r_{|P|-1})}{r_{|P|-1}},\\
a_{|P|-j\,1}=&\frac{r_{|P|}-r_{|P|-j}}{r_{|P|}-r_{|P|-1}}\cdot 
a_{|P|-1\,1},\qquad j=2,3,\cdots,|P|-1.\end{aligned}$$
Or better,
$$a_{|P|-j\,1}=\frac{r_1\cdot (r_{|P|}-r_{|P|-j})}{r_{|P|-1}},
\qquad j=1,2,\cdots,|P|-1.$$

More generally, we see that the conditions satisfied by $a_{ji}$ 
can be written as
$$R_i
=\begin{pmatrix} 1&0&0&\cdots&0&\frac{r_1}{d_{|P|}}\\
-1&1&0&\cdots&0&\frac{r_2-r_1}{d_{|P|}}\\
0&-1&1&\cdots&0&\frac{r_3-r_2}{d_{|P|}}\\
\cdots&\cdots&\cdots&\cdots&\cdots&\cdots\\
0&0&0&\cdots&1&\frac{r_{|P|-2}-r_{|P|-3}}{d_{|P|}}\\
0&0&0&\cdots&-1&1+\frac{r_{|P|-1}-r_{|P|-2}}{d_{|P|}}\end{pmatrix}\cdot A_i.$$
where $R_i=\Big(0,0,\cdots,r_i,\cdots,0\Big)^t$ with $r_i$ in the 
$i$-th position, and 
$A_i=\Big(a_{1i},a_{2i}\cdots, a_{|P|-1\,i}\Big)^t$.
By an easy calculation, 
the inverse matrix of $$\begin{pmatrix} 1&0&0&\cdots&0&\frac{r_1}{d_{|P|}}\\
-1&1&0&\cdots&0&\frac{r_2-r_1}{d_{|P|}}\\
0&-1&1&\cdots&0&\frac{r_3-r_2}{d_{|P|}}\\
\cdots&\cdots&\cdots&\cdots&\cdots&\cdots\\
0&0&0&\cdots&1&\frac{r_{|P|-2}-r_{|P|-3}}{d_{|P|}}\\
0&0&0&\cdots&-1&1+\frac{r_{|P|-1}-r_{|P|-2}}{d_{|P|}}\end{pmatrix}$$ 
is given by
$$\begin{pmatrix} 1-\frac{r_1}{r}&-\frac{r_1}{r}&-\frac{r_1}{r}&
\cdots&-\frac{r_1}{r}&-\frac{r_1}{r}\\
1-\frac{r_2}{r}&1-\frac{r_2}{r}&-\frac{r_2}{r}&\cdots&-\frac{r_2}{r}&
-\frac{r_2}{r}\\
1-\frac{r_3}{r}&1-\frac{r_3}{r}&1-\frac{r_3}{r}&\cdots&-\frac{r_3}{r}&
-\frac{r_3}{r}\\
\cdots&\cdots&\cdots&\cdots&\cdots&\cdots\\
1-\frac{r_{|P|-2}}{r}&1-\frac{r_{|P|-2}}{r}&1-\frac{r_{|P|-2}}{r}&
\cdots&1-\frac{r_{|P|-2}}{r}&-\frac{r_{|P|-2}}{r}\\
\frac{r-r_{|P|-1}}{r}&\frac{r-r_{|P|-1}}{r}&\frac{r-r_{|P|-1}}{r}&
\cdots&\frac{r-r_{|P|-1}}{r}&
\frac{r-r_{|P|-1}}{r}\end{pmatrix}.$$
As such, since $R_1,R_2,\cdots, R_{|P|-1}$ are linearly independent, so 
the linear forms $$\sum_{j=1}^{|P|-1}\Big(r_j+a_{ji}\Big)\alpha_{r_j},
\qquad i=1,\cdots,|P|-1$$
are linearly independent. Consequently, the open subset defined by
the linear inequalities 
$$p_P^g=d_1H_1+d_2H_2+\cdots+d_iH_i
=\bigg(\sum_{j=1}^{|P|-1}\Big(r_j+a_{ji}\Big)\alpha_{r_j}\bigg)(H)>p(r_i),
\qquad i=1,\cdots,|P|-1$$
is a cone in the sense of JLR. That is to say, we have shown the following

\noindent
{\bf Lemma.} With the same notation as above, the characteristic function 
${\bold 1}\Big(p_P^g>_Pp\Big)$
defines a cone in $\frak a_P$.

Consequently, all the discussion of JLR as recalled in \S6.2
about regularized integrations over cones works for 
${\bold 1}\Big(p_P^g>_Pp\Big)$ as well.

\subsection{Open Problems}
Motivated by our discussion on Arthur's analytic truncation and 
Jacquet-Lapid-Rogawski's
investigations on Arthur's period, for the modified analytic truncation 
$\Lambda_T$ using positive chambers 
and the geometrically oriented truncation $\Lambda_p$ using stability, 
we apply the regularization process of the previous
subsection to give the following discussion.

First recall that by definition, for an automorphic form on 
$G(F)\backslash G(\mathbb A)$,
$$\Lambda_T\phi(g):=\sum_P(-1)^{d(P)-d(G)}
\sum_{\delta\in P(F)\backslash G(F)}\phi_P(\delta g)
\cdot\tau_P\Big(H_0(\delta g)-T\Big).$$ By the fact that $\phi_P$ 
is again of the type we assume in the definition of regularized 
integration, and that $\tau_P$ is a characteristic function of 
positice chamber, we have  well-defined regularized integrations:
$$\int_{P(F)\backslash G(\mathbb A)^1}^*\phi_P(g)
\cdot \tau_P\Big(H_0(\delta g)-T\Big)\,dg.$$  Accordingly, 
it seems to be quite reasonable for us to introduce the following
regularized integration
$$\int_{G(F)\backslash G(\mathbb A)^1}^*\Lambda_T\phi(g)\,dg:=
\sum_P(-1)^{d(P)-d(G)}\int_{P(F)\backslash G(F)}^*\phi_P(g)
\cdot\tau_P\Big(H_0(g)-T\Big)\,dg$$ for any $G$-level automorphic form $\phi$, 
and any suitably regular element $-T\in\frak a_0^+$. As such, it appears that 
the relation 
$$\int_{G(F)\backslash G(\mathbb A)^1}^*\phi(g)\,dg
=\sum_P\int_{P(F)\backslash G(\mathbb A)^1}^*
\Lambda_{T,P}\phi(g)\cdot\widehat
\tau_P\Big(H_0(g)-T\Big)$$ also can be easily deduced from an 
induction on the rank since we have the inversion formula:
$$\phi(g)=\sum_P\sum_{\delta\in P(F)\backslash G(F)}
\Lambda_{T,P}\phi(\delta g)\cdot\widehat
\tau_P\Big(H_0(\delta g)-T\Big).$$

However, this process can be hardly justified, due to the fact that
ususally $\Lambda_{T,P}\phi$ is not of repidly decreasing, or better,
due to the fact that the sets defined by $\Lambda_{T,P}{\bold 1}$ are 
not compact -- Consequently,
the inner summation $$\sum_{\delta\in P(F)\backslash G(F)}\phi_P(\delta g)
\cdot\tau_P\Big(H_0(\delta g)-T\Big)$$ in \lq defining' 
$\int_{G(F)\backslash G(\mathbb A)^1}^*\phi(g)\,dg$ normally is an 
infinite one.

\vskip 0.30cm
Next, let us trun back to our main interests, the geometrically 
oriented truncation $\Lambda_p$. Then we expect the story to be 
quite different (from that for $\Lambda_{T,P}$), with our experiences
from Arthur's analytic truncations and periods. Recall that, by definition,
$$\Lambda_p\phi(g):=\sum_P(-1)^{d(P)-d(G)}
\sum_{\delta\in P(F)\backslash G(F)}\phi_P(\delta g)\cdot{\bold 1}
\Big(p_P^{\delta g}>_Pp\Big).$$ Similarly,  by the fact that $\phi_P$ 
is again of the type we assume in the definition of regularized 
integration, and that ${\bold 1}\Big(p_P^*>_Pp\Big)$ is a characteristic 
function of a certain cone as shown in the previous subsection, we have then  
well-defined regularized integrations 
$$\int_{P(F)\backslash G(\mathbb A)^1}^*\phi_P(g)\cdot {\bold 1}
\Big(p_P^{g}>_Pp\Big)\,dg.$$ 
Accordingly, we introduce the following

\noindent
{\bf Definition.} For any $G$-level automorphic form $\phi$, 
and any normalized convex polygon $p:[0,r]\to \mathbb R$, set
$$\int_{G(F)\backslash G(\mathbb A)^1}^*\Lambda_p\phi(g)\,dg:=
\sum_P(-1)^{d(P)-d(G)}\int_{P(F)\backslash G(F)}^*\phi_P(g)\cdot{\bold 1}
\Big(p_P^{g}>_Pp\Big)\,dg.$$

Motivated by the above discussion, concerning the truncation
$\Lambda_p\phi(g)\,dg$ and the regularized integration
 $\int_{G(F)\backslash G(\mathbb A)^1}^*
\Lambda_p\phi(g)\,dg$, we propose the following

\noindent
{\bf Open Problems.}  For any $G$-level automorphic form $\phi$, and any 
normalized convex polygon $p:[0,r]\to \mathbb R$,

\noindent
(1) $\Lambda_p\phi(g)$ is integrable. Hence the ordinary integration 
$\int_{G(F)\backslash G(\mathbb A)^1}\Lambda_p\phi(g)\,dg$ makes sense; 
Moreover, 

\noindent
(2) $\displaystyle{\int_{G(F)\backslash G(\mathbb A)^1}\Lambda_p\phi(g)\,dg
=\int_{G(F)\backslash G(\mathbb A)^1}^*\Lambda_p\phi(g)\,dg}.$
\vskip 0.30cm
As such, then the abelian parts of our non-abelian $L$-functions can be 
understood completely: With the help from
Bernstein's Principle, they are 
precisely given  by the formula stated in  Section 6.1, where it is 
obtained from a formal calculation.
\vskip 0.50cm
\noindent
{\bf {\large Acknowledgements:}} I would like to thank my wife, Shiao Ying and 
two kids, LerLer and AnAn for their patience, supports and understanding, 
and Ch. Deninger, I. Fesenko, I. Nakamura, T. Oda and my colleagues 
at KyuDai for their encouragement, helps and supports. Special thanks also 
due to Henry Kim for his keen interests, whose kind invitation to Toronto
is a primitive driving force for me to write down the details involved here,
and due to J. Arthur for helping me to understand analytic truncations. 

\noindent
This research is partially supported by JSPS.
\vfill
\eject
\vskip 0.5cm
\noindent
\centerline {\Large REFERENCES} 
\vskip 0.45cm
\noindent
[Ar1] Arthur, J.  Eisenstein series and the trace formula. {\it Automorphic 
forms, 
representations and $L$-functions}  pp. 253--274, Proc. Sympos. Pure Math., 
XXXIII, 
Amer. Math. Soc., Providence, R.I., 1979.
\vskip 0.20cm
\noindent
[Ar2] Arthur, J. A trace formula for reductive groups. I. Terms associated to 
classes 
in $G({\mathbb Q})$. Duke Math. J. 45  (1978), no. 4, 911--952
\vskip 0.20cm
\noindent
[Ar3] Arthur, J. A trace formula for reductive groups. II. Applications of a 
truncation 
operator. Compositio Math. 40 (1980), no. 1, 87--121.
\vskip 0.20cm
\noindent
[Ar4] Arthur, J. On the inner product of truncated Eisenstein series. Duke 
Math. J. 49 (1982), no. 1, 35--70.  
\vskip 0.20cm
\noindent
[Ar5] Arthur, J. A measure on the unipotent variety, Canadian J of Math, 
Vol. {\bf 37}, No. 6, pp.1237-1274, 1985
\vskip 0.20cm
\noindent
[Ar6] Arthur, J. {\it An Introduction to the Trace Formula}, to appear
\vskip 0.20cm
\noindent
[Bo1]  Borel, A. Some finiteness properties of adele groups over
number fields, Publ. Math., IHES, {\bf 16} (1963) 5-30
\vskip 0.20cm
\noindent
[Bo2] Borel, A. {\it Introduction aux groupes arithmetictiques}, Hermann, 1969
\vskip 0.20cm
\noindent
[Bu] Bump, D. The Rankin-Selberg method: a survey. {\it Number theory, 
trace formulas 
and discrete groups} (Oslo, 1987), 49--109, Academic Press, Boston, MA, 1989.
\vskip 0.20cm
\noindent 
[FT] Fr\"ohlich, A. \& Taylor, M.J. {\it Algebraic Number Theory}, Cambridge 
studies in advanced mathematics {\bf 27}, Cambridge Univ. Press (1991)
\vskip 0.20cm
\noindent
[H] Humphreys, J. {\it Introduction to Lie algebras and representation theory}, GTM {\bf 9}, 1972
\vskip 0.20cm
\noindent 
[J] Jacquet, H. On the residual spectrum of $GL(n)$, Springer LNM 1041,
185-208
\vskip 0.20cm
\noindent 
[JLR] Jacquet, H.; Lapid, E.\& Rogawski, J.  Periods of automorphic forms. 
J. Amer. Math. Soc. 12 (1999), no. 1, 173--240.
\vskip 0.20cm
\noindent 
[K] Kim, H. The residual spectrum of $G_2$, Can. J. Math., 48(6) (1996),
1245-1272
\vskip 0.20cm
\noindent 
[L1]  Lang, S. {\it Algebraic Number Theory}, 
Springer-Verlag, 1986
\vskip 0.20cm
\noindent 
[L2] Lang, S. {\it Fundamentals on Diophantine Geometry}, 
Springer-Verlag, 1983
\vskip 0.20cm
\noindent 
[L3] Lang, S. {\it Introduction to Arekelov Theory}, 
Springer-Verlag, 1988
\vskip 0.20cm
\noindent
[Laf] Lafforgue, L. {\it Chtoucas de Drinfeld et conjecture de 
Ramanujan-Petersson}. 
Asterisque No. 243 (1997)
\vskip 0.20cm
\noindent
[La1] Langlands, R. Eisenstein series. {\it Algebraic Groups and 
Discontinuous Subgroups}
 Proc. Sympos. Pure Math., IX, Amer. Math. Soc., Providence, R.I., 
1979,  pp. 235--252.
\vskip 0.20cm
\noindent
[La2] Langlands, R. {\it On the functional equations satisfied by 
Eisenstein series}, 
Springer LNM {\bf 544}, 1976
\vskip 0.20cm
\noindent
[MW1] Moeglin, C. \& Waldspurger, J.-L. Le spectre residuel de GL(n), 
Ann. Sci. Ec. Norm. Sup. {\bf 22}, (1989), 605-674
\vskip 0.20cm
\noindent
[MW2] Moeglin, C. \& Waldspurger, J.-L. {\it Spectral decomposition 
and Eisenstein series}. 
Cambridge Tracts in Mathematics, {\bf 113}. Cambridge University Press, 
Cambridge, 1995.
\vskip 0.20cm
\noindent 
[Neu] Neukirch, {\it Algebraic Number Theory}, Grundlehren der
Math. Wissenschaften, Vol. {\bf 322}, Springer-Verlag, 1999 
\vskip 0.20cm
\noindent 
[We1] Weng, L.  A Program for Geometric Arithmetic, not yet ready for a formal 
publication, but available at http://xxx.lanl.gov/abs/math.AG/0111241, 
\vskip 0.20cm
\noindent
[We2] Weng, L. Non-Abelian $L$-Functions for Function Fields, 
{\it Amer. J. Math},
to appear (see also http://xxx.lanl.gov/abs/math.NT/0412007)
\vskip 0.20cm
\noindent
[We3] Weng, L. Non-Abelian $L$-Functions for Number Fields, submitted, 
available at http://xxx.lanl.gov/abs/math.NT/0412008
\vskip 0.20cm
\noindent 
[We4] Weng, L. {\it Rank Two Non-Abelian Zeta and Its Zeros},  submitted, 
available at http://xxx.lanl.gov/abs/math.NT/0412009
\vskip 0.20cm
\noindent
[Z] Zagier, D.  The Rankin-Selberg method for automorphic functions 
which are not of rapid decay. 
J. Fac. Sci. Univ. Tokyo Sect. IA Math. 28 (1981), no. 3, 415--437 (1982).

\vskip 4.0cm
\noindent
{\bf {\large Lin WENG}}

\noindent
Graduate School of Mathematics, Kyushu University, Fukuoka 812-8581, Japan

\noindent
E-Mail Address: weng@math.kyushu-u.ac.jp
\enddocument
\vfill
\eject
\centerline {\large Appendix:}
\vskip 0.50cm
\centerline {\bf {\Large L Academy}} 
\vskip 0.30cm
\centerline {{\bf \large Lin ONG}}  
\centerline {\bf Institute of Fundamental Research, The L Academy}
\vskip 1.0cm
\noindent
Due to the unprecedented public demand, and on the behalf 
of the President Leonhard Euler of the $L$ Academy, I now publish this 
personal note on some of the historic events happened in the $L$ Academy. 
We hence sincerely hope that this helps to  bring the normal life back to the 
Academy.   
\vskip 0.50cm
\centerline {\bf {\large File 1. Epstein and Riemann}}
\vskip 0.30cm
\noindent
Epstein was once a junior faculty member of the Zeta Institute, one of the 
oldest institutes 
in the $L$ Academy. At that time, Epstein's position was not permanent. 
Rather, he held a 
temporary one, which in current terminology means a tenure-track position. 
So at some point, 
the $L$ Academy 
started a reviewing process on Epstein in order to decide whether to offer 
him a tenure with 
a promotion to the rank of Associate Professorship.

As usual, this evaluation was carried out by two committees, an internal one
and an external one, for the purpose of avoiding \lq possible ignorance, 
prejudice 
and abuse of powers of special interests groups'.

The internal evaluation process went extremely badly for Epstein. Riemann, 
the chairman of 
the committee, openly told everyone he met that Epstein should not be awarded 
with a tenure 
in the Zeta Institute because of several obvious reasons. In fact Riemann even 
went so far to
 say that everyone would reach such a conclusion by simply looking at this 
young 
man's \lq strange appearance\rq. Since Riemann, a distinguished member of 
the $L$ Academy, 
was also the Director of the Zeta Institute, practically speaking, 
Epstein's chance to be 
reminded at the Zeta Institute became very very slim.

However, the external committee, led by the CEO of Lufthansa, reached a 
totally different 
conclusion. They claimed that Epstein deserved a tenure position at the 
Zeta Institute: 
Not only this talent young man was no different from any others in the 
Zeta Institute, 
said the external committee, Epstein also admitted \lq Meromorphic 
Extension\rq$\ $ 
and satisfied \lq Functional Equation\rq$\ $, two advanced criterions 
for senior staffs 
in the Zeta Institute. This committee also pointed out that even though 
it was considered 
to be a top secrete, as
far as they understood, the \lq Meromorphic Extension\rq$\ $ and the 
\lq Functional 
Equation\rq$\ $ were the only qualifications for the Zeta Institute. And 
it was recorded 
that the external committee's evaluation of Epstein was SA, i.e., Super A.

As such, the $L$ Academy was in crisis: Never before in the history of the 
$L$ Academy, 
two committees had reached such totally opposite conclusions. Say, 
it was recorded that even in the case of Artin, when the $L$ Academy 
decided whether to 
offer him a tenure immediately after his Thesis, the internal committee's 
evaluation was
 $O^-$, i.e., not even the Ordinary, the bearly good enough grade to be 
admitted, while the 
external committee's evaluation was SA, which save Artin's case finally. 

Two finally reports were sent to Euler, who as the President of 
the $L$ Academy, was in 
a position to make the final decision for Epstein's case.

As usual, politically very skillful Euler first asked his Assistant to 
investigate what was 
really going on. In fact Euler told his Assistant that 
as Artin's case proved, despite the fact that normally
the internal committee was in a right position to give a proper valuation 
academically, 
occasionally  he was more intended
to accept the conclusion of the external committee, in particular, when the 
person 
under review was doing something \lq revolutionary'. Euler 
went on to say that even it was not proper 
to use the word revolutionary to describe
Epstein's work and that he was a bit \lq suspicious' about the external 
committee: the 
words such as  \lq Meromorphic Extension\rq$\ $ and
\lq Functional Equation\rq$\ $  are \lq too difficult to be understood by 
these noble 
gentlemen', the CEOs of multinational companies who had no special academic
 training in 
topics studied in the $L$ Academy, and it was \lq
not natural' for the external committee to know the secrete criterion of 
the Zeta Institute. 
\lq Find out the true story for me' said Euler to his Assistant. 

The later events proved that Euler was right about almost everything.
 The Assistant, 
after talking with various people involved both within and outside the 
$L$ Academy and 
analyzing
the collected evidences carefully, 
revealed the following picture to Euler: While Riemann was the one who 
openly talked 
about the decision of his committee, his view was shared by all members in 
his Institute.
 By contrast, not only the external committee strongly support Epstein, 
their conclusion 
was also supported by all members of the School of
Eisenstein in the Academy. As a matter of fact, there was a close 
connection between the 
external committee and the School of Eisenstein, say, Lufthansa 
once used the parabolic technique of the School of Eisenstein to 
improve the efficiency of 
its jet fleets \lq quite significantly'... Because of this, the 
Assistant, suggested to 
Euler that while the internal committee made a right academic 
conclusion, the external 
committee, strongly supported by the School of Eisenstein, should 
not be discarded completely. 
As we all know, Euler after bacame the President of the $L$ Academy was 
extremely poor in 
health: He lost his eye sights completely.  While he did have the final 
say, the daily 
matters were in fact operated mainly through two to three of his 
Assistants. The practice 
for Epstein's case was no different...
\vskip 0.30cm
A few days later, in an official poster of the $L$ Academy, we read:
\vskip 0.30cm
\centerline {\bf $L$ Academy Notice No: 23}
\vskip 0.30cm
\hskip 8.0cm Sep. 4, $****^*$
\vskip 0.30cm
{\it We are extremely please to announce that Epstein has been promoted 
to the rank of 
Associate Professor. With the effect from this date of announcement, 
Epstein is transfered 
from the Zeta Institute to the School of Eisenstein of this Academy. 
We thank both external 
and internal committees for their remarkable works.}
\vskip 0.30cm
\hskip 6.0cm the President: Leonardvs EVLER
\vskip 0.30cm
As such the crisis in the $L$ Academy around Epstein's case was solved 
smoothly. 
It was said that Epstein happily worked at the School of Eisenstein 
since after while Riemann 
was also extremely pleased: He  told his members that he was finally able 
to \lq get rid of that fellow'.

This episode of Riemann and Epstein does not stop here.
In the library of the $L$ Academy, from an official recorded, we read:
\lq Riemann and Epstein had their conciliation later on. It was good 
for the Academy -- 
the Non-Commutative Division in the Zeta Institute was created and 
Epstein was awarded 
the Integration Prize by Riemann.'  
\vskip 0.30cm
\noindent  
{\it Remark.} The exact date of the Notice No. 23 is not known. This
remains as one of the hot topics in the study of the history of the Academy. 
The best result so far has been that the last number is
supposed to be 3.
\vskip 0.50cm
\noindent
{\bf Background:}
\vskip 0.30cm
\noindent
(1) Epstein Zeta Functions, despite of their commonly accepted name, 
are indeed a special kind 
of Eisenstein series;
\vskip 0.30cm
\noindent
(2) Even though they are not zeta functions, Epsteins may be used in the 
construction of 
non-abelian zeta functions
for number fields. For example, the rank $r$ non-abelian zeta function for 
${\mathbb Q}$ is 
defined to be
$$\xi_{{\Bbb Q},r}(s):=\int_{{\mathcal M}_{{\mathbb Q},r}}
\widehat E(\Lambda,s)\,
d\mu(\Lambda),\qquad \mathrm{Re}(s)>1,$$ where $\widehat E$ denotes the 
(completed) Epstein, 
and ${\mathcal M}_{{\mathbb Q},r}$ denotes the moduli space of semi-stable 
lattices of rank $r$ 
over ${\mathbb Q}$.
\vfill
\eject
\centerline {\bf {\large File 2: School of Eisenstein}}
\vskip 0.30cm
For a certain period of time after the death of Eisenstein, a well-respected
figure in the Academy, Euler, the President, received enormous amounts of 
appealing letters 
from various interests groups of
relatives, students and follows of Eisenstein. Despite the fact that these 
groups normally 
were at odd with each other, this time, they were quite unified. They all 
claimed that even 
though the official reason for the earlier death of Eisenstein was 
tuberculosis, recently they
 all found that the true 
story was very very different. Their arguments, supported by the world 
famous detective 
agencies such as the Holmes and the Paulo, were that Eisenstein's death 
was mainly due to 
a plot of the L Academy -- the Academy sent Terminators from 20th century 
back to Eisenstein's 
time to take the soul of Eisenstein away around his 20's so that at the 
age of 29, 
Eisenstein died. Indeed, this also explained why as a 
young man, Eisenstein complained of nervous ailments, they claimed.

Due to the fact that the detective agencies were involved, the letters 
were full with details. 
Here were a summary of what they said.
\vskip 0.20cm  
\noindent
(1) Who
\vskip 0.20cm
This was mainly the concern of the students of Eisenstein. In many letters, 
they described 
how Eisenstein's soul was taken away by three different generations of 
Terminators, namely, 
Terminator I, Terminator II and Terminator III. Say, in letter A, it 
claimed that mainly due 
to a certain basic problem in the grand
design, Terminator I was quite primitive and
indeed consisted of two, named Rankin and Selberg respectively; and that
practically it was impossible to distinguish who was who, and they 
could not rule out 
the possibility that there might be more than two Terminator I.

However the mistake in the basic design of Terminators
was corrected later added in the letters B and C: Terminator II named Langlands
was created and sent via the so-called maximal parabolic technique and this 
was further improved when Terminator III named Arthur was created. 
Terminator III, 
manufactured using the so-called total parabolic technique
was so powerful that no further version of Terminators after the mission of
Terminator III was needed. Moreover, in one letter, it claimed that Terminator
III took the time machine at a place called the Earth Dynamic not far away from
Arthur's Seat.
\vskip 0.20cm  
\noindent
(2) Why 
\vskip 0.20cm
This was mainly the concern of relatives of Eisenstein. In many 
letters, they claimed that 
the sole reason for the $L$ Academy to sent these Terminators was that 
afterwards the $L$ 
Academy would  freely use the vast land  donated by Eisensetin. These 
relatives went further 
to claim that it was not fair that after the death of Eisenstein, his 
name was seldom quoted: 
Before Kronecker, only Hurwitz cited Eisensetin's name as a footnote 
in his Thesis, claimed by some  students of Eisenstein. 
\vskip 0.20cm  
\noindent
(3) How 
\vskip 0.20cm
This was mainly the concern of the follows of Eisensetin. In many of 
their letters, they 
described how Terminators completed their missions in every possible detail. 
Many letters were
not pleasant to read when they described how Terminators used  weapons 
called Fourier 
expansions to take the soul of Eisenstein away. However, if we put these 
details aside, 
the most serious
claim in these letters was the follows: Termninator II indeed could clone 
Eisenstein. 
According to them, Terminator II, heavily influenced by the philosophy of 
ancient Turks, believed that the soul of human beings determined every 
aspect of the related 
individuals, their past, their present and their future. They said that 
after his 
understanding of the soul of Eisensetin, Terminator II almost successfully  
cloned  Eisesnetin, 
if this was not prohibited by the law of 20th century.  As a matter of fact, 
asserted in several 
letters, Terminator II ever wrote a book about his study later on, 
but it was such a 
difficult one even some of the most promising young researchers of 
Eisetstein felt that one needed several years to
understand it. To support their claims, one example was provided: 
One researcher, working at a university
belonging to the Ever Green Group, even went so far to show others a 
photo of an infant with the book of Terminator II, to indicate that
his present understanding of the book was no better than that baby.
 (However, there was another extreme:  It is said that Dr. Lion
a legend, spent only about 4 and half months to get a good master of 
Terminator II's work: Without any knowledge of Eisenstein, 
he started reading Terminator II after his guest
Deninger, a leading figure working on the infinite dimensional cohomology 
at the Zeta Institute, left at the beginning of March, 
and was able to give  a series of intensive lectures about the  
applications of Terminator II's work to the non-commutative L functions 
at Kyodai. We also reminder the reader that this story
was disputed by some other letters: After all, Dr. Lion  
was helped by a besutiful manuscript written by Moeglin and Waldspurger of 
the College der Eiffel Tower 
after their seminar on Ternimator II's work, despite the fact that 
Dr. Lion did not really think that 
the second part of the manuscript was a success.) 
\vskip 0.20cm  
\noindent
(4) Demand 
\vskip 0.20cm 
In these letters, they all appealed to the $L$ Academy that due 
to the fact that Eisenstein 
donated all his enormous fortunes to the $L$ Academy, which as the 
result made the expansion of
 the $L$ Academy in an \lq unprecedented scale and enormous pace', 
the $L$ Academy should set a 
new institute named after Eisenstein.
\vskip 0.20cm
The wish of these appeal letters was mainly fulfilled: Soon after, 
the $L$ Academy sent letters 
to all related bodies to inform them that
a new school, the School of Eisesnetin, was founded at the Academy, 
with Langlands as the dean.
\vskip 0.30cm
\noindent
{\bf Background.} Eisenstain series are not integrable over fundamental 
domains, due to the 
fact that the constant terms obtained from the Fourier expansion are 
only of slow growth. 
It was Rankin and Selberg who first used the truncation of Eisenstein 
series to make the 
integration meaningful. This method was generalized by Langlands and 
Arthur later on. Among 
all, clearly, 
Arthur's fundamental work is the most general and hence powerful one.

On the other hand, the constant term, the cut off part in the truncation, 
is really the 
heart of Eisenstein series -- the theory for Eisenstein series, such as 
the meromorphic 
continuation and functional equations, is indeed the one for 
constant terms. By 
analyzing these constant terms, Langlands obtained one of the most deepest 
and most 
beautiful theory in mathematics: Eisenstein series and spectral decompositions.
\vfill
\eject
\centerline {\bf {\large File 3: Siegel's Garden}}
\vskip 0.40cm
Due to the exceptional contribution of Siegel to the $L$ Academy, 
a public garden originally 
donated by the Minkowskis was named the Siegel Garden. At 
the same time the $L$ Academy also announced that the Siegel Garden 
was to be managed by
Borel, a promising young star in the $L$ Academy.

Shortly after the garden was named, and learnt that  most of Siegel's 
ideals in research were 
formulated during his everyday walk,
Borel made a special arrangement so as to let Siegel alone have a walk 
in the garden to 
\lq publicize it'. Siegel accepted the invitation with great pleasure.

The stroll went smoothly. The weather was extremely good. Nice sunshine, 
beautiful scenery, 
vivid wind -- everything was just right for a walk. Siegel enjoyed it 
very much. As usual, 
he started his \lq thinking when walking'. This time was
about a calculation of Riemann he read recently from the manuscript of 
Riemann in the library 
of the University of Gottingen. Without any notice, along the bank 
of a lake, Siegel walked about 
a half day. He became a bit tired. At one point, Siegel realized that 
something possibly 
went wrong as 
he found that the lake became larger and larger, and smelled heavily 
salty breeze. \lq I lost 
my way' said Siegel to himself. So he decided to ask someone for 
help. Unfortunately, there was nobody around -- After learnt that 
Siegel tended to walk alone 
so as to minimize any possible disturbance of his \lq
research activity', the Siegel Garden was closed to public that day 
by the Administration... 
It was recorded that due to the complicated landform of the
Siegel Garden, the rescue team had a \lq great\rq$ $ difficulty to find Siegel.

It was quite embarrassing, and then became quite appearent that the 
Siegel Garden should do 
quite significant work to win the trust from the public and 
especially from Siegel. However 
it was proved that this was not an easy work: not only the financial 
resources was quite limited at the time, 
there were also some technical constrains.
After a series of serious discussions, the Administration decided that as the
first step to win the public trust, the Academy was going to build a bridge 
across the lake so as 
to cut the lake into the inner and outer parts. \lq In this way' they 
claimed, \lq at least, 
the public would  have a good chance to find their own way out when they were 
in the inner garden.'

The plan was confirmed by the board of trustee of the $L$ Academy 
as well. However due to the 
sluggish economical situation around the world, the board 
of trustee worried about the justification of using huge amount of
 money to build 
such a bridge during this difficult period of time. They demanded 
that the cost should be kept 
minimal, and recommended that the official purpose of building the 
bridge should be the 
one in promoting the research activity. \lq All in all, the project should 
be kept as a top secrete', said the Board of Trustee. 

As such the secret project got started. The Lang Design, a world 
famous design
company with Three Little Pigs as the trade mark, was chosen to 
design the bridge, and the Ong 
Construction was chosen to build the bridge. (The choice 
of the Ong Construction was not without objection since \lq the Ong 
Construction was not yet 
a world famous brand', but the Administration cited the following 
reasons to defend their 
decision: not only being a Chinese company,
the Ong Construction was very competitive due to the lower working cost, 
it had the advanced technique called GA truncation to build the bridge.)

The construction of the bridge was a big success: With the use of
GA truncation, it is very elegant --  various parabolic components beautifully
offered curved appearance of the bridge. Later, some who were saw 
the bridge claimed that 
it was not just the beautiful curves which made the bridge so elegant, 
it was the positions of the each parabolic components of the bridge 
which made the bridge 
so beautiful: If it was with other technique, say even the one used  
in creating 
the Terminator III, The result was definitely different, claimed by 
these people.) 

The bridge, full with surprise, was proved to be very important to 
the $L$ Academy as well: 
Due to very careful choice of the location of the bridge using GA 
truncation technique,
the Bridge made the direct connection among various institutes of the
 $L$ Academy possible. For instance,  
later on, people could easily use the bridge to commute between
the non-commutative division of the Zeta Institute and the School 
of Eisensetin.

The $L$ Academy held a big opening ceremony for the bridge. 
The Administration
invited Siegel again to have his walk in the Garden: \lq As 
the bridge was built, you 
certainly would have a quite extraordinary experience walking over it.'
However Siegel, still remembering his first encounter with the garden, 
was not fully convinced. \lq For safety 
and also for sharing my pleasure,' Siegel invited  his best 
friend Mordell to have \lq a small walk' together with him. 

Their walk went quite fine at the very beginning. The weather 
was excellent. The bridge was 
simply \lq beautiful'. Mordell and Siegel were very relax: 
Even they got to find a small secret 
about the bridge -- from \lq mirrow' image of the bridge in the 
lake, they read the words such as Mumford Stability. 
Mathematically it went also 
\lq superb' -- Mordell
told Siegel that the parabolic components (of the bridge)
reminded him hyperbolic 
curves. \lq You see', said, 
Mordell to Siegel, \lq  these hyperbolic curves are also very elegant -- 
only finitely many rational points we can spot there.' \lq That is 
great.' Said Siegel. \lq
Let us check this bridge more carefully. It might offer us more 
interesting ideals.' Siegel continued.

This latest attitude towards the Bridge costed them very badly. 
Their nice feelings did not last 
quite long.  Just about Mordell and Siegel started walking over the bridge, 
with their great care, they spot a tiny bronze mark of Three 
Little Pigs. Surely, they 
immediately realized this is the trademark of the Lang Design... 
It was said that never after
 Mordell and Siegel were seen to walk in the garden.
\vskip 0.40cm
\noindent
{\bf Background.} \vskip 0.30cm
Siegel domains are not bounded and a bit complicated for beginners 
due to the complication
of parabolic subgroups. However we may use geo-arithmetical truncation
associated with intersection stability,
to get  intrinsic compact subsets for them by throwing away cusp 
regions corresponding to 
parabolic subgroups. Over function fields, such a truncation
has its root in Momford's intersection stability via the so-called 
Harder-Narasimhan 
filtration consideration and is first systematically used by Lafforgue. 
As such the 
ill-defined integration of Eisensetin series over the Siegel domain may 
be modified to
 the one over truncated domain so as to obtain the non-abelian zeta and 
$L$ functions. 
\enddocument